\documentclass[12pt]{amsart}
\usepackage{amsmath, amsthm, amssymb,cite,enumerate,slashed}
\usepackage{fullpage}
\usepackage{color}
\usepackage[normalem]{ulem}
\usepackage{hyperref}
\usepackage[english]{babel}
\usepackage{framed}

\allowdisplaybreaks

\newtheorem{thm}{Theorem}
\newtheorem{lem}{Lemma}[section]
\newtheorem{exm}[lem]{Example}
\newtheorem{prop}[lem]{Proposition}
\newtheorem{cor}[lem]{Corollary}
\newtheorem{rmkn}{Remark}[thm]

\newtheorem{rmk}{Remark}[lem]

\newtheorem{defn}[lem]{Definition}
\numberwithin{equation}{section}

\newcommand{\R}{\mathbb{R}}
\newcommand{\D}{\mathcal{D}}

\newcommand{\Z}{\mathbb{Z}}

\newcommand{\eps}{\varepsilon}
\newcommand{\norm}[1]{\left\lVert#1\right\rVert}

\newcommand{\mcl}[1]{\mathcal{ #1 }}
\newcommand{\fm}[1]{\begin{align*} #1 \end{align*}}
\newcommand{\eq}[1]{\begin{equation}\begin{aligned} #1 \end{aligned}\end{equation}}
\newcommand{\lra}[1]{\langle #1 \rangle}
\newcommand{\abs}[1]{\left| #1 \right|}
\newcommand{\kh}[1]{\left( #1 \right)}

\newcommand{\wt}[1]{\widetilde{ #1 }}
\newcommand{\wh}[1]{\widehat{ #1 }}

\DeclareMathOperator{\supp}{supp}
\DeclareMathOperator{\sgn}{sgn}
\allowdisplaybreaks[1]

\begin{document}

\title{Nontrivial global solutions to some quasilinear wave equations in three space dimensions}

\author{Dongxiao Yu}
\address{Department of Mathematics, Vanderbilt University}
\thanks{}
\email{dongxiao.yu@vanderbilt.edu}

\begin{abstract}

In this paper, we seek to construct nontrivial global solutions to some quasilinear wave equations in three space dimensions. We first present a conditional result on the construction of nontrivial global solutions to a general system of quasilinear wave equations. Assuming that a global solution to the geometric reduced system exists and satisfies several well-chosen pointwise estimates, we find a matching exact global solution to the original wave equations. Such a conditional result is then applied to two types of equations which are of great interest. One is John's counterexamples $\Box u=u_t^2$ or $\Box u=u_t u_{tt}$, and the other is the 3D compressible Euler equations with no vorticity. We explicitly construct global solutions to the corresponding geometric reduced systems and show that these global solutions satisfy the required pointwise bounds. As a result, there exists a large family of nontrivial global solutions to these two types of equations.

\end{abstract}

\maketitle
\tableofcontents\addtocontents{toc}{\protect\setcounter{tocdepth}{1}}

\section{Introduction}
This paper is devoted to the study of global solutions\footnote{\label{footnote:futureglobal}By ``global'', we mean ``future global''. A solution is future global if it is defined for all $t\geq 0$. One can also discuss past global solutions, i.e.\ solutions defined for all $t\leq 0$. See Remark \ref{rmkfutureglobal}.} to a system of quasilinear wave
equations in $\R^{1+3}_{t,x}$, of the form
\eq{\label{qwe}g^{\alpha\beta}(u,\partial u)\partial_\alpha\partial_\beta u^{(I)}=f^{(I)}(u,\partial u),\qquad I=1,\dots,M.
}
Here the unknown function $u=(u^{(I)})_{I=1,\dots,M}$ is an $\R^M$-valued function. In \eqref{qwe} we use the Einstein summation convention with the sum taken over $\alpha,\beta=0,1,2,3$ with $\partial_0=\partial_t$ and $\partial_i=\partial_{x_i}$, $i=1,2,3$. We assume that the $g^{\alpha\beta}$  are smooth functions defined in $\R^{M+4M}$, such that $g^{\alpha\beta}=g^{\beta\alpha}$ and $g^{\alpha\beta}(0)=m^{\alpha\beta}$ where $(m^{\alpha\beta})=(m_{\alpha\beta})=\text{diag}(-1,1,1,1)$ is the Minkowski metric (so $g^{\alpha\beta}(0)\partial_\alpha\partial_\beta=\Box=-\partial_t^2+\Delta$). We also assume that  the $f^{(I)}$ are smooth functions defined in $\R^{M+4M}$ such that 
\eq{\label{qwefi}  f^{(I)}(u,\partial u)=f^{I,\alpha\beta}_{JK}\partial_\alpha u^{(J)}\partial_\beta u^{(K)}+O(|u|^4+|\partial u|(|\partial u|^2+|u|^2)),\qquad \text{whenever }|u|+|\partial u|\ll 1.}
In other words, each $f^{(I)}$ vanishes of second order at $0$ and does not contain terms like $u\cdot \partial u$, $u\cdot u$ or $u\cdot u\cdot u$ in its Taylor expansion at $0$.

The system \eqref{qwe} has attracted lots of attention because of its rich physical background. For example, the Einstein vacuum equations, when written in wave coordinates, are of the form \eqref{qwe}; see \cite{MR2134337,MR2680391,MR3638312}. Moreover, the $3$D compressible Euler equations with zero vorticity can be written as a system of quasilinear wave equations; see \cite{MR4109292,MR4011696}.

In general, there is no global existence result for \eqref{qwe}, not even if the initial data is small and localized. For example, John \cite{MR0600571,MR0808321} proved that any nonzero solution to $\Box u=u_t^2$ with $C_c^\infty$ data, or any nonzero solution to $\Box u=u_tu_{tt}$ whose data are $C_c^\infty$ and satisfy an integral inequality, blows up in finite time. To obtain a global existence result for small and localized initial data, we need extra assumptions on the coefficients $g^{\alpha\beta}$ and $f^{(I)}$. One such condition is the \emph{null condition}; see Klainerman \cite{MR0804771,MR0837683} and Christodoulou \cite{MR0820070}. The \emph{weak null condition}, which was first introduced by Lindblad and Rodnianski \cite{MR1994592}, is conjectured to be another sufficient condition for global existence, but this conjecture remains open at this moment\footnote{To be more precise, there are different versions of this conjecture. Some of them have been disproved, while others remain open. We will discuss it in detail in Section \ref{sec1.1}.}.

In this paper, we seek to construct nontrivial global solutions to a general system \eqref{qwe}. To achieve this goal, we consider a backward scattering problem for \eqref{qwe}. Here we identify a notion of asymptotic profile for \eqref{qwe} as a global solution to the geometric reduced system which was introduced in \cite{MR4232783,MR4315017}. The main result is a conditional one: assuming that an asymptotic profile exists and satisfies several well-chosen pointwise bounds, we show that there exists a global solution to \eqref{qwe} that matches this asymptotic profile at infinite time. While the proof here is essentially the same as those in \cite{MR4232783} and \cite[Chapter 3]{MR4315017} on the modified wave operators for a scalar quasilinear wave equation, in this paper we extend these previous results to a much larger class
of quasilinear wave equations and relax the assumptions on the asymptotic profile.

Note that the construction above does not necessarily generate a nontrivial global solution to an arbitrary system of quasilinear wave equations. This is because, for a general \eqref{qwe}, we do not know whether the associated geometric reduced system admits a nontrivial global solution satisfying the required bounds or not. Despite this,  our construction does apply to two types of equations which are of great interest. One is John's counterexamples $\Box u=u_t^2$ and $\Box u=u_t u_{tt}$.\footnote{We  apply the main result above directly to construct global solutions to $\Box u=u_t^2$. However, for $\Box u=u_tu_{tt}$, we apply the main result to a system of wave equations obtained by differentiating $\Box u=u_tu_{tt}$, and then show that these solutions to the differentiated equations induce global solutions to $\Box u=u_tu_{tt}$. See Sections \ref{sec1.4.1}, \ref{sec8.2}, and \ref{sec8.3} for more details. } The other is the 3D compressible Euler equations with no vorticity. For these two types of equations, we explicitly construct nonzero global solutions to the corresponding geometric reduced systems and show that these global solutions satisfy the required pointwise bounds in our main result above. Therefore, there exists a large family of nontrivial global solutions to both these types of equations. We remark that this existence result does not contradict John's finite time blowup results. The initial data of the global solutions constructed here are not localized enough (e.g.\ not in  $C_c^\infty$), while John assumed that the data are nonzero and $C_c^\infty$  in his finite time blowup results. We also remark that for each $C>0$, $\ln(t+C)$ is a global solution to $\Box u=u_t^2$ and that any polynomial of $(t,x)$ of degree at most one is a global solution to $\Box u=u_tu_{tt}$. The solutions constructed in this paper, however, have very different asymptotic behaviors when compared with those known solutions. We will discuss these differences later in the paper.


\subsection{Background}\label{sec1.1}

Let us consider  a generalization of the system \eqref{qwe} in~$\R^{1+3}_{t,x}$
\begin{equation}\label{nlw}
\square u^I=F^I(u,\partial u,\partial^2u),\qquad I=1,2,\dots,M.
\end{equation}
The nonlinear term is assumed to be smooth with the Taylor expansion
\begin{equation}\label{nonlinearity}F^I(u,\partial u,\partial^2u)=\sum a_{\alpha\beta,JK}^I\partial^\alpha u^J\partial^\beta u^K+O(|u|^3+|\partial u|^3+|\partial^2 u|^3).\end{equation}
The sum  is taken over all $1\leq J,K\leq M$ and all multiindices $\alpha,\beta$ with $|\alpha|\leq |\beta|\leq 2$, $|\beta|\geq1$ and $|\alpha|+|\beta|\leq 3$. Besides, the coefficients $a_{\alpha\beta,JK}^I$'s are all universal constants.

Since 1980's, several results on the lifespan of the solutions to the Cauchy problem \eqref{nlw} with sufficiently small $C_c^\infty$ initial data have been proved.  For example, John \cite{MR0600571,MR0808321} proved that \eqref{nlw} does not necessarily have a global solution; in fact, any nonzero solution to $\Box u=u_t^2$ with $C_c^\infty$ data, or any nonzero solution to $\Box u=u_tu_{tt}$ whose data are $C_c^\infty$ and satisfy an integral inequality, must blow up in finite time. In contrast, in $\R^{1+d}$ with $d\geq 4$, we do have a  small data global existence result for \eqref{nlw}. We refer our readers to \cite{MR1120284,MR2679723,MR4046191}. For arbitrary nonlinearities in three space dimensions, the best result on the lifespan is expected to be the almost global existence: the solution exists for $t\leq \exp(c/\eps)$ where $\eps\ll 1$ is the size of the initial data. This is indeed the case when $M=1$, or when the $F^I$'s do not depend on $u$ (i.e. when $F^I=F^I(\partial u,\partial^2u)$). See  \cite{MR1047332,MR0745325,MR0784477,MR0897781,MR1466700,MR2015331,MR1945285}. 
To the best of the author's knowledge,  the almost global existence result has not been proved when $M>1$ and when the $F^I$'s depend on $u$. However, we refer to a recent work by Metcalfe and Rhoads \cite{MR4565922} where the authors proved a result weaker than the almost global existence ($T_{\max}\geq \exp(c/\eps^{1/3})$).

Following John's results, there have been several results on the blowup mechanism for general quasilinear wave equations. For simplicity, we restrict our discussions to a scalar quasilinear wave equation:
\eq{\label{shockqwe3d} g^{\alpha\beta}(\partial u)\partial_\alpha\partial_\beta u=F(\partial u),\qquad\text{\rm in }\R^{1+3}.}
Here the $g^{\alpha\beta}$ and $F$ are smooth functions of $\partial u$ and satisfy the same assumptions as in \eqref{qwe}. In the case when $F\equiv 0$ and when the equation violates the classical null condition\footnote{We will define the null condition in the next paragraph. Here it is enough to realize that $\Box u=u_tu_{tt}$ violates the null condition.}, both Alinhac \cite{MR1867312} and Speck \cite{MR3561670} proved shock formation for \eqref{shockqwe3d}. The difference is that Alinhac's result holds for nontrivial initial data verifying some uniqueness and nondegeneracy conditions, while  Speck eliminated these uniqueness and nondegeneracy conditions. Moreover, Speck was also able to provide information about the maximal development of the data. Speck's work was motivated by Christodoulou \cite{MR2284927} who showed a remarkable result giving a detailed description of shock formation in solutions to
the relativistic Euler equations in three space dimensions. For an overview of these three results and for more discussions on shock formation, we refer to \cite{MR3474069} and the references therein. We also refer to Miao and P.\ Yu \cite{MR3595936} for shock formation of a quasilinear wave equation of the form \eqref{shockqwe3d} satisfying the null condition. Unlike all the aforementioned results where the initial data are assumed to be small, the blowup in \cite{MR3595936} only occurs in the large data case. Besides, Speck \cite{MR4047642} proved that the equation
\fm{-u_{tt}+\frac{\Delta u}{1+u_t}=-u_t^2} exhibits a distinct kind of blowup called the ODE-type blowup. Meanwhile, the blowup mechanism for the semilinear example, $\Box u=u_t^2$, remains poorly understood as of the present; see, e.g., the discussion in \cite[Section 1]{MR4047642}.

In contrast to the finite-time blowup in John's examples,  Klainerman \cite{MR0837683} and Christodoulou \cite{MR0820070} proved that the null condition is sufficient for small data global existence. The null condition, first introduced by Klainerman \cite{MR0804771}, states that for each $1\leq I,J,K\leq M$ and for each $0\leq m\leq n\leq 2$ with $m+n\leq 3$, we have\footnote{Since there is no term of the form $u\cdot u$ in the nonlinearities \eqref{nonlinearity}, the identity \eqref{nullcond} holds automatically for $m=n=0$.}
\begin{equation}\label{nullcond}A_{mn,JK}^I(\omega):=\sum_{|\alpha|=m,|\beta|=n}a_{\alpha\beta,JK}^I\widehat{\omega}^{\alpha}\widehat{\omega}^{\beta}=0,\hspace{1cm}\text{for all }\widehat{\omega}=(-1,\omega)\in\R\times\mathbb{S}^2.\end{equation}
Equivalently, we assume that $A_{mn,JK}^I\equiv 0$  on the null cone $\{\xi_0^2=\sum_{j=1}^3\xi_j^2\}$. The null condition leads to cancellations in the nonlinear terms \eqref{nonlinearity} so that the nonlinear effects of the equations are much weaker than the linear effects. 
We also remark that the null condition is not necessary for  small data global existence. One such example is the Einstein vacuum equations in wave coordinates; see \cite{MR1994592,MR2134337}.  We refer our readers to \cite{MR1957533} for a general introduction to the null condition. We also refer to \cite{MR2455195,MR2110542,MR1857981} and the references therein for a generalized version of the global existence results by Klainerman  and Christodoulou  for  nonlinear wave equations involving multiple wave speeds.

Later, Lindblad and Rodnianski \cite{MR1994592,MR2134337} introduced the weak null condition. To state this condition, we start with the asymptotic equations first introduced by H\"{o}rmander \cite{MR0897781,MR1466700,MR1120284}. We make the ansatz\begin{equation}\label{ansatz}
u^I(t,x)\approx \eps r^{-1}U^I(s,q,\omega),\hspace{1cm}r=|x|,\ \omega_i=x_i/r,\ s=\eps\ln(t),\ q=r-t,\ 1\leq I\leq M.
\end{equation}
Assuming that $t=r\to\infty$, we substitute this ansatz into \eqref{nlw} and compare the coefficients of  terms of order $\eps^2 t^{-2}$.  Nonrigorously, we can derive the following asymptotic PDE
\begin{equation}\label{asypde11}2\partial_s\partial_q U^I=\sum A_{mn,JK}^I(\omega)\partial_q^mU^J\partial_q^nU^K.\end{equation}
Here $A_{mn,JK}^I$ is defined in \eqref{nullcond} and the sum is taken over $1\leq J,K\leq M$ and $0\leq m\leq n\leq 2$ with $n\geq 1$ and $m+n\leq 3$. We say that the \emph{weak null condition} is satisfied if \eqref{asypde11} has a global solution for all $s\geq 0$ and if the solution and all its derivatives grow at most exponentially in $s$, provided that the initial data decay sufficiently fast in $q$. Note that under the null condition, \eqref{asypde11} becomes $\partial_s\partial_q U=0$. As a result, the weak null condition follows from and is indeed weaker than the null condition. We also remark that there are many examples of \eqref{nlw} satisfying the weak null condition and admitting small data global existence at the same time. One such example is the Einstein vacuum equations in wave coordinates; see \cite{MR2134337,MR2680391}.  There are also many examples violating the weak null condition and admitting finite-time blowup at the same time. Two such examples are $\square u=u_t u_{tt}$ and $\square u=u_t^2$: the corresponding asymptotic equations are  $(2\partial_s+U_q\partial_q)U_q=0$ (Burgers' equation) and $2\partial_s U_q=U_q^2$, respectively, whose solutions are known to blow up in finite time for $U|_{s=0}\in C_c^\infty\setminus\{0\}$.

It is natural to conjecture\footnote{While this conjecture is often attributed to Lindblad and Rodnianski, they did not make such a conjecture in \cite{MR1994592,MR2134337,MR2680391}. Instead, in \cite[Footnote 6]{MR2680391}, they mentioned that it was unclear whether this weak null condition is sufficient for small data global existence.} that the weak null condition is sufficient for small data global existence. However, if we consider a general system \eqref{nlw} with nonlinearities satisfying \eqref{nonlinearity} without adding further assumptions, this conjecture is false. Kadar \cite{kadar2023small} recently presented a system of semilinear wave equations satisfying the weak null condition but does not admit global existence:
\eq{\label{kadarwnceq}\left\{\begin{array}{l}\Box \phi_1=\phi_3^3,\\
\Box \phi_2=(\partial_t\phi_1)^2,\\
\Box \phi_3=(\partial_t\phi_2)^2+\phi_1^3.\end{array}\right.}
In his work, Kadar mentioned that this negative result is due to the undifferentiated terms ($\phi_3^3$ and $\phi_1^3$) in the semilinear nonlinearities. He also stated a more restricted version of the weak null conjecture by only considering those systems \eqref{nlw} of wave equations whose nonlinearities contain only differentiated semilinear terms; see \cite[Section 1.2]{kadar2023small} and \cite[Conjecture 5.3.2]{keir}. This restricted version, to the author's knowledge, is still open. On the positive side, Keir has recently made partial progress towards this conjecture. In \cite{keir}, he proved the small data global existence for a large class of quasilinear wave equations satisfying the weak null condition, significantly enlarging upon the class of equations for which global existence is known.  In \cite{keir2},  he also proved that if the solutions to the asymptotic system are bounded (given small initial data) and stable against rapidly decaying perturbations, then the corresponding system of nonlinear wave equations admits small data global existence.

We have a final remark for this subsection. When we discuss small data global existence, what we mean is that there exists a global solution for an  \emph{arbitrary}  initial data which is small and localized enough. If a certain system of nonlinear wave equations does not admit small data global existence, it is still possible that this system admits a global solution for \emph{some} small and localized data, or for  small and \emph{nonlocalized} data. The main goal of this paper is to construct examples that support this remark.

\subsection{The geometric reduced system}\label{s1introgrs}
The main tool to construct global solutions in this paper is the geometric reduced system introduced by the author. In \cite{MR4232783}, this reduced system is derived in the scalar case, and it is extended to the vector-valued case in Chapter 2 of the author's PhD dissertation \cite{MR4315017}. Our analysis starts as in H\"{o}rmander's derivation in \cite{MR0897781,MR1466700,MR1120284}, but diverges at a key point: the choice of $q$ is different. One may contend from the author's previous papers \cite{MR4315017,MR4232783,yu2021} that this new system is more accurate than \eqref{asypde11}, in that it both describes the long time evolution and contains full information about it.

To derive the new asymptotic equations, we still make the ansatz  \eqref{ansatz}, but now we replace $q=r-t$ with an optical function $q(t,x)$. That is,  $q(t,x)$ is a solution to the eikonal equation related to \eqref{qwe}\begin{equation}\label{eikeqn}g^{\alpha\beta}(u,\partial u)\partial_\alpha q\partial_\beta q=0.\end{equation} We remark that eikonal equations have been used in Alinhac \cite{MR2003417} and Lindblad \cite{MR2382144} on small data global existence for quasilinear wave equations; we refer to \cite[Section 1.2]{MR4232783} for a detailed discussion. We also refer to \cite{MR1316662,MR3638312} where eikonal equations were used to study the asymptotic behaviors of solutions to the Einstein vacuum equations. All these works suggest that eikonal equations play an important role in the study of the long time dynamics of \eqref{qwe} and of the Einstein vacuum equations.  Meanwhile, eikonal equations play an important role both in the proof of sharp local wellposedness for nonlinear wave equations and in the proof of the bounded $L^2$ curvature conjecture for the Einstein vacuum equations. We refer our readers to \cite{MR2178963,MR3402797} and the references therein.

Since $u$ is unknown, it is difficult to solve \eqref{eikeqn} directly. Instead,  we introduce a new auxiliary function $\mu=\mu(s,q,\omega)$ such that $q_t-q_r=\mu$. Note that $\mu$ measures the density of the stacking of the level sets of $q$; we refer to \cite{MR3561670} where Speck defined the foliation density (denoted by $1/\mu$ there) which is closely related to our definition of $\mu$ here. From \eqref{eikeqn}, we can express $q_t+q_r$ in terms of $\mu$ and $U$, and then solve for all partial derivatives of $q$, assuming that all the angular derivatives are negligible. Then from \eqref{qwe}, we can derive the following asymptotic equations for $\mu=\mu(s,q,\omega)$ and $U=(U^{(I)})(s,q,\omega)$:
\eq{\label{asypde2}\left\{
\begin{array}{l}\displaystyle
\partial_s(\mu U_q^{(I)})=-\frac{1}{4}F^{I}_{2,JK}(\omega) \mu^2U_q^{(J)} U_q^{(K)},\qquad I=1,\dots,M;\\[1em]\displaystyle
\partial_s\mu=\frac{1}{4}G_{2,J}(\omega) \mu^2U_q^{(J)}-\frac{1}{8}G_{3,J}(\omega)\mu^2\partial_q(\mu U_{q}^{(J)}).
\end{array}\right.}
In this system we use the Einstein summation convention. For each $\wh{\omega}=(-1,\omega)\in\R\times\mathbb{S}^2$ and each $1\leq I,J,K\leq M$, we set
\eq{G_{2,J}(\omega):=g^{\alpha\beta}_J\widehat{\omega}_\alpha\widehat{\omega}_\beta,\qquad G_{3,J}(\omega):=g^{\alpha\beta\lambda}_J\widehat{\omega}_\alpha\widehat{\omega}_\beta\widehat{\omega}_\lambda,\qquad F_{2,JK}^I(\omega):=f^{I,\alpha\beta}_{JK}\widehat{\omega}_\alpha\widehat{\omega}_\beta.}
Here the coefficients come from the  Taylor expansions of $g^{**}$ and $f^*$ in \eqref{qwe}
\eq{g^{\alpha\beta}(u,\partial u)&=m^{\alpha\beta}+g_J^{\alpha\beta}u^{(J)}+g_J^{\alpha\beta\lambda}\partial_\lambda u^{(J)}+O(|u|^2+|\partial u|^2),\\
f^{(I)}(u,\partial u)&=f^{I,\alpha\beta}_{JK}\partial_\alpha u^{(J)}\partial_\beta u^{(K)}+O(|u|^3+|\partial u|^3).}
We call this new system \eqref{asypde2} of asymptotic equations the \emph{geometric reduced system} for \eqref{qwe} since it is related to the geometry of the null cone with respect to the Lorentzian metric $(g_{\alpha\beta})=(g^{\alpha\beta}(u,\partial u))^{-1}$ instead of the Minkowski metric.  Heuristically,  if there exists a global solution to the system of quasilinear wave equations \eqref{qwe}, then one expects it to correspond to an approximate solution to this geometric reduced system, and to be well approximated by an exact solution to the geometric reduced system. For a derivation of \eqref{asypde2}, we refer our readers to  \cite[Section 3]{MR4232783}, or the author's PhD dissertation \cite[Chapter 2]{MR4315017}. For some examples of the geometric reduced system, we refer to Section \ref{sec1.4} and Section \ref{sec8} in this paper, or \cite[Example 2.3 and 2.4]{MR4315017}. 

Let us discuss the solvability of the reduced system \eqref{asypde2} in general. One can view \eqref{asypde2} as a system of quadratic ODE's coupled to a Riccati differential equation. In fact, from the first $M$ equations in \eqref{asypde2}, we notice that $(\mu U_q^{(J)})_{J=1,\dots,M}$ is a solution to a system of quadratic ODE's of the form
\eq{\label{s1ode}\dot{X}^J=X^TA^JX,\qquad J=1,\dots,M.}
Here the unknown $X=(X^J(s))\in\R^M$ is viewed as a column vector and each $A^J$ is an $M\times M$ matrix with real constant entries. This system determines $\mu U_q^{(*)}$. Once $\mu U_q^{(*)}$ is known, the last equation in \eqref{asypde2} becomes a Riccati differential equation for $\mu$. As a result, if we choose the initial data appropriately, we are able to solve \eqref{asypde2} at least locally in $s$ by applying Picard's theorem. However, given an arbitrary system \eqref{asypde2} and an arbitrary set of initial data, we are unlikely to determine whether the corresponding local solution can be extended globally for all $s\geq 0$. The reason is that, to the author's knowledge, the solvablity of a general system \eqref{s1ode} has not been completely understood\footnote{
Of course, in some special cases (such as when $M=1$ or $A^*\equiv 0$), we can solve \eqref{s1ode} explicitly. We also refer to \cite{MR1320184,MR132743} and the references therein for some interesting work on the connection between solutions to \eqref{s1ode} and the algebraic structures from the coefficients $A^*$.}. Despite this issue, it is still very interesting to study the connection between global solutions to \eqref{asypde2} and those to \eqref{qwe}.  As a result, in this paper, we will prove results based on the assumption that a global solution to \eqref{asypde2} exists.  

The geometric reduced system generates a new notion of asymptotic profile for nonlinear wave equations. As an application, it allows us to establish the modified scattering\footnote{The word ``modified scattering'' can have a different meaning in other papers. Here we mean that the scattering for \eqref{qwe2} is not linear. In general, a global solution to \eqref{qwe2} with small $C_c^\infty$ data does not scatter to a solution to the linear wave equation $\Box u=0$, due to the different pointwise decay rates in time. See \cite{MR2382144}.} for a scalar quasilinear wave equation satisfying the weak null condition:
\eq{\label{qwe2}g^{\alpha\beta}(u)\partial_\alpha\partial_\beta u=0\qquad\text{in }\R^{1+3}.}
It was proved by Lindblad \cite{MR2382144} that \eqref{qwe2} admits a global existence result for sufficiently small $C_c^\infty$ data. In modified scattering, we seek to answer the following two main problems:
\begin{enumerate}[1.]

\item \emph{Asymptotic completeness}. Given an exact global solution to the model equation, can we  find the corresponding asymptotic profile?

\item \emph{Existence of (modified) wave operators}. Given  an asymptotic profile, can we construct a unique exact global solution to the model equation which matches  the asymptotic profile at  infinite time?
\end{enumerate}
There have been only a few previous results on the (modified) scattering for general quasilinear wave equations and the Einstein equations.  In \cite{dafermos2013scattering}, Dafermos, Holzegel, and Rodnianski gave a scattering theory construction of nontrivial black hole solutions to the Einstein vacuum equations.
That is a backward scattering problem in General Relativity. In \cite{MR4591580}, Lindblad and Schlue  proved the existence of the wave operators for the semilinear models of Einstein's equations. In \cite{MR4078713}, Deng and Pusateri used the original H\"{o}rmander's asymptotic system \eqref{asypde11} to prove a partial scattering result for \eqref{qwe2}.  In their proof, they applied the spacetime resonance method; we refer to \cite{MR3084697,MR3202837} for some earlier applications of this method to the first order systems of wave equation. In these works, the authors made use of H\"{o}rmander's asymptotic equations associated to \eqref{qwe2}:
\eq{2\partial_s\partial_qU=G(\omega)U\partial_q^2U}
where $G(\omega):=g^{\alpha\beta}_0\wh{\omega}_\alpha\wh{\omega}_\beta$ and $g^{\alpha\beta}_0:=\frac{d}{du}g^{\alpha\beta}(u)|_{u=0}$.
Recently, in \cite{MR4232783,yu2021,MR4315017},   the author answered both of these problems and thus established the modified scattering for \eqref{qwe2} by making use  of the new geometric reduced system associated to \eqref{qwe2}:
\eq{\label{qwe2redsys}\left\{\begin{array}{l}
\displaystyle \partial_s(\mu U_q)=0,\\
\displaystyle \partial_s \mu=\frac{1}{4}G(\omega)\mu^2U_q.
\end{array}\right.}

In this paper, we will apply the method in \cite{MR4232783,MR4315017} to a larger family of quasilinear wave equations. If a global solution to the geometric reduced system for \eqref{qwe} exists and satisfies several good pointwise estimates, then we will construct a matching global solution to \eqref{qwe}. Therefore, the papers \cite{MR4232783,MR4315017} offer us a systematic way to construct global solutions to \eqref{qwe}. 

We finish this section by remarking that both H\"ormander's asymptotic equations and the geometric reduced system only provide good approximations of global solutions to \eqref{qwe} at the null infinity. That is, if $u$ is a solution to \eqref{qwe} and if $U$ is the corresponding solution to H\"ormander's asymptotic equations, then we only expect $u\approx\eps r^{-1}U$ for $(t,x)$ sufficiently close to the light cone (e.g.\ $|r-t|\ll t$). Similarly for the geometric reduced system. One can also see this from the estimates in Theorem \ref{mthm}. The reason is that, when we derive these asymptotic equations, we take $r=t$, send $t\to\infty$, and omit those lower order terms. Such lower order terms, however, may not be negligible when $|r-t|$ is large. It is also very interesting to study the asymptotic behaviors of the global solutions to \eqref{qwe} at the spacelike infinity (e.g.\ when $r/t>3/2$) and at the timelike infinity (e.g.\ when $r/t<1/2$), if such global solutions exist.

\subsection{The main result}
\label{sec1.3}

We now clarify what assumptions are needed in this paper.  Fix two constants  $\gamma_+>1$ and $\gamma_->2$.  Let   $(\mu,U)=(\mu,(U^{(I)}))(s,q,\omega)$ be a  solution to the geometric reduced system \eqref{asypde2} defined for all $s\geq -\delta_0$ and $(q,\omega)\in\R\times\mathbb{S}^2$ where   $0<\delta_0<1$ is a fixed constant. Suppose that at each $(s,q,\omega)\in[-\delta_0,\infty)\times\R\times\mathbb{S}^2$ we have
\eq{\label{def1.1a3}-C\exp(C s )\leq\mu\leq -C^{-1}\exp(-C s )<0;}
\eq{\label{def1.1a4}|\mu_q|\lesssim \lra{q}^{-1-\gamma_{\sgn (q)}}|s\mu|;}
\eq{\label{def1.1a5}|\mu U_q|\lesssim 1;}
\eq{\label{def1.1a51}\sum_{J,K}|G_{3,J}\mu\partial_q(\mu U_q^{(K)})|\lesssim 1;}
\eq{\label{def1.1a52}\sum_{J,K}|G_{2,J}\partial_q(\mu U_q^{(K)})|\lesssim \lra{q}^{-\gamma_{\sgn(q)}};}
\eq{\label{def1.1a6}|\partial_s^a\partial_q^b\partial_\omega^c(\mu+2)|\lesssim_{a,b,c}\exp(C_{a,b,c} s )\lra{q}^{-b-\gamma_{\sgn(q)}},\qquad \forall a,b,c\geq 0;}
\eq{\label{def1.1a7}|\partial_s^a\partial_\omega^cU|\lesssim_{a,c}\exp(C_{a,c} s )\lra{\max\{0,-q\}}^{1-\gamma_-},\qquad \forall a,c\geq 0;}
\eq{\label{def1.1a8}|\partial_s^a\partial_q^b\partial_\omega^cU_q|\lesssim_{a,b,c}\exp(C_{a,b,c} s )\lra{q}^{-b-\gamma_{\sgn(q)}},\qquad \forall a,b,c\geq 0.}
Here  we set $\gamma_{\sgn(q)}=\gamma_+$ if $q\geq 0$ and $\gamma_{\sgn(q)}=\gamma_-$ otherwise.
A global smooth solution $(\mu,U)$ to \eqref{asypde2} satisfying \eqref{def1.1a3}--\eqref{def1.1a8} is said to be  $(\gamma_+,\gamma_-)$-\emph{admissible}. 
Note that the estimates \eqref{def1.1a3}--\eqref{def1.1a8} are extracted from \cite{MR4232783,MR4315017}. In fact, for all $A=A(q,\omega)\in C_c^\infty(\R\times\mathbb{S}^2)$, we notice that
\eq{\left\{\begin{array}{l}
\displaystyle U(s,q,\omega)=\int_{-\infty}^qA(p,\omega)\exp(\frac{1}{2}G(\omega)A(p,\omega)s)\ dp,\\
\displaystyle \mu(s,q,\omega)=-2\exp(-\frac{1}{2}G(\omega)A(q,\omega)s),
\end{array}\right.}
is a smooth global solution to the associated geometric reduced system \eqref{qwe2redsys}. It has been proved in \cite{MR4232783,MR4315017} that this solution $(\mu,U)$ satisfies the bounds \eqref{def1.1a3}--\eqref{def1.1a8} for all $\gamma_+>1$ and $\gamma_->2$. Several bounds for the matching global solution $u$ to \eqref{qwe2} were derived from \eqref{def1.1a3}--\eqref{def1.1a8} in \cite{MR4232783,MR4315017}, and they were used to show the existence of modified wave operators for~\eqref{qwe2}. A more detailed discussion on the motivation behind these estimates can be found in Section~\ref{sec3.2assumptions}.

As mentioned in Section \ref{s1introgrs}, we do not know if there exists a nontrivial global solution to an arbitrary geometric reduced system \eqref{asypde2}, not to mention a nontrivial admissible global solution. However, we emphasize that for John's counterexamples and the 3D compressible Euler equations with no vorticity, we can explicitly solve the corresponding geometric reduced systems; see \eqref{exm141.sol} and \eqref{exm142.sol} in Section~\ref{sec1.4}. These explicit solutions are admissible and global if we choose the initial data appropriately. Thus, at least for the key examples studied in this paper, there does exist a large family of nontrivial admissible global solutions to the corresponding geometric reduced systems.

Starting with a $(\gamma_+,\gamma_-)$-admissible global solution  $(\mu,U)$, we can construct an approximate solution to \eqref{qwe} as follows. We first make a change of coordinates. For a small $\eps>0$,  we set $s=\eps\ln(t)-\delta$ with $0<\delta<\delta_0$. We remark that this choice of $s$ is related to the almost global existence, since now $s=0$ if and only if $t= e^{\delta/\eps}$. In fact, when $t\leq e^{\delta/\eps}$, we expect the solution to \eqref{qwe} behaves as a solution to $\Box u=0$, so our asymptotic equations play a role only when $t\geq e^{\delta/\eps}$. Set $\Omega:=\{t> 1/\eps,\ |x|> t/2\}$ and let $q=q(t,x)$ be the solution to \fm{q_t-q_r=\mu(\eps\ln(t)-\delta,q(t,r,\omega),\omega)\quad\text{in }\Omega;\qquad q(2|x|,x)=-|x|.} We can use the method of characteristics to solve this equation. Then, any function of $(s,q,\omega)$ induces a new function of $(t,x)$ in $\Omega$. With an abuse of notation, we set \fm{U^{(I)}(t,x)=U^{(I)}(\eps\ln(t)-\delta,q(t,x),\omega),\qquad I=1,\dots,M.} We can prove that, near the light cone $\{t=r\}$, $\eps r^{-1}U$ is an approximate solution to \eqref{qwe}, and $q(t,x)$ is an approximate optical function, i.e.\ an approximate solution to the eikonal equation corresponding with the metric $g^{\alpha\beta}(\eps r^{-1}U,\partial (\eps r^{-1}U))$.

We can now state the main theorem of this paper.  Let $Z^I$ denote a product of $|I|$  commuting vector fields: translations $\partial_\alpha$, scaling $S=t\partial_t+r\partial_r$, rotations $\Omega_{ij}=x_i\partial_j-x_j\partial_i$, and Lorentz boosts $\Omega_{0i}=x_i\partial_t+t\partial_i$.

\thm{\label{mthm}
Consider a system of quasilinear wave equations \eqref{qwe} satisfying \eqref{qwefi}.
Suppose that there exists a $(\gamma_+,\gamma_-)$-admissible global solution $(\mu,U)$ to the geometric reduced system \eqref{asypde2} of \eqref{qwe}  for all $s\geq- \delta_0$. Here $\gamma_+>1$, $\gamma_->2$ and $0<\delta_0<1$ are fixed constants.  Fix $\gamma_1\in(2,\min\{4,2(\gamma_--1)\})$, an integer $N\geq 2$ and a sufficiently small $\eps>0$ depending on  $\gamma_\pm,\gamma_1,(\mu,U),N$\footnote{In other words, there exists a small constant $\eps_0>0$ depending on $\gamma_\pm,\gamma_1,(\mu,U),N$ such that the following conclusions hold for all $\eps\in(0,\eps_0)$. We refer to Section \ref{sec2.1} for the convention on $\eps$ used in this paper.}. Let $q(t,x)$  and $U(t,x)$ be the associated approximate optical function and  asymptotic profile. Then, there is a $C^{N}$ solution $u$ to \eqref{qwe}  for $t\geq 0$ with the following properties:
\begin{enumerate}[{\normalfont (i)}]
\item The solution satisfies the energy bounds: for all $|I|\leq N-1$, we have
\fm{&\norm{w_0^{1/2}\partial Z^I(u-\eps r^{-1}U)(t)}_{L^2(\{x\in\R^3:\ |x|/t\in[1/2,3/2]\})}+\norm{w_0^{1/2}\partial Z^Iu(t)}_{L^2(\{x\in\R^3:\ |x|/t\notin[1/2,3/2]\})}\\&\lesssim_I \eps t^{-1/2+C_I\eps},\qquad \forall t\geq 1/\eps;}
\fm{&\norm{w_0^{1/2}\partial Z^Iu(t)}_{L^2(\R^3)}\lesssim_I \eps,\qquad \forall 0\leq t\leq 1/\eps.}
Here $w_0=1_{r\geq t}+\lra{r-t}^{\gamma_1}1_{r<t}$ where $\gamma_1$ is chosen above.
\item The solution satisfies the  pointwise bounds: for all $|I|\leq N-1$ and all $(t,x)$ with $t\geq 1/\eps$, \fm{
&|\partial Z^I(u-\eps r^{-1}U)(t,x)| 1_{|x|/t\in[1/2,3/2]}+|\partial Z^Iu(t,x)|1_{|x|/t\notin[1/2,3/2]}\\
&\lesssim_I \eps t^{-1/2+C_I\eps}\lra{t+r}^{-1}\lra{t-r}^{-1/2}w_0^{-1/2},}
\fm{&|Z^I(u-\eps r^{-1}U)(t,x)|1_{|x|/t\in[1/2,3/2]}+| Z^Iu(t,x)|1_{|x|/t\notin[1/2,3/2]}\\
&\lesssim_I \eps t^{-1/2+C_I\eps/2} (t^{-1}\lra{r-t}^{-(\gamma_1-1)/2}1_{r< t}+t^{-1}\lra{r-t}^{1/2}1_{t\leq r\leq 2t}+r^{-1/2}1_{r>2t}).}
Moreover, whenever $t\geq 1/\eps$ and $r/t\in[1/2,3/2]$,
\fm{|(\partial_t-\partial_r)u-\eps r^{-1}\mu U_q|\lesssim \eps t^{-3/2+C\eps}\lra{t-r}^{-1/2}w_0^{-1/2}.}
\end{enumerate}
Here all the constants $C_I$ and the implicit constants in $\lesssim_I$ are independent of the choice of~$\eps$. We also emphasize that in this paper, any constant without the subscript $\eps$ must be independent of $\eps$; see Section \ref{sec2.1}.
}

\rmkn{\label{rmkthmext}\rm This theorem is a generalized version of  \cite[Theorem 1]{MR4232783} or   \cite[Theorem 3.1]{MR4315017}. There are improvements in two aspects. First, we are working on a much larger class of quasilinear wave equations. Second, we relax the assumptions on the global solutions to the geometric reduced system. In \cite{MR4232783,MR4315017}, we assume that  $U_q$ vanishes whenever $q\leq -R$ for some constant $R>0$. Under such assumptions, we can show that the global solutions to the quasilinear wave equations vanish in the light cone $\{r-t<-R\}$. In this paper, we replace these vanishing assumptions with some weaker pointwise decaying assumptions such as \eqref{def1.1a8}.

For the scalar equation \eqref{qwe2}, our assumptions on the admissible global solution $(\mu,U)$ still seem too strong from the point of view of the forward problem. Suppose that $u$ is a global solution to \eqref{qwe2} with $C_c^\infty$ initial data of size $\eps\ll1$. In \cite{yu2021}, we find an exact solution $(\mu,U)$ to the geometric reduced system \eqref{qwe2redsys} matching the global solution $u$ at infinite time\footnote{This solution $(\mu,U)$ may depend on $\eps$. However, see Remark \ref{rmkmthmdepeps} below.}. However, while it has not been proved yet, we tend to believe that in general, this $(\mu,U)$ is not $(\gamma_-,\gamma_+)$-admissible for all $\gamma_+>1$ and $\gamma_->2$. For example, we proved that $|\mu U_q|\lesssim \lra{q}^{-1+C\eps}$ for $q<0$ in \cite{yu2021}, while we require $|\mu U_q|\lesssim \lra{q}^{-2-}$ for $q<0$ in the definition of admissibility. The author is currently working on a project to study whether the assumptions on the admissible global solutions in this paper can be further relaxed\footnote{For example, is it sufficient to assume that $\gamma_->1$ in the definition of admissibility instead of $\gamma_->2$?}.  }

\rmkn{\rm\label{mthmrmk3} Recall that the main goal of this paper is to construct nontrivial global solutions to \eqref{qwe}. One may guess that as long as we start with a nontrivial admissible global solution $(\mu,U)$ in Theorem \ref{mthm} (i.e.\ $U\not\equiv 0$), we will obtain a nontrivial global solution $u$. However, in this paper, we cannot prove such a strong result in general. Later we will apply part (ii) of Theorem \ref{mthm} to write $(u,u_t-u_r)$ as the sum of $(\eps r^{-1}U,\eps r^{-1}\mu U_q)$ and some remainder terms. To conclude that $u\not\equiv 0$, we need to assume that $\eps r^{-1}U$ or $\eps r^{-1}\mu U_q$ does not decay to zero faster than the remainder terms do as $s\to\infty$. 

Let $(\mu,U)$ be an admissible global solution to the geometric reduced system. For  $\eps\ll1$, let $u$ be the global solution constructed in Theorem~\ref{mthm}. Then, we have the following results.
\begin{enumerate}[(a)]
\item Suppose that there exist $(s^0,q^0,\omega^0)\in(1,\infty)\times\R\times\mathbb{S}^2$ such that \fm{|U(s,q^0,\omega^0)|+|\mu U_q(s,q^0,\omega^0)|\gtrsim \exp(-Cs),\qquad\forall s>s^0.} Then, for $\eps\ll1$, we have $u\not\equiv 0$.
\item Suppose that there exist $\kappa>1$ and  $(s^0,q^0,\omega^0)\in(1,\infty)\times\R\times\mathbb{S}^2$ such that \fm{|U(s,q,\omega^0)|+|\mu U_q(s,q,\omega^0)|\gtrsim \lra{q}^{-\kappa}\exp(-Cs),\qquad \forall s>s^0,\ q>q^0.}Then, for $\eps\ll1$, the initial data of $u$ cannot have compact support.
\end{enumerate}
Here the constants and the choice of $\eps$ depend on the admissible solution $(\mu,U)$. The accurate versions of these results can be found in Proposition \ref{mthmrmk3prop} in Section \ref{sec7.4pfprop1.1}. 

Note that an arbitrary admissible global solution $(\mu,U)$ does not necessarily satisfy the assumptions listed here in general and that these assumptions are not necessarily optimal. We choose them because they are simple and fit the examples studied in this paper well. 
We refer our readers to Section \ref{sec1.4} (especially Corollaries \ref{sec1:cor1.1:john} and \ref{sec1:cor1.3:euler}) and Sections \ref{sec8.2}--\ref{sec8.4} for discussions on how these results are applied to the examples studied in this paper.}

\rmkn{\rm\label{rmkthmunique} The solution in the main theorem is unique in the following sense. Suppose $N\geq 7$. Suppose $u_1,u_2$ are two $C^N$ solutions to \eqref{qwe} (given the same admissible solution $(\mu,U)$ and the same $\eps$) such that both of them satisfy the energy bounds and pointwise bounds in the main theorem. Then, we have $u_1=u_2$, assuming $\eps\ll 1$.  We refer our readers to Remark \ref{rmk712} and Section \ref{s7uni} for a sketch of the proof of this uniqueness result.

We also notice that when $(\mu,U)$ and $\eps$ are fixed, a solution $u$  satisfies all the bounds in this theorem if and only if it satisfies all the bounds (with possibly different constants) with the interval $[1/2,3/2]$ replaced by a smaller interval $[a,b]\subset[1/2,3/2]$ with   $1\in(a,b)$. The reason is that $Z^I(\eps r^{-1}U)\cdot 1_{r/t\in[1/2,3/2]\setminus[a,b]}$ can be controlled by the right sides of the estimates in the theorem; see Remark \ref{rmk7.1.3Uuapp}. As a result, the choice of the interval $[1/2,3/2]$ is not important in the main theorem. Following from this and the uniqueness in the previous paragraph, we conclude that, if $\wt{u}$ is another solution satisfying all the bounds in this theorem with the interval $[1/2,3/2]$ replaced by a smaller interval $[a,b]$, then we have $u=\wt{u}$ as long as $\eps\ll_{a,b}$1.}

\rmkn{\label{rmkmthm1.3}\rm In Theorem \ref{mthm}, we assume that there is no term of the form  $u\cdot u\cdot u$ in the Taylor expansions of the semilinear terms $f^{(*)}$ at the origin. See \eqref{qwefi}. As a result, we cannot apply Theorem \ref{mthm} to, e.g., Kadar's counterexample \eqref{kadarwnceq}. Fortunately, in John's counterexamples or the 3D compressible Euler equations with no vorticity, there are no terms of the form  $u\cdot u\cdot u$ in the semilinear terms, so Theorem \ref{mthm} is adequate. 
The reason why we cannot have terms of the form $u\cdot u\cdot u$ in this paper is explained after Lemma \ref{sec6.3l5}, and it is unclear to the author whether this assumption can be relaxed or not. Besides, though Kadar's work on \eqref{kadarwnceq} indicates that the undifferentiated cubic terms $u\cdot u\cdot u$ can lead to blowup in a forward Cauchy problem, it is unclear whether the same type of blowup will occur in a backward scattering problem. Thus, at least at this moment, we cannot say his results exclude the possibility to relax the assumptions on $u\cdot u\cdot u$ in this paper.
}

\rmkn{\label{rmkmthmdepeps}\rm In the statement of the main theorem, we fix an admissible global solution $(\mu,U)$ before we choose the value of $\eps$. In other words, we implicitly assume that $(\mu,U)$ is independent of $\eps$. However, our proofs in this paper rely not on the formulas of $(\mu,U)$ but on the bounds \eqref{def1.1a3}--\eqref{def1.1a8}. As long as $(\mu,U)$ satisfies these bounds where all the constants (such as the $C$ in \eqref{def1.1a3}) are independent of $\eps$, our proofs here will still work even if $(\mu,U)$ depends on $\eps$. If we compare this with the results proved in \cite{yu2021} where we studied the forward problem, we notice that this is indeed the case. That is, while the matching exact solution $(\mu,U)$ found in that paper depends on $\eps$, the bounds for $(\mu,U)$ there have constants independent of $\eps$.
}


\rmkn{\label{rmkinit}\rm There is an important question related to Theorem \ref{mthm}. Can we give an accurate description of the initial data of the global solutions constructed in this theorem?  If this is possible, then one may also attempt to show a global existence result for a forward Cauchy problem \eqref{qwe} with initial data
\fm{(u,u_t)|_{t=0}=(\eps u^0,\eps u^1),\qquad 0<\eps\ll1}where $u^0$ and $u^1$ satisfy certain assumptions arising from this accurate description. Such a description may be complicated in some cases. At the very least, we need to exclude the possibility that  $(u^0,u^1)\in C_c^\infty\setminus\{0\}$ when we study John's counterexamples.

Unfortunately, one can see from the estimates in Theorem \ref{mthm} that we do not have much information about what the initial data look like at $t=0$. The only estimates we have are
\eq{\label{rmkinitff}\norm{\partial Z^Iu(0)}_{L^2(\R^3)}\lesssim_I \eps,\qquad\forall I.} 
It is not even clear whether we have \fm{\norm{ Z^Iu(t)}_{L^2(\R^3)}<\infty,
\qquad t\geq 0.}
Even if we take our initial time as $t=1/\eps$, these estimates only tell us that  $u\approx\eps r^{-1}U$ when $|r-t|<t/2$. For $r\geq 3t/2$, we only have bounds like $|\partial Z^Iu|\lesssim \eps t^{-1/2+C\eps}r^{-3/2}$ and $|Z^Iu|\lesssim \eps t^{-1/2+C\eps}r^{-1/2}$ from part (ii). Clearly, we cannot exclude $C_c^\infty$ data by only assuming these bounds.

One possible way to give a description of the global solutions $u$ to \eqref{qwe} constructed in Theorem \ref{mthm} is to use the geometric reduced system. Recall that  $u$ matches a certain global solution $(\mu,U)$ to the geometric reduced system at infinite time. However, if we study a forward Cauchy problem \eqref{qwe} with data $(\eps u^0,\eps u^1)$, then to find the matching $(\mu,U)$, it seems that we must first solve \eqref{qwe} globally. This leads to a circulation. To resolve this issue, we make a compromise by using the Friedlander radiation field. Let $w$ be the global solution to the linear wave equation $\Box w=0$ with initial data $(u^0,u^1)$, and suppose that the radiation field $F_{0,0}=F_{0,0}(q,\omega)$ of $w_t$ exists. Since a solution to \eqref{qwe} is expected to  behave as a linear solution whenever $t\leq \exp(\delta/\eps)$, we expect the initial data $(\mu,U_q)|_{s=0}$ in the main theorem to be very close to $(-2,-F_{0,0})$. Thus, we may instead assume that the geometric reduced system \eqref{asypde2} with initial data\fm{(\mu,U_q)|_{s=0}=(-2,-F_{0,0})+\text{small perturbations}}  admits a global solution for all $s\geq 0$ with some pointwise estimates (such as \eqref{def1.1a3} -- \eqref{def1.1a8}). The author has to admit that such an assumption on the geometric reduced system is made naively, and currently there is no work  indicating that such a assumption guarantees the existence of a global solution with data $(\eps u^0,\eps u^1)$. However, on the positive side, for John's counterexample $\Box u=u_t^2$, this assumption does exclude the possibility that $(u^0,u^1)\in C_c^\infty$. In fact, one can see from the computations in Section \ref{sec8.2} that we  require $F_{0,0}\geq 0$ everywhere to guarantee a global solution to the geometric reduced system with data $(-2,-F_{0,0})$. However, if $(u^0,u^1)\in C_c^\infty$, then we must have $F_{0,0}<0$ at some point unless $(u^0,u^1)=0$.

The assumption above on the geometric reduced system can be viewed as a stability assumption: we consider a small perturbations of some fixed initial data. One can compare it with \cite{keir2} where Keir also added a  stability assumption on H\"ormander's asymptotic equations. Of course, his assumption is quite different from the assumption above. Instead of perturbing the data, he perturbed the asymptotic equations and considered all initial data.}

\rmkn{\label{rmkhyperboloid}\rm When we discuss a forward Cauchy problem \eqref{qwe} in Remark \ref{rmkinit}, we choose the initial data on a time slice $t=t_0$. Meanwhile, we notice that there have been some recent work (such as \cite{keir,keir2,dong2022generically}) where the authors studied the Cauchy problem \eqref{qwe} with initial data posed on a hyperboloid (e.g.\ $t^2-r^2=1$), or a combination of a time slice $r\leq r_0$, $t=t_0$ and a null cone $r-t=r_0-t_0$, $t\geq t_0$. We can thus ask similar questions: can we give an accurate description on the data of $u$ posed on a hyperboloid, or on a combination of a time slice and a null cone? Can we also prove a related global existence result under assumptions on radiation fields? While we can answer neither of them, it is interesting to notice that the estimates in Theorem \ref{mthm} fit such a setting of initial data quite well. In fact, since $\{t\geq 1/\eps, t^2-r^2=1\}\subset \{t\geq 1/\eps,r/t\in[1/2,3/2]\}$, we do have $u\approx \eps r^{-1}U$ on the hyperboloid for all large times.}

\rmkn{\label{rmkfutureglobal} \rm In Theorem \ref{mthm}, we construct global solutions, or more precisely, \emph{future global} solutions, to \eqref{qwe}. That is, they are solutions defined for all $t\geq 0$. It is natural to ask whether these future global solutions are also \emph{past global} in the sense that they are defined for all $t\leq 0$. Note that this question is closely related to the discussion in Remark \ref{rmkinit}. Suppose that we are given a future global solution $u$ to \eqref{qwe} constructed in Theorem \ref{mthm}. To determine whether $u$ is past global, we need to study a backward Cauchy problem\footnote{Such a backward Cauchy problem can be converted into a forward one if we define $w(t,x):=u(-t,x)$. Note that $w$ solves a (possibly different) system of quasilinear wave equations.} \eqref{qwe} with initial data given by $(u,u_t)|_{t=0}$ for $t\leq 0$. So far, the known global existence results (such as the ones in \cite{MR2382144,MR0837683,MR0820070}) and blowup results (such as the ones in \cite{MR0600571,MR0808321}) for \eqref{qwe} usually require the initial data to be sufficiently localized. However, as discussed in Remark \ref{rmkinit}, we do not have a precise description of the asymptotic behaviors of the initial data $(u,u_t)|_{t=0}$ as $|x|\to\infty$. As a result, we do not know whether $u$ is past global or not in general.

There, however, is a special case when the answer to the question above is known to be no. Later in this paper (Section \ref{sec8.2}), we will show that Theorem \ref{mthm} applies to $\Box u=u_t^2$, so we obtain a family of nontrivial future global solutions. It has been proved in \cite[Corollary 12]{bernhardt2024john} that these future global solutions are not past global. We will discuss this negative result in detail in Remark~\ref{rmkfutureglobal:john} in Section~\ref{sec1.4.1}.}

\rm\bigskip

Here we outline the main idea of the construction of $u$ in Theorem~\ref{mthm} which is very similar to that in \cite{MR4232783}. To construct a matching global solution, we solve a backward Cauchy problem with some initial data at $t=T$ and send $T$ to infinity. The initial data here are generated using the geometric reduced system. The same idea has appeared in, e.g.,\ \cite{MR4591580}, etc. Since our backward Cauchy problem is quasilinear, we follow the proof in \cite{MR2382144} to show the existence of the solution to this backward Cauchy problem. We use a continuity argument with the help of the adapted energy estimates and  Poincar$\acute{\rm e}$'s lemmas.

We now provide more detailed descriptions of the proof. We also ask our readers to check the beginning of each of Sections \ref{sec4}--\ref{sec7} for an overview of that specific section. First, we construct an approximate solution to \eqref{qwe}. This is done in Section \ref{sec4}. Let $q(t,x)$  and $U(t,x)$ be the approximate optical function and asymptotic profile associated to some admissible solution $(\mu,U)$. In Section \ref{sec4.2}, we prove several bounds for   $(q,U)(t,x)$. Most importantly, for $(t,x)$ close to the light cone, we show that $q$ is an approximate optical function in Proposition \ref{prop4.6} and that $\eps r^{-1}U$ is an approximate solution to \eqref{qwe} in Proposition \ref{prop4.8}. In Section \ref{sec4.3}, we define $u_{app}=(u_{app}^{(I)})$ by
\begin{equation}\label{uappdefn2}u_{app}(t,x)=\eps r^{-1}\psi(r/t)U(\eps\ln(t)-\delta,q(t,x),\omega)\end{equation} for all $t> 0$ and $x\in\R^3$. Here $\psi\equiv 1$ when $|r-t|\leq ct/2$, $\psi\equiv 0$ when $|r-t|\geq ct$, and $c\in(0,1/4)$ is a small parameter whose value will be chosen later. This $\psi$ is used to localize $\eps r^{-1}U$ near the light cone $\{r=t\}$. We can check that $u_{app}$ is a good approximate solution to \eqref{qwe} in the sense that \fm{g^{\alpha\beta}(u_{app})\partial_\alpha\partial_\beta u_{app}^{(I)}-f^{(I)}(u_{app},\partial u_{app})=O(\eps t^{-3+C\eps}),\qquad t\gg 1.}
In fact, whenever $|x|<t$,  a better estimate is expected; see Proposition \ref{mainprop4}.

In Section \ref{sec5}, we establish the energy estimates and Poincar$\acute{\rm e}$'s lemmas. These are the main tools of the proofs in the rest of the paper. Note that we do not choose the value of $c$ in the definition of $u_{app}$ until Section \ref{sec5.3}. Here we have an important remark. The estimates in  Section \ref{sec5} are very different from those in the author's previous work.  This is because in this paper our solutions are not supported outside some light cone $r-t=R$. We only know that they decay relatively fast inside the light cone $t=r$. To deal with this issue, we introduce a new weight function which is motivated by \cite{MR2134337,MR2680391}. This weight looks like a positive power of $\lra{r-t}$ when $|x|<t$ and looks like a constant when $|x|>t$. To overcome this difference, we need to assume that our solutions decay faster inside than outside. This explains why we have different decay assumptions on the admissible solution $(\mu,U)$ in \eqref{def1.1a3}--\eqref{def1.1a8}.

Next we seek to construct an exact solution matching $u_{app}$ at infinite time. Fixing a large time $T$, we consider the following equation 
\eq{\label{c3introeqn}
&g^{\alpha\beta}(u_{app}+v,\partial (u_{app}+v))\partial_\alpha\partial_\beta v^{(I)}\\&=-[ g^{\alpha\beta}(u_{app}+v,\partial (u_{app}+v))-g^{\alpha\beta}(u_{app},\partial u_{app})]\partial_\alpha\partial_\beta u_{app}^{(I)}\\
&\quad+[f^{(I)}(u_{app}+v,\partial (u_{app}+v))-f^{(I)}(u_{app},\partial u_{app})]\\
&\quad-\chi(t/T)[g^{\alpha\beta}(u_{app},\partial (u_{app}))\partial_\alpha\partial_\beta u_{app}^{(I)}-f^{(I)}(u_{app},\partial u_{app})]}
along with $v\equiv 0$ for $t\geq 2T$. Here $\chi\in C^\infty(\R)$ satisfies $\chi(t/T)=1$ for $t\leq T$ and $\chi(t/T)=0$ for $t\geq 2T$. Note that  $u_{app}+v$ is now an exact solution to \eqref{qwe} for $t\leq T$. In Section \ref{sec6} we prove that, for $\eps\ll1$, \eqref{c3introeqn} has a solution $v=v^T$ for all $t\geq 0$ that satisfies some decay in energy as $t\to\infty$; see Proposition \ref{prop6}. To prove this, we use a continuity argument. An overview of the proof of Proposition \ref{prop6} can be found right after its statement in Section \ref{sec6.1}.  The proof relies on the energy estimates and Poincar$\acute{\rm e}$'s lemmas, which are established in Section \ref{sec5}.

Finally we prove in Section \ref{sec7} that $v^T$ does converge to some $v^\infty$ in suitable function spaces, as $T\to\infty$. Thus we obtain a global solution $u_{app}+v^\infty$ to \eqref{qwe} for $t\geq 0$, such that it ``agrees with'' $u_{app}$ at infinite time, in the sense that the energy of $v^\infty$ tends to $0$ as $t\to\infty$. By the Klainerman-Sobolev inequality, we can derive the pointwise bounds in the main theorem from the estimates for the energy of $v^\infty$.

Note that to obtain a candidate for $v^\infty$, we have a more natural choice of PDE than \eqref{c3introeqn}. We may consider the Cauchy problem \eqref{qwe} for $t\leq T$ with initial data $(u_{app}(T),\partial_tu_{app}(T))$. The problem with such a choice is that for $u_{app}$ constructed above, $Z^I(u-u_{app})(T)$ does not seem to have a good decay in $T$ if  $Z^I$ only contains the scaling $S=t\partial_t+r\partial_r$ and Lorentz boosts $\Omega_{0i}=t\partial_i+x_i\partial_t$. This issue has been discussed in the author's previous paper; see the last paragraph of \cite[Section 1.4]{MR4232783}. It occurs even if we consider the simplest equation $\Box u=0$.  In the linear case, one possible way to deal with this difficulty is to consider more terms in the asymptotic expansion of the solutions, say take
\fm{u_{app}=\sum_{n=0}^N\frac{\eps}{r^{n+1}}F_n(r-t,\omega)}where $F_0$ is the usual Friedlander radiation field, and $F_n$ satisfies some PDE based on $F_{n-1}$.  However,  it does not seem to work in the quasilinear case, since we do not have such a good asymptotic expansion for a solution to  \eqref{qwe}. In this paper, we avoid such a difficulty by considering a variant \eqref{c3introeqn} of \eqref{qwe}. Such a difficulty does not appear in \eqref{c3introeqn}, since  $v\equiv 0$ for all $t\geq 2T$.

\subsection{Examples}\label{sec1.4}
We now show some applications of our main theorem. For a detailed discussion, we refer to Section \ref{sec8}.

\subsubsection{Global solutions to John's counterexamples}\label{sec1.4.1}

The first application of Theorem \ref{mthm} is the construction of a nontrivial global solution to
\eq{\label{exm141.1}\Box u=u_tu_{tt}\qquad\text{ in }\R^{1+3}.}
The asymptotic equation of H\"ormander's type associated to \eqref{exm141.1} is the inviscid Burgers' equation \eq{\label{exm141.2} 2\partial_sU_q+U_qU_{qq}=0.}
Since a nontrivial solution to  Burgers' equation with $C_c^\infty$ data is known to blow up in finite time, the equation \eqref{exm141.1} violates the weak null condition. We also remind our readers of the finite time blowup result for \eqref{exm141.1} by  John \cite{MR0600571}. One way to understand John's results is that nonzero $C_c^\infty$ data for \eqref{exm141.1} correspond
to nonzero $C_c^\infty$ data for \eqref{exm141.2}, which then leads to finite time blowup. Near the blowup time, it turns out the blowup mechanisms of both \eqref{exm141.1} and \eqref{exm141.2} are the same: a shock occurs in each case; see for example \cite{MR1466700,MR3561670,MR1867312}.

From these facts,  it seems that we can predict the long time dynamics of solutions to \eqref{exm141.1} by analyzing the long time dynamics of solutions to Burgers' equation. Thus, in order to construct a global solution to \eqref{exm141.1}, we can try starting with a global solution to \eqref{exm141.2}. Provided that $U_q|_{s=0}$ is $C^1$,  \eqref{exm141.2} admits a global solution if and only if $\partial_q^2U|_{s=0}\geq 0$ everywhere because of the method of characteristics. In general, we do not even need to assume that $U_q|_{s=0}$ is continuous. For example, we can take $U_q|_{s=0}=a_-\cdot 1_{q\leq 0}+a_+\cdot 1_{q>0}$ where $a_-<a_+$ are two fixed constants.  With this data,  Burgers' equation \eqref{exm141.2} admits a smooth global solution of a simple form. This solution is called a \emph{rarefaction wave solution}; we refer to  Example 2 in \cite[Section 3.4.2]{MR2597943}. However, there is a technical issue if we study \eqref{exm141.2} directly. In this paper, we prefer to work with smooth initial data for the asymptotic equations, but the only bounded initial data with $\partial_q^2U|_{s=0}\geq 0$ everywhere is a constant function, which is not general enough. The same issue occurs if we use the geometric reduced system. There is no nontrivial admissible global solution to the geometric reduced system. 

To resolve such a technical issue, we notice that if $u$ solves \eqref{exm141.1}, then  $\partial u=(\partial_\alpha u)_{\alpha=0,1,2,3}$ solves the system of quasilinear wave equations (here we set $v_{(\alpha)}:=\partial_\alpha u$)\footnote{Such a system of differentiated equations is also useful in the study of shock formation. For example, when Speck \cite{MR3561670} proved the shock formation for \eqref{shockqwe3d} with $F\equiv 0$, he converted the equation to a system of a system of covariant wave equations by differentiation. See \cite[Appendix A]{MR3561670}.}:
\eq{\label{exm141.qwe}\Box v_{(\sigma)}-v_{(0)}\Delta v_{(\sigma)}=\frac{1}{2}(\partial_\sigma v_{(0)}+\partial_t v_{(\sigma)})\cdot\partial_tv_{(0)},\qquad \sigma=0,1,2,3}
coupled to a system of constraint equations
\eq{\label{exm141.ceq}\partial_\alpha v_{(\beta)}=\partial_\beta v_{(\alpha)},\qquad \alpha,\beta=0,1,2,3.}
Meanwhile, if $v=(v_{(\alpha)})$ solves \eqref{exm141.qwe} and \eqref{exm141.ceq}, then we obtain a solution $u$ to \eqref{exm141.1} with $v_{(\alpha)}=\partial_\alpha u$. Thus, it suffices to construct a nontrivial global solution to these new equations. By easy computations, we notice that  the geometric reduced system associated to \eqref{exm141.qwe} is
\eq{\label{exm141.asy}
\left\{
\begin{array}{l}
\displaystyle \partial_s(\mu \partial_qU_{(\sigma)})=\frac{1}{8} (\wh{\omega}_\sigma \partial_q U_{(0)}- \partial_q U_{(\sigma)}) \partial_q U_{(0)}\mu^2,\quad \sigma=0,1,2,3;\\[1em]
\displaystyle \partial_s\mu=-\frac{1}{4}\mu^2 \partial_q U_{(0)}. 
\end{array}
\right.}Here the unknown $(\mu,U)=(\mu,(U_{(I)})_{I=0,\dots,3})$ is a function of $(s,q,\omega)$.
Inspired by the constraint equations \eqref{exm141.ceq}, we have an additional assumption
\eq{\label{exm141.asy2}\wh{\omega}_\alpha \partial_qU_{(\beta)}=\wh{\omega}_\beta \partial_qU_{(\alpha)},\qquad \alpha,\beta=0,1,2,3.}
Using  \eqref{exm141.asy2}, we can solve \eqref{exm141.asy} explicitly. In fact, if we set $(\mu,\partial_qU_{(0)})|_{s=0}=(-2,A)$ for some function $A=A(q,\omega)$, then we obtain a solution
\eq{\label{exm141.sol}\mu(s,q,\omega)&=\frac{4}{A(q,\omega)s-2},\qquad U_{(\alpha)}(s,q,\omega)=\int_{-\infty}^q-\wh{\omega}_\alpha A(p,\omega)\ dp.}
As long as  $0\leq -A(q,\omega)\leq C$,  $(\mu,U)$ is defined for all $s\geq -\delta_0$ where $0<\delta_0\ll_A1$ is a fixed small constant. 
If moreover we assume that $\partial_q^a\partial_\omega^cA=O(\lra{q}
^{-a-\gamma_{\sgn(q)}})$ for each $a,c\geq 0$ for some constants $\gamma_+>1,\gamma_->2$, then \eqref{exm141.sol} is indeed a $(\gamma_+,\gamma_-)$-admissible solution to the geometric reduced system; see Lemma \ref{lemexm3.1}. Thus, we can apply Theorem \ref{mthm} to construct a matching global solution $v=(v_{(\alpha)})_{\alpha=0,1,2,3}$.

Moreover, if we assume $A\not\equiv 0$, then the global solution $v$ constructed above is nonzero and its data is not compactly supported. In fact, since $A$ is nonpositive and continuous everywhere, and since $A(q^0,\omega^0)<0$ for some $(q^0,\omega^0)$, we have
\fm{U_{(0)}(s,q,\omega^0)=\int_{-\infty}^q A(p,\omega^0)\ dp\leq \int_{-\infty}^{q^0} A(p,\omega^0)\ dp<0,\qquad \forall s>0,\ q\geq q^0.}In other words, we have $|U(s,q,\omega^0)|\gtrsim_{(\mu,U)}1$. Now we apply part (b) in Remark~\ref{mthmrmk3}.

In summary, we have the following corollary.

\cor{\label{sec1:cor1.1:john}Fix $A\in C^\infty(\R\times\mathbb{S}^2)$, and fix two constants $\gamma_+>1,\gamma_->2$. Suppose that $A(q,\omega)\leq 0$ for all $(q,\omega)$ and that $\partial_q^a\partial_\omega^cA=O_{a,c}(\lra{q}^{-a-\gamma_{
\sgn(q)}})$ for each $a,c\geq 0$. For each $\eps\ll1$, there exists a global solution $v=(v_{(\alpha)})_{\alpha=0,1,2,3}$ for $t\geq 0$ to the system \eqref{exm141.qwe} matching $(\mu,U)$ defined by \eqref{exm141.sol} at infinite time of size $\eps$. Moreover, as long as $A\not\equiv 0$, the global solution $v$ is nontrivial, and its initial data are not compactly supported.}
\rmk{\label{rmk1.1.1}\rm The global solution $(\mu,(U_{(\alpha)}))$ constructed above induces a global solution to the geometric reduced system for \eqref{exm141.1}. In fact, if we set $W(s,q,\omega)$ by  \eq{\label{rmk1.1.1Wdefn}W(s,q,\omega)=\int_{-\infty}^q 2(\mu^{-1}U_{(0)})(s,p,\omega)\ dp} then this $(\mu,W)$ is a global solution to 
\fm{\left\{\begin{array}{l}
\displaystyle \partial_s(\mu W_q)=0,\\
\displaystyle \partial_s \mu=-\frac{1}{8}\mu^2\partial_q(\mu W_q).
\end{array}\right.}
Since $A=O(\lra{q}^{-2-})$ for $q<0$, we have $\mu^{-1}U_{(0)}=O(\lra{s}\lra{q}^{-1-})$ and thus $W=O(\lra{s}\lra{q}^{0-})$ is well defined. Later we will see in part (iii) of Remark \ref{rmkasyjohnquasi} that this $W$ offers an approximation for $u$ itself at least near the light cone.
However, this  $(\mu,W)$ is not an admissible solution because $\mu W_q=2U_{(0)}$ does not decay as $q\to\infty$ unless $A\equiv 0$. Thus, we cannot apply Theorem \ref{mthm} directly to \eqref{exm141.1}. This again explains why we study the differentiated equations \eqref{exm141.qwe}.}
\rmk{\rm \label{sec1:cor1.1:rmkene} As discussed in Remark \ref{rmkinit}, it is unclear whether we have $\norm{v(t)}_{L^2(\R^3)}<\infty$ or not. However, because of the explicit formulas for the solution $(\mu,U)$ in \eqref{exm141.sol}, we can apply a more quantitative version of part (b) of the results in Remark \ref{mthmrmk3} (see Proposition \ref{mthmrmk3prop}) to prove a lower bound for $\norm{v(t)}_{L^2(\R^3)}$: if $A\not\equiv 0$, we have
\fm{\norm{v(t)}_{L^2(\R^3)}\gtrsim \eps t^{1/2-},\qquad \forall t\geq e^{(1+\delta)/\eps}.}
The implicit constant here depends on $A$. See Remark \ref{corexm3:rmk}. Note that the left side of this estimate is not necessarily finite and that the lower bound is not necessarily optimal. Meanwhile, by applying the pointwise bounds in Section \ref{sec4}, we have
\fm{\norm{\partial(\eps r^{-1}U)(t)}_{L^2(\{x\in\R^3:|x|/t\in[1/2,3/2]\})}\lesssim\eps t^{-1+C\eps}\norm{\lra{r-t}^{-1}1_{r\sim t}}_{L^2(\R^3)}\lesssim \eps t^{C\eps},\qquad \forall t\geq 1/\eps.}
It follows from Theorem \ref{mthm} that
\fm{\norm{\partial v(t)}_{L^2(\R^3)}\lesssim \eps t^{C\eps},\qquad \forall t\geq 1/\eps.}
}
\rm

\bigskip

At this moment, we obtain a solution $v=(v_{(\alpha)})$ to \eqref{exm141.qwe}, but it is unknown whether there exists $u$ so that $v=\partial u$. We thus need to check the constraint equations \eqref{exm141.ceq}. The idea here is that if we set $w_{\alpha\beta}:=\partial_\alpha v_{(\beta)}-\partial_\beta v_{(\alpha)}$, then $w=(w_{\alpha\beta})$ is a solution to a certain system of general \emph{linear} wave equations with zero energy at infinite time. We then show that $w\equiv 0$ by applying the (backward) energy estimate. See Lemma \ref{lemexm3.2} -- Lemma \ref{lemexm3.4} for more details. In conclusion, we obtain the following result.

\cor{There exists a family of nontrivial global solutions to \eqref{exm141.1} for $t\geq 0$. Moreover, the initial data of these solutions are not compactly supported. }
\rmk{\label{rmkasyjohnquasi}\rm Since a solution $u$ constructed in this corollary behaves very differently from those constructed in the main theorem, we state the estimates that are satisfied by $u$.
\begin{enumerate}[{\normalfont (i)}]
\item The solution satisfies the energy bounds: for all $|I|\leq N-1$, we have
\fm{&\norm{w_0^{1/2}\partial Z^I(\partial_\alpha u-\eps r^{-1}U_{(\alpha)})(t)}_{L^2(\{x\in\R^3:\ |x|/t\in[1/2,3/2]\})}+\norm{w_0^{1/2}\partial^2 Z^Iu(t)}_{L^2(\{x\in\R^3:\ |x|/t\notin[1/2,3/2]\})}\\&\lesssim_I \eps t^{-1/2+C_I\eps},\qquad \forall t\geq 1/\eps;}
\fm{&\norm{w_0^{1/2}\partial^2 Z^Iu(t)}_{L^2(\R^3)}\lesssim_I \eps,\qquad \forall 0\leq t\leq 1/\eps.}
Here $w_0=1_{r\geq t}+\lra{r-t}^{\gamma_1}1_{r<t}$ for some constant $\gamma_1$ with $2<\gamma_1<\min\{4,2(\gamma_--1)\}$.
\item The solution satisfies the  pointwise bounds: for all $|I|\leq N-1$ and all $(t,x)$ with $t\geq 1/\eps$, \fm{
&|\partial Z^I(\partial_\alpha u-\eps r^{-1}U_{(\alpha)})(t,x)| 1_{|x|/t\in[1/2,3/2]}+|\partial Z^Iu(t,x)|1_{|x|/t\notin[1/2,3/2]}\\
&\lesssim_I \eps t^{-1/2+C_I\eps}\lra{t+r}^{-1}\lra{t-r}^{-1/2}w_0^{-1/2},}
\fm{&|Z^I(\partial_\alpha u-\eps r^{-1}U_{(\alpha)})(t,x)|1_{|x|/t\in[1/2,3/2]}+| Z^I\partial u(t,x)|1_{|x|/t\notin[1/2,3/2]}\\
&\lesssim_I \eps t^{-1/2+C_I\eps/2} (t^{-1}\lra{r-t}^{-(\gamma_1-1)/2}1_{r< t}+t^{-1}\lra{r-t}^{1/2}1_{t\leq r\leq 2t}+r^{-1/2}1_{r>2t}).}
Moreover, whenever $t\geq 1/\eps$ and $r/t\in[1/2,3/2]$,
\fm{|(\partial_t-\partial_r)\partial_\alpha u-\eps r^{-1}\mu \partial_q(U_{(\alpha)})|\lesssim \eps t^{-3/2+C\eps}\lra{t-r}^{-1/2}w_0^{-1/2}.}
\item The limit $u_{\infty,0}=\lim_{t\to\infty} u(t,0)$ exists. For all $t\geq 1/\eps$ and $x\in\R^3$, we have
\fm{&|u(t,x)-u_{\infty,0}-\frac{\eps}{|x|}W(t,x)1_{|x|/t\in[1/2,3/2]}|\lesssim
\left\{
\begin{array}{ll}
\eps t^{-\kappa+C\eps},& r\leq t;\\
\eps t^{-\kappa+C\eps}+\eps t^{-3/2+C\eps}\lra{r-t}^{3/2},&t< r\leq 2t;\\
\eps t^{-1/2+C\eps}r^{1/2},&r>2t.
\end{array}
\right.}
Here $\kappa>1$ is a constant depending only on $\gamma_1$, and we recall that $W$ is defined by \eqref{rmk1.1.1Wdefn} in Remark \ref{rmk1.1.1}. As a result, $\eps r^{-1}W+u_{\infty,0}$ provides a good approximation of $u$ whenever, e.g.,\ $r-t<t^{1/3-}$. One can also derive similar bounds for $Z^Iu$, but we omit them here for simplicity.
\end{enumerate}

Note that the estimates in (i) and (ii) are essentially the same as those in Theorem \ref{mthm}. In fact, to obtain the estimates here, we simply replace $u$ with $\partial u$ in Theorem \ref{mthm}. This is because we apply Theorem \ref{mthm} to the differentiated equations \eqref{exm141.qwe} and construct the solutions to \eqref{exm141.1} by using the condition \eqref{exm141.ceq}. In part (iii), we present pointwise bounds for $u$ itself obtained by integrating the pointwise bounds in (ii), so the bounds in part (iii) are worse than those in part (ii). We refer to Proposition \ref{propasymujohn} for a more accurate version of part (iii).

Meanwhile, as Remark \ref{rmkinit} suggests, it is not clear whether the $L^2$ norms of $\partial u$ and $u$ are finite or not. However, we refer to Remark \ref{sec1:cor1.1:rmkene} for a lower bound of the $L^2$ norm of $\partial u(t)$ when $A\not\equiv 0$.

Note that in part (ii) and part (iii), if we add a restriction that $t^2-r^2=1$ (or $t^2-r^2=1/\eps^2$ if we want to consider the whole hyperboloid), we obtain some estimates on the initial data of $u$ posed on the hyperboloid. In particular, we have $\partial u\approx\eps r^{-1}U$. Also see Remark \ref{rmkhyperboloid}.}

\rmk{\rm\label{rmk:ext:john:semilinear} The constructions above work not just for \eqref{exm141.1} but for a general quasilinear wave equation of the form
\fm{g^{\alpha\beta}(\partial u)\partial_\alpha\partial_\beta u=0.}
However, now we cannot say that the initial data of the global solutions constructed here is not compactly supported. We refer to Section \ref{sec8.3} (and Corollary \ref{corexm3}, in particular) for this general quasilinear case.

We  can also construct nontrivial global solutions  to another John's counterexample $\Box u=u_t^2$, or more generally
$\Box u=f(u,\partial u)$ where $f$ satisfies \eqref{qwefi}. The proof there is even simpler because we study \emph{not} the differentiated equations \emph{but} the equation $\Box u=u_t^2$ itself. We can use the geometric reduced system (which coincides with H\"ormander's asymptotic equations) associated to $\Box u=f(u,\partial u)$ directly.  See Section \ref{sec8.2} (and Corollary \ref{corexm2}, in particular). Similarly, the initial data of the solutions to $\Box u=u_t^2$ constructed here are not compactly supported, while we cannot make such a statement for the general equation $\Box u=f(u,\partial u)$.}

\rmk{\rm If we compare the constructions in the quasilinear case and the semilinear case, the following question arises naturally. Can we construct nontrivial global solutions to $\Box u=u_tu_{tt}$ by avoiding the differentiated equations and studying this scalar equation directly? More precisely, can we prove Theorem \ref{mthm} under a weaker notion of admissibility so that $(\mu,W)$ defined in Remark \ref{rmk1.1.1} is also admissible? It is unclear to the author how to achieve this goal at this moment because the assumptions \eqref{def1.1a3}--\eqref{def1.1a8} seem necessary in the proof in this paper. Moreover, Remark \ref{sec1:cor1.1:rmkene} suggests that the estimate for $\norm{\partial u(t)}_{L^2}$ should be worse than that for $\norm{\partial^2 u(t)}_{L^2}$, which brings difficulty in studying \eqref{exm141.1} directly via an energy method. Not to mention that we may even have $\norm{\partial u(t)}_{L^2}=\infty$. }

\rmk{\rm There is an interesting difference between the construction for $\Box u=u_t^2$ and that for $\Box u=u_tu_{tt}$ in this paper.  After we construct nontrivial global solutions to $\Box u=u_t^2$, we can conclude immediately that the initial data of these global solutions cannot have compact support by directly applying John's results in \cite{MR0600571}. In contrast, this is not the case for $\Box u=u_tu_{tt}$; we need to prove this result by ourselves instead of applying John's results directly. The reason is that, when John proved the finite time results for $\Box u=u_tu_{tt}$, he assumed that the initial data are $C_c^\infty$ and satisfy a certain integral inequality. In contrast, when he proved the results for $\Box u=u_t^2$, he only needed to assume that the data are $C_c^\infty$.
}

\rmk{\rm Note John's counterexamples admit other nonzero global solutions. For example, we notice that any polynomial of $(t,x)$ of degree at most $1$ (i.e.\ any linear combination of $t,x_1,x_2,x_3,1$ with real constant coefficients) is a global solution to $g^{\alpha\beta}(\partial u)\partial_\alpha\partial_\beta u=0$. A direct corollary of part (iii) in Remark \ref{rmkasyjohnquasi} is that the solutions constructed in this paper cannot be a nonzero polynomial of $(t,x)$ of degree at most $1$, because we have $|u(t,x)-\lim_{t\to\infty}u(t,0)|\lesssim \eps t^{-1-}$ whenever $t\geq 1/\eps$ and $|x|<t/2$. In other words, the asymptotics of the solutions constructed in this paper at the timelike infinity are very different from those of  nonzero polynomials of $(t,x)$ of degree at most $1$.

In the semilinear case, if we suppose that $u$ is independent of $x$, i.e.\ if $u(t,x)=u(t)$, then from the original wave equations we will get some ODE's. By solving them, we may also obtain nonzero global solutions to $\Box u=f(u,\partial u)$. In fact, by assuming $u(t,x)=u(t)$, we reduce $\Box u=u_t^2$ to $-u_{tt}=u_t^2$, and thus obtain a family of  global solutions $u_C:=\ln (t+C)$ with $C> 0$. We emphasize that the asymptotic behaviors of the global solutions $u$ constructed in this paper at the spacelike infinity are quite different from those of $u_C$. For example, by part (ii) of the main theorem, we have $|u(t,x)|\lesssim \eps t^{-1/2+C\eps}|x|^{-1/2}$ for $t\geq 1/\eps$ and $|x|>2t$, while $u_C(t,x)$ does not decay in $x$. In addition, we have $\norm{\partial  u(t)}_{L^2(\R^3)}<\infty$ while $\norm{\partial u_C(t)}_{L^2(\R^3)}=\infty$.}

\rmk{\rm\label{rmkfutureglobal:john}
In Remark \ref{rmkfutureglobal}, we mentioned that the solutions constructed in Theorem~\ref{mthm} are future global and that we do not know whether they are past global in general. However, there is one case when we know the answer is no. As discussed in Remark \ref{rmk:ext:john:semilinear}, we have constructed a large family of nontrivial (future) global solutions to John's counterexample $\Box u=u_t^2$ for all $t\geq 0$. Recently, it has been proved in \cite[Corollary 12]{bernhardt2024john} that such nontrivial future solutions are not past global. In other words, they must blow up in finite negative time. Such a blowup result is due to Bernhardt, the first author of \cite{bernhardt2024john}. In \cite[Theorem 11]{bernhardt2024john}, Bernhardt presented a novel blowup result for finite energy solutions to $\Box u=u_t^2$ satisfying a sign condition. This blowup result was then invoked to show that future solutions above blow up in the past. We refer to \cite[Section 4]{bernhardt2024john} for a detailed proof.

Note that the proof in \cite{bernhardt2024john} relies heavily on the structure of the semilinear wave equation $\Box u=u_t^2$. It is thus unclear to the author whether the blowup result in \cite{bernhardt2024john} can be extended to a different system \eqref{qwe} of quasilinear wave equations, such as the quasilinear example $\Box u=u_tu_{tt}$ and the 3D compressible Euler equations with no vorticity (which will be discussed in Section \ref{sec1.4.2}). In fact, it is even unknown whether the proof in \cite{bernhardt2024john} can be applied to a different semilinear equation $\Box u=f(u,\partial u)$ that violates the null condition.
}
\rm

\subsubsection{The 3D compressible Euler equations with no vorticity}\label{sec1.4.2}
Another application is the construction of a nontrivial global solution to the 3D compressible Euler equations with no vorticity: \eq{\label{exm142.1}\left\{
\begin{array}{l}
\rho_t+ \nabla\cdot(\rho v)=0,\\[1em]
\rho(v_t+v\cdot\nabla v)+\nabla p=0,\\[1em]
\partial_av^{b}=\partial_bv^a,\qquad a,b\in\{1,2,3\}.
\end{array}
\right.}
Here the unknowns are the velocity $v=(v^a)_{a=1,2,3}:\R^{1+3}\to\R^{3}$ and the density $\rho:\R^{1+3}\to (0,\infty)$. The pressure $p$ in the equations is a function of the density $p=p(\rho)$ where the map $\rho\mapsto p(\rho)$ is given. Here we assume that $dp/d\rho>0$ everywhere. 

In general,  we do not expect a global existence result for the 3D compressible Euler equations. For example, Sideris  \cite{MR0815196} proved that a $C^1$ solution with some well chosen initial data has a finite lifespan. There have also been several papers on the shock formation for the 3D compressible Euler equations with or without zero vorticity; see for example \cite{MR3288725,MR4011696,MR4109292,MR2284927,MR4612576}. 
Despite these negative results, there are grounds to believe that a nontrivial global solution might exist in a slightly different regime. In fact, under an arbitrary barotropic equation of state, \eqref{exm142.1} can be reduced to a scalar quasilinear wave equation 
\eq{\label{exm142tt}g^{\alpha\beta}(\partial\Phi)\partial_\alpha\partial_\beta \Phi=0\qquad\text{ in }\R^{1+3}}
where $\Phi$ is a potential function. While its associated H\"{o}rmander's asymptotic equation is Burgers' equation, we have successfully constructed nontrivial global solutions to \eqref{exm142tt} for all $t\geq 0$. We now show that we can also construct nontrivial global solutions to \eqref{exm142.1} by using the same method.

In our construction, we make use of a  system of covariant quasilinear wave equations derived in Luk-Speck \cite{MR4109292} and Speck \cite{MR4011696} for the 3D compressible Euler equations (with possibly nonzero vorticity). Their work was motivated by Christodoulou \cite{MR2284927} and Christodoulou-Miao \cite{MR3288725}, and one of their goals to introduce the new wave equation formulation is to prove the shock formation result for the 3D compressible Euler equations. It also turns out that the system they derived is also very useful in the study of  local wellposdeness; see \cite{MR4387234,MR4370374}.

We now formulate the system of covariant quasilinear wave equations used in this paper.
Since $\rho>0$ everywhere, we can set $\varrho:=\ln(\rho/\rho_0)$ where $\rho_0>0$ is a constant background density. Set the sound of speed $c_{\rm s}:=\sqrt{dp/d\rho}$ and view it as a function of $\varrho$, so $c_{\rm s}=c_{\rm s}(\varrho)$ and $c_{\rm s}'(\varrho)=\frac{d}{d\varrho}c_{\rm s}(\varrho)$. By setting $B=\partial_t+v^a\partial_a$,  we reduce \eqref{exm142.1} to 
\eq{\label{exm142.2}\left\{
\begin{array}{l}
B\varrho+ \partial_c v^c=0,\\[1em]
Bv+c_{\rm s}^2(\varrho)\nabla\varrho=0,\\[1em]
\partial_av^{b}=\partial_bv^a,\qquad a,b\in\{1,2,3\}.
\end{array}
\right.}
Now, we define the acoustical metric $g=(g_{\alpha\beta})$ by 
\eq{g:=-dt\otimes dt+c_{\rm s}^{-2}\sum_{a=1}^3 (dx^a-v^adt)\otimes (dx^a-v^adt).}
The inverse acoustical metric $g^{-1}=(g^{\alpha\beta})$ is defined by
\eq{g^{-1}:=-B\otimes B+c_{\rm s}^2\sum_{a=1}^3 \partial_a\otimes \partial_a.}
It is obvious that both $g$ and $g^{-1}$ are functions of $(v,\varrho)$. Without loss of generality, we also assume that $c_{\rm s}(0)=1$, so $g$ and $g^{-1}$ are both the Minkowski metrics when $(v,\varrho)=0$. With all these notations, it has been proved in \cite{MR4011696,MR4109292} that if $(v,\rho)$ solves \eqref{exm142.1}, then $(v,\varrho)$ solves the following system:
\eq{\label{exm142.qwe}\left\{
\begin{array}{l}
\displaystyle \Box_g v^a=-(1+c_{\rm s}^{-1}c_{\rm s}')g^{\alpha\beta}\partial_\alpha\varrho\partial_\beta v^a,\qquad a=1,2,3;\\[1em]
\displaystyle \Box_g \varrho=-3c_{\rm s}^{-1}c_{\rm s}'g^{\alpha\beta}\partial_\alpha\varrho\partial_\beta \varrho+2\sum_{1\leq a<b\leq 3}(\partial_av^a\partial_bv^b-\partial_av^b\partial_bv^a).
\end{array}
\right.}
For convenience, we  write $u^0=\varrho$ and $u^a=v^a$.
Here $\Box_g$ is a covariant wave operator defined by 
\eq{\Box_g\phi:=\frac{1}{\sqrt{|\det g|}}\partial_\alpha(\sqrt{|\det g|}g^{\alpha\beta}\partial_\beta \phi).}
After some computations, we obtain the following geometric reduced system\footnote{Note that the subscript $\rm s$ of $c_{\rm s}'(0)$ is different from and irrelevant to the variable $s$ in this equation.}
\eq{\label{exm142.asy}\left\{
\begin{array}{l}
\displaystyle \partial_s(\mu U_q^{a})=\frac{1}{2}(\omega_b U_q^b+c_{\rm s}'(0) U_q^0)\mu^2U_q^a,\qquad a=1,2,3;\\[1em]
\displaystyle \partial_s(\mu U_q^{0})=\frac{1}{2}(\omega_aU_q^a+c_{\rm s}'(0) U_q^0 )\mu^2U_q^0;\\[1em]
\displaystyle \partial_s \mu=\frac{1}{2}(\omega_a U^a_q+c_{\rm s}'(0) U_q^0)\mu^2.
\end{array}
\right.}
Here our unknown $(\mu,U)=(\mu,U^0,(U^a)_{a=1,2,3})$ is a function of $(s,q,\omega)$. Moreover, from \eqref{exm142.1},  we obtain  additional constraint equations for $U$:
\eq{\label{exm142.asy2}\left\{
\begin{array}{l}
\displaystyle U_q^0=\omega_cU_q^c;\\[1em]
\displaystyle U_q^a=\omega_aU_q^0,\qquad a=1,2,3;\\[1em]
\omega_aU_q^b=\omega_bU_q^a,\qquad a,b=1,2,3.
\end{array}
\right.}

We now solve the \eqref{exm142.asy} explicitly. Again, set $(\mu,U^0)|_{s=0}=(-2,A)$ for some smooth function $A=A(q,\omega)$. Then, we obtain a solution
\eq{\label{exm142.sol}\mu(s,q,\omega)=\frac{-2}{(
1+c_{\rm s}'(0))A(q,\omega)s+1},\qquad (U^0,U^a)(s,q,\omega)=\int_{-\infty}^q(A,\omega_aA)(p,\omega)\ dp.}
As long as $0\leq (
1+c_{\rm s}'(0))A(q,\omega)s\leq C$, $(\mu,U)$ is a defined for all $s\geq -\delta_0$ where $0<\delta_0\ll_A1$ is a fixed small constant. If moreover we assume that $\partial_q^a\partial_\omega^cA=O(\lra{q}
^{-a-\gamma_{\sgn(q)}})$ for each $a,c\geq 0$ for some constants $\gamma_+>1,\gamma_->2$, then \eqref{exm142.sol} is indeed a $(\gamma_+,\gamma_-)$-admissible solution to the geometric reduced system; see Lemma \ref{lemexm4.1}. Thus, we can apply Theorem~\ref{mthm}. Now, if $A\not\equiv 0$ and $c_{\rm s}'(0)\neq -1$, we can follow the proof in Section \ref{sec1.4.1} to show that the solution constructed here is nontrivial and its data are not compactly supported. In the case when $A\not\equiv 0$ but $c_{\rm s}'(0)=-1$, we notice that now $\mu U^0_q=-2A$, so we conclude that $u\not\equiv 0$ by applying part (a) in Remark \ref{mthmrmk3}, but we cannot make comments on the support of the initial data. In summary, we have the following corollary.

\cor{\label{sec1:cor1.3:euler}Fix $A\in C^\infty(\R\times\mathbb{S}^2)$, and fix two constants $\gamma_+>1,\gamma_->2$. Suppose that $(1+c_{\rm s}'(0))A(q,\omega)\geq 0$ for all $(q,\omega)$ and that $\partial_q^a\partial_\omega^cA=O_{a,c}(\lra{q}^{-a-\gamma_{
\sgn(q)}})$ for each $a,c\geq 0$. For each $\eps\ll1$, there exists a global solution $u=(u^{\alpha})_{\alpha=0,1,2,3}$ for $t\geq 0$ to the system \eqref{exm142.qwe}  matching $(\mu,U)$ defined by \eqref{exm142.sol} at infinite time of size $\eps$.

Moreover, as long as $A\not\equiv 0$, the solution $u$  is nontrivial. If furthermore we have $c_{\rm s}'(0)\neq -1$,  then the data of $u$ are not compactly supported.}
\rmk{\rm We have a remark similar to Remark \ref{sec1:cor1.1:rmkene}. As discussed in Remark \ref{rmkinit}, it is unclear whether we have $\norm{u(t)}_{L^2(\R^3)}<\infty$ or not. However, because of the explicit formulas for the solution $(\mu,U)$ in \eqref{exm142.sol}, we can apply Proposition \ref{mthmrmk3prop} to prove a lower bound for $\norm{u(t)}_{L^2(\R^3)}$: if $(1+c_{\rm s}'(0))A\not\equiv 0$, we have
\fm{\norm{u(t)}_{L^2(\R^3)}\gtrsim \eps t^{1/2-},\qquad \forall t\geq e^{(1+\delta)/\eps}.}
The implicit constant here depends on $A$. See Remark \ref{corexm4:rmk}. Again, the left side of this estimate is not necessarily finite and the lower bound is not necessarily optimal. Meanwhile, from Theorem \ref{mthm} we also have
\fm{\norm{\partial u(t)}_{L^2(\R^3)}\lesssim \eps t^{C\eps},\qquad \forall t\geq 1/\eps.}}
\rm

\bigskip

It remains to check that this $u=(\varrho,v)$ is indeed a solution to the original equations \eqref{exm142.2}. The idea of the proof is very similar to that in the previous example. We set
\fm{w_{00}&=Bu^0+\partial_cu^c,\qquad w_{0a}=-w_{a0}=c_{\rm s}^2\partial_au^0+Bu^a,\qquad w_{ab}=\partial_au^b-\partial_bu^a,}
and then show that this $w$ is a solution to a certain system of general linear wave equations with zero energy at infinite time. Again, we conclude $w\equiv 0$  by applying the (backward) energy estimate. We refer our readers to Section \ref{sec8.4} for the proof. Our final result is as follows.

\cor{ There exists a family of nontrivial global solutions to \eqref{exm142.2} for $t\geq 0$. If moreover we have $c_{\rm s}'(0)\neq -1$, then the initial data of all the solutions constructed here are not compactly supported.}
\rmk{\rm Unlike the previous example, here the global solutions in this corollary do satisfy the properties listed in Theorem \ref{mthm}. This is because here we apply the constructions directly to the system \eqref{exm142.qwe} of wave equations while we need to take antiderivatives when we study the previous example. However, again, we do not know whether the $L^2$ norm of $u$ itself is finite or not.}\rm
\bigskip 

\subsection{Acknowledgement}
The author would like to thank Herbert Koch, Sung-Jin Oh, and Daniel Tataru for many helpful discussions. The author would also like to
thank the anonymous reviewers for their valuable comments and suggestions on this paper. 
The author has been funded by the Deutsche Forschungsgemeinschaft (DFG, German Research Foundation)  through the Hausdorff Center for Mathematics under
Germany's Excellence Strategy - GZ 2047/1, Projekt-ID 390685813, by a Simons Investigator grant of Daniel Tataru from the Simons Foundation, and by a VandyGRAF Fellowship from Vanderbilt University.
\section{Preliminaries}
\subsection{Notation}\label{sec2.1}
In this paper, we always assume that $\eps\ll 1$ which means $0<\eps<\eps_0$ for some sufficiently small constant $\eps_0<1$. We write $\eps\ll_v 1$ or $\eps \ll1$ depending on $v$ if we want to emphasize that this constant $\eps_0$ depends on a parameter $v$. 

We use $C$ to denote universal positive constants. We write $A\lesssim B$, $B\gtrsim A$, or $A=O(B)$ if $|A|\leq CB$ for some  $C>1$. The values of all constants in this paper may vary from line to line. We write $A\sim B$ if $A\lesssim B$ and $B\lesssim A$. We use $C_{v}$, $\gtrsim_v$, or $\lesssim_v$ if we want to emphasize that the constant depends on a parameter $v$. Moreover, if we make a statement such as ``$|A|\leq CB$ for all $v$ in a certain set $V$'', then it implies that the constant $C$ can be chosen to be uniform for all $v\in V$. 

The convention above allows us to write a constant depending on $v$ as $C$ (i.e.\ without the subscript $v$). However, here we make a special convention that the statement above is not true for $\eps$. That is, in this paper, if a constant depends on the choice of $\eps$, then we must put $\eps$ in its subscript. In other words, any constant without the subscript $\eps$ must be independent of the choice of $\eps$. Similarly for $\lesssim_\eps$, $\gtrsim_\eps$, or $\ll_\eps$. 

Unless specified otherwise, we always assume that the Latin indices $i,j,l$ take values in $\{1,2,3\}$ and the Greek indices $\alpha,\beta$ take values in $\{0,1,2,3\}$.   We use subscript to denote partial derivatives unless specified otherwise. For example, $u_{\alpha\beta}=\partial_\alpha\partial_\beta u$, $q_r=\partial_rq=\sum_i \omega_i\partial_iq$, $A_q=\partial_qA$ and etc. For a fixed integer $k\geq 0$, we  use $\partial^k$ to denote either a specific partial derivative  of order $k$, or the collection of partial derivatives  of order $k$.

Given a function $f=f(\omega)$ defined on $\mathbb{S}^2$, in order to define angular derivatives $\partial_\omega$, we first extend $f$ to $\R^3\setminus 0$ by setting $f(x):=f(x/|x|)$ and then set $\partial_{\omega_i}f:=\partial_{x_i}f|_{\mathbb{S}^2}$. To prevent confusion, we will only use $\partial_\omega$ to denote the angular derivatives  under the coordinate $(s,q,\omega)$, and will never use it under the coordinate $(t,r,\omega)$. For a fixed integer $k\geq 0$, we will use $\partial_\omega^k$ to denote  either a specific angular derivative  of order $k$, or the collection of all angular derivatives of order $k$.

\subsection{Commuting vector fields}

Let $Z$ be any of the following vector fields:
\begin{equation}\label{c1vf} \partial_\alpha,\ \alpha=0,1,2,3;\ S=t\partial t+r\partial_r;\ \Omega_{ij}=x_i\partial_j-x_j\partial_i,\ 1\leq i<j\leq 3;\ \Omega_{0i}=x_i\partial_t+t\partial_i,\ i=1,2,3.\end{equation}We write these vector fields as $Z_1,Z_2,\dots,Z_{11}$, respectively. For any multiindex $I=(i_1,\dots,i_m)$ with length $m=|I|$ such that $1\leq i_*\leq 11$, we set $Z^I=Z_{i_1}Z_{i_2}\cdots Z_{i_m}$. Then we have Leibniz's rule
\begin{equation}Z^I(fg)=\sum_{|J|+|K|=|I|}C^I_{JK}Z^JfZ^Kg,\qquad\text{where $C_{JK}^I$ are constants.}\end{equation}

We have the following commutation properties.
\begin{equation}[S,\square]=-2\square,\qquad[Z,\square]=0\text{ for other $Z$};\end{equation}
\begin{equation}[Z_1,Z_2]=\sum_{|I|=1} C_{Z_1,Z_2,I}Z^I,\qquad\text{where $C_{Z_1,Z_2,I}$ are constants};\end{equation}
\begin{equation}\label{comf3}[Z,\partial_\alpha]=\sum_\beta C_{Z,\alpha\beta}\partial_\beta,\hspace{2em} \text{where $C_{Z,\alpha\beta}$ are constants}.\end{equation}

In this paper, we need the following lemma related to the commuting vector fields. For simplicity, we introduce the following notation. We use $f_0$ to denote an arbitrary  polynomial of $\{Z^I\omega\}$: in other words, $f_0$ is a finite sum of terms of the form\fm{C\cdot Z^{I_1}\omega_{i_1}\cdots Z^{I_p}\omega_{i_p}}
where $I_*$ are multiindices and $i_*\in\{1,2,3\}$. We allow $f_0$ to vary from line to line, so it is easy to check that $Z^If_0=f_0$ for each $I$ and that $f_0\cdot f_0=f_0$. We also remark that while the definition of $f_0$ will be modified in the rest of this paper, an arbitrary  polynomial of $\{Z^I\omega\}$ could always be denoted as $f_0$.

\lem{\label{c1lemtrcom}For each multiindex $I$ and each function $F$,  we have
\eq{\label{lemtrcomf1}(\partial_t-\partial_r)Z^IF&=Z^I(F_t-F_r)+\sum_{|J|<|I|}[f_0Z^J(F_t-F_r)+\sum_if_0(\partial_i+\omega_i\partial_t)Z^JF].}
Note that in $\sum_i(\dots)$,  the sum is taken over all $i=1,2,3$.}
\begin{proof}First, note that $[\partial_t-\partial_r,Z]=f_0\cdot \partial$ and $\partial=f_0(\partial_t-\partial_r)+\sum_i f_0(\partial_i+\omega_i\partial_t)$. We now prove \eqref{lemtrcomf1} by induction on $|I|$. If $|I|=0$, there is nothing to prove. Now suppose we have proved \eqref{lemtrcomf1} for each $|I|<n$. Now we fix a multiindex $I$ with $|I|=n>0$. Then, by writing $Z^I=ZZ^{I'}$, we have \fm{&\hspace{1.5em}(\partial_t-\partial_r)Z^IF=[\partial_t-\partial_r,Z]Z^{I'}F+Z((\partial_t-\partial_r)Z^{I'}F)\\
&=f_0\cdot\partial Z^{I'}F+Z\left(Z^{I'}(F_t-F_r)+\sum_{|J|<n-1}[f_0Z^J(F_t-F_r)+\sum_if_0(\partial_i+\omega_i\partial_t)Z^JF]\right)\\
&=f_0(f_0(\partial_t-\partial_r)+\sum_jf_0(\partial_j+\omega_j\partial_t))Z^{I'}F+Z^I(F_t-F_r)\\
&\hspace{1em}+\sum_{|J|<n-1}Z[f_0Z^J(F_t-F_r)+\sum_if_0(\partial_i+\omega_i\partial_t)Z^JF]\\
&=f_0(\partial_t-\partial_r)Z^{I'}F+\sum_jf_0(\partial_j+\omega_j\partial_t)Z^{I'}F+Z^I(F_t-F_r)\\
&\hspace{1em}+\sum_{|J|<n-1}[(Zf_0)Z^J(F_t-F_r)+\sum_i(Zf_0)(\partial_i+\omega_i\partial_t)Z^JF]\\&\hspace{1em}+\sum_{|J|<n-1}[f_0ZZ^J(F_t-F_r)+\sum_if_0Z(\partial_i+\omega_i\partial_t)Z^JF].}
In the second  equality, we can apply \eqref{lemtrcomf1} by the induction hypotheses. Moreover, we note that $[\partial_i+\omega_i\partial_t,Z]=f_0\cdot\partial$, so
\fm{Z(\partial_i+\omega_i\partial_t)Z^JF&=(\partial_i+\omega_i\partial_t)ZZ^JF+f_0\cdot\partial Z^JF\\
&=(\partial_i+\omega_i\partial_t)ZZ^JF+f_0(\partial_t-\partial_r) Z^JF+\sum_jf_0(\partial_j+\omega_j\partial_t)Z^JF.}
Now \eqref{lemtrcomf1} follows from the induction hypotheses and the computations above.
\end{proof}\rm

\subsection{Several pointwise bounds}
We have the pointwise estimates for partial derivatives.  
\begin{lem}\label{c1l2.1}
For any function $\phi$, we have 
\begin{equation}\label{c1l2.1f1}|\partial^k\phi|\leq C\lra{t-r}^{-k}\sum_{|I|\leq k}|Z^I\phi|,\qquad\forall k\geq 0,\end{equation}
and
\begin{equation}|(\partial_t+\partial_r)\phi|+|(\partial_i-\omega_i\partial_r)\phi|\leq C\lra{t+r}^{-1}|Z \phi|.\end{equation}
Here, for each $x\in\R$, we define the Japanese bracket $\lra{x}:=\sqrt{1+|x|^2}$. We also define $|Z\phi|:=\sum_{|I|=1}|Z^I\phi|$.
\end{lem}\rm\bigskip

In addition, we have the Klainerman-Sobolev inequality.  
\begin{prop} For $\phi\in C^\infty(\R^{1+3})$ which vanishes for large $|x|$, we have
\begin{equation}(1+t+|x|)(1+|t-|x||)^{1/2}|\phi(t,x)|\leq C\sum_{|I|\leq 2}\norm{Z^I\phi(t,\cdot)}_{L^2(\R^3)}.\end{equation}
\end{prop}\rm

We  have the following corollary which is essentially Proposition 14.1 in \cite{MR2680391}.
\cor{\label{corkl}
Fix a constant $\eta>0$. Fix a smooth function $w\in C^\infty(\R)$ such that $w\geq 1$, $w|_{[0,\infty)}\equiv 1$ and $w|_{(-\infty,-1]}= 2(1-q)^{\eta}$. Then, for any function $\phi\in C^\infty(\R^{1+3})$ vanishing for large $|x|$, whenever $t\geq 0$ and $x\in\R^3$, we have
\eq{\label{corklc}(1+t+|x|)[(1+|t-|x||)w(r-t)]^{1/2}|\phi(t,x)|\leq C_\eta\sum_{|I|\leq 2}\norm{w(|\cdot|-t)^{1/2}Z^I\phi(t,\cdot)}_{L^2(\R^3)}.}
The estimate \eqref{corklc} still holds if  $w$ is replaced by $1_{q\geq 0}+(1-q)^\eta 1_{q<0}$. 
}
\begin{proof}
If $t\leq 2$, then we have $r-t\geq -2$ so $w(r-t)\sim 1$ whenever $t\leq 2$. In this case, we can simply apply the classic Klainerman-Sobolev inequality. So, from now on we assume $t>2$.

We first prove \eqref{corklc} under the assumption that $\phi\equiv0$ whenever $r/t\leq 1/3$.  Note that 
\fm{\lra{r-t}|(w^{1/2})'(r-t)|+\lra{r-t}^2|(w^{1/2})''(r-t)|&\lesssim w^{1/2}.}
Also note that $Z(r-t)=O(\lra{r-t})$ and $Z(Z(r-t))=O(r^{-1}(1+r+t)\lra{r-t})=O(\lra{r-t})$ whenever $r/t\geq 1/3$ and $t>2$. We conclude that
\fm{\sum_{|I|\leq 2}|Z^I(w(r-t)^{1/2}\phi)|\lesssim w(r-t)^{1/2}\sum_{|I|\leq 2}|Z^I\phi|.}
Now \eqref{corklc} follows if we apply the classic Klainerman-Sobolev inequality to $w(r-t)^{1/2}\phi$.

Next we prove \eqref{corklc} under the assumption $\phi\equiv0$ whenever $r/t\geq 3t/5$. In this case, we have $w(r-t)\sim t^\eta$ in the support of $\phi$. Thus, \eqref{corklc} is equivalent to 
\fm{(1+t+|x|)[(1+|t-|x||)t^\eta]^{1/2}|\phi(t,x)|\lesssim \sum_{|I|\leq 2}\norm{t^{\eta/2}Z^I\phi(t,\cdot)}_{L^2(\R^3)}.}
Divide both sides by $t^{\eta/2}$ and we notice that this is the classic Klainerman-Sobolev inequality.

In general, we fix a cutoff function $\chi\in C_c^\infty(\R)$ such that $\chi|_{[-1/4,\infty)}\equiv 0$, $\chi|_{(-\infty,-1/2]}\equiv 1$ and $0\leq \chi\leq 1$. Set 
\fm{\phi_1:=\chi(\frac{r-t}{r+t})\phi,\qquad \phi_2:=\phi-\phi_1.}
For each multiindex $I$ with $|I|>0$, we can write $Z^I(\chi(\frac{r-t}{r+t}))$ as a linear combination of terms of the form
\fm{\chi^{(p)}(\frac{r-t}{r+t})\cdot \prod_{i=1}^pZ^{I_i}(\frac{r-t}{r+t}),\qquad \sum |I_*|=|I|,\ p\geq 1,\ |I_*|>0.}
Since $p>0$, we have $\chi^{(p)}(\frac{r-t}{r+t})=0$ unless $\frac{r-t}{r+t}\in[-1/2,-1/4]$, or equivalently $r/t\in [1/3,3/5]$.
If $r\sim t\geq 2$, we have $|Z^J(\frac{r-t}{r+t})|\lesssim \lra{r-t}(r+t)^{-1}\lesssim 1$ for each $J$. As a result, we have $Z^I(\chi(\frac{r-t}{r+t}))=O(1)$ whenever $t>2$. It follows that
\fm{\sum_{|I|\leq 2}|Z^I\phi_j|&\lesssim \sum_{|J_1|+|J_2|\leq 2}(1+|Z^{J_1}(\chi(\frac{r-t}{r+t}))|)|Z^{J_2}\phi|\lesssim \sum_{|I|\leq 2}|Z^I\phi|,\qquad j=1,2.}
Then, \eqref{corklc} follows from combining the two cases discussed above.
\end{proof}
\rm

\bigskip

We also state  Gronwall's inequality.
\prop{Suppose $A,E,r$ are bounded functions from $[a,b]$ to $[0,\infty)$. Suppose that $E$ is increasing. If \fm{A(t)\leq E(t)+\int_a^br(s)A(s)\ ds,\qquad \forall t\in[a,b],}then\fm{A(t)\leq E(t)\exp(\int_a^t r(s)\ ds),\qquad \forall t\in[a,b].}}\rm
The proofs of these results are standard. See, for example, \cite{MR2382144,MR2455195,MR1466700,MR4591580} for the proofs.

We also need the following lemma, which can be viewed as the estimates for Taylor's series adapted to $Z$ vector fields.

\lem{\label{l24}Fix an integer $k\geq 0$ and a multiindex $I$. Suppose there are two $\R^d$-valued functions $u,v$ on $(t,x)$ such that \eq{\label{l24aa1}\sum_{k'\leq k,\ |I'|\leq |I|\atop k'+|I'|\leq (k+ |I|)/2}(|\partial^{k'}Z^{I'}u|+|\partial^{k'}Z^{I'}v|)\leq 1,\qquad \forall(t,x).} Suppose for some integer $l_0>0$ we have $f\in C^\infty(\R^d)$ with $(\partial^lf)(0)=0$ for all $l\leq l_0$. Then, for all $(t,x)$, we have 
\begin{equation}\label{l242}\begin{aligned}&|\partial^kZ^I(f(u+v)-f(u))|\\&\lesssim_{f,k,I} \sum_{k'+k''\leq k\atop |I'|+|I''|\leq |I|}(|\partial^{k'}Z^{I'}v|(|\partial^{k''}Z^{I''}u|+|\partial^{k''}Z^{I''}v|))\cdot \left(\sum_{k'\leq k,\ |I'|\leq |I|\atop k'+|I'|\leq (k+|I|)/2}(|\partial^{k''}Z^{I''}u|+|\partial^{k''}Z^{I''}v|)\right)^{l_0-1}.\end{aligned}\end{equation}
}

\begin{proof} By the chain rule  Leibniz's rule, $\partial^kZ^I(f(u))$ can be written as a sum of terms of the form
\fm{(\partial^lf)(u) \prod_{j=1}^l\partial^{k_j}Z^{I_j}u^{i_j}}where $l\leq k+|I|$, $k_i+|I_i|>0$ for each $i$ and $\sum_ik_i= k$, $\sum_iI_i=I$.  Thus, $\partial^kZ^I(f(u+v)-f(u))$ can be written as a sum of terms of the form
\eq{\label{l24fff0}&(\partial^lf)(u+v) \prod_{j=1}^{l}\partial^{k_j}Z^{I_j}(u^{i_j}+v^{i_j}) -(\partial^lf)(u) \prod_{j=1}^{l}\partial^{k_j}Z^{I_j}u^{i_j}
\\&=((\partial^lf)(u+v)-(\partial^lf)(u))\prod_{j=1}^{l}\partial^{k_j}Z^{I_j}(u^{i_j}+v^{i_j})\\&\quad+ (\partial^lf)(u)\sum_{p=1}^l[\prod_{j=1}^{p-1}\partial^{k_j}Z^{I_j}u^{i_j}\cdot\partial^{k_p}Z^{I_p}v^{i_p}\cdot \prod_{j=p+1}^{l}\partial^{k_{j}}Z^{I_{j}}(u^{i_j}+v^{i_j}) ]}
where $k_i+|I_i|>0$ for each $i$ and $\sum_ik_i=k$, $\sum_iI_i=I$. In particular, if the term \eqref{l24fff0} with $l=0$ appears, we must have $k+|I|=0$; if the term \eqref{l24fff0} with $l=1$ appears, then \eqref{l24fff0} must be of the form $((\partial f)(u+v)-(\partial f)(u))\partial^kZ^Iu+(\partial f)(u)\partial^kZ^Iv$.

Moreover, we claim that for each $l\leq l_0$,
\eq{\label{l24fff1}|(\partial^lf)(u+v)-(\partial^lf)(u)|\lesssim_{l,l_0} (|u|+|v|)^{l_0-l}|v|,\qquad \forall u,v\in B(0,1).}
Here we set $B(0,1)=\{z\in\R^d:\ |z|\leq 1\}$. In fact, we have
\fm{|(\partial^{l_0}f)(u+v)-(\partial^{l_0}f)(u)|\leq \sup_{\beta\in[0,1]}|(\partial^{l_0+1}f)(u+\beta v)||v|\leq \norm{\partial^{l_0+1}f}_{L^\infty(B(0,1))}|v|.}
So \eqref{l24fff1} holds for $l=l_0$. Next, if we have proved \eqref{l24fff1} for some $1\leq l\leq l_0$, then we have
\fm{|(\partial^{l-1}f)(u+v)-(\partial^{l-1}f)(u)|\leq \sup_{\beta\in[0,1]}|(\partial^{l}f)(u+\beta v)||v|.}
By replacing $(u,v)$ with $(0,u+\beta v)$ in \eqref{l24fff1}, we have
\fm{\sup_{\beta\in[0,1]}|(\partial^{l}f)(u+\beta v)|\lesssim \sup_{\beta\in[0,1]}|u+\beta v|^{l_0-l+1}\lesssim (|u|+|v|)^{l_0-l+1}.}
This implies \eqref{l24fff1} with $l$ replaced by $l-1$. 

Now let us prove \eqref{l242}. If $l=0$ in \eqref{l24fff0}, we must have $k=|I|=0$, so \eqref{l242} follows directly from \eqref{l24fff1}. If $l\geq 1$  in \eqref{l24fff0}, then there exists at most one index $j$ such that $k_j+|I_j|>(k+|I|)/2$; for all other terms, we can apply our pointwise assumptions \eqref{l24aa1}.
Thus, if $1\leq l\leq l_0$ in \eqref{l24fff0},  by \eqref{l24fff1} we have
\fm{&|((\partial^lf)(u+v)-(\partial^lf)(u))\prod_{j=1}^{l}\partial^{k_j}Z^{I_j}(u^{i_j}+v^{i_j})|\\&\lesssim (|u|+|v|)^{l_0-l}|v|\cdot 
\left(\sum_{k'\leq k,\ |I'|\leq|I|\atop k'+|I'|\leq (k+|I|)/2 } (|\partial^{k'}Z^{I'}u|+|\partial^{k'}Z^{I'}v|)\right)^{l-1}\cdot\sum_{k'\leq k \atop |I'|\leq|I| } (|\partial^{k'}Z^{I'}u|+|\partial^{k'}Z^{I'}v|)\\
&\lesssim \left(\sum_{k'\leq k,\ |I'|\leq|I|\atop k'+|I'|\leq (k+|I|)/2 } (|\partial^{k'}Z^{I'}u|+|\partial^{k'}Z^{I'}v|)\right)^{l_0-1}\cdot\sum_{k'\leq k \atop |I'|\leq|I| } |v|(|\partial^{k'}Z^{I'}u|+|\partial^{k'}Z^{I'}v|).}
If $l>l_0$, then
\fm{&|((\partial^lf)(u+v)-(\partial^lf)(u))\prod_{j=1}^{l}\partial^{k_j}Z^{I_j}(u^{i_j}+v^{i_j})|\\&\lesssim |v|\cdot \left(\sum_{k'\leq k,\ |I'|\leq|I|\atop k'+|I'|\leq (k+|I|)/2 } (|\partial^{k'}Z^{I'}u|+|\partial^{k'}Z^{I'}v|)\right)^{l-1}\cdot \sum_{k'\leq k\atop |I'|\leq |I|} (|\partial^{k'}Z^{I'}u|+|\partial^{k'}Z^{I'}v|)\\
&\lesssim \left(\sum_{k'\leq k,\ |I'|\leq|I|\atop k'+|I'|\leq (k+|I|)/2 } (|\partial^{k'}Z^{I'}u|+|\partial^{k'}Z^{I'}v|)\right)^{l_0-1}\cdot \sum_{k'\leq k\atop |I'|\leq |I|} |v|(|\partial^{k'}Z^{I'}u|+|\partial^{k'}Z^{I'}v|).}

Moreover, if $l=1$, then the second term in \eqref{l24fff0} is of the form $(\partial f)(u)\cdot \partial^k Z^I v$. By replacing $(u,v)$ with $(0,u)$ in \eqref{l24fff1}, we have
\fm{|(\partial f)(u)\cdot \partial^k Z^I v|&\lesssim |u|^{l_0}| \partial^k Z^I v|.} If $2\leq l\leq l_0$, by replacing $(u,v)$ with $(0,u)$ in \eqref{l24fff1}, we have 
\fm{&|(\partial^lf)(u)[\prod_{j=1}^{p-1}\partial^{k_j}Z^{I_j}u^{i_j}\cdot\partial^{k_p}Z^{I_p}v^{i_p}\cdot \prod_{j=p+1}^{l}\partial^{k_{j}}Z^{I_{j}}(u^{i_j}+v^{i_j}) |\\
&\lesssim|u|^{l_0-l+1} \sum_{k'+k''\leq k\atop |I'|+|I''|\leq |I|}(|\partial^{k'}Z^{I'}v|(|\partial^{k''}Z^{I''}u|+|\partial^{k''}Z^{I''}v|))\cdot \left(\sum_{k'\leq k,\ |I'|\leq |I|\atop k'+|I'|\leq (k+|I|)/2}(|\partial^{k''}Z^{I''}u|+|\partial^{k''}Z^{I''}v|)\right)^{l-2}\\
&\lesssim \sum_{k'+k''\leq k\atop |I'|+|I''|\leq |I|}(|\partial^{k'}Z^{I'}v|(|\partial^{k''}Z^{I''}u|+|\partial^{k''}Z^{I''}v|))\cdot \left(\sum_{k'\leq k,\ |I'|\leq |I|\atop k'+|I'|\leq (k+|I|)/2}(|\partial^{k''}Z^{I''}u|+|\partial^{k''}Z^{I''}v|)\right)^{l_0-1}.}
If $l>l_0$, we have $l-2\geq l_0-1$ and thus
\fm{&|(\partial^lf)(u)[\prod_{j=1}^{p-1}\partial^{k_j}Z^{I_j}u^{i_j}\cdot\partial^{k_p}Z^{I_p}v^{i_p}\cdot \prod_{j=p+1}^{l}\partial^{k_{j}}Z^{I_{j}}(u^{i_j}+v^{i_j}) |\\
&\lesssim \sum_{k'+k''\leq k\atop |I'|+|I''|\leq |I|}(|\partial^{k'}Z^{I'}v|(|\partial^{k''}Z^{I''}u|+|\partial^{k''}Z^{I''}v|))\cdot \left(\sum_{k'\leq k,\ |I'|\leq |I|\atop k'+|I'|\leq (k+|I|)/2}(|\partial^{k''}Z^{I''}u|+|\partial^{k''}Z^{I''}v|)\right)^{l-2}\\
&\lesssim \sum_{k'+k''\leq k\atop |I'|+|I''|\leq |I|}(|\partial^{k'}Z^{I'}v|(|\partial^{k''}Z^{I''}u|+|\partial^{k''}Z^{I''}v|))\cdot \left(\sum_{k'\leq k,\ |I'|\leq |I|\atop k'+|I'|\leq (k+|I|)/2}(|\partial^{k''}Z^{I''}u|+|\partial^{k''}Z^{I''}v|)\right)^{l_0-1}.}

We finish our proof by combining all the estimates above.
\end{proof}

\rm

\subsection{A function space}\label{c1sec1.6.4}
Fix an open set $\mcl{D}\subset\R_{t,x}^{1+3}$ which may depend on the parameter $\eps$. Suppose that in $\mcl{D}$ we have $t\geq 2C$ and $r/t\in[1/C,C]$ for some constant $C>1$ which is independent of $\eps$. We make the following definition.

\defn{\label{c1defn1.5}\rm Fix $n,s,p\in\R$.  We say that a function $F=F_\eps(t,x)$ is in $\eps^nS^{s,p}=\eps^nS^{s,p}_\D$, if for each fixed integer $N\geq 1$, for all $\eps\ll_{n,s,p,N,\D}1$, we have $F\in C^N(\D)$ and 
\eq{\label{c1def1.5}\sum_{|I|\leq N}|Z^IF(t,x)|\lesssim_{N,\D} \eps^nt^{s+C_{N,\D}\eps}\lra{r-t}^{p},\hspace{2em}\forall (t,x)\in\D.}
Here $F$ is allowed to depend on $\eps$, but all the constants in \eqref{c1def1.5} must be independent of $\eps$.

If $n=0$, we write $\eps^0S^{s,p}$ as $S^{s,p}$ for simplicity.
}\rm\bigskip

We have the following  key lemma.

\lem{\label{c1lem6} We have the following two properties.

{\rm (a)} For any $F_1\in\eps^{n_1} S^{s_1,p_1}$ and $F_2\in \eps^{n_2}S^{s_2,p_2}$, we have \fm{F_1+ F_2\in \eps^{\min\{n_1,n_2\}}S^{\max\{s_1,s_2\},\max\{p_1,p_2\}},\hspace{2em}F_1F_2\in\eps^{n_1+n_2} S^{s_1+s_2,p_1+p_2}.}

{\rm (b)} For any $F\in\eps^n S^{s,p}$,  we have $ZF\in \eps^nS^{s,p}$, $\partial F\in \eps^nS^{s,p-1}$ and $(\partial_i+\omega_i\partial_t)F\in \eps^nS^{s-1,p}$.}
\begin{proof} Note that (a) follows directly from the definition and  Leibniz's rule. In (b), if $F\in \eps^nS^{s,p}$, then  $ZF\in \eps^nS^{s,p}$ follows directly from the definition. Next, we fix an arbitrary integer $N\geq 1$. Since $F\in \eps^nS^{s,p}$, for all $\eps\ll_{n,s,p,N+1}1$ we have $F\in C^{N+1}(\mcl{D})$ and \fm{\sum_{|I|\leq N+1}|Z^IF(t,x)|\lesssim \eps^nt^{s+C\eps}\lra{r-t}^{p},\hspace{2em}\forall (t,x)\in\D.}
Thus, $\partial F\in C^N(\mcl{D})$. Moreover, by \eqref{comf3} and Lemma \ref{c1l2.1}, in $\mcl{D}$ we have
\fm{\sum_{|I|\leq N}|Z^I\partial F|&\lesssim \sum_{|I|\leq N}|\partial Z^I F|\lesssim \sum_{|J|\leq N+1}\lra{r-t}^{-1}|Z^J F|\lesssim \eps^nt^{s+C\eps}\lra{r-t}^{p-1}.}
In conclusion, $\partial F\in \eps^nS^{s,p-1}$.

Next, we note that
\fm{(\partial_i+\omega_i\partial_t)F&=r^{-1}\Omega_{0i}F+(r+t)^{-1}r^{-1}\omega_i(rSF-\sum_j t\omega_j\Omega_{0j}F)+(r-t)r^{-2}\sum_j\omega_j\Omega_{ji}F.}
By the definition, we can easily show  $t^m, r^{m},(r+t)^{m}\in S^{m,0}$ for each $m\in\R$,  $r-t\in S^{0,1}$ and  $\partial^m\omega_i\in S^{-m,0}$ for each integer $m\geq 0$. And since $ZF\in S^{s,p}$, by part (a) we conclude that
\fm{(\partial_i+\omega_i\partial_t)F\in \eps^nS^{s-1,p}+\eps^nS^{s-2,p-1}=\eps^nS^{s-1,p}.}
Here we have $\eps^nS^{s-2,p-1}\subset \eps^nS^{s-1,p}$. In fact, recall that $t\geq 2C$ and $C^{-1}\leq r/t\leq C$ for some constant $C$ independent of $\eps$. Thus, in $\D$ we have
\fm{\lra{r-t}/t= \sqrt{t^{-2}+(r/t-1)^2}\leq \sqrt{1/(4C^2)+C^2}\lesssim 1.}
In summary, in $\D$ we have \fm{\eps^nt^{s-2+C\eps}\lra{r-t}^{p-1}\lesssim\eps^nt^{s-1+C\eps}\lra{r-t}^{p},} so $S^{s-2,p-1}\subset S^{s-1,p}$. This finishes the proof.
\end{proof}
\exm{\rm\label{c1exm1.7} We have \fm{t^m, r^{m},(r+t)^{m}\in S^{m,0},\ \forall m\in\R;\hspace{2em}  r-t\in S^{0,1};\hspace{2em} \partial^m\omega_i\in S^{-m,0}\ \forall m\geq 0,\ m\in\Z.}
It also follows from $\frac{d}{ds}\lra{s}=s/\lra{s}$, the chain rule and Lemma \ref{c1lem6} that
\fm{(r-t)^m,\lra{r-t}^m\in S^{0,m},\ \forall m\in\R.}}\rm\bigskip

In addition, we have the following lemma which is relevant to the Taylor expansion of a function.
\lem{\label{lem2.9} Suppose that $t\sim r$ and $t\gg 1$ in $\D$. Suppose $f\in C^\infty(\R^d)$ and let $u=(u^j)\in \eps^n S^{s,p}$ for some $n>0$, $s<0$ and $p\leq 0$. Suppose that we have the Taylor expansion of $f$ at $0$:
\fm{f(X)=a_0+\sum a_iX^i+\sum a_{ij}X^iX^j+O(|X|^4),\qquad X\in\R^d.} Then, we have \fm{f(u)-a_0-\sum a_iu^i-\sum a_{ij}u^iu^j\in \eps^{3n}S^{3s,3p}.}
As a result, we have
\fm{f(u)-a_0\in \eps^nS^{s,p},\quad f(u)-a_0-\sum a_iu^i\in \eps^{2n}S^{s,p}.}}
\begin{proof}
It suffices to prove the result for $f$ with $a_i=a_{ij}=0$ for all $1\leq i,j\leq d$. From now on we assume that $a_i=a_{ij}=0$ for all $1\leq i,j\leq d$.

Since $u\in \eps^nS^{s,p}$, $t\gtrsim1$,  $s<0$ and $p\leq 0$, by choosing $\eps\ll_{n,s} 1$, we have \fm{|u|\leq C\eps^n t^{s+C\eps}\lra{r-t}^p\leq C\eps^n\leq 1.}
In this estimate, we can choose $\eps\ll_{n,m,s,p}1$ so that $s+C\eps<0$. Then, by  Taylor's theorem, for each $l=0,1,2$ we have
\eq{\label{lem2.9ff}|(\partial^lf)(u)|& \lesssim\norm{\partial^3f}_{L^\infty(B(0,1))} |u|^{3-l}\lesssim \eps^{(3-l)n} t^{(3-l)s+C\eps}\lra{r-t}^{(3-l)p}.}
Here $B(0,1):=\{X\in\R^d:\ |X|\leq 1\}$. If $l\geq 3$, then we have \fm{|(\partial^lf)(u)|\lesssim 1.}

We now prove $Z^I(f(u))=O(\eps^{3n} t^{3s+C\eps}\lra{r-t}^{3p})$ for each $|I|$. The case $|I|=0$ follows from \eqref{lem2.9ff} with $l=0$. In general, we fix a multiindex $I$ with $|I|=m>0$.  By  Leibinz's rule and the chain rule, we can write $Z^I(f(u))$ as  a linear combination of terms of the form
\fm{(\partial^lf)(u)\cdot\prod_{j=1}^lZ^{I_j}u^{k_j},\hspace{2em}\text{ where }1\leq l\leq m,\ \sum |I_j|=m,\ |I_j|>0\text{ for each }j.}
We choose $\eps\ll_{n,m,s,p}1$ so that $Z^Ju$ exists and satisfies the estimate \eqref{c1def1.5} with $N$ replaced by $m$ and $F$ replaced by $m$.
If $l=1,2$, each of these terms  is controlled by \fm{\eps^{(3-l)n}t^{(3-l)s+C\eps}\lra{r-t}^{(3-l)p}\cdot (\eps^nt^{s+C\eps}\lra{p})^l\lesssim\eps^{3n}t^{3s+C\eps}\lra{r-t}^{3p}.}If $l\geq 3$, each of these terms  is controlled by \fm{1\cdot (\eps^nt^{s+C\eps}\lra{r-t}^{p})^l\lesssim\eps^{3n}t^{3s+C\eps}\lra{r-t}^{3p}.}
As a result, we conclude that $f(u)\in\eps^{3n}S^{3s,3p}$.
\end{proof}\rm

\section{Global solutions to the geometric reduced system}

\subsection{The geometric reduced system}
In this paper, we study a general system of quasilinear wave equations \eqref{qwe}:
\fm{g^{\alpha\beta}(u,\partial u)\partial_\alpha\partial_\beta u^{(I)}=f^{(I)}(u,\partial u),\qquad I=1,\dots,M.}
Here our unknown function is  vector-valued. That is, we have $u=(u^{(I)}):\R^{1+3}_{t,x}\to\R^M$ for some positive integer $M$. In addition, we assume that  $(g^{\alpha\beta})$ are smooth, symmetric and independent of the index $I$ and that $g^{\alpha\beta}(0,0)=m^{\alpha\beta}$. Moreover, we assume that the $f^{(I)}$ are all smooth functions such that  $f^{(I)}(0,0)=0$ and $d f^{(I)}(0,0)=0$.

Assume that we have the  Taylor expansions
\eq{\label{sec3tay}g^{\alpha\beta}(u,\partial u)&=m^{\alpha\beta}+g_J^{\alpha\beta}u^{(J)}+g_J^{\alpha\beta\lambda}\partial_\lambda u^{(J)}+O(|u|^2+|\partial u|^2),\\
f^{(I)}(u,\partial u)&=f^{I,\alpha\beta}_{JK}\partial_\alpha u^{(J)}\partial_\beta u^{(K)}+O(|u|^4+|\partial u|(|\partial u|^2+|u|^2)).}
Here  $m^{\alpha\beta},g^{*}_*,f^{*}_{*}$ are all real constants, and we use the Einstein summation convention. In particular, note that there is no term of order $u^3$ in the expansion of the $f^I$'s. The reason for such a setting is explained in Remark \ref{rmkmthm1.3}.

For each $\wh{\omega}=(-1,\omega)\in\R\times\mathbb{S}^2$ and each $1\leq I,J,K\leq M$, we set
\eq{G_{2,J}(\omega):=g^{\alpha\beta}_J\widehat{\omega}_\alpha\widehat{\omega}_\beta,\qquad G_{3,J}(\omega):=g^{\alpha\beta\lambda}_J\widehat{\omega}_\alpha\widehat{\omega}_\beta\widehat{\omega}_\lambda,\qquad F_{2,JK}^I(\omega):=f^{I,\alpha\beta}_{JK}\widehat{\omega}_\alpha\widehat{\omega}_\beta.}
Then, the geomteric reduced system for $(\mu,(U^{(I)})_{I=1}^M)(s,q,\omega)$ is defined by
\eq{\label{asy}\left\{
\begin{array}{l}\displaystyle
\partial_s(\mu U_q^{(I)})=-\frac{1}{4}F^{I}_{2,JK}(\omega) \mu^2U_q^{(J)} U_q^{(K)},\qquad I=1,\dots,M;\\[1em]\displaystyle
\partial_s\mu=\frac{1}{4}G_{2,J}(\omega) \mu^2U_q^{(J)}-\frac{1}{8}G_{3,J}(\omega)\mu^2\partial_q(\mu U_{q}^{(J)}).
\end{array}\right.}
In this system we use the Einstein summation convention with respect to $J,K=1,\dots,M$. For a derivation of this reduced system, we refer our readers to   \cite[Section 2.2]{MR4315017} and \cite[Section 3.2]{MR4232783}.  The solvability of the geometric reduced system has been discussed in Section \ref{s1introgrs}.

\subsection{Assumptions}\label{sec3.2assumptions}
As explained in the introduction, it is not a realistic goal to solve a general geometric reduced system  explicitly. Instead, in this paper we  assume that a global solution $(\mu,U)(s,q,\omega)$ satisfying several assumptions exists. For convenience, we make the following definition.
\defn{\label{def3.1}\rm  Fix two constants  $\gamma_+>1$ and $\gamma_->2$.  A smooth solution $(\mu,U)=(\mu,(U^{(I)}))(s,q,\omega)$ to the geometric reduced system \eqref{asy} is said to be $(\gamma_+,\gamma_-)$-\emph{admissible}, if there exists  $0<\delta_0<1$, such that $(\mu,U)$ is defined for all $(s,q,\omega)\in[-\delta_0,\infty)\times\R\times\mathbb{S}^2$ and if for each $(s,q,\omega) $ we have
\eq{\label{def3.1a3}-C\exp(C s )\leq\mu\leq -C^{-1}\exp(-C s )<0;}
\eq{\label{def3.1a4}|\mu_q|\lesssim \lra{q}^{-1-\gamma_{\sgn (q)}}|s\mu|;}
\eq{\label{def3.1a5}|\mu U_q|\lesssim 1;}
\eq{\label{def3.1a51}\sum_{J,K}|G_{3,J}\mu\partial_q(\mu U_q^{(K)})|\lesssim 1;}
\eq{\label{def3.1a52}\sum_{J,K}|G_{2,J}\partial_q(\mu U_q^{(K)})|\lesssim \lra{q}^{-\gamma_{\sgn(q)}};}
\eq{\label{def3.1a6}|\partial_s^a\partial_q^b\partial_\omega^c(\mu+2)|\lesssim_{a,b,c}\exp(C_{a,b,c} s )\lra{q}^{-b-\gamma_{\sgn(q)}},\qquad \forall a,b,c\geq 0;}
\eq{\label{def3.1a7}|\partial_s^a\partial_\omega^cU|\lesssim_{a,c}\exp(C_{a,c} s )\lra{\max\{0,-q\}}^{1-\gamma_-},\qquad \forall a,c\geq 0;}
\eq{\label{def3.1a8}|\partial_s^a\partial_q^b\partial_\omega^cU_q|\lesssim_{a,b,c}\exp(C_{a,b,c} s )\lra{q}^{-b-\gamma_{\sgn(q)}},\qquad \forall a,b,c\geq 0.}
Here  we set $\gamma_{\sgn(q)}=\gamma_+$ if $q\geq 0$ and $\gamma_{\sgn(q)}=\gamma_-$ otherwise.

For convenience, we set $\gamma:=\min\{\gamma_+,\gamma_-\}>1$. Then, the $\gamma_{\sgn(q)}$ in \eqref{def3.1a4}, \eqref{def3.1a6} and \eqref{def3.1a8} can be replaced by $\gamma$.

As explained in Remark \ref{rmkmthmdepeps}, here we implicitly assume that $(\mu,U)$ is independent of $\eps$. Even if $(\mu,U)$ depends on $\eps$, our proofs in this paper will work as long as all the constants in the bounds above are uniform in $\eps$.}

\rmk{\label{sec3:rmk:defnadm}\rm One way to understand the estimates \eqref{def3.1a3}--\eqref{def3.1a8} is as follows. Recall that our main goal is to construct a global solution $u$ to \eqref{qwe} matching a given admissible solution $(\mu,U)$. At the same time, we construct an approximate optical function $q$ by solving $q_t-q_r=\mu$. In the proof of Theorem \ref{mthm}, we require that the $q,u$, and their derivatives satisfy several bounds. Since the asymptotic behaviors of $q,u$, and their derivatives (at least near the light cone $|x|=t$) are described by this $(\mu,U)$, those required bounds can be translated into bounds for $(\mu,U)$. This is how we choose the estimates \eqref{def3.1a3}--\eqref{def3.1a8} in the definition.

For example, later we will show the following approximate identities whenever $||x|-t|\lesssim t^{C\eps}$ and $t\gtrsim_\eps 1$:
\fm{q_{\alpha}\approx -\frac{1}{2}\mu\wh{\omega}_\alpha,\ q_{\alpha\beta}\approx \frac{1}{4}\mu_q\mu\wh{\omega}_\alpha\wh{\omega}_\beta,\ 
u_\alpha\approx -\frac{1}{2}\eps r^{-1}\mu U_q\wh{\omega}_\alpha,\   u_{\alpha\beta}\approx \frac{1}{4}\eps r^{-1}\mu\partial_q(\mu U_q)\wh{\omega}_\alpha\wh{\omega}_\beta.}
Here recall that $\wh{\omega}_0=-1$ and $\omega_j=x_j/|x|$ for $j=1,2,3$. Now, \eqref{def3.1a3} implies that $t^{-C\eps}\lesssim |\partial q|\lesssim t^{C\eps}$. Another way to understand \eqref{def3.1a3} is that the foliation determined by the level sets of $q$ can be neither too dense nor too sparse. Recall from the introduction that $\mu$ measures the density of the stacking of the level sets of $q$. The estimate \eqref{def3.1a4} gives a control of  $|\partial^2q|$ in terms of $t,q$, and $|\partial q|$. The estimate \eqref{def3.1a5} implies that $|\partial u|\lesssim \eps t^{-1}$. The estimates \eqref{def3.1a51} and \eqref{def3.1a52} give controls of some linear combinations of $\partial^2u$. Similarly, we can express the higher derivatives of $q$ and $u$ in terms of the higher derivatives of $(\mu,U)$, so the bounds \eqref{def3.1a6}--\eqref{def3.1a8} give controls of the higher derivatives of $q$ and $u$. For example, we have $|Z^Iu|\lesssim \eps t^{-1+C\eps}$ for each multiindex $I$.

All the bounds for $q,u$, and their derivatives obtained above will be necessary for the proof of Theorem \ref{mthm}. For example, the estimates $|\partial u|\lesssim \eps t^{-1}$ and $|Z^Iu|\lesssim \eps t^{-1+C\eps}$ are necessary for the proof of the energy estimates; see \eqref{sec5a1} in Section \ref{sec5.2}. The estimates derived from \eqref{def3.1a51} and \eqref{def3.1a52} are used in the proof of Lemmas \ref{lem4.11} and \ref{lem4.10}.
}

\rmk{\label{def3.1rmk1}\rm The assumption \eqref{def3.1a6} can be interpreted as follows: we have \fm{|\partial_s^a\partial_q^b\partial_\omega^c(\mu+2)|\lesssim\exp(C s )\lra{q}^{-b-\gamma_+}}
whenever $q\geq 0$, and 
\fm{|\partial_s^a\partial_q^b\partial_\omega^c(\mu+2)|\lesssim\exp(C s )\lra{q}^{-b-\gamma_-}}
whenever $q\leq 0$. That is, we allow $\mu+2$ and its derivatives to satisfy different estimates for $q\geq 0$ and for $q\leq 0$. Similarly for \eqref{def3.1a4} and \eqref{def3.1a8}.

The assumption \eqref{def3.1a7} can be interpreted as follows: we have
\fm{|\partial_s^a\partial_\omega^cU|\lesssim \exp(C s )\lra{q}^{1-\gamma_-}}
whenever $q\leq 0$, and 
\fm{|\partial_s^a\partial_\omega^cU|\lesssim \exp(C s )}
whenever $q\geq 0$.

Note that we assume that $\gamma_->2$ and $\gamma_+>1$. That is, we have different requirements on the decay rates of $(\mu,U)$ and the derivatives (as well as $u$ and its derivatives) in $q$  for $q<0$ and for $q>0$. Such a difference is caused by the weights in Poincar$\acute{\rm e}$'s estimates in Section \ref{sec5.3}. For example, later we will apply Lemma \ref{lp1} to the solution $u$. In that lemma, we estimate a weighted $L^2$ norm of $u$ with the weight $\lra{r-t}^{-2}\cdot 1_{r>t}+\lra{r-t}^{-1+}\cdot 1_{r<t}$. The different powers of $\lra{r-t}$ for $r<t$ and for $r>t$ explain why we need different decay rates of $(\mu,U)$ and the derivatives for $q<0$ and for $q>0$.}

\rmk{\rm We only need to assume \eqref{def3.1a6}-\eqref{def3.1a8} hold for $a=0$. The estimates with $a>0$ actually follow from these estimates with $a=0$ and the geometric reduced system \eqref{asy}. In addition, if we assume that $\lim_{q\to-\infty}\partial_s^a\partial_\omega^cU=0$ for each $(s,\omega)$, then \eqref{def3.1a7} is implied by \eqref{def3.1a8} and the following  integral inequality
\eq{\label{intineq}\int_{-\infty}^q\lra{\rho}^{-\gamma_{\sgn(\rho)}}\ d\rho\lesssim \lra{\max\{0,-q\}}^{1-\gamma_-}.}
This inequality turns out to be very useful in the next section.}

\rmk{\rm By the product rule, we can prove \eqref{def3.1a8} with $U_q$ replaced by $\mu U_q$. In addition, since $(-\mu)\exp(C s )\gtrsim 1$, we can prove \eqref{def3.1a6}, \eqref{def3.1a7} and \eqref{def3.1a8} with an additional factor $(-\mu)$ on the right hand side.}

\rm

\section{Construction of an approximate solution}\label{sec4}

Our main goal in this section is to construct an approximate solution $u_{app}$ to \eqref{qwe}. Fix  $\gamma_+>1$ and $\gamma_->2$. Fix a $(\gamma_+,\gamma_-)$-admissible global solution $(\mu,U)(s,q,\omega)$ to the geometric reduced system \eqref{asy}. By Definition \ref{def3.1}, we obtain  $0<\delta_0<1$ such that $(\mu,U)$ is defined whenever $s\geq -\delta_0$. We fix a small constant $0<\delta<\delta_0$ whose value will be chosen later. Then, for all sufficiently small $\eps\ll1$ depending on $(\mu,U)$, $\gamma_\pm$ and $\delta$, we define a function $q(t,x)$ in the region
\fm{\Omega:=\{(t,x)\in\R^{1+3}:\ t>T_\eps=1/\eps,\ |x|>t/2\}} 
by solving   the PDE \fm{(\partial_t-\partial_r)q(t,x)=\mu(\eps\ln(t)-\delta,q(t,x),\omega);\qquad q(2|x|,x)=-|x|.}
With the function $q(t,x)$, we also define \fm{U^{(I)}(t,x):=U^{(I)}(\eps\ln(t)-\delta,q(t,x),\omega),\qquad I=1,\dots,M.}
In this section, we will show that in a conic neighborhood of the light cone $\{t=r\}$, $\eps r^{-1}U(t,x)$  and $q(t,x)$ are an approximate solution to \eqref{qwe} and an approximate optical function, respectively,  in the sense that for all $(t,x)\in\Omega$ with $|r-t|< t/2$, we have
\fm{g^{\alpha\beta}(\eps r^{-1}U,\partial(\eps r^{-1}U))\partial_\alpha\partial_\beta(\eps r^{-1}U)-f(\eps r^{-1}U,\partial(\eps r^{-1}U))=O(\eps t^{-3+C\eps}),}
\fm{g^{\alpha\beta}(\eps r^{-1}U,\partial(\eps r^{-1}U))q_\alpha q_\beta=O( t^{-2+C\eps}\lra{r-t}).}
In fact, we expect better estimates whenever  $q(t,x)\leq 0$, or whenever $|x|\leq t$.

To obtain an approximate solution for all $(t,x)$ with $t\geq  T_\eps$, we multiply $\eps r^{-1}U$ by a cutoff function.  Let $0<c<1/2$ be a  small constant. We set \fm{u_{app}(t,x)=\eps r^{-1}\psi(r/t)U(\eps\ln(t)-\delta,q(t,r,\omega),\omega).}
Here $\psi\equiv 1$ when $|r-t|\leq ct/2$ and $\psi\equiv 0$ when $|r-t|\geq ct$, which is used to localize $\eps r^{-1}U$ near the light cone $\{r=t\}$. The value of $c$ is not important in the current section, but it is closely related to Poincar$\acute{\rm e}$'s lemmas in the next section.

Our main proposition in this section is the following:

\prop{\label{mainprop4} Fix  $0<c<1/2$,  $\gamma_+>1$ and $\gamma_->2$. Let $(\mu,U)(s,q,\omega)$ be a $(\gamma_+,\gamma_-)$-admissible solution  to \eqref{asy} defined for all $s\geq -\delta_0$ with $0<\delta_0<1$, and fix $0<\delta<\delta_0$.  For all  sufficiently small $\eps\ll1$ depending on $(\mu,U)$, $\gamma_\pm$, $c$ and $\delta$,  we define $q=q(t,x)$ and $u_{app}=u_{app}(t,x)$ as above.   Then, for all $(t,x)$ with $t\geq T_\eps=1/\eps$, we have \fm{|\partial u_{app}(t,x)|\lesssim \eps t^{-1}.} Moreover, for all multiindices $I$ and for all $(t,x)$ with $t\geq T_\eps$, we have
\fm{|Z^I u_{app}(t,x)|\lesssim_I \eps t^{-1+C_I\eps}\lra{r-t}^{1-\gamma_-}1_{q\leq 0}+\eps t^{-1+C_I\eps}1_{q>0},}
\fm{|Z^I(g^{\alpha\beta}(u_{app})\partial_\alpha\partial_\beta u_{app}-f(u_{app},\partial u_{app}))(t,x)|\lesssim_I \eps t^{-3+C_I\eps}\lra{r-t}^{1-\gamma_-}1_{q\leq 0}+\eps t^{-3+C_I\eps}1_{q>0}.}}
\rmk{\rm The statement holds for all $\eps\in(0,\eps_0)$, where $\eps_0$ is a constant depending on $(\mu,U),\gamma_\pm,c,\delta$. Though the value of $\eps_0$ depends on $\delta$, the implicit constants in the estimates listed in this proposition can be chosen to be independent of $\delta$ as long as we have $\delta\in(0,1)$.

In contrast, all the implicit constants in this proposition depend on $c$. However,   $c$ is not introduced until the definition of $u_{app}$, so all the estimates  in Section \ref{sec4.1} and Section  \ref{sec4.2} below are independent of $c$.}
\rmk{\rm We can replace $1_{q\leq 0}$ and $1_{q>0}$ in the proposition with $1_{|x|\leq t}$ and $1_{|x|>t}$, respectively. See Remark \ref{rmk4.1.1} below. We also note that we have different estimates for $q\geq 0$ and $q<0$ because of \eqref{intineq}. }
\rm\bigskip

This proposition is proved in three steps. First, in Section \ref{sec4.1}, we construct $q(t,x)$ and $U(t,x)$ for all $(t,x)\in\Omega$, by solving the transport equation $q_t-q_r=\mu$. This can be done by applying the method of characteristics. Next, in Section \ref{sec4.2}, we prove that $\eps r^{-1}U(t,x)$ is an approximate solution to \eqref{qwe} near the light cone $\{t=r\}$ when $t\geq T_\eps$. To achieve this goal we prove several estimates for $q$ and $U$ in the region $t\sim r$. Finally, in Section \ref{sec4.3}, we define $u_{app}$ and prove the pointwise bounds listed in the proposition. To define $u_{app}$, we use a  cutoff function to localize $\eps r^{-1}U$ in a conical neighborhood of $\{t=r\}$.

In Section \ref{sec4.4}, we prove several extra pointwise estimates for $u_{app}$. These estimates come from the assumptions \eqref{def3.1a51} and \eqref{def3.1a52}, and they turn out to be very useful in the rest of this paper.

\subsection{Setup}\label{sec4.1}
Fix  a $(\gamma_+,\gamma_-)$-admissible solution $(\mu,U)$ and a sufficiently small constant $0<\delta<1$.  Here we assume that  $\gamma_+>1$, $\gamma_->2$ and $0<\delta<\delta_0<1$. As commented in Definition \ref{def3.1}, we set $\gamma:=\min\{\gamma_+,\gamma_-\}>1$ for convenience. Fix $\eps\ll 1$ and set $T_\eps=1/\eps$.  Now we  define a new function $q=q(t,x)$ in \eq{\Omega=\Omega_{\eps}:=\{(t,x)\in\R^{1+3}:\ t>T_\eps,\ |x|>t/2\}} by solving
\begin{equation}\label{qeqn}
(\partial_t-\partial_r)q(t,x)=\mu(\eps\ln(t)-\delta,q(t,x),\omega);\qquad q(2|x|,x)=-|x|.\end{equation}
Note that the right hand side of the first equation in \eqref{qeqn} is defined everywhere in $\Omega$ because \fm{\eps\ln T_\eps-\delta=\eps\ln\eps^{-1}-\delta> -\delta_0,\qquad \text{for }\eps\ll1.}
Also note that the second equation  implies that $q=r-t$ on the cone $\{r=t/2\}$. We can rewrite \eqref{qeqn} as an integral equation
\eq{\label{qeqn2}q(t,x)&=-\frac{r+t}{3}-\int_t^{2(r+t)/3}\mu(\eps\ln \tau-\delta,q(\tau,(r+t-\tau)\omega),\omega)\ d\tau\\
&=r-t-\int_t^{2(r+t)/3}[2+\mu(\eps\ln \tau-\delta,q(\tau,(r+t-\tau)\omega),\omega)]\ d\tau.}

This PDE can be solved by the method of characteristics. For each fixed $(t,x)\in\Omega $, we set $s(\tau)=\tau$ and $z(\tau)=q(\tau,(r+t-\tau)\omega)$. Then the characteristic ODE's associated to \eqref{qeqn} is 
\eq{\label{charode}\left\{\begin{array}{l}\displaystyle\dot{z}(\tau)=\mu(\eps\ln(s(\tau))-\delta,z(\tau),\omega),\\[1em]
\displaystyle\dot{s}(\tau)=1,\\[1em]
(z,s)|_{\tau=2(r+t)/3}=(-(r+t)/3,2(r+t)/3).\end{array}\right.}
To show that these ODE's admit a solution for all $\tau\geq T_\eps$, we apply  Picard's theorem. 
Recall that $|\mu(s,q,\omega)|\leq C \exp(C s )$ for each $(s,q,\omega)$ where $C$  is independent of $(q,\omega)$.  Thus, for each $\tau_0\geq T_\eps$, we have
\fm{|z(\tau_0)|&\leq \frac{r+t}{3}+\int_{\tau_0}^{2(r+t)/3} |\mu(\eps\ln \tau-\delta,z(\tau),\omega)|\ d\tau\\
&\lesssim r+t+\int_{\tau_0}^{2(r+t)/3}  \exp(C|\eps\ln \tau-\delta|)\ d\tau\\
&\leq r+t+\int_{\tau_0}^{2(r+t)/3}   \tau^{C\eps}\ d\tau.}
Thus, $|z(\tau)|$ cannot blow up when $\tau\geq  T_\eps$. By applying  Picard's theorem (e.g.\ Theorem 1.7 in \cite{MR2233925}), we conclude that \eqref{charode} has a solution for all $\tau\geq T_\eps$.

With the function $q=q(t,x)$, every function of $(s,q,\omega)$ induces a function of $(t,x)$. For a given function $F=F(s,q,\omega)$, by abuse of notation we  set \fm{F(t,x):=F(\eps\ln t-\delta,q(t,x),\omega).}So $F$ can denote both a function of $(s,q,\omega)$ and a function of  $(t,x)$. This allows us to define $U=(U^{(I)})(t,x)$ from $U(s,q,\omega)$. In addition,  for a given  function $F=F(t,x)$ and for each fixed point $(t,x)$,  by abuse of notation we write \fm{F(\tau):=F(\tau,(r+t-\tau)\omega),\qquad \tau\geq T_\eps.}So on the right hand side of \eqref{qeqn2}, we can simply write the integrand  as $\mu(\tau)$.

\subsection{Estimates for $(q,U)$}\label{sec4.2}
Let $0<\delta<\delta_0$ and $\eps\ll1$ be two small constants. Note that the choice of $\eps\ll1$ can depend on $\delta$ (in the sense that $\eps\in(0,\eps_0)$ where $\eps_0$ can depend on $\delta$) but not vice versa. We set
\eq{\label{c3defnD}\D:=\{(t,x)\in\Omega_{\eps}:\ |r-t|\leq t/2\}} and recall Definition \ref{c1defn1.5}  in Section \ref{c1sec1.6.4}. Our main goal now is to prove that  $\eps r^{-1}U$ has some good pointwise bounds and is an approximate solution whenever $t\geq T_\eps$ and $t\sim r$.

We start with a lemma for $q(t,x)$. Note that in this lemma we do not require $(t,x)\in\D$.
\lem{\label{lem4.1} Fix $(t,x)\in\Omega$. Then, for $\eps\ll1$, we have
\eq{\label{lem4.1c1} |q-(r-t)|\lesssim(t+r)^{C\eps}\lra{\max\{0,-q\}}^{1-\gamma_-};}
\eq{\label{lem4.1c3} (t+r)^{-C\eps}\lesssim\lra{q}/\lra{r-t}\lesssim (t+r)^{C\eps};}
for each fixed $\kappa\in(0,1)$, as long as $\eps\ll_{\kappa} 1$, we have
\eq{\label{lem4.1c4} (t,x)\in\D\qquad\text{whenever }|q(t,x)|\lesssim t^\kappa.}}
\begin{proof} Fix $(t,x)\in\Omega$ so we have $t\geq T_\eps$ and $r\geq t/2$. For any $\tau\in[t,2(r+t)/3]$, we have \fm{-Ce^{C s }\leq \dot{z}(\tau)=\mu(\tau)\leq -C^{-1}e^{-C s }<0,} \fm{|\mu(\tau)+2|\lesssim \exp(C|\eps\ln \tau-\delta|)\lra{z(\tau)}^{-\gamma_{\sgn(z(\tau))}}\lesssim  \tau^{C\eps}\lra{z(\tau)}^{-\gamma_{\sgn(z(\tau))}}} by  \eqref{def3.1a3} and \eqref{def3.1a6}. Here recall that $\gamma=\min\{\gamma_+,\gamma_-\}>1$. It follows from \eqref{qeqn2} that
\fm{|q-(r-t)|&\leq\int_t^{2(r+t)/3}|\mu(\tau)+2|\ d\tau\lesssim \int_t^{2(r+t)/3}\tau^{C\eps}(-\dot{z}(\tau))\cdot\tau^{C\eps}\lra{z(\tau)}^{-\gamma_{\sgn(z(\tau))}}\ d\tau\\
&\lesssim (r+t)^{C\eps}\int_t^{2(r+t)/3}(-\dot{z}(\tau))\lra{z(\tau)}^{-\gamma_{\sgn(z(\tau))}}\ d\tau\\
&\lesssim(r+t)^{C\eps}\int_{-((r+t)/3)}^{z(t)} \lra{\rho}^{-\gamma_{\sgn(\rho)}}\ d\rho\lesssim (r+t)^{C\eps}\lra{\max\{0,-q(t,x)\}}^{1-\gamma_-}.}
Here the integral is taken along the characteristic $(\tau,(r+t-\tau)\omega)$ for $\tau\in[t,2(r+t)/3]$. In the last step, we use the fact that $\gamma_+,\gamma_->1$ and the integral inequality \eqref{intineq}. So we obtain \eqref{lem4.1c1}. 

By \eqref{lem4.1c1}, we have
\fm{\lra{r-t}&\lesssim 1+|r-t|\leq 1+|q|+|q-(r-t)|\lesssim \lra{q}+(r+t)^{C\eps}\lesssim \lra{q}(r+t)^{C\eps},\\
\lra{q}&\lesssim 1+|q|\leq 1+|r-t|+|q-(r-t)|\lesssim \lra{r-t}+(r+t)^{C\eps}\lesssim \lra{r-t}(r+t)^{C\eps}.}
This gives \eqref{lem4.1c3}. 

Finally, if $|q|\lesssim t^\kappa$, then 
\fm{|r-t|\leq |q-(r-t)|+|q|\lesssim (t+r)^{C\eps}+t^\kappa.}
By choosing $\eps\ll 1$, we have $C\eps<1/2$. It follows that
\fm{r/t-1\leq Ct^{-1}(t+r)^{1/2}+Ct^{\kappa-1}=Ct^{-1/2}(r/t+1)^{1/2}+Ct^{\kappa-1}.}
By choosing $\eps\ll_\kappa 1$, we have $Ct^{-1/2}+Ct^{\kappa-1}\leq 1$ for all $t\geq T_\eps$. As a result, we have
\fm{0<r/t\leq \sqrt{1+r/t}+2\Longrightarrow r/t\lesssim1}
and thus
\fm{|r-t|\leq Ct^{1/2}+Ct^\kappa=\frac{1}{2}t(2Ct^{-1/2}+2Ct^{\kappa-1})\leq \frac{1}{2}t(2CT_\eps^{-1/2}+2CT_\eps^{\kappa-1}).}
By choosing $\eps\ll_\kappa 1$, we obtain \eqref{lem4.1c4}.
\end{proof}
\rmk{\label{rmk4.1.1}\rm Here is an important corollary of this lemma. For any two fixed real constants $\gamma_1,\gamma_2$, we have
\eq{\label{rmk4.1.1c2}t^{-C\eps}\lesssim\frac{\lra{q}^{\gamma_1}1_{q\geq 0}+\lra{q}^{\gamma_2}1_{q< 0}}{\lra{r-t}^{\gamma_1}1_{r-t\geq 0}+\lra{r-t}^{\gamma_2}1_{r-t< 0}}\lesssim t^{C\eps},\qquad\forall(t,x)\in\Omega.}
All the constants here  could depend on $\gamma_1,\gamma_2$.

To prove \eqref{rmk4.1.1c2}, we start to prove a weaker estimate
\eq{\label{rmk4.1.1c1}(r+t)^{-C\eps}\lesssim\frac{\lra{q}^{\gamma_1}1_{q\geq 0}+\lra{q}^{\gamma_2}1_{q< 0}}{\lra{r-t}^{\gamma_1}1_{r-t\geq 0}+\lra{r-t}^{\gamma_2}1_{r-t< 0}}\lesssim (r+t)^{C\eps},\qquad\forall(t,x)\in\Omega.}
In fact, if $q\cdot(r-t)\geq 0$, then the quotient above is equal to either $(\lra{q}/\lra{r-t})^{\gamma_1}$ or $(\lra{q}/\lra{r-t})^{\gamma_2}$, so we can directly  apply \eqref{lem4.1c3}. If $q\cdot(r-t)< 0$, then $|q|+|r-t|=|q-r+t|\lesssim (r+t)^{C\eps}$, so the estimate follows.

To end the proof, we claim that $r\sim q$ whenever $r>2t$ and $t\geq 1/\eps$. In fact, if $r>2t$, then $|q-r+t|\lesssim (r+t)^{C\eps}\lesssim r^{1/2}$ as long as $\eps\ll1$. It follows that
\fm{|q/r-1|\leq |q-r+t|/r+t/r\leq Cr^{-1/2}+\frac{1}{2}\leq CT_\eps^{-1/2}+\frac{1}{2}.}
By choosing $\eps\ll1$, we have $|q/r-1|<3/4$ and thus $q\sim r$.
This finishes the proof of \eqref{rmk4.1.1c2}.
}
\rm

\bigskip

We now move on to estimates for $\partial q$.  In Lemma \ref{lem4.2} and Lemma \ref{lem4.4}, we give the pointwise bounds for $\nu=q_t+q_r$, $\lambda_i=q_i-\omega_i q_r$ and $\nu_r$. In Lemma \ref{lem4.3}, we find the first terms in the asymptotic expansion of $\nu$ and $\nu_q$ when $1/\eps\leq t\sim r$.

\lem{\label{lem4.2}Set $\nu(t,x):=(\partial_t+\partial_r)q$ and $\lambda_i(t,x):=(\partial_i-\omega_i\partial_r)q$ for each $i=1,2,3$. For $(t,x)\in\Omega$, we have
\eq{\label{lem4.2c1}|\nu(t,x)|&\lesssim (r+t)^{-\gamma_{-}+C\eps}1_{r\leq 3t}+(r+t)^{-\gamma_{+}+C\eps}1_{r>3t}+\eps (r+t)^{-1+C\eps}\lra{\max\{0,-q\}}^{1-\gamma_-}\\&\lesssim  (r+t)^{-1+C\eps}\lra{\max\{0,-q\}}^{1-\gamma_-},}
\begin{equation}\label{lem4.2c2}|\lambda_i(t,x)|\lesssim (r+t)^{-1+C\eps}\lra{\max\{0,-q\}}^{1-\gamma_-}.\end{equation}
Note that we do not need to assume that $(t,x)\in\D$ in this lemma.}
\begin{proof}
Fix $(t,x)\in\Omega$. We have
\eq{\label{nueqn}(\partial_t-\partial_r)\nu&=(\partial_t+\partial_r)\mu=\mu_q\nu+\eps t^{-1}\mu_s.}
Recall from \eqref{def3.1a4} that $|\mu_q|\lesssim \lra{q}^{-1-\gamma}|s\mu|$ where $\gamma=\min\{\gamma_+,\gamma_-\}$, so \eq{\label{lem4.2f1}\int_t^{2(r+t)/3}|\mu_q(\tau)|\ d\tau&\lesssim\int_t^{2(r+t)/3}\lra{z(\tau)}^{-1-\gamma}|(\eps\ln \tau-\delta)\mu(\tau)|\ d\tau\\
&\lesssim (\eps\ln (r+t)+1)\int_t^{2(r+t)/3}\lra{z(\tau)}^{-1-\gamma}(-\mu(\tau))\ d\tau\\
&\lesssim(\eps\ln (r+t)+1)\int^{z(t)}_{z(2(r+t)/3)}\lra{\rho}^{-1-\gamma}\ d\rho\lesssim\eps\ln (r+t)+1.}
Here the integral is taken along the characteristic; see the notation introduced at the end of Section \ref{sec4.1}. Moreover, we have 
\fm{\int_t^{2(r+t)/3}|\eps\tau^{-1}\mu_s(\tau)|\ d\tau&\lesssim \int_t^{2(r+t)/3} \eps \tau^{-1}\cdot\tau^{C\eps}|\mu(\tau)|\lra{z(\tau)}^{-\gamma_{\sgn(z(\tau))}}\ d\tau\\&\lesssim \eps t^{-1+C\eps}\int^{z(t)}_{-(r+t)/3} \lra{\rho}^{-\gamma_{\sgn(\rho)}}\ d\rho\lesssim  \eps t^{-1+C\eps}\lra{\max\{0,-q\}}^{1-\gamma_-}.}
Since $q(2|x|,x)=-|x|$, by \eqref{def3.1a6} we have 
\fm{\nu|_{\tau=2(r+t)/3}&=\frac{2}{3}(2\partial_t+\partial_r)q(2(r+t)/3,(r+t)\omega/3)-\frac{1}{3}\mu|_{\tau=2(r+t)/3}\\
&=-\frac{2}{3}(2+\mu)|_{(s,q,\omega)=(\eps\ln(2(r+t)/3)-\delta,-(r+t)/3,\omega)}=O(\lra{r+t}^{-\gamma_-+C\eps}).}
Here we have $(2\partial_t+\partial_r)q(2(r+t)/3,(r+t)\omega/3)=-1$ because $2\partial_t+\partial_r$ is tangent to the surface $|x|=t/2$.
Thus, by  Gronwall's inequality, we conclude that
\fm{|\nu(t,x)|&\lesssim\exp(C\eps\ln(r+t)+C)\cdot(\lra{r+t}^{-\gamma_-+C\eps}+\eps t^{-1+C\eps}\lra{\max\{0,-q\}}^{1-\gamma_-})\\
&\lesssim (r+t)^{-\gamma_-+C\eps}+\eps t^{-1+C\eps}(r+t)^{C\eps}\lra{\max\{0,-q\}}^{1-\gamma_-}\\
&\lesssim t^{-1+C\eps}(r+t)^{C\eps}\lra{\max\{0,-q\}}^{1-\gamma_-}.}
The last estimate holds because $\lra{q}\lesssim \lra{r-t}(r+t)^{C\eps}\lesssim (r+t)^{1+C\eps}$ by \eqref{lem4.1c3}.
That is, we have proved \eqref{lem4.2c1} for $r\leq 3t$.

If $r\geq 2t$, we have proved in Remark \ref{rmk4.1.1} that $r\sim q>0$ as long as $\eps\ll1$, so for all  $\tau\in[t,(r+t)/3]$ we have \fm{0<z(\tau)\sim(r+t-2\tau)\sim (r+t).}
In particular, we have $\lra{\max\{0,-z(\tau)\}}^{1-\gamma_-}=1$.
And since $|\mu(\tau)|\lesssim \tau^{C\eps}$, we have
\fm{\int_t^{(r+t)/3}|\mu_q(\tau)|\ d\tau&\lesssim\int_t^{(r+t)/3}\lra{z(\tau)}^{-1-\gamma}|(\eps\ln \tau-\delta)\mu(\tau)|\ d\tau\\
&\lesssim (\eps\ln (r+t)+1)\int_t^{(r+t)/3}(r+t)^{-(1+\gamma)}\cdot\tau^{C\eps}\ d\tau\\
&\lesssim(r+t)^{-(1+\gamma)+C\eps}\cdot(r+t)^{1+C\eps}.}
For $\eps\ll 1$, we have
\fm{-(1+\gamma)+C\eps+1+C\eps\leq -\gamma+C\eps<-1.}
Thus, by choosing $\eps\ll1$, we have
\fm{\int_t^{(r+t)/3}|\mu_q(\tau)|\ d\tau&\leq  C T_\eps^{-1}\leq 1,\qquad \forall t\geq T_\eps.}
Similarly, for $\tau\in[t,(r+t)/3]$ we have
\fm{\int_t^{(r+t)/3}|\eps\tau^{-1}\mu_s(\tau)|\ d\tau&\lesssim \int_t^{(r+t)/3} \eps \tau^{-1}\cdot\tau^{C\eps}|\mu(\tau)|\lra{z(\tau)}^{-\gamma_+}\ d\tau\\
&\lesssim \int_t^{(r+t)/3}\eps \tau^{-1+C\eps}(r+t)^{-\gamma_+}\ d\tau\lesssim (r+t)^{-\gamma_++C\eps}.}
At $\tau=(r+t)/3$, we have proved that \fm{|\nu(\tau)|\lesssim (r+t)^{-\gamma_-+C\eps}+ \eps(r+t)^{-1+C\eps}\lesssim \eps(r+t)^{-1+C\eps}} because  $\gamma_->2$ and $r+t\geq \eps^{-1}$. Thus, by  Gronwall's inequality, we have
\fm{|\nu(t,x)|&\lesssim \exp(C)((r+t)^{-\gamma_++C\eps}+\eps(r+t)^{-1+C\eps})\lesssim (r+t)^{-1+C\eps}.}
Here we use $\gamma_+>1$.
This finishes the proof of \eqref{lem4.2c1}.

To prove \eqref{lem4.2c2}, we note that
\eq{\label{lameqn}(\partial_t-\partial_r)\lambda_i&=(\partial_i-\omega_i\partial_r)\mu+r^{-1}\lambda_i=(\mu_q+r^{-1})\lambda_i+\sum_{l=1}^3 \mu_{\omega_l}\partial_i\omega_l.}
By \eqref{lem4.2f1} we have
\fm{\int_t^{2(r+t)/3}|\mu_q(\tau)|+(r+t-\tau)^{-1}\ d\tau&\lesssim \eps\ln(t+r)+1+\int_t^{2(r+t)/3}(r+t-2(r+t)/3)^{-1}\ d\tau\\&\lesssim \eps\ln(t+r)+1.}
Besides, note that  by \eqref{def3.1a6}
\fm{\int_t^{2(r+t)/3}|\mu_{\omega_l}\partial_i\omega_l|\ d\tau&\lesssim \int_t^{2(r+t)/3}(r+t-\tau)^{-1}\lra{z(\tau)}^{-\gamma_{\sgn(z(\tau))}}\tau^{C\eps}\ d\tau\\
&\lesssim (r+t)^{-1+C\eps}\int_t^{2(r+t)/3}\lra{z(\tau)}^{-\gamma_{\sgn(z(\tau))}}(-\dot{z}(\tau))\ d\tau\\
&\lesssim  (r+t)^{-1+C\eps}\lra{\max\{0,-q\}}^{1-\gamma_-}.}
Since $\lambda_i(2(r+t)/3,(r+t)\omega/3)=0$, we conclude that $|\lambda_i|\lesssim (r+t)^{-1+C\eps}\lra{\max\{0,-q\}}^{1-\gamma_-}$ by applying  Gronwall's inequality. 
\end{proof}
\rmk{\label{rmk4.2.1}\rm By choosing $\eps\ll 1$, for all $(t,x)\in\Omega$   we have \fm{q_r&=\frac{-\mu+\nu}{2}\geq C^{-1}t^{-C\eps}-C  (t+r)^{-1+C\eps}\geq (2C)^{-1}t^{-C\eps},\\
q_t&=\frac{\mu+\nu}{2}\leq -C^{-1}t^{-C\eps}+C(t+r) ^{-1+C\eps}\leq -(2C)^{-1}t^{-C\eps}.}
As a result, for fixed $t\geq T_\eps$ and $\omega$, the function $r\mapsto q(t,r\omega)$ is a strictly increasing function, and we have $\lim_{r\to\infty}q(t,r\omega)=\infty$. This implies that for each $t\geq  T_\eps$ and $q^0\geq -t/2$, there exists a unique $r$ such that $q(t,r\omega)=q^0$. So \fm{\Omega\ni(t,x)\mapsto (\eps\ln(t)-\delta,q(t,x),\omega)} has an inverse map $(s,q,\omega)\mapsto (e^{(s+\delta)/\eps},r(s,q,\omega)\omega)$. By the inverse function theorem, we can check that the map $(t,x)\mapsto(s,q,\omega)$ is a diffeomorphism.

From now on, any function $V$ can be written as both $V(t,x)$ and $V(s,q,\omega)$ at the same time. Thus, for any function $V$ on $(t,x)$, we can define $\partial_s^a\partial_q^b\partial_{\omega}^{c}V$ using the chain rule and Leibniz's rule. Note that in this paper, $\partial_\omega$ will only be used under the coordinate $(s,q,\omega)$ and will never be used under the coordinate $(t,x)$.}

\lem{\label{lem4.3} For $(t,x)\in\Omega$ with $r\lesssim t$,
\eq{\label{lem4.3c1}&\nu-\frac{\eps G_{2,J}(\omega)}{4 r}\mu U^{(J)}+\frac{\eps G_{3,J}(\omega)}{8r}\mu^2U_q^{(J)}\\&=O(t^{-\gamma_-+C\eps}+\eps t^{-2+C\eps}\lra{q}\lra{\max\{0,-q\}}^{1-\gamma_-})=O( t^{-2+C\eps}\lra{q}\lra{\max\{0,-q\}}^{1-\gamma_-}),}
\eq{\label{lem4.3c2}&\nu_q-\frac{\eps G_{2,J}(\omega)}{4 r}\partial_q(\mu U^{(J)})+\frac{\eps G_{3,J}(\omega)}{8r}\partial_q(\mu^2U_q^{(J)})\\&=O(t^{-1-\gamma_-+C\eps}+(t^{-\gamma_-+C\eps}+\eps  t^{-2+C\eps})\lra{\max\{0,-q\}}^{1-\gamma_-})=O( t^{-2+C\eps}\lra{\max\{0,-q\}}^{1-\gamma_-}).} }
\begin{proof}
Set \fm{V:=\nu-\frac{\eps G_{2,J}(\omega)}{4 r}\mu U^{(J)}+\frac{\eps G_{3,J}(\omega)}{8r}\mu^2U_q^{(J)}.}
Then, by  \eqref{nueqn} and \eqref{asy}, at each $(t,x)\in\Omega$ with $r\lesssim t$, we have
\eq{\label{lem4.3f1}&(\partial_t-\partial_r)V\\&=\mu_q\nu+\eps t^{-1}\mu_s+\frac{\eps G_{2,J}(\omega)}{4 r^2}\mu U^{(J)}-\frac{\eps G_{3,J}(\omega)}{8r^2} \mu^2U_q^{(J)} \\
&\quad-\frac{\eps G_{2,J}(\omega)}{4 r}[\mu\partial_q(\mu U^{(J)})+\eps t^{-1}\partial_s(\mu U^{(J)})]+\frac{\eps G_{3,J}(\omega)}{8r}[\mu\partial_q(\mu^2U_q^{(J)})+\eps t^{-1}\partial_s(\mu^2U_q^{(J)})]\\
&=\mu_qV+\eps t^{-1}\mu_s+\frac{\eps G_{2,J}(\omega)}{4 r^2}\mu U^{(J)}-\frac{\eps G_{3,J}(\omega)}{8r^2} \mu^2U_q^{(J)} \\
&\quad-\frac{\eps G_{2,J}(\omega)}{4 r}[\mu^2U_q^{(J)}+\eps t^{-1}\partial_s(\mu U^{(J)})]+\frac{\eps G_{3,J}(\omega)}{8r}[\eps t^{-1}\partial_s(\mu^2U_q^{(J)})+\mu^2\partial_q(\mu U_q^{(J)})]\\
&=\mu_qV+\eps (tr)^{-1}(r-t)\mu_s+\frac{\eps G_{2,J}(\omega)}{4 r^2}\mu U^{(J)}-\frac{\eps G_{3,J}(\omega)}{8r^2} \mu^2U_q^{(J)} \\
&\quad-\frac{\eps^2 G_{2,J}(\omega)}{4 rt}\partial_s(\mu U^{(J)})+\frac{\eps^2 G_{3,J}(\omega)}{8rt}\partial_s(\mu^2U_q^{(J)})\\
&=\mu_qV+O(\eps t^{-2}\lra{r-t}|\mu_s|+\eps t^{-2}(|\mu U|+|\mu^2U_q|+|\partial_s(\mu U)|+|\partial_s(\mu^2U_q)|)).}
By \eqref{asy}, \eqref{def3.1a6}--\eqref{def3.1a8} and the product rule, we have 
\fm{|\mu U|+|\mu U_s|\lesssim t^{C\eps}\lra{\max\{0,-q\}}^{1-\gamma_-}|\mu|,}
\fm{|\mu_s|&\lesssim |\mu^2U_q|+|\mu^2\partial_q(\mu U_q)|\lesssim t^{C\eps}\lra{q}^{-\gamma_{\sgn (q)}}|\mu|,
}
\fm{|\partial_s(\mu U)|+|\partial_s(\mu^2U_q)|&\lesssim |\mu_sU|+|\mu U_s|+|\mu_s\mu U_q|+|\mu\cdot \mu^2U_q\cdot U_q|\\
&\lesssim t^{C\eps}\lra{q}^{-\gamma_{\sgn (q)}}|\mu U|+|\mu U_s|+t^{C\eps}\lra{q}^{-2\gamma_{\sgn (q)}}|\mu|\\
&\lesssim t^{C\eps}\lra{\max\{0,-q\}}^{1-\gamma_-}|\mu|+t^{C\eps}\lra{q}^{-2\gamma_{\sgn (q)}}|\mu|.}
Thus, in the region where $t\geq T_\eps$ and $r\sim t$, we have
\fm{&|V_t-V_r-\mu_qV|\\
&\lesssim \eps t^{-2+C\eps}\lra{r-t}\lra{q}^{-\gamma_{\sgn (q)}}|\mu|+\eps t^{-2+C\eps}\lra{\max\{0,-q\}}^{1-\gamma_-}|\mu|+\eps t^{-2+C\eps}\lra{q}^{-\gamma_{\sgn (q)}}|\mu|
\\
&\lesssim \eps t^{-2+C\eps}\lra{q}^{1-\gamma_{\sgn (q)}}|\mu|+\eps t^{-2+C\eps}\lra{\max\{0,-q\}}^{1-\gamma_-}|\mu|.}
In the last estimate,  we apply \eqref{lem4.1c3}. And since  $(r+t-\tau)\sim\tau$ whenever $t\leq \tau\leq 2(r+t)/3$ and $r\lesssim t$, we conclude that
\eq{\label{lem4.3f2}&\int_t^{2(r+t)/3} |V_t-V_r-\mu_qV|\ d\tau\\
&\lesssim\eps t^{-2+C\eps}\int_t^{2(r+t)/3} [\lra{z(\tau)}^{1-\gamma_{\sgn (z(\tau))}}+\lra{\max\{0,-z(\tau)\}}^{1-\gamma_-}](-\dot{z}(\tau))\ d\tau\\
&\lesssim\eps t^{-2+C\eps}\int_{z(2(r+t)/3)}^{z(t)} [\lra{\rho}^{1-\gamma_{\sgn (\rho)}}+\lra{\max\{0,-\rho\}}^{1-\gamma_-}]\ d\rho\\
&\lesssim \eps t^{-2+C\eps}(\lra{q}^{2-\gamma_-}1_{q\leq 0}+\lra{q}1_{q>0}).}
Here we  explain how to obtain the last estimate. Note that $z(2(r+t)/3)=-(r+t)/3<0$. Thus, if $z(t)=q(t,x)\leq 0$, we have
\fm{\int_{z(2(r+t)/3)}^{z(t)} [\lra{\rho}^{1-\gamma_{\sgn (\rho)}}+\lra{\max\{0,-\rho\}}^{1-\gamma_-}]\ d\rho&=\int_{z(2(r+t)/3)}^{z(t)} \lra{\rho}^{1-\gamma_-}\ d\rho\lesssim\lra{q}^{2-\gamma_-}.}
Here we use $\gamma_->2$. If $z(t)=q(t,x)\geq 0$, then
\fm{&\int_{z(2(r+t)/3)}^{z(t)} [\lra{\rho}^{1-\gamma_{\sgn (\rho)}}+\lra{\max\{0,-\rho\}}^{1-\gamma_-}]\ d\rho\\
&=\int_{z(2(r+t)/3)}^{0} \lra{\rho}^{1-\gamma_-}\ d\rho+\int_{0}^{z(t)} (\lra{\rho}^{1-\gamma_+}+1)\ d\rho\lesssim 1+q(t,x)\sim\lra{q}.}
This finishes the proof of \eqref{lem4.3f2}.

Now, recall from the proof of Lemma \ref{lem4.2} that $\nu|_{\tau=2(r+t)/3}=O(\lra{r+t}^{-\gamma_-+C\eps})$. And since $z(2(r+t)/3)<0$, at  $(2(r+t)/3,(r+t)\omega/3)$ we have
\fm{|\mu U|+|\mu^2U_q|\lesssim (r+t)^{C\eps}\lra{-(r+t)/3}^{1-\gamma_-}+(r+t)^{C\eps}\lra{-(r+t)/3}^{-\gamma_-}\lesssim (r+t)^{1-\gamma_-+C\eps}.}
As a result, we have $V|_{\tau=2(r+t)/3}=O(t^{-\gamma_-+C\eps})$. By  Gronwall's inequality, we conclude that whenever $t\geq  T_\eps$ and $t/2\leq r\lesssim t$, \fm{|V(t,x)|&\lesssim \exp(C\eps\ln(r+t)+C)(t^{-\gamma_-+C\eps}+\eps t^{-2+C\eps}(\lra{q}^{2-\gamma_-}1_{q\leq 0}+\lra{q}1_{q>0}))\\
&\lesssim t^{-\gamma_-+C\eps}+\eps t^{-2+C\eps}(\lra{q}^{2-\gamma_-}1_{q\leq 0}+\lra{q}1_{q>0})\\
&\lesssim  t^{-2+C\eps}(\lra{q}^{2-\gamma_-}1_{q\leq 0}+\lra{q}1_{q>0}).}
In the last estimate, we note that $t^{2-\gamma_--C\eps}\lesssim \lra{q}^{2-\gamma_-}$ because of \eqref{lem4.1c3}. This finishes the proof of \eqref{lem4.3c1}.

By differentiating \eqref{lem4.3f1}, we have
\eq{&(\partial_t-\partial_r)V_r\\&=\mu_qV_r+\mu_{qq}q_rV+\eps r^{-2}\mu_s+\eps (rt)^{-1}(r-t)\mu_{sq}q_r -\frac{\eps G_{2,J}(\omega)}{2 r^3}\mu U^{(J)}+\frac{\eps G_{3,J}(\omega)}{2r^3} \mu^2U_q^{(J)}\\
&\quad+\frac{\eps G_{2,J}(\omega)}{4 r^2}q_r\partial_q(\mu U^{(J)})-\frac{\eps G_{3,J}(\omega)}{8r^2} q_r\partial_q(\mu^2U_q^{(J)}) +\frac{\eps^2 G_{2,J}(\omega)}{4 r^2t}\partial_s(\mu U^{(J)})-\frac{\eps^2 G_{3,J}(\omega)}{8r^2t}\partial_s(\mu^2U_q^{(J)})\\
&\quad-\frac{\eps^2 G_{2,J}(\omega)}{4 rt}q_r\partial_q\partial_s(\mu U^{(J)})+\frac{\eps^2 G_{3,J}(\omega)}{8rt}q_r\partial_q\partial_s(\mu^2U_q^{(J)}).}
By \eqref{def3.1a3} and \eqref{lem4.2c1}, we have $q_r=(\nu-\mu)/2=-\mu/2+O(\eps t^{-1+C\eps})=|\mu|/2+O(\eps t^{-1+C\eps}|\mu|)$, so by choosing  $\eps\ll 1$, we have $q_r\sim (-\mu)$. By applying \eqref{lem4.3c1} and \eqref{def3.1a6}-\eqref{def3.1a8}, for $r\sim t$ we have
\fm{&|(\partial_t-\partial_r)V_r-\mu_qV_r|\\
&\lesssim t^{C\eps}\cdot\lra{q}^{-2-\gamma_{\sgn(q)}}|\mu|\cdot(t^{-\gamma_-+C\eps}+ \eps t^{-2+C\eps}\lra{q}^{2-\gamma_-}1_{q\leq 0}+ \eps t^{-2+C\eps}\lra{q}1_{q>0})\\&\quad+\eps t^{-2+C\eps}\lra{q}^{-\gamma_{\sgn(q)}}|\mu|+\eps t^{-3+C\eps}\lra{\max\{0,-q\}}^{1-\gamma_-}|\mu|.}
It follows that for $t/2\leq r\lesssim t$, 
\fm{&\int_t^{2(r+t)/3}|(\partial_t-\partial_r)V_r-\mu_qV_r|\ d\tau\\
&\lesssim \int_{-(r+t)/3}^{q(t,x)} t^{-\gamma_-+C\eps}\lra{\rho}^{-2-\gamma_{\sgn(\rho)}}+ \eps t^{-2+C\eps}\lra{\rho}^{-2\gamma_-}1_{\rho\leq 0}+ \eps t^{-2+C\eps}\lra{\rho}^{-1-\gamma_+}1_{\rho> 0}\\&\qquad\qquad+\eps t^{-2+C\eps}\lra{\rho}^{-\gamma_{\sgn(\rho)}}+\eps t^{-3+C\eps}(\lra{\rho}^{1-\gamma_-}1_{\rho\leq 0}+1_{\rho>0})\ d\rho\\
&\lesssim ( t^{-\gamma_-+C\eps}\lra{q}^{-1-\gamma_-}+\eps t^{-2+C\eps}\lra{q}^{1-\gamma_-}+\eps t^{-3+C\eps}\lra{q}^{2-\gamma_-})1_{q\leq 0}\\&\quad+(t^{-\gamma_-+C\eps}+\eps  t^{-2+C\eps}+\eps t^{-3+C\eps}q)1_{q>0}\\
&\lesssim  ( t^{-\gamma_-+C\eps}\lra{q}^{-1-\gamma_-}+\eps t^{-2+C\eps}\lra{q}^{1-\gamma_-})1_{q\leq 0}+(t^{-\gamma_-+C\eps}+\eps  t^{-2+C\eps})1_{q>0}.}

It remains to estimate $V_r$ at $\tau =2(r+t)/3$. Recall that $z(2(r+t)/3)=-(r+t)/3<0$ and that $r\sim t$, so $|\mu_q|\lesssim \lra{q}^{-1-\gamma_{\sgn(q)}}t^{C\eps}\lesssim t^{-1-\gamma_-+C\eps}$ at $\tau =2(r+t)/3$. Similarly, we can derive bounds for $\mu_s,U_q,$ etc., at $\tau=2(r+t)/3$. It follows that at $\tau=2(r+t)/3$,
\fm{-\frac{\eps G_{2,J}(\omega)}{4 r}\partial_q(\mu U^{(J)})+\frac{\eps G_{3,J}(\omega)}{8r}\partial_q(\mu^2U_q^{(J)})=O(\eps t^{-1-\gamma_-+C\eps}).}
Since $q=r-t$ on the cone $r=t/2$, we have $(2\partial_t+\partial_r)(2\partial_t+\partial_r)q=0$ on this cone.  Thus, at $\tau=2(r+t)/3$, we have
\fm{0&=(2\partial_t+\partial_r)(\frac{1}{2}\mu+\frac{3}{2}\nu)=\frac{3}{2}(2\partial_t+\partial_r)\nu+\frac{1}{2}(\mu_q(2q_t+q_r)+2\eps t^{-1}\mu_s)\\
&=\frac{3}{2}(2\partial_t+\partial_r)\nu+\frac{1}{2}(-\mu_q+2\eps t^{-1}\mu_s)=\frac{3}{2}(2\partial_t+\partial_r)\nu+O(t^{-1-\gamma_-+C\eps}).}
Meanwhile, since $\nu|_{\tau=2(r+t)/3}=O(t^{-\gamma_-+C\eps})$, by \eqref{lem4.2f1} we have
\fm{|(\partial_t-\partial_r)\nu|&\lesssim |\mu_q\nu|+\eps t^{-1}|\mu_s|\lesssim t^{-1-2\gamma_-+C\eps}+\eps t^{-1-\gamma_-+C\eps}\lesssim t^{-1-\gamma_-+C\eps}.}
at $\tau=2(r+t)/3$. As a result, we have $\nu_r|_{\tau=2(r+t)/3}=O( t^{-1-\gamma_-+C\eps})$. And since $\partial_r=q_r\partial_q$, we have $V_r|_{\tau=2(r+t)/3}=O(t^{-1-\gamma_-+C\eps})$. It now follows from  Gronwall's inequality that
\fm{|V_r(t,x)|&\lesssim   (t^{-1-\gamma_-+C\eps}+ t^{-\gamma_-+C\eps}\lra{q}^{-1-\gamma_-}+\eps t^{-2+C\eps}\lra{q}^{1-\gamma_-})1_{q\leq 0}+(t^{-\gamma_-+C\eps}+\eps  t^{-2+C\eps})1_{q>0}\\&\lesssim t^{-2+C\eps}\lra{\max\{0,-q\}}^{1-\gamma_-}.}
To finish the proof of \eqref{lem4.3c2}, we simply note that $q_r\geq C^{-1}t^{-C\eps}$ and that $\partial_q=q_r^{-1}\partial_r$.
\end{proof}

\lem{\label{lem4.4} Fix $(t,x)\in\Omega$. If $r\leq 3t$ we have
\begin{equation}\label{lem4.4c1}\nu_r=O(\eps t^{-1+C\eps}\lra{q}^{-\gamma_{\sgn(q)}}+t^{-1-\gamma_-+C\eps}+(t^{-\gamma_-+C\eps}+\eps  t^{-2+C\eps})\lra{\max\{0,-q\}}^{1-\gamma_-}).\end{equation}
If $r>2t$ we have
\begin{equation}\label{lem4.4c2}\nu_r=O( (r+t)^{-\gamma_-+C\eps}+(r+t)^{-1-\gamma_++C\eps}+\eps (r+t)^{-2+C\eps}).\end{equation}
The same estimates hold if we replace $\nu_r$ with $\nu_q=q_r^{-1}\nu_r$.}
\begin{proof} Fix $(t,x)\in\Omega$.  First assume that $r\leq 3t$. By \eqref{lem4.3c2}, we have
\fm{\nu_r&=q_r(\frac{\eps G_{2,J}(\omega)}{4 r}\partial_q(\mu U^{(J)})-\frac{\eps G_{3,J}(\omega)}{8r}\partial_q(\mu^2U_q^{(J)}))\\&\quad+O( (t^{-1-\gamma_-+C\eps}+ t^{-\gamma_-+C\eps}\lra{q}^{-1-\gamma_-}+\eps t^{-2+C\eps}\lra{q}^{1-\gamma_-})1_{q\leq 0}+(t^{-\gamma_-+C\eps}+\eps  t^{-2+C\eps})1_{q>0})\\&=O(\eps t^{-1+C\eps}\lra{q}^{-\gamma_{\sgn(q)}}+t^{-1-\gamma_-+C\eps}+(t^{-\gamma_-+C\eps}+\eps  t^{-2+C\eps})\lra{\max\{0,-q\}}^{1-\gamma_-}).}

Next we assume $r>2t$. Here we have proved that $0<q\sim r$ as long as $\eps\ll1$. It follows from Lemma \ref{lem4.2} that $|\nu|\lesssim (r+t)^{-\gamma+C\eps}+ \eps (t+r)^{-1+C\eps}$. By differentiating \eqref{nueqn}, whenever $r>2t$ and $t\geq  T_\eps$ we have
\eq{\label{nuqeqn} (\partial_t-\partial_r)\nu_r&=\partial_r(\partial_t+\partial_r)\mu=\mu_q\nu_r+\mu_{qq}q_r\nu+\eps t^{-1}\mu_{sq}q_r\\
&=\mu_q\nu_r+O(((r+t)^{-\gamma+C\eps}+ \eps (t+r)^{-1+C\eps})\lra{q}^{-2-\gamma_+}+\eps t^{-1+C\eps}\lra{q}^{-1-\gamma_+})\\
&=\mu_q\nu_r+O(\eps t^{-1+C\eps}(r+t)^{-1-\gamma_+}).}
To obtain the last estimate, we use $\gamma_+>1$ and $r+t\geq \eps^{-1}$. 
 Then,
\fm{&\int_{t}^{(r+t)/3}|(\partial_t-\partial_r)\nu_r-\mu_q\nu_r|\ d\tau\lesssim \int_{t}^{(r+t)/3}\eps\tau^{-1+C\eps}(r+t)^{-(1+\gamma_+)}\ d\tau\lesssim (r+t)^{-1-\gamma_++C\eps}.}
Also recall that in Lemma \ref{lem4.2} we proved that 
\fm{\int_t^{(r+t)/3}|\mu_q(\tau)|\ d\tau\lesssim1.}
Moreover, we have $z|_{\tau=(r+t)/3}>0$, so by \eqref{lem4.4c1} we have \fm{|\nu_r((r+t)/3,2(r+t)\omega/3)|&\lesssim\eps (r+t)^{-1-\gamma_++C\eps}+(r+t)^{-\gamma_-+C\eps}+ \eps (r+t)^{-2+C\eps}\\&\lesssim  (r+t)^{-\gamma_-+C\eps}+\eps (r+t)^{-2+C\eps}.}
By applying  Gronwall's inequality for $\tau\in[t,(r+t)/3]$, we conclude \eqref{lem4.4c2}.
\end{proof}\rm

\bigskip

Most of the estimates in the previous  lemmas  will still hold, if $Z^I$ is applied to the left hand sides for each multiindex $I$. To prove these results, we need some preliminary lemmas.

\lem{\label{lem4.5pre} Fix an integer $k\geq 0$  and a real constant $\sigma>0$. Fix $(t,x)\in\D$. Suppose that $Z^Iq=O(\lra{q}t^{C\eps})$ at $(t,x)$ whenever $|I|\leq k$. Also suppose that $F=F(s,q,\omega)$ is a smooth function defined for all $(s,q,\omega)\in[-\delta_0,\infty)\times\R\times\mathbb{S}^2$ such that at $(s,q,\omega)=(\eps\ln t-\delta,q(t,x),\omega)$, we have $\partial_s^a\partial_q^b\partial_\omega^c F=O_{a,b,c}(\lra{q}^{-b-\sigma}\exp(C_{a,b,c} s ))$ for each $a,b,c\geq 0$. Then, at $(t,x)$ we have
\fm{\sum_{|I|\leq k}|Z^IF(t,x)|\lesssim \lra{q}^{-\sigma}t^{C\eps}.}
On the right hand side, $\lra{q}$ can be replaced by $\lra{r-t}$.
}
\begin{proof}
It is clear that $|F|\lesssim \lra{q}^{-\sigma}t^{C\eps}$. Now, fix a multiindex $I$ with $0<|I|\leq k$. By the chain rule and the product rule, we  express $Z^IF$ as a linear combination of terms of the form
\eq{\label{lem4.5pref1} \partial_s^a\partial_q^b\partial_\omega^cF(\eps\ln t-\delta,q,\omega)\cdot\prod_{j=1}^aZ^{I_j}(\eps\ln t-\delta)\cdot \prod_{j=1}^bZ^{J_j}q\cdot\prod_{j=1}^cZ^{K_j}\omega.}
Here $a+b+c\leq k$, the sum of all $|I_*|,|J_*|,|K_*|$ is $|I|$, and each of $|I_*|,|J_*|,|K_*|$ is nonzero. By the assumptions, we conclude that each term of the form \eqref{lem4.5pref1} is controlled by \fm{\lra{q}^{-\sigma-b}t^{C\eps}\cdot \eps^a\cdot \lra{q}^bt^{C\eps}\lesssim \lra{q}^{-\sigma}t^{C\eps}.}
Here we use the fact that $r\sim t\gg1$ in $\D$. We conclude that $Z^IF=O(\lra{q}^{-\sigma}t^{C\eps})$. To replace $\lra{q}$ with $\lra{r-t}$, we apply \eqref{lem4.1c3} in Lemma \ref{lem4.1}.
\end{proof}

\lem{\label{lem4.5bou} For $(t,x)$ with $|x|=t/2$ and $t\geq  T_\eps$, we have $|Z^I(q-r+t)|\lesssim_I t^{1-\gamma_-+C_I\eps}$ and $|Z^I\nu|+\sum_i|Z^I\lambda_i|\lesssim t^{-\gamma_-+C_I\eps}$ for each multiindex $I$.}
\begin{proof}
Set $H:=\{(t,x):\ |x|=t/2,\ t\geq  T_\eps\}$. We first remark that once we prove $Z^I(q-r+t)=O(t^{1-\gamma_-+C\eps})$ for each $I$ on $H$, we can follow the proof of part (b) in Lemma \ref{c1lem6} to conclude that
\fm{|Z^I(\partial_t+\partial_r)q|+\sum_i|Z^I(\partial_i-\omega_i\partial_r)q|\lesssim t^{-\gamma_-+C\eps}.}
Besides, on $H$ we have
\eq{\label{lem4.5bouf1}|Z^I(q-r+t)|\lesssim \sum_{k\leq |I|}(r+t)^{k}|\partial^k(q-r+t)|\lesssim \sum_{k\leq |I|}t^{k}|\partial^k(q-r+t)|.}
So it suffices to prove $\partial^k(q-r+t)=O(t^{1-\gamma_--k+C\eps})$ for each $k\geq 0$ on $H$. 

Set $V_0:=\partial_t-\partial_r$ and $V_i=\partial_i+2\omega_i\partial_t$, so $V_i$ is tangent to $H$. For simplicity we use $V^k$ to denote an arbitrary product of $k$ vector fields $V_*$. We first claim the following inequality: for each function $F=F(t,x)$ and each integer $k\geq 0$, on $H$ we have
\eq{\label{lem4.5bouclaim}\sum_{l=0}^kt^{l}|V^lF|\sim\sum_{l=0}^kt^{l}|\partial^lF|\qquad \text{on }H.}
We prove by induction on $k$. If $k=0$, there is nothing to prove. If we have proved \eqref{lem4.5bouclaim} for all $k<k_0$, then on $H$ we have
\fm{\sum_{l=0}^{k_0}t^{l}|V^lF|&=|F|+\sum_{\alpha=0}^3\sum_{l=0}^{k_0-1}t^{l+1}|V^{l}V_{\alpha}F|\lesssim |F|+\sum_{\alpha=0}^3\sum_{l=0}^{k_0-1}t^{l+1}|\partial^{l}V_{\alpha}F|\\
&\lesssim |F|+\sum_{l=0}^{k_0-1}t^{l+1}(|\partial^{l+1}F|+|\partial^{l}(\omega\cdot\partial F)|)\\
&\lesssim |F|+\sum_{l=0}^{k_0-1}t^{l+1}(|\partial^{l+1}F|+\sum_{l'=0}^l|\partial^{l-l'}\omega\cdot\partial^{l'+1} F|)\\&\lesssim |F|+\sum_{l=0}^{k_0-1}t^{l+1}(|\partial^{l+1}F|+\sum_{l'=0}^lr^{l'-l}|\partial^{l'+1} F|)\lesssim \sum_{l=0}^{k_0}t^{l}|\partial^l F|.}
The other direction of this inequality can be proved similarly; we simply notice that \fm{\partial_\alpha=f_0(\partial_t-\partial_r)+\sum_if_0(\partial_i+2\omega_i\partial_t)}
where $f_0$ denotes an arbitrary polynomial of $\omega$. This finishes the proof of the claim.

Next we claim that $V^k(q-r+t)=O(t^{1-\gamma_--k+C\eps})$ on $H$ for each $k\geq 0$. We prove this by induction. If $k=0$, there is nothing to prove since $q=r-t$ on $H$. Now suppose $k>0$ and suppose we have proved $V^l(q-r+t)=O(t^{1-\gamma_--l+C\eps})$ on $H$ for each $l<k$. By \eqref{lem4.5bouclaim} and the induction hypotheses, on $H$ we have \fm{\sum_{l=0}^{k-1}t^{l}|\partial^l(q-r+t)|\lesssim \sum_{l=0}^{k-1}t^{l}\cdot t^{1-\gamma_--l+C\eps}\lesssim t^{1-\gamma_-+C\eps}.}
Note that this and \eqref{lem4.5bouf1} imply that $|Z^J(q-r+t)|\lesssim t^{1-\gamma_-+C\eps}$ and thus $|Z^Jq|\lesssim \lra{r-t}t^{C\eps}\lesssim \lra{q}t^{C\eps}$ on $H$ whenever $|J|<k$, so we are ready to apply Lemma \ref{lem4.5pre} near $H$.
To prove the claim, we also need to induct on $n_0(V^k)$ where $n_0(V^k)$ is the number of $V_0$ in $V^k$. If $n_0(V^k)=0$, then $V^k$ only contains vector fields which are tangent to $H$. This implies that $V^k(q-r+t)=V^k(0)=0$ on $H$. In general, if $n_0(V^k)>0$, then we can write $V^k=V^{k_1}V_0V^{k_2}$ where $k_1+k_2=k-1$ and $n_0(V^{k_2})=0$. Since $[V_0,V_i]=\partial\omega\cdot\partial=f_0\partial \omega\cdot V$ for $i=1,2,3$, we have
\fm{|V^k(q-r+t)|&\lesssim|V^{k_1}V^{k_2}V_0(q-r+t)|+\sum_{V^{k_2}=V^{l_1}V_jV^{l_2}}|V^{k_1}V^{l_1}[V_0,V_j]V^{l_2}(q-r+t)|\\
&\lesssim |V^{k_1}V^{k_2}(\mu+2)|+\sum_{V^{k_2}=V^{l_1}V_jV^{l_2}}|V^{k_1}V^{l_1}(f_0\partial \omega \cdot V V^{l_2}(q-r+t))|\\
&\lesssim|V^{k_1}V^{k_2}(\mu+2)|+\sum_{ l\leq k-2}|V^{l}(f_0\partial \omega)  V^{k-1-l}(q-r+t))|. }
Note that $V_0(q-r+t)=\mu+2$ everywhere in $\Omega$ (not only on $H$). Then, by \eqref{lem4.5bouclaim}, Lemma \ref{lem4.5pre}, Lemma \ref{c1l2.1} and \eqref{def3.1a6}, on $H$ we have
\fm{|V^{k_1}V^{k_2}(\mu+2)|&\lesssim t^{1-k}\sum_{l=0}^{k-1}t^l|\partial^l(\mu+2)|\lesssim t^{1-k}\sum_{|J|<k}t^{|J|}\lra{r-t}^{-|J|}|Z^J(\mu+2)|\\
&\lesssim t^{1-k}\lra{r-t}^{-\gamma_-}t^{C\eps}\lesssim t^{1-\gamma_--k+C\eps}.}
Moreover, by the induction hypotheses, \eqref{lem4.5bouclaim} and  Leibniz's rule, on $H$ we have
\fm{|V^l(f_0\partial\omega)V^{k-l-1}(q-r+t)|&\lesssim t^{-l}\sum_{l'\leq l}t^{l'}|\partial^{l'}(f_0\partial\omega)|\cdot t^{1-\gamma_--(k-l-1)+C\eps}\\
&\lesssim t^{-l}\sum_{l'\leq l}t^{l'}\cdot t^{-1-l'}\cdot t^{1-\gamma_--(k-l-1)+C\eps}\lesssim t^{1-\gamma_--k+C\eps}.}
We thus finish the induction  and conclude that $V^k(q-r+t)=O(t^{1-\gamma_--k+C\eps})$ on $H$ for each $k\geq 0$. By \eqref{lem4.5bouclaim}, we have $\partial^k(q-r+t)=O(t^{1-\gamma_--k+C\eps})$ on $H$ for each $k\geq 0$.
\end{proof}
\rm

The next lemma is a key lemma in this subsection. We recall the notation $\eps^nS^{m,p}=\eps^nS^{m,p}_\D$ from  Definition \ref{c1defn1.5} with $\D$ defined by \eqref{c3defnD}. Since the estimates in the region where $q\geq 0$ are usually expected to be different from those in the region where $q\leq 0$, we also define
\eq{\eps^nS^{m,p}_+=\eps^nS^{m,p}_{\D\cap\{q\geq 0\}},\qquad \eps^nS^{m,p}_-=\eps^nS^{m,p}_{\D\cap\{q\leq 0\}}.}
Because of Remark \ref{rmk4.1.1}, it is okay to replace $q\geq 0$ (or $q\leq 0$) with $r-t\geq 0$ (or $r-t\leq 0$) in this definition.

\lem{\label{lem4.5} {\rm (a)} We have $q-r+t\in S^{0,1-\gamma_-}_-\cap S^{0,0}_+$.  That is, for all $(t,x)\in\D$ and for all multiindices $I$, we have
\begin{equation}\label{lq31}|Z^I(q-r+t)|\lesssim_I t^{C_I\eps}\lra{\max\{0,t-r\}}^{1-\gamma_-},\end{equation}
As a result, we have $q\in S^{0,1}$, $\partial q\in S^{0,0}$ and $\Omega_{kk'}q\in S^{0,1-\gamma_-}_-\cap S^{0,0}_+$ for each $1\leq k,k'\leq 3$.

{\rm (b)}  We have $\mu\in S^{0,0}$, $\partial_s^a\partial_q^b\partial_\omega^c(\mu+2)\in S^{0,-b-\gamma_-}_-\cap S^{0,-b-\gamma_+}_+$ for $a,b,c\geq 0$; $\partial_s^b\partial_\omega^cU\in S^{0,1-\gamma_-}_-\cap S^{0,0}_+$ for all $b,c\geq 0$, and $\partial_s^a\partial_q^b\partial_\omega^cU_q\in S^{0,-b-\gamma_-}_-\cap S^{0,-b-\gamma_+}_+$ for all $a,b,c\geq 0$. Here all the functions are of $(s,q,\omega)=(\eps\ln t-\delta,q(t,x),\omega)$ defined in $\D$. 

{\rm (c)} We have $\nu\in  S^{-1,1-\gamma_-}_-\cap  S^{-1,0}_+$, $\nu_q\in  S^{-1,-\gamma_-}_-\cap( S^{-1,-\gamma_+}_++ S^{-2,0}_+)$, $\lambda_i\in S^{-1,1-\gamma_-}_-\cap S^{-1,0}_+$, and 
\eq{\label{lem4.5c3}\nu-\frac{\eps G_{2,J}(\omega)}{4 r}\mu U^{(J)}+\frac{\eps G_{3,J}(\omega)}{8r}\mu^2U_q^{(J)}&\in S^{-2,2-\gamma_-}_-\cap  S^{-2,1}_+,}
\eq{\label{lem4.5c4}\nu_q-\frac{\eps G_{2,J}(\omega)}{4 r}\partial_q(\mu U^{(J)})+\frac{\eps G_{3,J}(\omega)}{8r}\partial_q(\mu^2U_q^{(J)})&\in S^{-2,1-\gamma_-}_-\cap  S^{-2,0}_+.}
Here all the functions are of $(s,q,\omega)=(\eps\ln t-\delta,q(t,x),\omega)$ defined in $\D$.}
\begin{proof}
In this proof, we always assume that $|r-t|\leq t/2$ and $t\geq  T_\eps$. As a result, we have estimates such as $|Z^I\omega|\lesssim 1$, $(r+t)\sim t$, etc.

(a) First, we recall from Remark \ref{rmk4.1.1} that \eq{\label{lem4.5claim}t^{-C\eps}\lra{\max\{0,-q\}}\lesssim \lra{\max\{0,t-r\}}\lesssim t^{C\eps}\lra{\max\{0,-q\}},\qquad \text{in }\D;}
\eq{\label{lem4.5claim2}t^{-C\eps}\lra{\max\{0,q\}}\lesssim \lra{\max\{0,r-t\}}\lesssim t^{C\eps}\lra{\max\{0,q\}},\qquad \text{in }\D.}

We now prove \eqref{lq31} by induction on $|I|$. The case $|I|=0$ has been proved in  Lemma \ref{lem4.1}.  In general, we fix an integer $k\geq 0$ and  suppose that \eqref{lq31} holds for all $|I|\leq k$. Since $Z^I(r-t)=O(\lra{r-t})$, we also have $|Z^Iq|\lesssim \lra{r-t}t^{C\eps}\lesssim \lra{q}t^{C\eps}$ for all $|I|\leq k$. It then follows  from Lemma \ref{lem4.5pre} that \fm{|Z^I(\mu+2)|\lesssim \lra{q}^{-\gamma_{\sgn(q)}}t^{C\eps},\qquad\forall (t,x)\in\D,\ |I|\leq k.}  Now fix a multiindex $I$ with $|I|=k+1$. By the chain rule and Leibniz's rule, we express $Z^I \mu $ as a linear combination of terms of the form\eq{\label{lem4.5f1}  (\partial_s^a\partial_q^b\partial_\omega^c\mu)\cdot\prod_{j=1}^a Z^{I_j}(\eps\ln t-\delta)\cdot \prod_{j=1}^bZ^{J_j}q\cdot \prod_{j=1}^cZ^{K_j}\omega}
where $a+b+c>0$, $|I_*|,|J_*|,|K_{*}|$ are all nonzero, and the sum of all these multiindices is $k+1$. Note that the only term with some $|J_*|>k$ is $\mu_qZ^Iq$. All the other terms can be controlled in the exact same way as in the proof of Lemma \ref{lem4.5pre}. We thus conclude that \fm{Z^I\mu=\mu_qZ^Iq+O(\lra{q}^{-\gamma_{\sgn(q)}}t^{C\eps})=\mu_qZ^Iq+O(\lra{q}^{-\gamma_{\sgn(q)}}t^{C\eps}|\mu|).}

We now apply \eqref{lemtrcomf1} in Lemma \ref{c1lemtrcom} to $F=q-r+t$. If we use $f_0$ to denote an arbitrary polynomial of $\{Z^J\omega\}$, then 
\eq{\label{lem4.5f2}&|(\partial_t-\partial_r)Z^I(q-r+t)|\\
&\leq |Z^I(\mu+2)|+\sum_{|J|\leq k}|f_0Z^J(\mu+2)+\sum_i f_0(\partial_i+\omega_i\partial_t)Z^J(q-r+t)|\\
&\lesssim |\mu_qZ^Iq|+ \lra{q}^{-\gamma_{\sgn(q)}}t^{C\eps}|\mu| +\sum_{|J|\leq k} \sum_i |(\partial_i-\omega_i\partial_r)Z^J(q-r+t)+\omega_i(\partial_t+\partial_r)Z^J(q-r+t)|\\
&\lesssim |\mu_qZ^Iq|+ \lra{q}^{-\gamma_{\sgn(q)}}t^{C\eps}|\mu| +(r+t)^{-1}\sum_{1\leq |J|\leq k+1}  |Z^J(q-r+t)|\\
&\lesssim|\mu_qZ^Iq|+ \lra{q}^{-\gamma_{\sgn(q)}}t^{C\eps}|\mu|+(r+t)^{-1}\lra{\max\{0,-q\}}^{1-\gamma_-}t^{C\eps} +(r+t)^{-1}\sum_{ |J|=k+1}  |Z^J(q-r+t)|.}
In the second estimate, we notice that $Z^I(2)=0$ because $|I|=k+1>0$. In the second last estimate, we apply Lemma \ref{c1l2.1}. In the last estimate, we apply the induction hypotheses and  \eqref{lem4.5claim}.

Now we fix $(t,x)\in\D$.   By integrating \eqref{lem4.5f2} along the characteristic $(\tau,(r+t-\tau)\omega)$ for $t\leq \tau\leq 2(r+t)/3$ and summing over $|I|=k+1$, we have
\fm{&\sum_{|I|=k+1}|(Z^I(q-r+t))|_{\tau=2(r+t)/3}-(Z^I(q-r+t))|_{\tau=t}|\\
&\lesssim\int_t^{2(r+t)/3} \sum_{|I|=k+1}|\mu_qZ^Iq|(\tau)+\lra{z(\tau)}^{-\gamma_{\sgn(z(\tau))}}\tau^{C\eps}(-\dot{z}(\tau))\\&\qquad+(r+t)^{-1}\lra{\max\{0,-z(\tau)\}}^{1-\gamma_-}\tau^{C\eps}\cdot(-\tau^{C\eps}\dot{z}(\tau))+(r+t)^{-1}\sum_{|J|=k+1}|Z^J(q-r+t)|(\tau)\ d\tau\\
&\lesssim \int_t^{2(r+t)/3} \sum_{|I|=k+1}|\mu_qZ^I(q-r+t)|(\tau)+(r+t)^{-1}\sum_{|J|=k+1}|Z^J(q-r+t)|(\tau)\ d\tau\\&\quad+\int_{-(r+t)/3}^{q(t,x)} (r+t)^{C\eps}\lra{\rho}^{-\gamma_{\sgn(\rho)}}+(r+t)^{-1+C\eps}\lra{\max\{0,-\rho\}}^{1-\gamma_-}\ d\rho\\
&\lesssim\int_t^{2(r+t)/3} \sum_{|I|=k+1}|\mu_qZ^I(q-r+t)|(\tau)+(r+t)^{-1}\sum_{|J|=k+1}|Z^J(q-r+t)|(\tau)\ d\tau\\&\quad+(r+t)^{C\eps}\lra{\max\{0,-q\}}^{1-\gamma_-}+ (r+t)^{-1+C\eps}(\lra{q}1_{q>0}+\lra{q}^{2-\gamma_-}1_{q\leq 0}).}
In the last estimate, we use $\gamma_->2$. Also note that $\lra{q}^{-1}\gtrsim \lra{r-t}^{-1}t^{-C\eps}\gtrsim (r+t)^{-1}t^{-C\eps}$, so $(r+t)^{-1+C\eps}(\lra{q}1_{q> 0}+\lra{q}^{2-\gamma_-}1_{q\leq 0})\lesssim (r+t)^{C\eps}\lra{\max\{0,-q\}}^{1-\gamma_-}$.
By Lemma \ref{lem4.5bou}, we have $\sum_{|I|=k+1}|(Z^I(q-r+t))|_{\tau=2(r+t)/3}|\lesssim (r+t)^{1-\gamma_-+C\eps}$. By \eqref{lem4.2f1} we have 
\fm{\int_t^{2(r+t)/3}(|\mu_q|+(r+t)^{-1})\ d\tau\lesssim \eps\ln (t+r)+1,}
so we conclude by  Gronwall's inequality that 
\fm{\sum_{|I|=k+1}|(Z^I(q-r+t))(t,x)|&\lesssim (r+t)^{1-\gamma_-+C\eps}+(r+t)^{C\eps}\lra{\max\{0,-q\}}^{1-\gamma_-}\\&\lesssim t^{C\eps}\lra{\max\{0,-q\}}^{1-\gamma_-}\lesssim t^{C\eps}\lra{\max\{0,t-r\}}^{1-\gamma_-}.}
This finishes the proof of \eqref{lq31}.  Since  $r-t\in S^{0,1}$ (see Example \ref{c1exm1.7}) and $q-r+t\in S^{0,1-\gamma_-}_-\cap S^{0,0}_+$, we have $q\in S^{0,1}$. Then it follows from Lemma \ref{c1lem6} that $\partial q\in S^{0,0}$ and that  $\Omega_{kk'}q=\Omega_{kk'}(q-r+t)\in S^{0,1-\gamma_-}_-\cap S^{0,0}_+$.

(b) All these results  follow from the assumptions \eqref{def3.1a6}--\eqref{def3.1a8} on $(\mu,U)$ and Lemma \ref{lem4.5pre}. 

(c) In (a) we have proved that $\Omega_{kk'}q\in S^{0,1-\gamma_-}_-\cap S^{0,0}_+$ for each $1\leq k<k'\leq 3$. Thus,
\fm{\lambda_i&=r^{-1}\sum_{j}\omega_j\Omega_{ji}q\in S^{-1,0}\cdot S^{0,0}\cdot (S^{0,1-\gamma_-}_-\cap S^{0,0}_+)\subset S^{-1,1-\gamma_-}_-\cap S^{-1,0}_+.}
Here we recall from Example \ref{c1exm1.7} that $r^{-1}\in S^{-1,0}$ and $\omega\in S^{0,0}$.

Next, we set \fm{V:=\nu-\frac{\eps G_{2,J}(\omega)}{4 r}\mu U^{(J)}+\frac{\eps G_{3,J}(\omega)}{8r}\mu^2U_q^{(J)}.}
By \eqref{lem4.3f1} and part (b), we have
\fm{V_t-V_r&=\mu_qV+\eps (tr)^{-1}(r-t)\mu_s+\frac{\eps G_{2,J}(\omega)}{4 r^2}\mu U^{(J)}-\frac{\eps G_{3,J}(\omega)}{8r^2} \mu^2U_q^{(J)} \\
&\quad-\frac{\eps^2 G_{2,J}(\omega)}{4 rt}\partial_s(\mu U^{(J)})+\frac{\eps^2 G_{3,J}(\omega)}{8rt}\partial_s(\mu^2U_q^{(J)})\\
&=\mu_qV+\eps S^{-2,1}\cdot (S^{0,-\gamma_-}_-\cap S^{0,-\gamma_+}_+)+\eps S^{-2,0}\cdot S^{0,0}\cdot(S^{0,1-\gamma_-}_-\cap S^{0,0}_+)\\&\quad+\eps S^{-2,0}\cdot S^{0,0}\cdot (S^{0,-\gamma_-}_-\cap S^{0,-\gamma_+}_+)\\
&=\mu_qV+\eps S^{-2,1-\gamma_-}_-\cap \eps S^{-2,0}_+.}
To prove \eqref{lem4.5c3}, we need to prove \fm{|Z^IV|\lesssim   t^{-2+C\eps}\lra{r-t}\lra{\max\{0,t-r\}}^{1-\gamma_-}.}
Here $\lra{r-t}$ and $\max\{0,t-r\}$ can be replaced  with $\lra{q}$ and $\max\{0,-q\}$ because of \eqref{lem4.1c3} and \eqref{lem4.5claim}. We induct on $|I|$. The base case is proved in   Lemma \ref{lem4.3}. Next we assume  $|I|>0$ and that this estimate for $Z^JV$ has been proved for each $|J|<|I|$. 
Then, by \eqref{lemtrcomf1} and our induction hypotheses, in $\D$ we have
\fm{&|(\partial_t-\partial_r)Z^IV|\\&\lesssim |Z^I(V_t-V_r)|+\sum_{|J|<|I|}[|Z^J(V_t-V_r)|+\sum_i |(\partial_i+\omega_i\partial_t)Z^JV|]\\&\lesssim |Z^I(\mu_qV)|+\sum_{|J|<|I|}|Z^J(\mu_qV)|+\sum_{1\leq |J|\leq |I|} (r+t)^{-1}|Z^JV|]+\eps t^{-2+C\eps}\lra{\max\{0,t-r\}}^{1-\gamma_-}\\
&\lesssim |\mu_qZ^IV|+\sum_{|J|=|I|}(r+t)^{-1}|Z^JV|+[\sum_{|J|\leq |I|}|Z^J\mu_q|]\cdot[\sum_{|J|<|I|}|Z^JV|]\\&\quad+ t^{-3+C\eps}\lra{r-t}\lra{\max\{0,t-r\}}^{1-\gamma_-}+\eps t^{-2+C\eps}\lra{\max\{0,t-r\}}^{1-\gamma_-}\\
&\lesssim |\mu_qZ^IV|+\sum_{|J|=|I|}(r+t)^{-1}|Z^JV|+t^{C\eps}\lra{q}^{-1-\gamma_{\sgn(q)}}\cdot t^{-2+C\eps}\lra{r-t}\lra{\max\{0,t-r\}}^{1-\gamma_-}\\&\quad+  t^{-2+C\eps}\lra{\max\{0,t-r\}}^{1-\gamma_-}\\
&\lesssim |\mu_qZ^IV|+\sum_{|J|=|I|}(r+t)^{-1}|Z^JV|+ t^{-2+C\eps}\lra{\max\{0,-q\}}^{1-\gamma_-}|\mu|.}
In the second last step, we apply $\mu_q\in S^{0,-1-\gamma_-}_-\cap S^{0,-1-\gamma_+}_+$ and \eqref{lem4.1c3}. In the last step, we apply \eqref{lem4.5claim}, \eqref{lem4.5claim2}, and the estimate $1\lesssim t^{C\eps}|\mu|$.   We now integrate this estimate along the characteristic $(\tau,(r+t-\tau)\omega)$ for $t\leq\tau\leq 2(r+t)/3$. Note that 
\fm{&\int_{t}^{2(r+t)/3} \tau^{-2+C\eps}\lra{\max\{0,-z(\tau)\}}^{1-\gamma_-}|\mu|\ d\tau\\
&\lesssim  \int_{-(r+t)/3}^{q(t,x)}t^{-2+C\eps}\lra{\max\{0,-\rho\}}^{1-\gamma_-}\ d\rho\lesssim t^{-2+C\eps}\lra{q}\lra{\max\{0,-q\}}^{1-\gamma_-},
}
that
\fm{\int_t^{2(r+t)/3}|\mu_q|\ d\tau\lesssim \eps\ln (t+r)+1,}
and that on $H$ we have \fm{|Z^IV|&\lesssim|Z^I\nu|+|Z^I(\eps S^{-1,1-\gamma_-}_-+\eps S^{-1,-\gamma_-}_-)|\lesssim t^{-\gamma_-+C\eps}+\eps t^{-1+C\eps}\lra{r-t}^{1-\gamma_-}\lesssim  t^{-\gamma_-+C\eps}} by Lemma \ref{lem4.5bou}. Thus, by  Gronwall's inequality, we have
\fm{|Z^IV|\lesssim  t^{-2+C\eps}\lra{q}\lra{\max\{0,-q\}}^{1-\gamma_-}+t^{-\gamma_-+C\eps}\lesssim t^{-2+C\eps}\lra{q}\lra{\max\{0,-q\}}^{1-\gamma_-},\qquad \text{in }\D.}
We thus conclude \eqref{lem4.5c3}.

Next we have
\fm{&(\partial_t-\partial_r)V_r\\&=\mu_qV_r+\mu_{qq}q_rV+\eps r^{-2}\mu_s+\eps (rt)^{-1}(r-t)\mu_{sq}q_r -\frac{\eps G_{2,J}(\omega)}{2 r^3}\mu U^{(J)}+\frac{\eps G_{3,J}(\omega)}{2r^3} \mu^2U_q^{(J)}\\
&\quad+\frac{\eps G_{2,J}(\omega)}{4 r^2}q_r\partial_q(\mu U^{(J)})-\frac{\eps G_{3,J}(\omega)}{8r^2} q_r\partial_q(\mu^2U_q^{(J)}) +\frac{\eps^2 G_{2,J}(\omega)}{4 r^2t}\partial_s(\mu U^{(J)})-\frac{\eps^2 G_{3,J}(\omega)}{8r^2t}\partial_s(\mu^2U_q^{(J)})\\
&\quad-\frac{\eps^2 G_{2,J}(\omega)}{4 rt}q_r\partial_q\partial_s(\mu U^{(J)})+\frac{\eps^2 G_{3,J}(\omega)}{8rt}q_r\partial_q\partial_s(\mu^2U_q^{(J)})\\
&=\mu_qV_r+(S^{0,-2-\gamma_-}_-\cap S^{0,-2-\gamma_+}_+)\cdot ( S^{-2,2-\gamma_-}_-\cap  S^{-2,1}_+)+\eps S^{-2,-\gamma_-}_-\cap \eps S^{-2,-\gamma_+}_+\\
&\quad+\eps S^{-3,1-\gamma_-}_-\cap \eps S^{-3,0}_+\\
&=\mu_qV_r+S^{-2,-\gamma_-}_-\cap (S^{-2,-\gamma_+}_++ S^{-3,0}_+).}
To prove \eqref{lem4.5c4}, we first prove it with $\partial_q$ replaced by $\partial_r$. That is, we need to prove 
\fm{|Z^IV_r|\lesssim  t^{-2+C\eps}\lra{\max\{0,t-r\}}^{1-\gamma_-}.}
Similarly we can also replace $\lra{\max\{0,t-r\}}$ with $\lra{\max\{0,-q\}}$. The base case is proved in Lemma \ref{lem4.3}. Next we assume  $|I|>0$ and that this estimate for $Z^JV_r$ has been proved for each $|J|<|I|$. 
Then, by \eqref{lemtrcomf1} and our induction hypotheses, in $\D$ we have
\fm{&|(\partial_t-\partial_r)Z^IV_r|\\&\lesssim |Z^I(\partial_t-\partial_r)V_r|+\sum_{|J|<|I|}[|Z^J(\partial_t-\partial_r)V_r|+\sum_i |(\partial_i+\omega_i\partial_t)Z^JV_r|]\\&\lesssim |Z^I(\mu_qV_r)|+\sum_{|J|<|I|}|Z^J(\mu_qV_r)|+\sum_{1\leq |J|\leq |I|} (r+t)^{-1}|Z^JV_r|]\\&\quad+ t^{-2+C\eps}\lra{r-t}^{-\gamma_{\sgn(q)}}+ t^{-3+C\eps}1_{q\geq 0}\\
&\lesssim |\mu_qZ^IV_r|+\sum_{|J|=|I|}(r+t)^{-1}|Z^JV_r|+[\sum_{|J|\leq |I|}|Z^J\mu_q|]\cdot[\sum_{|J|<|I|}|Z^JV_r|]\\&\quad+ t^{-2+C\eps}\lra{r-t}^{-\gamma_{\sgn(q)}}+ t^{-3+C\eps}1_{q\geq 0}+t^{-3+C\eps}\lra{\max\{0,-q\}}^{1-\gamma_-} \\
&\lesssim |\mu_qZ^IV|+\sum_{|J|=|I|}(r+t)^{-1}|Z^JV|+t^{C\eps}\lra{q}^{-1-\gamma_{\sgn(q)}}\cdot t^{-2+C\eps}\lra{\max\{0,t-r\}}^{1-\gamma_-}\\&\quad+ t^{-2+C\eps}\lra{r-t}^{-\gamma_{\sgn(q)}}+ t^{-3+C\eps}1_{q\geq 0}\\
&\lesssim |\mu_qZ^IV|+\sum_{|J|=|I|}(r+t)^{-1}|Z^JV|+ t^{-2+C\eps}\lra{r-t}^{-\gamma_{\sgn(q)}}+ t^{-3+C\eps}1_{q\geq 0}.}
Here $\lra{r-t}^{-\gamma_{\sgn(q)}}$ can be replaced by $\lra{q}^{-\gamma_{\sgn(q)}}|\mu|t^{C\eps}$.
Now integrate this estimate along the characteristic for $\tau\in[t,2(r+t)/3]$. We have
\fm{&\int_t^{2(r+t)/3} \tau^{-2+C\eps}\lra{z(\tau)}^{-\gamma_{\sgn(z(\tau))}}|\mu|+ \tau^{-3+C\eps}1_{z(\tau)\geq 0} \ d\tau\\
&\lesssim \int_{-(r+t)/3}^{q(t,x)}( t^{-2+C\eps}\lra{\rho}^{-\gamma_{\sgn(\rho)}}+\eps t^{-3+C\eps}1_{\rho\geq 0})\ d\rho\\
&\lesssim  t^{-2+C\eps}\lra{\max\{0,-q\}}^{1-\gamma_-}+ t^{-3+C\eps}\lra{q}1_{q\geq 0}\lesssim t^{-2+C\eps}\lra{\max\{0,-q\}}^{1-\gamma_-}.}
On $H$, we have $|Z^IV_r|\lesssim t^{-1-\gamma_-+C\eps}$ by Lemma \ref{lem4.5bou} and \eqref{c1l2.1}. By applying  Gronwall's inequality, we conclude that $V_r\in S^{-\gamma_--1,0}+( S_-^{-2,1-\gamma_-}\cap  S^{-2,0}_+)\subset S_-^{-2,1-\gamma_-}\cap S^{-2,0}_+$.  Note that
\fm{V_r&=\nu_r-\frac{\eps G_{2,J}}{4r}\partial_r(\mu U^{(J)})+\frac{\eps G_{3,J}}{8r}\partial_r(\mu^2U_q^{(J)})+\frac{\eps G_{2,J}}{4r^2}\mu U^{(J)}-\frac{\eps G_{3,J}}{8r^2}\mu^2 U_q^{(J)}\\
&=q_r(\nu_q-\frac{\eps G_{2,J}}{4r}\partial_q(\mu U^{(J)})+\frac{\eps G_{3,J}}{8r}\partial_q(\mu^2U_q^{(J)}))\mod \eps S^{-2,1-\gamma_-}_-\cap \eps S^{-2,0}_+,}
so 
\fm{q_r(\nu_q-\frac{\eps G_{2,J}}{4r}\partial_q(\mu U^{(J)})+\frac{\eps G_{3,J}}{8r}\partial_q(\mu^2U_q^{(J)}))\in S_-^{-2,1-\gamma_-}\cap S^{-2,0}_+.}
Now, $q_r\in S^{0,0}$. Since $q_r\geq C^{-1}t^{-C\eps}$, we also have $q_r^{-1}\in S^{0,0}$. This easily follows from the fact that $Z^I(q_r^{-1})$ can be written as a linear combination of terms of the form
\fm{q_r^{-1-m}Z^{I_1}q_r\cdots Z^{I_m}q_r,\hspace{2em}\sum|I_*|=|I|.}  Then \eqref{lem4.5c4} follows.

Finally, by applying the results in (b), we have 
\fm{\frac{\eps G_{2,J}(\omega)}{4 r}\mu U^{(J)}-\frac{\eps G_{3,J}(\omega)}{8r}\mu^2U_q^{(J)}&\in \eps S^{-1,1-\gamma_-}_-\cap \eps S^{-1,0}_+,\\
\frac{\eps G_{2,J}(\omega)}{4 r}\partial_q(\mu U^{(J)})-\frac{\eps G_{3,J}(\omega)}{8r}\partial_q(\mu^2U_q^{(J)})&\in \eps S^{-1,-\gamma_-}_-\cap\eps S^{-1,-\gamma_+}_+.}
As a result, \fm{\nu&\in(\eps S^{-1,1-\gamma_-}_-+ S^{-2,2-\gamma_-}_-)\cap(\eps S^{-1,0}_++ S^{-2,1}_+)\subset  S^{-1,1-\gamma_-}_-\cap  S^{-1,0}_+ ,\\ 
\nu_q&\in (\eps S^{-1,-\gamma_-}_-+ S^{-2,1-\gamma_-}_-)\cap(\eps S^{-1,-\gamma_+}_++ S^{-2,0}_+)\\
&\subset  S^{-1,-\gamma_-}_-\cap( S^{-1,-\gamma_+}_++ S^{-2,0}_+)} Here we notice that $\lra{r-t}\lesssim  t$ in $\D$, so $\eps S^{m,p}\subset \eps S^{m+1,p-1}$.
\end{proof}\rm

The following proposition states that $q$ is an approximate optical function.

\prop{\label{prop4.6} We have $g^{\alpha\beta}(\eps r^{-1}U,\partial(\eps r^{-1}U))q_\alpha q_\beta\in S^{-2,2-\gamma_-}_-\cap S_+^{-2,1}$.}
\begin{proof}Fix a constant symmetric matrix $(c_0^{\alpha\beta})\in\R^{4\times 4}$. Since $\nu,\lambda_i\in S^{-1,1-\gamma_-}_-\cap S^{-1,0}_+
$, 
\eq{\label{lem4.6f1}&c_0^{\alpha\beta}q_\alpha q_\beta\\&=c^{00}_0(\frac{\mu+\nu}{2})^2+2c_0^{0i}(\frac{\mu+\nu}{2})(\lambda_i+\frac{\omega_i(\nu-\mu)}{2})+c_0^{ij}(\lambda_i+\frac{\omega_i(\nu-\mu)}{2})(\lambda_j+\frac{\omega_j(\nu-\mu)}{2})\\&=\frac{1}{4}c_0^{\alpha\beta}\wh{\omega}_\alpha\wh{\omega}_\beta\mu^2+\frac{1}{2}(c_0^{00}-c_0^{ij}\omega_i\omega_j)\mu\nu+(c_0^{0i}-c_0^{ij}\omega_j)\mu\lambda_i+\frac{1}{4}c^{00}_0\nu^2+\frac{1}{2}c_0^{0i}\nu(2\lambda_i+\omega_i\nu)\\&\hspace{1em}+\frac{1}{4}g_0^{ij}(2\lambda_i+\omega_i\nu)(2\lambda_j+\omega_j\nu)\\
&=\frac{1}{4}c_0^{\alpha\beta}\wh{\omega}_\alpha\wh{\omega}_\beta\mu^2+\frac{1}{2}(c_0^{00}-c_0^{ij}\omega_i\omega_j)\mu\nu+(c_0^{0i}-c_0^{ij}\omega_j)\mu\lambda_i\mod S^{-2,2-2\gamma_-}_-\cap S^{-2,0}_+\\
&=\frac{1}{4}c_0^{\alpha\beta}\wh{\omega}_\alpha\wh{\omega}_\beta\mu^2\mod S^{-1,1-\gamma_-}_-\cap S^{-1,0}_+.}
By taking $(c_0^{\alpha\beta})=(m^{\alpha\beta})$, we have
\eq{\label{lem4.6f2}m^{\alpha\beta}q_\alpha q_\beta=-\mu\nu+\mu\sum_i\omega_i\lambda_i\mod S^{-2,2-2\gamma_-}_-\cap S^{-2,0}_+=-\mu\nu\mod S^{-2,2-2\gamma_-}_-\cap S^{-2,0}_+.}
Here we note that $\sum_i\omega_i\lambda_i=\sum_i\omega_i(q_i-\omega_iq_r)=0$.

Note that $\eps r^{-1}U,\partial (\eps r^{-1}U)\in \eps S^{-1,1-\gamma_-}_-\cap \eps S^{-1,0}_+$. By Lemma \ref{lem2.9}, we have 
\eq{\label{lem4.6f3}g^{\alpha\beta}(\eps r^{-1}U,\partial(\eps r^{-1}U))-m^{\alpha\beta}-g^{\alpha\beta}_J\cdot\eps r^{-1}U^{(J)}-g^{\alpha\beta\lambda}_{J} \partial_\lambda(\eps r^{-1}U^{(J)})\in \eps^2 S^{-2,2-2\gamma_-}_-\cap \eps^2 S^{-2,0}_+.}
And since $\partial q\in S^{0,0}$, it suffices to prove 
\fm{\ [m^{\alpha\beta}+g^{\alpha\beta}_J\cdot\eps r^{-1}U^{(J)}+g^{\alpha\beta\lambda}_{J} \partial_\lambda(\eps r^{-1}U^{(J)})]q_\alpha q_\beta\in S^{-2,2-\gamma_-}_-\cap S_+^{-2,1}. }
By \eqref{lem4.6f1} and \eqref{lem4.6f2} we have
\fm{&[m^{\alpha\beta}+g^{\alpha\beta}_J\cdot\eps r^{-1}U^{(J)}+g^{\alpha\beta\lambda}_{J} \partial_\lambda(\eps r^{-1}U^{(J)})]q_\alpha q_\beta\\
&=-\mu\nu+\frac{1}{4}G_{2,J}\mu^2\cdot\eps r^{-1}U^{(J)}+\frac{1}{4}g^{\alpha\beta\lambda}_J\wh{\omega}_\alpha\wh{\omega}_\beta \mu^2\cdot\partial_\lambda(\eps r^{-1}U^{(J)})\mod S^{-2,2-2\gamma_-}_-\cap S^{-2,0}_+\\
&=-\mu\nu+\frac{1}{4}G_{2,J}\mu^2\cdot\eps r^{-1}U^{(J)}+\frac{1}{4}g^{\alpha\beta\lambda}_J\wh{\omega}_\alpha\wh{\omega}_\beta \mu^2(-\eps r^{-2}(\partial_\lambda r)U^{(J)}+\eps r^{-1}\partial_\lambda U^{(J)})\mod S^{-2,2-2\gamma_-}_-\cap S^{-2,0}_+\\
&=-\mu\nu+\frac{\eps}{4r}G_{2,J}\mu^2U^{(J)}+\frac{\eps}{4r}g^{\alpha\beta\lambda}_J\wh{\omega}_\alpha\wh{\omega}_\beta \mu^2 \partial_\lambda U^{(J)}+\eps S^{-2,1-\gamma_-}_-\cap\eps S^{-2,0}_+\mod S^{-2,2-2\gamma_-}_-\cap S^{-2,0}_+\\
&=-\mu\nu+\frac{\eps}{4r}G_{2,J}\mu^2U^{(J)}+\frac{\eps}{4r}g^{\alpha\beta\lambda}_J\wh{\omega}_\alpha\wh{\omega}_\beta \mu^2 \partial_\lambda U^{(J)}\mod S^{-2,1-\gamma_-}_-\cap S^{-2,0}_+.}
Moreover, \fm{\partial_tU&=q_tU_q+\eps t^{-1}U_s=\frac{\mu+\nu}{2}U_q\mod \eps S^{-1,1-\gamma_-}_-\cap \eps S^{-1,0}_+=\frac{1}{2}\mu U_q\mod  S^{-1,1-\gamma_-}_-\cap  S^{-1,0}_+,\\
\partial_iU&=q_iU_q+\sum_lU_{\omega_l}\partial_i\omega_l=(\frac{(\nu-\mu)\omega_i}{2}+\lambda_i)U_q\mod  S^{-1,1-\gamma_-}_-\cap S^{-1,0}_+\\
&=-\frac{1}{2}\mu\omega_iU_q\mod S^{-1,1-\gamma_-}_-\cap S^{-1,0}_+.}
In conclusion, we have
\eq{\label{prop4.6f1} \partial_\alpha U=-\frac{1}{2}\mu\wh{\omega}_\alpha U_q\mod S^{-1,1-\gamma_-}_-\cap S^{-1,0}_+.}
Thus,
\fm{&[m^{\alpha\beta}+g^{\alpha\beta}_J\cdot\eps r^{-1}U^{(J)}+g^{\alpha\beta\lambda}_{J} \partial_\lambda(\eps r^{-1}U^{(J)})]q_\alpha q_\beta\\
&=-\mu\nu+\frac{\eps}{4r}G_{2,J}\mu^2U^{(J)}+\frac{\eps}{4r}g^{\alpha\beta\lambda}_J\wh{\omega}_\alpha\wh{\omega}_\beta \mu^2\cdot\frac{-\wh{\omega}_\lambda}{2}\mu  U_q^{(J)}\mod S^{-2,1-\gamma_-}_-\cap S^{-2,0}_+\\
&=-\mu(\nu-\frac{\eps}{4r}G_{2,J}\mu U^{(J)}+\frac{\eps}{8r}G_{3,J}\mu^2  U_q^{(J)})\mod S^{-2,1-\gamma_-}_-\cap S^{-2,0}_+\\
&\in   S^{-2,2-\gamma_-}_-\cap  S^{-2,1}_++S^{-2,1-\gamma_-}_-\cap S^{-2,0}_+\subset S^{-2,2-\gamma_-}_-\cap S_+^{-2,1}.}
\end{proof}\rm

In order to prove that $\eps r^{-1}U$ is an approximate solution to \eqref{qwe}, we need the following lemma.

\lem{\label{lem4.7}We have \fm{g^{\alpha\beta}(\eps r^{-1}U,\partial(\eps r^{-1}U))q_{\alpha\beta}=-\eps r^{-1}\mu_s-\mu r^{-1}\mod S^{-2,1-\gamma_-}_-\cap  S^{-2,0}_+.}}
\begin{proof}
We first note that \fm{\eps t^{-1}\nu_s&=\nu_t-q_t\nu_q=(\nu_t+\nu_r)-\nu\nu_q,\\
\sum_l\partial_i\omega_l\nu_{\omega_l}&=\nu_i-q_i\nu_q=(\nu_i-\omega_i\nu_r)-\lambda_i\nu_q.}
By Lemma \ref{c1lem6} and since $\nu\in  S^{-1,1-\gamma_-}_- \cap  S^{-1,0}_+$, we conclude that $\nu_t+\nu_r,\nu_i-\omega_i\nu_r\in S^{-2,1-\gamma_-}_- \cap  S^{-2,0}_+$.  By Lemma \ref{lem4.5}, we have $\lambda_i\nu_q,\nu\nu_q\in S^{-2,1-2\gamma_-}\cap ( S^{-2,-\gamma_+}_++ S^{-3,0}_+)$.  Thus, we have $\eps t^{-1}\nu_s,\sum_l\partial_i\omega_l\nu_{\omega_l}\in  S^{-2,1-\gamma_-}_-\cap S^{-2,0}_+$.  Moreover, we have $\partial\lambda_i\in S^{-1,-\gamma_-}_-\cap S^{-1,-1}_+$ by Lemma \ref{c1lem6} and Lemma \ref{lem4.5}. It follows that
\fm{q_{tt}&=\partial_t(\frac{\mu+\nu}{2})=\frac{1}{4}\mu_q(\mu+\nu)+\frac{\eps}{2 t}\mu_s+\frac{1}{4}\nu_q(\mu+\nu)+\frac{\eps}{2 t}\nu_s\\&=\frac{1}{4}\mu_q\mu+\frac{1}{4}\mu_q\nu+\frac{\eps}{2 t}\mu_s+\frac{1}{4}\nu_q\mu\mod  S^{-2,1-\gamma_-}_- \cap  S^{-2,0}_+\\
&=\frac{1}{4}\mu_q\mu\mod S^{-1,-\gamma_-}_- \cap( S^{-1,-\gamma_+}_++ S^{-2,0}_+),}
\fm{q_{ti}&=\partial_i(\frac{\mu+\nu}{2})=\frac{1}{2}(\mu_q+\nu_q)(\lambda_i+\frac{\omega_i(\nu-\mu)}{2})+\frac{1}{2}\sum_l(\mu_{\omega_l}+\nu_{\omega_l})\partial_i\omega_l\\&=-\frac{1}{4}\omega_i\mu_q\mu\mod  S^{-1,-\gamma_-}_- \cap( S^{-1,-\gamma_+}_++ S^{-2,0}_+),}
\fm{q_{ij}&=\partial_i(\lambda_j+\frac{\omega_j(\nu-\mu)}{2})\\&=\partial_i\lambda_j+\frac{1}{2}\partial_i\omega_j (\nu-\mu)+\frac{1}{2}\omega_j (\nu_q-\mu_q)(\lambda_i+\frac{\omega_i(\nu-\mu)}{2})+\frac{1}{2}\omega_j\sum_l(\mu_{\omega_l}+\nu_{\omega_l})\partial_i\omega_l\\&=\frac{1}{4}\omega_i\omega_j\mu\mu_q+\partial_i\lambda_j-\frac{1}{2}\mu\partial_i\omega_j -\frac{1}{4}\omega_j \mu_q(2\lambda_i+\omega_i\nu)\\&\hspace{1em}-\frac{1}{4}\omega_j \nu_q\omega_i\mu+\frac{1}{2}\omega_j\sum_l\mu_{\omega_l}\partial_i\omega_l\mod  S^{-2,1-\gamma_-}_-\cap  S^{-2,0}_+\\
&=\frac{1}{4}\omega_i\omega_j\mu\mu_q\mod S^{-1,0}_-\cap S^{-1,0}_+.}
Thus, we have \eq{\label{lem4.7f1}q_{\alpha\beta}=\frac{1}{4}\wh{\omega}_\alpha\wh{\omega}_\beta\mu\mu_q\mod S_-^{-1,0}\cap S^{-1,0}_+\in (S^{0,-1-\gamma_-}_-+S^{-1,0}_-)\cap (S^{0,-1-\gamma_+}_++S^{-1,0}_+).} Moreover, we have
\fm{\Box q&=-(\frac{1}{4}\mu_q\mu+\frac{1}{4}\mu_q\nu+\frac{\eps}{2 t}\mu_s+\frac{1}{4}\nu_q\mu)+\frac{1}{4}\mu\mu_q-\frac{1}{4}\mu\nu_q-\frac{1}{4}\mu_q\nu\\&\quad+\sum_i[\partial_i\lambda_i-\frac{1}{2}\mu\partial_i\omega_i -\frac{1}{2}\omega_i \mu_q\lambda_i+\frac{1}{2}\omega_i\sum_l\mu_{\omega_l}\partial_i\omega_l] \mod S^{-2,1-\gamma_-}_-\cap  S^{-2,0}_+\\
&=- \frac{\eps}{2 t}\mu_s-\frac{1}{2}\mu\nu_q-\frac{1}{2}\mu_q\nu-\mu r^{-1} \mod   S^{-2,1-\gamma_-}_-\cap S^{-2,0}_+.}
Here we note that $\sum_i\partial_i\omega_i=2/r$, $\sum_i\omega_i\partial_i\omega_l=0$ and $\sum_i\omega_i\lambda_i=0$. Moreover, we have $\sum_i\omega_i\partial_r\lambda_i=\partial_r(\sum_i\omega_i\lambda_i)=0$ and $(\partial_i-\omega_i\partial_r)\lambda_i\in S^{-2,1-\gamma_-}_-\cap S^{-2,0}_+$, so
\fm{\sum_i\partial_i\lambda_i&=\sum_i\omega_i\partial_r\lambda_i+\sum_i(\partial_i-\omega_i\partial_r)\lambda_i\in S^{-2,1-\gamma_-}_-\cap S^{-2,0}_+.}

Recall that we have \eqref{lem4.6f3} and \eqref{prop4.6f1}. So to finish the proof, we need to prove
\fm{\ [m^{\alpha\beta}+\eps r^{-1}g^{\alpha\beta}_J U^{(J)}-\frac{1}{2}\eps r^{-1}g^{\alpha\beta\lambda}_{J} \wh{\omega}_\lambda \mu U_q^{(J)})]q_{\alpha\beta}=-\eps r^{-1}\mu_s-\mu r^{-1}\mod S^{-2,1-\gamma_-}_-\cap S^{-2,0}_+. }
Note that
\fm{&[m^{\alpha\beta}+\eps r^{-1}g^{\alpha\beta}_J U^{(J)}-\frac{1}{2}\eps r^{-1}g^{\alpha\beta\lambda}_{J} \wh{\omega}_\lambda \mu U_q^{(J)})]q_{\alpha\beta}\\
&=- \frac{\eps}{2 t}\mu_s-\frac{1}{2}\mu\nu_q-\frac{1}{2}\mu_q\nu-\mu r^{-1}+ \frac{1}{4}\eps r^{-1}G_{2,J} \mu \mu_qU^{(J)}-\frac{1}{8}\eps r^{-1}G_{3,J} \mu^2 \mu_qU_q^{(J)}\mod  S^{-2,1-\gamma_-}_-\cap S^{-2,0}_+\\
&=\frac{t-r}{2rt}\mu_s- \frac{\eps}{2 r}(\frac{1}{4}G_{2,J} \mu^2U_q^{(J)}-\frac{1}{8}G_{3,J}(\mu^3U_{qq}^{(J)}+\mu^2\mu_q U_q^{(J)}))\\&\quad-\frac{1}{2}\mu\nu_q-\frac{1}{2}\mu_q\nu-\mu r^{-1}+ \frac{1}{4}\eps r^{-1}G_{2,J} \mu \mu_qU^{(J)}-\frac{1}{8}\eps r^{-1}G_{3,J} \mu^2 \mu_qU_q^{(J)}\mod  S^{-2,1-\gamma_-}_-\cap S^{-2,0}_+.
}
By \eqref{lem4.5c3} and \eqref{lem4.5c4} we have
\fm{&\mu\nu_q+\nu\mu_q\\&=\mu(\frac{\eps G_{2,J}}{4 r}\partial_q(\mu U^{(J)})-\frac{\eps G_{3,J}}{8r}\partial_q(\mu^2U_q^{(J)}))+ S^{-2,1-\gamma_-}_-\cap  S^{-2,0}_+\\&\quad+(\frac{\eps G_{2,J}}{4 r}\mu U^{(J)}-\frac{\eps G_{3,J}}{8r}\mu^2U_q^{(J)})\mu_q+ S^{-2,1-\gamma_-}_-\cap  S^{-2,0}_+\\
&=\frac{\eps G_{2,J}}{2r}\mu\mu_q U^{(J)}+\frac{\eps G_{2,J}}{4r}\mu^2U_q^{(J)}-\frac{\eps G_{3,J}}{8r}(3\mu^2\mu_qU_q^{(J)}+\mu^3U_{qq}^{(J)})\mod S^{-2,1-\gamma_-}_-\cap  S^{-2,0}_+.}
As a result,
\fm{&[m^{\alpha\beta}+\eps r^{-1}g^{\alpha\beta}_J U^{(J)}-\frac{1}{2}\eps r^{-1}g^{\alpha\beta\lambda}_{J} \wh{\omega}_\lambda \mu U_q^{(J)})]q_{\alpha\beta}\\
&=S^{-2,1}\cdot (S^{0,-\gamma_-}_-\cap S^{0,-\gamma_+}_+)- \frac{\eps}{2 r}(\frac{1}{4}G_{2,J} \mu^2U_q^{(J)}-\frac{1}{8}G_{3,J}(\mu^3U_{qq}^{(J)}+\mu^2\mu_q U_q^{(J)}))\\&\quad-[\frac{\eps G_{2,J}}{4r}\mu\mu_q U^{(J)}+\frac{\eps G_{2,J}}{8r}\mu^2U_q^{(J)}-\frac{\eps G_{3,J}}{16r}(3\mu^2\mu_qU_q^{(J)}+\mu^3U_{qq}^{(J)})]\\&\quad-\mu r^{-1}+ \frac{1}{4}\eps r^{-1}G_{2,J} \mu \mu_qU^{(J)}-\frac{1}{8}\eps r^{-1}G_{3,J} \mu^2 \mu_qU_q^{(J)}\mod  S^{-2,1-\gamma_-}_-\cap  S^{-2,0}_+\\
&=\frac{\eps G_{3,J}}{8r}\mu^2\partial_q(\mu U_{q}^{(J)})-\frac{\eps G_{2,J}}{4r}\mu^2U_q^{(J)}-\mu r^{-1}\mod  S^{-2,1-\gamma_-}_-\cap  S^{-2,0}_+\\&=-\eps r^{-1}\mu_s-\mu r^{-1}\mod S^{-2,1-\gamma_-}_-\cap  S^{-2,0}_+.}

\end{proof}\rm

Finally we prove that $\eps r^{-1}U$ has good pointwise bounds and is an approximate solution to \eqref{qwe} in $\D$.
\prop{\label{prop4.8}We have $\eps r^{-1}U\in \eps S^{-1,1-\gamma_-}_-\cap\eps S^{-1,0}_+$,
\fm{ g^{\alpha\beta}(\eps r^{-1}U,\partial(\eps r^{-1}U))\partial_\alpha\partial_\beta(\eps r^{-1}U)-f(\eps r^{-1}U,\partial(\eps r^{-1}U))\in \eps S^{-3,1-\gamma_-}_-\cap\eps S^{-3,0}_+.}
In other word, for $(t,x)\in\D$ and each $I$,
\fm{|Z^I(\eps r^{-1}U)|\lesssim_I \eps t^{-1+C_I\eps}\lra{\max\{0,t-r\}}^{1-\gamma_-},}
\fm{|Z^I(g^{\alpha\beta}(\eps r^{-1}U,\partial(\eps r^{-1}U))\partial_\alpha\partial_\beta(\eps r^{-1}U)-f(\eps r^{-1}U,\partial(\eps r^{-1}U)))|\lesssim_I\eps t^{-3+C_I\eps}\lra{\max\{0,t-r\}}^{1-\gamma_-}.}
Note that we  have a better bound for $\partial (\eps r^{-1}U)$: for all $(t,x)\in\D$, \fm{|\partial(\eps r^{-1}U)|\lesssim \eps t^{-1}.}
}
\begin{proof}Since $r^{-1}\in S^{-1,0}$ and  $U\in S^{0,1-\gamma_-}_-\cap S^{0,0}_+$, we have $\eps r^{-1}U\in \eps S^{-1,1-\gamma_-}_-\cap\eps S^{-1,0}_+$. By \eqref{prop4.6f1}, we have $\partial_\alpha U=-\frac{1}{2}\mu\wh{\omega}_\alpha U_q\bmod S^{-1,1-\gamma_-}_-\cap  S^{-1,0}_+$, so we have 
\eq{\label{prop4.9eru1}\partial_\alpha(\eps r^{-1}U)&=-\eps (\partial_\alpha r)r^{-2}U-\frac{1}{2}\eps r^{-1}\wh{\omega}_\alpha \mu U_q\mod \eps S^{-2,1-\gamma_-}_-\cap  \eps S^{-2,0}_+\\&=-\frac{1}{2}\eps r^{-1}\wh{\omega}_\alpha \mu U_q\mod \eps S^{-2,1-\gamma_-}_-\cap  \eps S^{-2,0}_+.}
By \eqref{def3.1a5}, we have $|\mu U_q|\lesssim 1$, so
\fm{|\partial(\eps r^{-1}U)|\lesssim \eps r^{-1}+\eps t^{-2+C\eps}\lesssim \eps t^{-1}.}

We have 
\fm{(\eps r^{-1}U)_{tt}&=\eps r^{-1}(-U_s\eps t^{-2}+2U_{sq}q_t\eps t^{-1}+U_{ss}\eps^2 t^{-2}+q_{tt}U_q+q_t^2U_{qq})\\
&=\eps r^{-1}(2U_{sq}q_t\eps t^{-1}+q_{tt}U_q+q_t^2U_{qq})\mod \eps S^{-3,1-\gamma_-}_-\cap \eps S^{-3,0}_+\\
&=\eps r^{-1}(q_{tt}U_q+q_t^2U_{qq})\mod \eps S^{-2,-\gamma_-}_-\cap (\eps S^{-2,-\gamma_+}_++\eps S^{-3,0}_+),}
\fm{(\eps r^{-1}U)_{ti}&=\eps r^{-1}(U_{qq}q_tq_i+\sum_lU_{\omega_lq}q_t\partial_i\omega_l+U_qq_{it}+U_{sq}q_i\eps t^{-1}+\sum_lU_{s\omega_l}\partial_i\omega_l\eps t^{-1})\\
&\quad-\eps r^{-2}\omega_i(U_qq_t+U_s\eps t^{-1})\\
&=\eps r^{-1}(U_{qq}q_tq_i+U_qq_{it}) \mod \eps S^{-2,-\gamma_-}_-\cap (\eps S^{-2,-\gamma_+}_++\eps S^{-3,0}_+),}
\fm{(\eps r^{-1}U)_{ij}&=\eps r^{-1}(U_{qq}q_iq_j+\sum_lU_{q\omega_l}(q_i\partial_j\omega_l+q_j\partial_i\omega_l)+U_qq_{ij}+\sum_{l,l^\prime}U_{\omega_l\omega_{l^\prime}}\partial_i\omega_l\partial_j\omega_{l^\prime})\\
&\quad-\eps r^{-2}\omega_i(U_qq_j+\sum_lU_{\omega_l}\partial_j\omega_l)-\eps r^{-2}\omega_j(U_qq_i+\sum_lU_{\omega_l}\partial_i\omega_l)+\eps \partial_j(r^{-2}\omega_i)U\\
&=\eps r^{-1}(U_{qq}q_iq_j+\sum_lU_{q\omega_l}(q_i\partial_j\omega_l+q_j\partial_i\omega_l)+U_qq_{ij})\\&\quad-\eps r^{-2}U_q(\omega_iq_j+\omega_jq_i) \mod \eps S^{-3,1-\gamma_-}_-\cap\eps S^{-3,0}_+ \\
&=\eps r^{-1}(U_{qq}q_iq_j+U_qq_{ij}) \mod \eps S^{-2,-\gamma_-}_-\cap (\eps S^{-2,-\gamma_+}_++\eps S^{-3,0}_+).}
In summary, we have
\eq{\label{prop4.9eru}\partial_{\alpha}\partial_\beta(\eps r^{-1}U)=\eps r^{-1}(U_{qq}q_\alpha q_\beta+U_{q}q_{\alpha\beta})\mod \eps S^{-2,-\gamma_-}_-\cap (\eps S^{-2,-\gamma_+}_++\eps S^{-3,0}_+),}
\fm{&\Box (\eps r^{-1}U)\\&=-2\eps^2 (tr)^{-1}U_{sq}q_t+\eps r^{-1} \sum_{i,l}2U_{q\omega_l}q_i\partial_i\omega_l-2\eps r^{-2}U_q q_r\\
&\quad+\eps r^{-1}(\Box q)U_q+\eps r^{-1}(m^{\alpha\beta}q_\alpha q_\beta)U_{qq}\mod \eps S^{-3,1-\gamma_-}_-\cap\eps S^{-3,0}_+\\
&=-\eps^2 (tr)^{-1}U_{sq}\mu+\eps r^{-2}\mu U_q +\eps r^{-1}(\Box q)U_q+\eps r^{-1}(m^{\alpha\beta}q_\alpha q_\beta)U_{qq}\mod \eps S^{-3,1-\gamma_-}_-\cap\eps S^{-3,0}_+.}
Here note that
\fm{\eps r^{-1} \sum_{i,l}2U_{q\omega_l}q_i\partial_i\omega_l&=\eps r^{-1} \sum_{i,l}2U_{q\omega_l}\lambda_i\partial_i\omega_l+\eps r^{-1} \sum_{i,l}2U_{q\omega_l}q_r\omega_i\partial_i\omega_l\\
&=\eps r^{-2}\sum_{i,l}2U_{q\omega_l}\lambda_i(\delta_{il}-\omega_i\omega_l)+0\in \eps S^{-3,1-\gamma_-}_-\cap\eps S^{-3,0}_+.}
Moreover, by \eqref{lem4.6f3}, we have
\fm{g^{\alpha\beta}(\eps r^{-1}U,\partial(\eps r^{-1}U))-m^{\alpha\beta}&=\eps S^{-1,1-\gamma_-}_-\cap\eps S^{-1,0}_++\eps^2S^{-2,2-2\gamma_-}_-\cap \eps^2S^{-2,0}_+\\&\subset \eps S^{-1,1-\gamma_-}_-\cap\eps S^{-1,0}_+.}
As a result, by applying Proposition \ref{prop4.6} and Lemma \ref{lem4.7}, we have
\fm{&g^{\alpha\beta}(\eps r^{-1}U,\partial(\eps r^{-1}U))\partial_\alpha\partial_\beta(\eps r^{-1}U)\\
&=\Box (\eps r^{-1}U)+[g^{\alpha\beta}(\eps r^{-1}U,\partial(\eps r^{-1}U))-m^{\alpha\beta}]\cdot \eps r^{-1}(U_{qq}q_\alpha q_\beta+U_qq_{\alpha\beta})\\&\quad\mod\eps^2 S^{-3,1-2\gamma_-}_-\cap (\eps S^{-3,-\gamma_+}_++\eps S^{-4,0}_+)\\
&=-\eps^2(tr)^{-1}\mu U_{sq}+\eps r^{-2}\mu U_q+\eps r^{-1}g^{\alpha\beta}(\eps r^{-1}U,\partial(\eps r^{-1}U))q_\alpha q_\beta U_{qq}\\
&\quad+\eps r^{-1}g^{\alpha\beta}(\eps r^{-1}U,\partial(\eps r^{-1}U))q_{\alpha\beta}U_q\mod \eps S^{-3,1-\gamma_-}_-\cap\eps S^{-3,0}_+\\
&=-\eps^2(tr)^{-1}\mu U_{sq}+\eps r^{-2}\mu U_q+\eps S^{-3,1-2\gamma_-}_-\cap \eps S^{-3,-\gamma_+}_+\\
&\quad+\eps r^{-1}(-\eps r^{-1}\mu_s-\mu r^{-1})U_q+\eps S^{-3,1-2\gamma_-}_-\cap \eps S^{-3,-\gamma_+}_+\mod \eps S^{-3,1-\gamma_-}_-\cap\eps S^{-3,0}_+\\
&=\eps^2r^{-2}t^{-1}(t-r)\mu U_{sq}-\eps^2 r^{-2}\partial_s(\mu U_q) \mod\eps S^{-3,1-\gamma_-}_-\cap\eps S^{-3,0}_+\\
&=-\eps^2 r^{-2}\partial_s(\mu U_q) \mod \eps S^{-3,1-\gamma_-}_-\cap\eps S^{-3,0}_+.}
In the last step, we note that $\eps^2r^{-2}t^{-1}(t-r)\mu U_{sq}\in \eps S^{-3,1}\cdot S^{0,0}\cdot (S^{0,-\gamma_-}_-\cap S^{0,-\gamma_+}_+)$.

Moreover, since $f(0,0)=0$ and $df(0,0)=0$, by Lemma \ref{lem2.9} we have
\fm{&f^{(I)}(\eps r^{-1}U,\partial(\eps r^{-1}U))\\&=f^{I,\alpha\beta}_{JK}\cdot\partial_\alpha (\eps r^{-1}U^{(J)})\cdot\partial_\beta (\eps r^{-1}U^{(K)})\mod \eps^3S^{-3,3-3\gamma_-}_-\cap \eps^3S^{-3,0}_+\\
&=f^{I,\alpha\beta}_{JK}\cdot(-\frac{1}{2}\mu\wh{\omega}_\alpha\eps r^{-1}U_q^{(J)}+\eps S^{-2,1-\gamma_-}_-\cap\eps S^{-2,0}_+)\cdot(-\frac{1}{2}\mu\wh{\omega}_\beta\eps r^{-1}U_q^{(K)}+\eps S^{-2,1-\gamma_-}_-\cap\eps S^{-2,0}_+)\\&\qquad\qquad\qquad\mod \eps^3S^{-3,3-3\gamma_-}_-\cap \eps^3S^{-3,0}_+\\
&=\frac{\eps^2}{4r^2}F_{2,JK}^I\mu^2U_q^{(J)}U^{(K)}_q \mod \eps^2 S^{-3,1-2\gamma_-}_-\cap\eps^2 S^{-3,0}_+\\
&=- \eps^2r^{-2}\partial_s(\mu U^{(I)}_q)\mod \eps^2 S^{-3,1-\gamma_-}_-\cap\eps^2 S^{-3,0}_+.}
\end{proof}\rm

\subsection{An approximate solution $u_{app}$}\label{sec4.3} For each small $\delta\in(0,1)$, $c\in(0,1/2)$ and $\eps\ll1$, we have obtained a function $\eps r^{-1}U$ defined in $\D$. Choose $\psi\in C_c^\infty(\R)$ such that $\psi\equiv 1$ in $[1-c/2,1+c/2]$ and $\psi\equiv 0$ in $\R\setminus(1-c,1+c)$.
 We now define  $u_{app}$ by
\begin{equation}\label{uappdef}u_{app}(t,x):=\eps r^{-1}\psi(r/t)U(\eps\ln(t)-\delta,q(t,x),\omega).\end{equation}
Note that $u_{app}(t,x)$ is an $\R^M$-valued function defined for all $(t,x)\in[ T_\eps,\infty)\times\R^{3}$.  In fact,  since $\psi\equiv 0$ when $|t-r|\geq ct$, we have $u_{app}\equiv 0$ unless $(t,x)\in\D$; since $\psi\equiv 1$ when $|t-r|\leq ct/2$, we have $u_{app}=\eps r^{-1}U$ whenever $|t-r|\leq ct/2$ and $t\geq  T_\eps$.

We now prove the  estimates for $u_{app}$ in Proposition \ref{mainprop4}.  The estimates are in fact the same as those in Proposition \ref{prop4.8}. However, the estimates in Proposition \ref{prop4.8} hold for $(t,x)\in\D$ while in Proposition \ref{mainprop4} we only assume $t\geq  T_\eps$.

\begin{proof}[Proof of Proposition \ref{mainprop4}]
For $|r-t|\leq ct/2$, all the estimates  follow directly from Proposition \ref{prop4.8}, we have $u_{app}=\eps r^{-1}U$, so all the estimates in Proposition \ref{mainprop4} follow from Proposition \ref{prop4.8}. For $|r-t|\geq ct$, we have $u_{app}=0$, so there is nothing to prove.

Now suppose that $ct/2\leq |r-t|\leq ct$.  Note that now we have $|r-t|\sim t$ and $(t,x)\in\D$. Since $\partial^kZ^I(r/t)=O(t^{-k})$ for $t\sim r$, we have $\partial^kZ^I(\psi(r/t))=O(t^{-k})$, so $\partial^kZ^I(\psi(r/t))\in S^{-k,0}$ for each $k\geq 0$. As a result, we have
\fm{u_{app}=\psi(r/t)\cdot \eps r^{-1}U\in \eps S^{-1,1-\gamma_-}_-\cap \eps S^{-1,0}_+.}
By Lemma \ref{c1lem6}, we have 
\fm{\partial^ku_{app}\in\eps S^{-1,1-\gamma_--k}_-\cap \eps S^{-1,-k}_+,\qquad \forall k\geq 0.}
As a result, at each $(t,x)\in\D$ with $|r-t|\sim t$, for each $k\geq 0$ we have 
\eq{\label{c3pfff1}|Z^I\partial^k u_{app}|&\lesssim\eps t^{-1+C\eps}\lra{r-t}^{1-\gamma_--k}1_{q\leq 0}+\eps t^{-1+C\eps}\lra{r-t}^{-k}1_{q>0}\\
&\lesssim \eps t^{-\gamma_--k+C\eps}1_{q\leq 0}+\eps t^{-1-k+C\eps}1_{q>0}.}
In other words, by setting $\D'_+=\{(t,x)\in\D:\ r\geq (1+c/2)t\}$ and $\D'_-=\{(t,x)\in\D:\ r\leq (1-c/2)t\}$, we have $\partial^lu_{app}\in \eps S_{\D'_-}^{-\gamma_--l,0}\cap \eps S^{-1-l,0}_{\D'_+}$. Here we remind our readers of Definition \ref{c1defn1.5} with $\D$ replaced by $\D'_\pm$. We also notice that $r\geq (1+c/2)t$ implies $q>0$ and that $r\leq (1-c/2)t$ implies $q<0$. By applying Lemma \ref{lem2.9}, we have
\fm{g^{\alpha\beta}(u_{app},\partial u_{app})=m^{\alpha\beta}+\eps S_{\D'_-}^{-\gamma_-,0}\cap \eps S_{\D'_+}^{-1,0}\in S^{0,0}}
and since $\gamma_->2$,
\fm{f(u_{app},\partial u_{app})&=C\cdot \partial u_{app}\cdot\partial u_{app}+\eps^3 S_{\D'_-}^{-3\gamma_-,0}\cap \eps^3 S_{\D'_+}^{-3,0}\\&\in\eps^2 S_{\D'_-}^{-2\gamma_--2,0}\cap \eps^2 S_{\D'_+}^{-4,0}+\eps^3 S_{\D'_-}^{-3\gamma_-,0}\cap \eps^3 S_{\D'_+}^{-3,0}\subset \eps^2 S_{\D'_-}^{-2\gamma_--2,0}\cap \eps^2 S_{\D'_+}^{-3,0}.}
We conclude that
\fm{g^{\alpha\beta}(u_{app},\partial u_{app})\partial_\alpha\partial_\beta u_{app}-f(u_{app},\partial u_{app})\in \eps S_{\D'_-}^{-\gamma_--2,0}\cap \eps S_{\D'_+}^{-3,0},}
and that
\fm{|Z^I(g^{\alpha\beta}(u_{app},\partial u_{app})\partial_\alpha\partial_\beta u_{app}-f(u_{app},\partial u_{app}))|\lesssim \eps t^{-2-\gamma_-+C\eps}1_{q\leq 0}+\eps t^{-3+C\eps}1_{q> 0}.}
\end{proof}\rm

\subsection{More pointwise estimates for $u_{app}$}\label{sec4.4}
We end this section with more pointwise estimates for $u_{app}$.

\lem{\label{lem4.11} For $t\geq  T_\eps$ with $\eps\ll1$, we have \eq{\label{lem4.11c1}\partial_\alpha\partial_\beta u_{app}^{(I)}=\frac{\eps}{4 r}\psi(r/t)\wh{\omega}_\alpha\wh{\omega}_\beta \mu \partial_q(\mu U_q^{(I)})+O(\eps t^{-2+C\eps}\lra{q}^{-\gamma_{\sgn(q)}}+\eps t^{-3+C\eps}).}
As a result, we have
\eq{\label{lem4.11c2}g^{\alpha\beta}_J\partial_\alpha\partial_\beta u_{app}^{(K)}&=\frac{\eps}{4 r}\psi(r/t)G_{2,J}(\omega)\mu \partial_q(\mu U^{(K)})+O(\eps t^{-2+C\eps}\lra{q}^{-\gamma_{\sgn(q)}}+\eps t^{-3+C\eps}).}
}
\begin{proof}For $(t,x)\notin\D$, we have $\psi\equiv 0$ and $u_{app}\equiv 0$, so there is nothing to prove. Now assume that $(t,x)\in\D$. Since $\partial^k(\psi(r/t))\in S^{-k,0}$ for each $k\geq 0$, by the definition \eqref{uappdef}, we have
\fm{&\partial_\alpha\partial_\beta u_{app}^{(I)}\\&=\psi(r/t)\partial_\alpha\partial_\beta(\eps r^{-1}U^{(I)})+\partial(\psi(r/t))\cdot\partial(\eps r^{-1}U^{(I)})+\partial^2(\psi(r/t))\eps r^{-1}U^{(I)}\\
&=\psi(r/t)\partial_\alpha\partial_\beta(\eps r^{-1}U^{(I)})+S^{-1,0}\cdot [\eps S^{-1,-\gamma_-}_-\cap(\eps S^{-1,-\gamma_+}_++\eps S^{-2,0}_+)]+S^{-2,0}\cdot (\eps S^{-1,1-\gamma_-}_-\cap \eps S^{-1,0}_+).}
Recall from \eqref{prop4.9eru1} that
\fm{\partial (\eps r^{-1}U)&=-\frac{1}{2}\eps r^{-1}\wh{\omega} \mu U_q\mod \eps S^{-2,1-\gamma_-}_-\cap  \eps S^{-2,0}_+=\eps S^{-1,-\gamma_-}_-\cap(\eps S^{-1,-\gamma_+}_++\eps S^{-2,0}_+).} By applying \eqref{prop4.9eru}, \eqref{lem4.6f1} and \eqref{lem4.7f1}, we have
\fm{\partial_\alpha\partial_\beta u_{app}^{(I)}&=\psi(r/t)\cdot \eps r^{-1}(U_{qq}^{(I)}q_\alpha q_\beta+U_q^{(I)}q_{\alpha\beta})\mod \eps S^{-2,-\gamma_-}_-\cap (\eps S^{-2,-\gamma_+}_++\eps S^{-3,0}_+)\\
&=\psi(r/t)\cdot \eps r^{-1}(U_{qq}^{(I)}\cdot \frac{1}{4}\wh{\omega}_\alpha \wh{\omega}_\beta\mu^2+U_q^{(I)}\cdot \frac{1}{4}\wh{\omega}_\alpha \wh{\omega}_\beta\mu\mu_q)\\&\quad+\eps S^{-1,0}\cdot ((\eps S^{0,-1-\gamma_-}_-\cap \eps S^{0,-1-\gamma_+}_+)\cdot (S^{-1,1-\gamma_-}_-\cap S^{-1,0}_+)+(\eps S^{0,-\gamma_-}_-\cap \eps S^{0,-\gamma_+}_+)\cdot S^{-1,0}) \\&\quad\mod \eps S^{-2,-\gamma_-}_-\cap (\eps S^{-2,-\gamma_+}_++\eps S^{-3,0}_+)\\
&=\frac{\eps}{4 r}\psi(r/t)\wh{\omega}_\alpha\wh{\omega}_\beta\mu\partial_q(\mu U_q^{(I)})\mod \eps S^{-2,-\gamma_-}_-\cap (\eps S^{-2,-\gamma_+}_++\eps S^{-3,0}_+).}
Since $\lra{r-t}^{-1}\lesssim t^{C\eps}\lra{q}^{-1}$, we obtain \eqref{lem4.11c1}. The estimate \eqref{lem4.11c2} then follows from \eqref{lem4.11c1} and \eqref{def3.1a52}.
\end{proof}\rm

To state the next lemma, we define vector fields $T_\alpha=\partial_\alpha+\wh{\omega}_\alpha \partial_t$ and $\wt{T}_\alpha=q_t\partial_\alpha-q_\alpha \partial_t$ for $\alpha=0,1,2,3$. 

\lem{\label{lem4.10} Let $\phi=\phi(t,x)$ be a smooth function defined for $t\geq  T_\eps$ with $\eps\ll1$. Then, for each $(t,x)$ with $t\geq  T_\eps$, we have
\eq{\label{lem4.10c0}&|g_J^{\alpha\beta\lambda}(\partial_\lambda\phi)\partial_\alpha\partial_\beta u_{app}^{(K)}|\lesssim\eps t^{-1}|\partial\phi|+\eps t^{-1+C\eps}|T\phi|\lra{q}^{-1-\gamma_{\sgn(q)}};}
\eq{\label{lem4.10c}&|g_J^{\alpha\beta\lambda}(\partial_\alpha\phi)(\partial_\beta\phi)\partial_t\partial_\lambda u_{app}^{(K)}|\lesssim\eps t^{-1}|\partial\phi|^2+\eps t^{-1+C\eps}|T\phi||\partial\phi|\lra{q}^{-1-\gamma_{\sgn(q)}};}
\eq{\label{lem4.10c1}&|g_J^{\alpha\beta\lambda}(\partial_\alpha\partial_\beta \phi)\partial_\lambda\partial_{\sigma} u_{app}^{(K)}|\lesssim\eps t^{-1} |\partial^2\phi|+\eps t^{-1+C\eps}|T\partial\phi|\lra{q}^{-1-\gamma_{\sgn(q)}}.}
We can also replace $|T\phi|$ with $|\wt{T}\phi|$ in \eqref{lem4.10c0} and \eqref{lem4.10c}, and replace $|T\partial \phi|$ with $|\wt{T}\partial \phi|$ in \eqref{lem4.10c1}.
}
\begin{proof}There is nothing to prove for $(t,x)\notin\D$, so let us fix $(t,x)\in\D$. The idea of the proof is to use the previous lemma in order to handle $\partial^2u_{app}$, and to use $\partial_\lambda \phi=T_\lambda\phi-\wh{\omega}_\lambda\partial_t\phi$ to replace $\partial\phi$ with $\wh{\omega}\partial_t\phi$. By \eqref{lem4.11c1} we have
\fm{&g_J^{\alpha\beta\lambda}\partial_\lambda\phi\partial_\alpha\partial_\beta u_{app}^{(K)}\\&=\frac{\eps}{4r}\psi(r/t)g_J^{\alpha\beta\lambda}(T_\lambda\phi-\wh{\omega}_\lambda\partial_t\phi)\wh{\omega}_\alpha\wh{\omega}_\beta\mu\partial_q(\mu U_q^{(K)})+O((\eps t^{-2+C\eps}\lra{q}^{-\gamma_{\sgn(q)}}+\eps t^{-3+C\eps})|\partial \phi|)\\
&=-\frac{\eps}{4r}\psi(r/t)G_{3,J}\phi_t\mu\partial_q(\mu U_q^{(K)})+O(\eps t^{-1+C\eps}|T\phi|\lra{q}^{-1-\gamma_{\sgn(q)}}+(\eps t^{-2+C\eps}\lra{q}^{-\gamma_{\sgn(q)}}+\eps t^{-3+C\eps})|\partial \phi|)\\
&=O(\eps t^{-1}|\partial\phi|+\eps t^{-1+C\eps}|T\phi|\lra{q}^{-1-\gamma_{\sgn(q)}}+(\eps t^{-2+C\eps}\lra{q}^{-\gamma_{\sgn(q)}}+\eps t^{-3+C\eps})|\partial \phi|)).}
In the second last estimate we apply \eqref{def3.1a51}. In addition, since
\fm{g^{\alpha\beta\lambda}_J(\partial_\alpha\phi)( \partial_\beta \phi)&=g^{\alpha\beta\lambda}_J(T_\alpha\phi-\wh{\omega}_\alpha\phi_t)\partial_\beta\phi=-g^{\alpha\beta\lambda}_J \wh{\omega}_\alpha\phi_t(T_\beta\phi-\wh{\omega}_\beta\phi_t)+O(|T\phi||\partial\phi|) \\
&=g^{\alpha\beta\lambda}_J \wh{\omega}_\alpha\phi_t\cdot  \wh{\omega}_\beta\phi_t+O(|T\phi||\partial\phi|),}
we have
\fm{&g_J^{\alpha\beta\lambda}(\partial_\alpha\phi)(\partial_\beta\phi)\partial_t\partial_\lambda u_{app}^{(K)}\\
&=\frac{\eps}{4r}\psi(r/t)\cdot(g_J^{\alpha\beta\lambda}\phi_\alpha \phi_\beta)\cdot\wh{\omega}_\lambda\wh{\omega}_0\mu\partial_q(\mu U_q^{(K)})+O((\eps t^{-2+C\eps}\lra{q}^{-\gamma_{\sgn(q)}}+\eps t^{-3+C\eps})|\partial \phi|^2)\\
&=\frac{\eps}{4r}\psi(r/t)g_J^{\alpha\beta\lambda}\wh{\omega}_\alpha\wh{\omega}_\beta\phi_t^2  \wh{\omega}_\lambda\mu\partial_q(\mu U_q^{(K)})\\&\quad+O(\eps t^{-1+C\eps}|T\phi||\partial\phi|\lra{q}^{-1-\gamma_{\sgn(q)}}+(\eps t^{-2+C\eps}\lra{q}^{-\gamma_{\sgn(q)}}+\eps t^{-3+C\eps})|\partial \phi|^2)\\
&=-\frac{\eps}{4r}G_{3,J}\psi(r/t) \phi_t^2 \mu\partial_q(\mu U_q^{(K)})\\&\quad+O(\eps t^{-1+C\eps}|T\phi||\partial\phi|\lra{q}^{-1-\gamma_{\sgn(q)}}+(\eps t^{-2+C\eps}\lra{q}^{-\gamma_{\sgn(q)}}+\eps t^{-3+C\eps})|\partial \phi|^2)\\
&=O(\eps t^{-1}|\partial\phi|^2+\eps t^{-1+C\eps}|T\phi||\partial\phi|\lra{q}^{-1-\gamma_{\sgn(q)}}+(\eps t^{-2+C\eps}\lra{q}^{-\gamma_{\sgn(q)}}+\eps t^{-3+C\eps})|\partial \phi|^2).}
Again, we apply \eqref{def3.1a51} in the last estimate. 

Next, since
\fm{\partial_\alpha\partial_\beta\phi&=(T_\alpha-\wh{\omega}_\alpha\partial_t)\partial_\beta \phi=-\wh{\omega}_\alpha(T_\beta-\wh{\omega}_\beta\partial_t)\phi_t+O(|T\partial\phi|)=\wh{\omega}_\alpha \wh{\omega}_\beta\phi_{tt}+O(|T\partial\phi|),}
we have
\fm{&g_J^{\alpha\beta\lambda}(\partial_\alpha\partial_\beta \phi)\partial_\lambda\partial_{\sigma} u_{app}^{(K)}\\
&=\frac{\eps}{4r}\psi(r/t)g_J^{\alpha\beta\lambda}\wh{\omega}_\lambda\wh{\omega}_\sigma \mu \partial_q(\mu U_q^{(K)})\cdot \partial_\alpha\partial_\beta \phi+O((\eps t^{-2+C\eps}\lra{q}^{-\gamma_{\sgn(q)}}+\eps t^{-3+C\eps})|\partial^2 \phi|)\\
&=\frac{\eps}{4r}\psi(r/t)g_J^{\alpha\beta\lambda}\wh{\omega}_\lambda\wh{\omega}_\sigma \mu \partial_q(\mu U_q^{(K)})\cdot \wh{\omega}_\alpha\wh{\omega}_\beta\phi_{tt}\\&\quad+O(\eps t^{-1+C\eps}|T\partial\phi|\lra{q}^{-1-\gamma_{\sgn(q)}}+(\eps t^{-2+C\eps}\lra{q}^{-\gamma_{\sgn(q)}}+\eps t^{-3+C\eps})|\partial^2 \phi|)\\
&=O(\eps t^{-1} |\partial^2\phi|+\eps t^{-1+C\eps}|T\partial\phi|\lra{q}^{-1-\gamma_{\sgn(q)}}+(\eps t^{-2+C\eps}\lra{q}^{-\gamma_{\sgn(q)}}+\eps t^{-3+C\eps})|\partial^2 \phi|).}

To end the proof, we notice that $\eps t^{-2+C\eps}\lra{q}^{-\gamma_{\sgn(q)}}+\eps t^{-3+C\eps}\lesssim \eps t^{-1}$ and that
\fm{|q_tT_\alpha \phi-\wt{T}_\alpha\phi|&=|(q_\alpha+q_t\wh{\omega}_\alpha)\phi_t|\lesssim (|\nu|+\sum_i|\lambda_i|)|\phi_t|\lesssim t^{-1+C\eps}|\partial\phi|.}
It follows that
\fm{\eps t^{-1+C\eps}\lra{q}^{-1-\gamma_{\sgn(q)}}|T\phi|&\lesssim \eps t^{-1+C\eps}\lra{q}^{-1-\gamma_{\sgn(q)}}|q_t^{-1}\wt{T}\phi|+\eps t^{-2+C\eps}\lra{q}^{-1-\gamma_{\sgn(q)}}|\partial\phi|\\
&\lesssim\eps t^{-1+C\eps}\lra{q}^{-1-\gamma_{\sgn(q)}}|\wt{T}\phi|+\eps t^{-1}|\partial\phi|.}
Similarly, we obtain
\fm{\eps t^{-1+C\eps}\lra{q}^{-1-\gamma_{\sgn(q)}}|T\partial \phi|&\lesssim \eps t^{-1+C\eps}\lra{q}^{-1-\gamma_{\sgn(q)}}|\wt{T}\partial\phi|+\eps t^{-1}|\partial^2\phi|.}
This finishes our proof.
\end{proof}\rm

\section{Energy estimates and Poincare's lemmas}\label{sec5}

We now derive the energy estimates and Poincar$\acute{\rm e}$'s lemmas, which are the main tools in the proof of our main theorem. For most results in this section, we assume $t\geq T_\eps=1/\eps$. 

In Section \ref{sec4}, our approximate optical function $q=q(t,x)$ is defined in $\Omega=\{t\geq T_\eps,\ |x|\geq t/2\}$. From now on, we extend it to all $(t,x)$ with $t\geq T_\eps$ by setting $q(t,x)=|x|-t$ whenever $|x|< t/2$. Note that $q$ is  continuous everywhere, and that it is  differentiable everywhere except on the cone $|x|= t/2$.

We now summarize the results in this section. In Sections \ref{sec5.2} and \ref{sec5.1}, we prove energy estimates. The one in Section \ref{sec5.2} holds for $t\geq 1/\eps$, and the one in Section \ref{sec5.1} holds for $t\leq 1/\eps$. See Propositions \ref{prop5.1} and \ref{prop5.1st}. In the proof of Proposition \ref{prop5.1}, we define a weight $w$ by \eqref{sec5wei} and an energy $E_u(\phi)$ by \eqref{sec5ene}. Then we differentiate the energy with respect to time and apply integration by parts. The proof of Proposition \ref{prop5.1st} is similar but simpler since we can use a weight involving $r-t$ instead of the approximate optical function $q$.

In Section \ref{sec5.3}, we prove Poincar$\acute{\rm e}$'s lemmas. There are three versions: Lemmas \ref{lp1}--\ref{lp2}. These lemmas are very different from the ones proved in \cite{MR4232783,MR4315017}. For example, there is an integral of $|\phi|^2$ on the right hand side of the estimate in Lemma \ref{lp1}, while such an integral does not appear in \cite[Lemma 5.4]{MR4232783}. This is because our solutions here do not necessarily vanish inside a light cone $r-t=\text{constant}$. The most important result in Section \ref{sec5.3} is Lemma \ref{lp2}. In its proof, we need to use the asymptotics of the approximate optical function $q$, i.e.\ Lemma \ref{lp2.l1}. With the help of Lemma \ref{lp2}, we obtain Lemma \ref{lp3}. Note that Lemma \ref{lp3} is the reason why we want to include the approximate optical function in Poincar$\acute{\rm e}$'s lemmas. If one uses Lemma \ref{lp1} or \ref{lp1st} instead of Lemma \ref{lp2}, then the factor $\eps t^{-1}$ on the right hand side of the estimate in Lemma \ref{lp3} needs to be replaced by $\eps t^{-1+C\eps}$. This causes a serious problem because in this case we cannot end our continuity argument in Section \ref{sec6}.

\subsection{Energy estimates when $t\geq T_\eps=1/\eps$}\label{sec5.2}
Assume that $u=u(t,x)$ is a function  satisfying the pointwise estimates: for all $t\geq T_\eps$ with $\eps\ll 1$ we have
\begin{equation}\label{sec5a1}\sum_{|I|\leq 2}|Z^Iu|\leq B \eps t^{-1+B\eps},\qquad |\partial u(t,x)|\leq B \eps t^{-1};\end{equation}
if $|q(t,x)|\leq t^{1/4}$ and $t\geq T_\eps$, we have 
\begin{equation}\label{sec5a2}\sum_{l\leq 2}|\partial^l (u-\eps r^{-1}U)|\leq B\eps t^{-5/4+B\eps}.\end{equation}Recall that $U=U(t,x)$ is defined in Section \ref{sec4}. In Section \ref{sec6} we will check these estimates when we apply the energy estimates in the continuity argument.

Note that we do not need to assume that $u$ is a solution to \eqref{qwe}. In fact, in this section, the semilinear terms $f^{(I)}$ in \eqref{qwe} never appear. We also emphasize that, while $B$ denoted a vector field when we discussed the Euler equations in Section \ref{sec1.4.2}, it is always a constant in Section \ref{sec5}. We do  not use $C$ to denote the constants in \eqref{sec5a1} and \eqref{sec5a2} because in the main results of this section, it is important to determine whether each constant depends on the constants in \eqref{sec5a1} and \eqref{sec5a2} or not. For example, in Proposition \ref{prop5.1}, we want to show that the implicit constant does not depend on $B$.

Fix a smooth function $\phi(t,x)$ such that $\phi(t)\in C_c^\infty(\R^3)$ for each $t\geq T_\eps$. Define the energy \begin{equation}\label{sec5ene}\begin{aligned}
E_u(\phi)(t)&=\int_{\R^3}w(t,x)(-2g^{0\alpha}(u,\partial u)\phi_t\phi_\alpha+g^{\alpha\beta}(u,\partial u)\phi_\alpha\phi_\beta)(t,x)\ dx\\&=\int_{\R^3}w(t,x)(|\partial \phi|^2-2(g^{0\alpha}(u,\partial u)-m^{0\alpha})\phi_t\phi_\alpha-(g^{\alpha\beta}(u,\partial u)-m^{\alpha\beta})\phi_\alpha\phi_\beta)(t,x)\ dx.
\end{aligned}\end{equation}
The weight function $w$ is defined by\begin{equation}\label{sec5wei}w=w_{\gamma_1,\gamma_2}(t,x):=m(q(t,x))\cdot\exp\kh{c_0 \eps\ln t\cdot \sigma(q(t,x))}\end{equation}with
\eq{\label{sec5wei2}m(q)=m_{\gamma_1}(q):=1_{q\geq 0}+(1-q)^{\gamma_1} 1_{q<0},}
\eq{\label{sec5wei3}\sigma(q)=\sigma_{\gamma_2}(q):=(1+q)^{-\gamma_2} 1_{q\geq 0}+(2-(1-q)^{-\gamma_2})1_{q<0}.}
Here $q(t,x)$ is originally defined in Section \ref{sec4} and  extended at the beginning of  Section \ref{sec5}, so $w$ is  defined whenever $t\geq T_\eps$. Moreover,  $c_0\gg  1$ is a large constant to be chosen later, and  $\gamma_1,\gamma_2$ are two arbitrary constants  such that \eq{\label{sec5wei4}2<\gamma_1<\min\{2(\gamma_--1),4\},\qquad 0<\gamma_2<\min\{\gamma_--2,\gamma_+-1,1/2\}} where $\gamma_\pm$ comes from Definition \ref{def3.1}. Since $\gamma_->2$ and $\gamma_+>1$, there exist an infinite number of pairs of $(\gamma_1,\gamma_2)$ satisfying \eqref{sec5wei4}. We emphasize that $c_0,\gamma_1,\gamma_2$ are independent of the constant $B$ in our assumptions \eqref{sec5a1} and \eqref{sec5a2}.

We  remark that this type of  weight $w$ is motivated by several previous papers on the long time dynamics of quasilinear wave equations such as  \cite{MR3419884,MR2382144,MR2134337,MR2680391,MR2003417,MR2666888}, etc.

Our goal  is to prove the following energy estimates. For convenience, we set
\eq{\label{sec5eneq}E_q(\phi)(t)&:= \sum_{\alpha=0}^3\int_{\R^3}((1-q)^{-1}1_{q<0}+\eps\ln\eps^{-1}(1+|q|)^{-1-\gamma_2})|q_t|^{-1}|(q_t\partial_\alpha-q_\alpha \partial_t)\phi(t)|^2 w\ dx.}
From now on, we use $C_B$ to denote a constant which depends on $B$, and $C$ to denote a constant which is independent of $B$. Similarly for $\lesssim_B$ and $\lesssim$.

\prop{\label{prop5.1}For $T_\eps\leq t\leq T$ with $\eps\ll_{c_0,\gamma_1,\gamma_2,\delta} 1$, and for all smooth functions $\phi=\phi(t,x)$ with $\phi(t)\in C_c^\infty(\R^3)$ for each $t$, we have \eq{\label{energymain}
&E_u(\phi)(t)+\int_t^TE_q(\phi)(\tau) d\tau\lesssim_{\gamma_1,\gamma_2} E_u(\phi)(T)\\
&
\qquad+\int_t^T\norm{g^{\alpha\beta}(u,\partial u)\partial_\alpha\partial_\beta\phi(\tau)}_{L^2(w)}\norm{\partial\phi(\tau)}_{L^2(w)}+(C_B\eps \tau^{-1}+\tau^{-17/16})\norm{\partial\phi}^2_{L^2(w)}\ d\tau.}
Here $\norm{f}_{L^2(w)}^2:=\int_{\R^3}|f|^2w\ dx$.
Note that the implicit constant in \eqref{energymain}  depends on $\gamma_1,\gamma_2$, but not on $c_0,\delta,B,C_B$. 
}\rm\bigskip

The proof starts with a computation of $\frac{d}{dt}E_u(\phi)(t)$. For simplicity, we write $ g^{\alpha\beta}= g^{\alpha\beta}(u,\partial u)$. Then, by applying integration by parts, we have 
\fm{&\frac{d}{dt}E_u(\phi)(t)\\&=\int_{\R^3}w_t(-2 g^{0\alpha}\phi_t\phi_\alpha+ g^{\alpha\beta}\phi_\alpha\phi_\beta)\\&\hspace{3em}+w(-2 g^{0\alpha}\phi_{tt}\phi_\alpha-2 g^{0\alpha}\phi_{t}\phi_{\alpha t}-2\partial_t g^{0\alpha}\phi_t\phi_\alpha+2 g^{\alpha\beta}\phi_{\alpha t}\phi_\beta+\partial_t g^{\alpha\beta}\phi_\alpha\phi_\beta)\ dx\\&=\int_{\R^3}w_t(-2 g^{0\alpha}\phi_t\phi_\alpha+ g^{\alpha\beta}\phi_\alpha\phi_\beta)+w(-2 g^{0\alpha}\phi_{\alpha t}\phi_t+2 g^{i\beta}\phi_{it}\phi_{\beta}-2\partial_t g^{0\alpha}\phi_t\phi_\alpha+\partial_t g^{\alpha\beta}\phi_\alpha\phi_\beta)\ dx\\&=\int_{\R^3}w_t(-2 g^{0\alpha}\phi_t\phi_\alpha+ g^{\alpha\beta}\phi_\alpha\phi_\beta)-2w_i g^{i\beta}\phi_t\phi_\beta\\&\hspace{3em}+w(-2 g^{0\alpha}\phi_{\alpha t}\phi_t-2 g^{i\beta}\phi_{t}\phi_{i\beta}-2\partial_t g^{0\alpha}\phi_t\phi_\alpha-2\partial_i g^{i\beta}\phi_t\phi_\beta+\partial_t g^{\alpha\beta}\phi_\alpha\phi_\beta)\ dx\\&=\int_{\R^3}w_t g^{\alpha\beta}\phi_\alpha\phi_\beta+w(-2 g^{\alpha\beta}\phi_{\alpha \beta}\phi_t-2\partial_\alpha g^{\alpha\beta}\phi_t\phi_\beta+\partial_t g^{\alpha\beta}\phi_\alpha\phi_\beta)-2w_\alpha  g^{\alpha\beta}\phi_t\phi_\beta\ dx.}

By setting $\wt{T}_\alpha:=q_t\partial_\alpha-q_\alpha\partial_t$, we have $\phi_\alpha=q_t^{-1}(\wt{T}_\alpha \phi+q_\alpha\phi_t)$. We also have
\fm{w_t&=(\frac{m'(q)}{m(q)}+c_0\eps\ln t\cdot\sigma'(q)) q_tw+  c_0\eps t^{-1}\sigma(q)w,\\
w_i&=(\frac{m'(q)}{m(q)}+c_0\eps\ln t\cdot\sigma'(q)) q_iw.} Thus,
\fm{g^{\alpha\beta}\phi_\alpha\phi_\beta q_t-2 g^{\alpha\beta}\phi_t\phi_\beta q_\alpha&= g^{\alpha\beta}q_t^{-1}(\wt{T}_\alpha\phi+q_\alpha\phi_t)(\wt{T}_\beta\phi+q_\beta\phi_t)-2 g^{\alpha\beta}q_\alpha\phi_tq_t^{-1}(\wt{T}_\beta\phi+q_\beta\phi_t)\\&= g^{\alpha\beta}q_t^{-1}\wt{T}_\alpha\phi \wt{T}_\beta\phi- g^{\alpha\beta}q_t^{-1}q_\alpha q_\beta\phi_t^2}
and 
\fm{&w_t g^{\alpha\beta}\phi_\alpha\phi_\beta-2w_\alpha  g^{\alpha\beta}\phi_t\phi_\beta\\&=(\frac{m'(q)}{m(q)}+c_0\eps\ln t\cdot\sigma'(q)) w( g^{\alpha\beta}q_t^{-1}\wt{T}_\alpha\phi \wt{T}_\beta\phi- g^{\alpha\beta}q_t^{-1}q_\alpha q_\beta\phi_t^2)+c_0\eps t^{-1}\sigma(q)w(-g^{00}\phi_t^2+ g^{ij}\phi_i\phi_j).}
Since $\wt{T}_0=0$ and that $(g^{ij})=(\delta_{ij}+O(B\eps t^{-1+B\eps}))$, we have $g^{\alpha\beta}\wt{T}_\alpha \phi\wt{T}_\beta\phi\geq (1-O(\eps))|\wt{T}\phi|^2$ as long as $\eps\ll_B 1$ and $t\geq T_\eps$. By Lemma \ref{lem4.2} and Remark \ref{rmk4.2.1} we have \fm{q_t=(\mu+\nu)/2\leq -C^{-1}t^{-C\eps}+C\eps (t+r)^{-1+C\eps}<0.}  Moreover, we have  $ t\geq 1/\eps$ and thus
\fm{\frac{m'(q)}{m(q)}+c_0\eps \ln t\cdot \sigma'(q)&=-\gamma_1(1-q)^{-1}1_{q<0}-c_0\eps\ln t\cdot \gamma_2(1+|q|)^{-1-\gamma_2}\\&\leq -\gamma_1(1-q)^{-1}1_{q<0}-c_0\eps\ln\eps^{-1}\cdot \gamma_2(1+|q|)^{-1-\gamma_2}\leq 0.}
We  conclude that \fm{&(\frac{m'(q)}{m(q)}+c_0\eps \ln t\cdot \sigma'(q))w g^{\alpha\beta}q_t^{-1}\wt{T}_\alpha\phi \wt{T}_\beta\phi\\
&\geq [\gamma_1(1-q)^{-1}1_{q<0}+c_0\eps\ln\eps^{-1}\cdot \gamma_2(1+|q|)^{-1-\gamma_2}]|q_t|^{-1} w |\wt{T}\phi|^2.}

In addition, we have the following  lemma. 

\lem{For all $t\geq T_\eps$ with $\eps\ll_{c_0,B}1$, we have \fm{ |(\frac{m'(q)}{m(q)}+c_0\eps \ln t\cdot \sigma'(q)) w g^{\alpha\beta}q_t^{-1}q_\alpha q_\beta\phi_t^2|\lesssim   t^{-17/16}\sigma(q)w\phi_t^2.}
The implicit constant here is independent of $B$.}
\begin{proof}
First, we recall that
\fm{&|\frac{m'(q)}{m(q)}+c_0\eps \ln t\cdot \sigma'(q)|\\&=\gamma_1(1-q)^{-1}1_{q<0}+c_0\eps\ln t\cdot \gamma_2(1+|q|)^{-1-\gamma_2}\\
&\leq [(\gamma_1(1-q)^{-1}+c_0\eps\ln t\cdot \gamma_2(1-q)^{-1-\gamma_2})1_{q<0}+c_0\gamma_2\eps\ln t\cdot (1+q)^{-1}1_{q\geq 0}]\sigma(q)\\
&\lesssim_{\gamma_1,\gamma_2} (1_{q<0}+c_0\eps\ln t)\lra{q}^{-1}\sigma(q).}
Here we recall that $\sigma(q)\in[1,2]$ whenever $q<0$. 

Now we suppose $|q(t,x)|\leq t^{1/4}$. By \eqref{lem4.1c4} in Lemma \ref{lem4.1}, we have $(t,x)\in\D$ as long as $\eps\ll 1$.   By Proposition \ref{prop4.6} and \eqref{sec5a2}, we have
\fm{&|g^{\alpha\beta}(u,\partial u)q_\alpha q_\beta|\\&\leq |(g^{\alpha\beta}(u,\partial u)-g^{\alpha\beta}(\eps r^{-1}U,\partial (\eps r^{-1}U)))q_\alpha q_\beta|+|g^{\alpha\beta}(\eps r^{-1}U,\partial (\eps r^{-1}U))q_\alpha q_\beta|\\
&\lesssim (|u-\eps r^{-1}U|+|\partial (u-\eps r^{-1}U)|)\cdot|\partial q|^2+t^{-2+C\eps}\lra{r-t}\\&\lesssim  C_B\eps t^{-5/4+C_B\eps}+t^{-2+C\eps}\lra{q}\lesssim C_B\eps t^{-5/4+C_B\eps}+t^{-7/4+C\eps}.}
Also recall from Lemma \ref{lem4.5} that $\partial q\in S^{0,0}$. Thus,
\fm{ &|(\frac{m'(q)}{m(q)}+c_0\eps \ln t\cdot \sigma'(q))  g^{\alpha\beta}q_t^{-1}q_\alpha q_\beta|\\
&\lesssim (1_{q<0}+c_0\eps\ln t)\lra{q}^{-1}\sigma(q)\cdot t^{C\eps}\cdot(C_B\eps t^{-5/4+C_B\eps}+t^{-7/4+C\eps})\\
&\lesssim (1+c_0\eps t^{1/8})(C_B\eps t^{-5/4+C_B\eps}+t^{-7/4+C\eps})\sigma(q)\lesssim  t^{-9/8+C_B\eps}\sigma(q)\lesssim t^{-17/16}\sigma(q).}
Here we note that $\ln (t)\lesssim t^{1/8}$ and that $|q_t|\gtrsim t^{-C\eps}$ by Remark \ref{rmk4.2.1}. In the last two steps, we choose $\eps\ll_{B,c_0}1$ so that $c_0\eps\leq 1$ and $C_B\eps<1/16$.

Next we suppose $|q(t,x)|\geq t^{1/4}$. If $|x|\geq t/2$, by  Lemma \ref{lem4.2} (instead of Lemma \ref{lem4.5}) we have $\partial q=O(t^{C\eps})$ and $|q_t|\gtrsim t^{-C\eps}$. Thus, by Lemma \ref{lem4.2} and  \eqref{lem4.6f2}
\fm{|g^{\alpha\beta}q_\alpha q_\beta|&\lesssim  |m^{\alpha\beta}q_\alpha q_\beta|+(|u|+|\partial u|)|\partial q|^2\lesssim |\mu\nu|+\sum_{i}|\lambda_i|^2+B\eps t^{-1+(B+C)\eps}\\
&\lesssim t^{-1+C\eps}\lra{\max\{0,-q\}}^{1-\gamma_-}+B\eps t^{-1+C_B\eps}\lesssim t^{-1+C_B\eps}.}
If $|x|<t/2$, then  $m^{\alpha\beta}q_\alpha q_\beta =m^{\alpha\beta}\wh{\omega}_\alpha \wh{\omega}_\beta=0$, so \fm{|g^{\alpha\beta}q_\alpha q_\beta|&\lesssim  |m^{\alpha\beta}q_\alpha q_\beta|+(|u|+|\partial u|)|\partial q|^2\lesssim 0+B\eps t^{-1+B\eps}\lesssim t^{-1+B\eps}.}
For $\eps\ll_B1$, it follows that
\fm{&|(\frac{m'(q)}{m(q)}+c_0\eps \ln t\cdot \sigma'(q))  g^{\alpha\beta}q_t^{-1}q_\alpha q_\beta|\\
&\lesssim  (1_{q<0}+c_0\eps\ln t)\lra{q}^{-1}\sigma(q)\cdot  t^{-1+C_B\eps}\lesssim (1+c_0\eps t^{1/8})t^{-1/4}\cdot  t^{-1+C_B\eps}\sigma(q)\\&\lesssim t^{-9/8+C_B\eps}\sigma(q) \lesssim t^{-17/16} \sigma(q).}
In the   last two steps, we choose $\eps_{B,c_0}1$ so that $C_B\eps<1/16$, $c_0\eps\leq 1$ and  $\ln t\lesssim  t^{1/8}$.

\end{proof}\rm\bigskip

Moreover, we have \fm{-g^{00}\phi_t^2+g^{ij}\phi_i\phi_j&=|\partial\phi|^2+O(|u||\partial\phi|^2)\sim |\partial\phi|^2.}
In summary,  we have
\eq{\label{lem5.3f1}&w_t g^{\alpha\beta}\phi_\alpha\phi_\beta-2w_\alpha  g^{\alpha\beta}\phi_t\phi_\beta\\&\geq[\gamma_1(1-q)^{-1}1_{q<0}+c_0\eps\ln\eps^{-1}\cdot \gamma_2(1+|q|)^{-1-\gamma_2}]|q_t|^{-1} w |\wt{T}\phi|^2\\&\quad+(C^{-1}c_0\eps t^{-1}-Ct^{-17/16})\sigma(q)w|\partial \phi|^2.}

To finish the proof of Proposition \ref{prop5.1}, we need to control $w(-2\partial_\alpha g^{\alpha\beta}\phi_t\phi_\beta+\partial_tg^{\alpha\beta}\phi_\alpha\phi_\beta)$. This is done in the next lemma. We remark that  the assumption \eqref{def3.1a51} is introduced in Definition \ref{def3.1} because of this lemma.
\lem{\label{lem5.4}We have \fm{w(-2\partial_\alpha g^{\alpha\beta}\phi_t\phi_\beta+\partial_tg^{\alpha\beta}\phi_\alpha\phi_\beta)\geq -C\eps t^{-1+C\eps}\lra{q}^{-2-2\gamma_{\sgn(q)}}w|\wt{T}\phi|^2-C_B\eps t^{-1}w|\partial\phi|^2.}
}
\begin{proof}
We have
\fm{\partial_\alpha (g^{\beta\beta'}(u,\partial u))&=g^{\beta\beta'}_J\partial_\alpha u^{(J)}+g^{\beta\beta'\lambda}_J\partial_\alpha\partial_\lambda u^{(J)}+O(C_B\eps^2 t^{-2+C_B\eps})=g^{\beta\beta'\lambda}_J\partial_\alpha\partial_\lambda u^{(J)}+O(C_B\eps  t^{-1}).}

If $|q|\geq t^{1/4}$, then  we  have either $|q|=|r-t|\geq t^{1/4}$ or $\lra{r-t}\gtrsim \lra{q}t^{-C\eps}\geq t^{1/4-C\eps}$ by \eqref{lem4.1c3}. It follows that
\fm{|\partial^2u|\lesssim\lra{r-t}^{-2}\sum_{|J|\leq 2}|Z^Ju|\lesssim  t^{-1/2+C\eps}\cdot C_B\eps t^{-1+C_B\eps}\lesssim C_B\eps t^{-3/2+C_B\eps}\lesssim C_B\eps t^{-1}.}
We thus have $\partial g^{**}=O(C_B\eps t^{-1})$ in this case.

Now we assume that  $|q|\leq t^{1/4}$, so $(t,x)\in\D$ by \eqref{lem4.1c4}. By \eqref{sec5a2}, \eqref{lem4.10c0} and \eqref{lem4.10c},  we have
\fm{&-2\partial_\alpha g^{\alpha\beta}\phi_t\phi_\beta+\partial_tg^{\alpha\beta}\phi_\alpha\phi_\beta\\
&=-2 g_J^{\alpha\beta\lambda}\partial_\alpha\partial_\lambda u^{(J)}\phi_t\phi_\beta+g^{\alpha\beta\lambda}_J\partial_t\partial_\lambda u^{(J)}\phi_\alpha\phi_\beta+O((|\partial u|+|\partial^2u|^2)|\partial\phi|^2)
\\&=-2 g_J^{\alpha\beta\lambda}\partial_\alpha\partial_\lambda u_{app}^{(J)}\phi_t\phi_\beta+g^{\alpha\beta\lambda}_J\partial_t\partial_\lambda u_{app}^{(J)}\phi_\alpha\phi_\beta+O(|\partial^2(u-u_{app})||\partial\phi|^2+C_B\eps t^{-1}|\partial\phi|^2)\\
&=O(C_B\eps t^{-1}|\partial\phi|^2+\eps t^{-1+C\eps}|\wt{T}\phi||\partial\phi|\lra{q}^{-1-\gamma_{\sgn(q)}}).}
In the last step we apply Lemma \ref{lem4.10} where the constant $B$ is not involved.
We end the proof by noticing that
\fm{\eps t^{-1+C\eps}|\wt{T}\phi||\partial\phi|\lra{q}^{-1-\gamma_{\sgn(q)}}\lesssim \eps t^{-1}|\partial\phi|^2+\eps t^{-1+C\eps}\lra{q}^{-2-2\gamma_{\sgn(q)}}|\wt{T}\phi|^2.}
\end{proof}\rm

\bigskip

In conclusion, 
\fm{&\frac{d}{dt}E_u(\phi)(t)\\&\geq \int_{\R^3}-2w\phi_tg^{\alpha\beta}\partial_\alpha\partial_\beta \phi-C\eps t^{-1+C\eps}\lra{q}^{-2-2\gamma_{\sgn(q)}}w|\wt{T}\phi|^2-C_B\eps t^{-1}w|\partial\phi|^2\\&\qquad+[\gamma_1(1-q)^{-1}1_{q<0}+c_0\eps\ln\eps^{-1}\cdot \gamma_2\lra{q}^{-1-\gamma_2}]|q_t|^{-1} w |\wt{T}\phi|^2-Ct^{-17/16}\eps  t^{-1}\sigma (q)w|\partial\phi|^2\ dx\\&\geq -2\norm{\partial\phi(t)}_{L^2(w)}\norm{g^{\alpha\beta}\partial_\alpha\partial_\beta \phi}_{L^2(w)}-(C_B\eps t^{-1}+Ct^{-17/16})\norm{\partial\phi}_{L^2(w)}^2\\&\quad+\int_{\R^3}[\gamma_1(1-q)^{-1}1_{q<0}+ \gamma_2\eps\ln\eps^{-1}\cdot\lra{q}^{-1-\gamma_2}]|q_t|^{-1} w |\wt{T}\phi|^2\ dx.}
To obtain the last estimate, we compute the coefficients of $w|\wt{T}\phi|^2$:
\fm{&-C\eps t^{-1+C\eps}\lra{q}^{-2-2\gamma_{\sgn(q)}}+[\gamma_1(1-q)^{-1}1_{q<0}+c_0\eps\ln \eps^{-1}\cdot \gamma_2\lra{q}^{-1-\gamma_2}]|q_t|^{-1} \\
&\geq[-C\eps t^{-1+C\eps}\lra{q}^{-2-2\gamma_{\sgn(q)}}+\gamma_1(1-q)^{-1}1_{q<0}+c_0\eps\ln \eps^{-1}\cdot \gamma_2\lra{q}^{-1-\gamma_2}]|q_t|^{-1}\\
&\geq [\gamma_1(1-q)^{-1}1_{q<0}+(-C\eps  +c_0\eps\ln \eps^{-1}\cdot \gamma_2)\lra{q}^{-1-\gamma_2}]|q_t|^{-1}.}
Now we choose $c_0\geq 2$ and then choose $\eps\ll1$ so that $C\leq \gamma_2\ln\eps^{-1}$.

Finally, we integrate the inequality for $dE_u(\phi)/dt$ with respect to $t$ on $[t,T]$ and we conclude~\eqref{energymain}.

\subsection{Energy estimates when $0\leq t\leq  T_\eps$}\label{sec5.1}
Assume that $u=u(t,x)$ is a function satisfying the pointwise estimates: for all $0\leq t\leq T_\eps$ with $\eps\ll 1$ we have
\begin{equation}\label{sec5a3}\sum_{|I|\leq 2}|Z^Iu|\leq B \eps\lra{t}^{-1/2} .\end{equation}
Again, we do not need to assume that $u$ is a solution to \eqref{qwe}.
In Section \ref{sec6} we will check this  estimate when we apply the energy estimates.

Recall that $m=m_{\gamma_1}(q)$ is defined by \eqref{sec5wei2}. Then, we have the following energy estimates.

\prop{\label{prop5.1st} For $0\leq t\leq t'\leq T_\eps$ and for all smooth functions $\phi=\phi(t,x)$ with $\phi(t)\in C_c^\infty(\R^3)$ for each $t$, we have \eq{\label{energymain2}&\norm{\partial\phi(t)}_{L^2(m)}^2\\&\lesssim_{\gamma_1}\norm{\partial\phi(t')}_{L^2(m)}^2+\int_t^{t'}\norm{g^{\alpha\beta}(u,\partial u)\partial_\alpha\partial_\beta\phi(\tau)}_{L^2(m)}\norm{\partial\phi(\tau)}_{L^2(m)}+C_B\eps \lra{\tau}^{-1/2}\norm{\partial\phi(\tau)}^2_{L^2(m)}\ d\tau.}
Here $\norm{f}_{L^2(m)}^2:=\int_{\R^3}|f|^2m(r-t)\ dx$.
Note that the implicit constant here depends on $\gamma_1$ but not on $B$. }\rm

\bigskip 

The proof of \eqref{energymain2} is similar to and simpler than that of \eqref{energymain}. We first define an energy
\eq{\label{sec5.4wei} E_{m,u}(\phi)(t):=\int_{\R^3}m(|x|-t)(-2g^{0\alpha}(u,\partial u)\phi_t\phi_\alpha+g^{\alpha\beta}(u,\partial u)\phi_\alpha\phi_\beta)(t,x)\ dx.}
By \eqref{sec5a3}, we have $E_{m,u}(\phi)\sim \norm{\partial\phi}_{L^2(m)}^2$. Note that \fm{&\frac{d}{dt}E_{m,u}(\phi)(t)\\&=\int_{\R^3}-m'(|x|-t)(-2 g^{0\alpha}\phi_t\phi_\alpha+ g^{\alpha\beta}\phi_\alpha\phi_\beta)\\&\hspace{3em}+m(|x|-t)(-2 g^{0\alpha}\phi_{tt}\phi_\alpha-2 g^{0\alpha}\phi_{t}\phi_{\alpha t}-2\partial_t g^{0\alpha}\phi_t\phi_\alpha+2 g^{\alpha\beta}\phi_{\alpha t}\phi_\beta+\partial_t g^{\alpha\beta}\phi_\alpha\phi_\beta)\ dx\\&=\int_{\R^3}-m'(|x|-t)(-2 g^{0\alpha}\phi_t\phi_\alpha+ g^{\alpha\beta}\phi_\alpha\phi_\beta)\\&\hspace{3em}+m(|x|-t)(-2 g^{0\alpha}\phi_{\alpha t}\phi_t+2 g^{i\beta}\phi_{it}\phi_{\beta}-2\partial_t g^{0\alpha}\phi_t\phi_\alpha+\partial_t g^{\alpha\beta}\phi_\alpha\phi_\beta)\ dx\\&=\int_{\R^3}-m'(|x|-t)(-2 g^{0\alpha}\phi_t\phi_\alpha+ g^{\alpha\beta}\phi_\alpha\phi_\beta)-2m'(|x|-t)\omega_i g^{i\beta}\phi_t\phi_\beta\\&\hspace{3em}+m(|x|-t)(-2 g^{0\alpha}\phi_{\alpha t}\phi_t-2 g^{i\beta}\phi_{t}\phi_{i\beta}-2\partial_t g^{0\alpha}\phi_t\phi_\alpha-2\partial_i g^{i\beta}\phi_t\phi_\beta+\partial_t g^{\alpha\beta}\phi_\alpha\phi_\beta)\ dx\\&=\int_{\R^3}-m'(|x|-t) g^{\alpha\beta}\phi_\alpha\phi_\beta+m(|x|-t)(-2 g^{\alpha\beta}\phi_{\alpha \beta}\phi_t-2\partial_\alpha g^{\alpha\beta}\phi_t\phi_\beta+\partial_t g^{\alpha\beta}\phi_\alpha\phi_\beta)\\&\hspace{3em}-2m'(|x|-t)\wh{\omega}_\alpha g^{\alpha\beta}\phi_t\phi_\beta\ dx.}
By setting $T_\alpha=\partial_\alpha+\wh{\omega}_\alpha\partial_t$, we have
\fm{&-m'(|x|-t) g^{\alpha\beta}\phi_\alpha\phi_\beta-2m'(|x|-t)\wh{\omega}_\alpha g^{\alpha\beta}\phi_t\phi_\beta\\
&=m'(|x|-t)(-g^{\alpha\beta}T_\alpha\phi T_\beta\phi+g^{\alpha\beta}\wh{\omega}_\alpha\wh{\omega}_\beta \phi^2_t).}
Here we have 
\fm{g^{\alpha\beta}\wh{\omega}_\alpha\wh{\omega}_\beta&=m^{\alpha\beta}\wh{\omega}_\alpha\wh{\omega}_\beta+O((|u|+|\partial u|))=O(B\eps\lra{t}^{-1/2}).}
Since $m'(|x|-t)\leq 0$ and since $g^{\alpha\beta}T_\alpha\phi T_\beta\phi=g^{ij}T_i\phi T_j\phi\geq (1-C\eps)|T\phi|^2$, we conclude that
\fm{-m'(|x|-t) g^{\alpha\beta}\phi_\alpha\phi_\beta-2m'(|x|-t)\wh{\omega}_\alpha g^{\alpha\beta}\phi_t\phi_\beta\geq -CB\eps\lra{t}^{-1/2}|m'(|x|-t)|\phi_t^2.}
Since $m'(q)=-\gamma_1(1-q)^{\gamma_1-1}1_{q<0}=O(\lra{q}^{-1}m(q))$ and since 
\fm{|-2\partial_\alpha g^{\alpha\beta}\phi_t\phi_\beta+\partial_t g^{\alpha\beta}\phi_\alpha\phi_\beta|\lesssim(|\partial u|+|\partial^2u|)|\partial\phi|^2\lesssim B\eps \lra{t}^{-1/2}|\partial\phi|^2,}
we conclude that
\fm{\frac{d}{dt}E_{m,u}(\phi)(t)&\geq \int_{\R^3}-2m(|x|-t) g^{\alpha\beta}\phi_{\alpha \beta}\phi_t-CB\eps \lra{t}^{-1/2}|\partial\phi|^2m(|x|-t)\ dx\\
&\geq -2\norm{ g^{\alpha\beta}\phi_{\alpha \beta}(\tau)}_{L^2(m)}\norm{\phi_t(\tau)}_{L^2(m)}-CB\eps \lra{t}^{-1/2}\norm{\partial\phi(\tau)}^2_{L^2(m)}.}
Integrate $dE_{m,u}(\phi)/dt$ with respect to $t$ and we conclude \eqref{energymain2}.

\subsection{Poincar$\acute{\bf e}$'s lemmas}\label{sec5.3}

Our first result in this section is the following version of  Poincar$\acute{\rm e}$'s lemmas. We emphasize that here we assume $t\geq 100(\eta+1)^2/\eta^2$ and that there is no $\eps$ involved. Note that our assumptions here are weaker than those of Lemma 5.4 in \cite{MR4232783}, so our conclusion here is also weaker than that in   \cite{MR4232783}.

\lem{\label{lp1} Fix $\eta>0$. Then, for all $c\ll_\eta 1$ and all smooth functions $\phi=\phi(t,x)$ with $\phi(t)\in C^\infty_c(\R^3)$,  we have\fm{
&\int_{x\in\R^3:\ |x|\geq t}\lra{t-r}^{-2}|\phi|^2\ dx+\int_{x\in\R^3:\ |x|< t}\lra{t-r}^{\eta-1}|\phi|^21_{r/t\geq 1-c/2}\ dx\\
&\lesssim_{\eta,c} \int_{x\in\R^3:\ |x|\geq t}|\partial\phi|^2\ dx+\int_{x\in\R^3:\ |x|< t}\lra{t-r}^{\eta+1}(|\partial\phi|^2+r^{-2}|\phi|^21_{r/t\in[1-c,1-c/2]})\ dx}
whenever $t\geq 100(\eta+1)^2/\eta^2$.}
\begin{proof} In this proof, for each $q\in\R$ we set 
\eq{\kappa(q):=-(1+q)^{-1}\cdot 1_{q\geq 0}-(1-q)^{\eta}\cdot 1_{q<0}.}
It is clear that
\eq{\label{lp1f2}\kappa'(q)&=(1+q)^{-2}\cdot 1_{q> 0}+\eta(1-q)^{\eta-1}\cdot 1_{q<0}.}
It follows that 
\fm{\partial_r(r^2\kappa(r-t)\phi^2)&=2r\kappa(r-t)\phi^2+r^2\kappa'(r-t)\phi^2+2r^2\kappa(r-t)\phi\phi_r.}
Because $\phi(t)$ has a compact support, we have
\fm{0&=\int_{\mathbb{S}^2}\int_0^\infty2r\kappa(r-t)\phi^2+r^2\kappa'(r-t)\phi^2+2r^2\kappa(r-t)\phi\phi_r\ drdS_\omega\\
&=\int_{\R^3}\kappa'(r-t)\phi^2\ dx+\int_{\R^3}2\kappa(r-t)\phi(r^{-1}\phi+\phi_r )\ dx.}
By the Cauchy-Schwarz inequality, we have
\fm{\int_{\R^3}\kappa'(r-t)\phi^2\ dx&= \left|\int_{\R^3}2\kappa(r-t)\phi(r^{-1}\phi+\phi_r )\ dx\right|\\
&\leq 2\kh{\int_{\R^3}\kappa'(r-t)\phi^2\ dx}^{1/2}\kh{\int_{\R^3}\frac{\kappa^2(r-t)}{\kappa'(r-t)}|r^{-1}\phi+\phi_r |^2\ dx}^{1/2}.}
As a result, we have
\eq{\label{lp1f3}\int_{\R^3}\kappa'(r-t)\phi^2\ dx&\leq4\int_{\R^3}\frac{\kappa^2(r-t)}{\kappa'(r-t)}|r^{-1}\phi+\phi_r |^2\ dx.}

To continue, we set \fm{\theta(q)=1_{q\geq 0}+(1-q)^{\eta+1}1_{q<0}.}
We can thus replace $\kappa^2/\kappa'$ with $(1+1/\eta)\theta$ in \eqref{lp1f3}. Note that $\kappa^2/\kappa'$ is not  continuous at $q=0$ and that $\theta$ is continuous for all $q$, so we prefer to use $\theta$. Now, since
\fm{\partial_r(r\theta(r-t)\phi^2)&=\theta(r-t)\phi^2+2r\theta(r-t)\phi\phi_r+r\theta'(r-t)\phi^2}
and since
\fm{\theta'&=-(1+\eta)(1-q)^{\eta}\cdot 1_{q<0},}
we have
\fm{0&=\int_{\R^3}\theta(r-t)(r^{-2}\phi^2+2r^{-1}\phi\phi_r)+r^{-1}\theta'(r-t)\phi^2\ dx\\
&=\int_{\R^3}\theta(r-t)(|r^{-1}\phi+\phi_r|^2-|\phi_r|^2)\ dx-\int_{\R^3:\ |x|<t}(\eta+1)(1+t-r)^{\eta}r^{-1}\phi^2\ dx\\
&=\int_{\R^3}\theta(r-t)(|r^{-1}\phi+\phi_r|^2-|\phi_r|^2)\ dx-\int_{\R^3:\ |x|<t}\frac{(\eta+1)(1+t-r)}{\eta r}\kappa'(r-t)\phi^2\ dx.}

To continue, we fix a small constant $0<c<1/2$ and temporarily assume that $\phi$ is supported in the region where $r>(1-c)t$. Then, whenever $r<t$ and $\phi(t,x)\neq 0$, we have
\fm{0\leq \frac{(\eta+1)(1+t-r)}{\eta r}\leq \frac{\eta+1}{\eta}(r^{-1}+\frac{c}{1-c})\leq \frac{\eta+1}{\eta}(r^{-1}+2c)\leq \frac{\eta}{50(\eta+1)}+\frac{2c(\eta+1)}{\eta}.}
To obtain the last estimate, we notice  that $r>(1-c)t\geq 50(\eta+1)^2/\eta^2$. By \eqref{lp1f3}, we have
\fm{\int_{\R^3}\kappa'(r-t)\phi^2\ dx&\leq 4(1+1/\eta)\int_{\R^3}\theta(r-t)(r^{-1}\phi+\phi_r)^2\ dx\\
&\leq 4(1+1/\eta)\int_{\R^3}\theta(r-t)\phi_r^2\ dx+(\frac{2}{25}+\frac{8c(\eta+1)^2}{\eta^2})\int_{x\in \R^3:\ |x|<t}\kappa'(r-t)\phi^2\ dx.}
We choose $c\ll_\eta1$ so that $\frac{8c(\eta+1)^2}{\eta^2}\leq 1/10$. Then we  conclude that
\fm{
&\int_{x\in\R^3:\ |x|\geq t}\lra{t-r}^{-2}|\phi|^2\ dx+\int_{x\in\R^3:\ |x|< t}\lra{t-r}^{\eta-1}|\phi|^2\ dx\\
&\lesssim_{\eta} \int_{x\in\R^3:\ |x|\geq t}|\partial\phi|^2\ dx+\int_{x\in\R^3:\ |x|< t}\lra{t-r}^{\eta+1}|\partial\phi|^2\ dx}
in this special case.

In general, if we do not add any assumption on the support of $\phi$, we choose a cutoff function $\chi\in C^\infty(\R)$ such that $0\leq \chi\leq 1$, $\chi|_{[1-c/2,\infty)}\equiv 1$ and $\chi|_{(-\infty,1-c]}\equiv 0$. Now, $\chi(r/t)\phi$ satisfies the support condition stated in the previous paragraph, so we obtain
\fm{&\int_{x:\in \R^3:\ |x|\geq t}\lra{r-t}^{-2}|\phi|^2\ dx+\int_{x\in\R^3:\ |x|<t}\lra{r-t}^{\eta-1}|\chi(r/t)\phi|^2 \ dx\\
&\lesssim \int_{x:\in \R^3:\ |x|\geq t}|\partial\phi|^2\ dx+\int_{x\in\R^3:\ |x|<t}\lra{r-t}^{\eta+1}|\partial(\chi(r/t)\phi)|^2 \ dx\\
&\lesssim \int_{x:\in \R^3:\ |x|\geq t}|\partial\phi|^2\ dx+\int_{x\in\R^3:\ |x|<t}\lra{r-t}^{\eta+1}(|\partial\phi|^2 +r^{-2}|\chi'(r/t)\phi|^2)\ dx.}
In the last estimate, we notice that $\chi'(r/t)=0$ unless $r/t\in[1-c,1-c/2]$. We also notice that $\chi(r/t)\geq 1_{r/t\geq 1-c/2}$, so we finish the proof.
\end{proof}
\rm

If $0\leq t\leq 100(\eta+1)^2/\eta^2$, we have the following version of Poincar$\acute{\rm e}$'s lemmas. 

\lem{\label{lp1st} Fix $\eta>0$. Then, for all smooth functions $\phi=\phi(t,x)$ with $\phi(t)\in C^\infty_c(\R^3)$,  we have\fm{
&\int_{x\in\R^3:\ |x|\geq t}\lra{t-r}^{-2}|\phi|^2\ dx+\int_{x\in\R^3:\ |x|< t}\lra{t-r}^{\eta-1}|\phi|^2\ dx\lesssim_{\eta} \int_{\R^3}|\partial\phi|^2\ dx}
whenever $0\leq t\leq 100(\eta+1)^2/\eta^2$.}
\begin{proof}
Set $\kappa=\kappa(q)$ as in the previous proof. Since the proof of \eqref{lp1f3} does not rely on the assumption $t\geq 100(\eta+1)^2/\eta^2$, we have \eqref{lp1f3} holds for all $t\geq 0$. Note that $\kappa^2/\kappa'=1_{q\geq 0}+\frac{1}{\eta}(1-q)^{\eta+1}1_{q<0}$ and that $\lra{r-t}\leq 1+t\lesssim_\eta 1$ whenever $0\leq t\leq 100(\eta+1)^2/\eta^2$ and $r\leq t$. Thus, we have
\fm{\int_{\R^3}\kappa'(r-t)\phi^2\ dx\lesssim_\eta \int_{\R^3}|r^{-1}\phi+\phi_r|^2\ dx\lesssim \int_{\R^3}|\partial\phi|\ dx.}
\end{proof}

\rm

\bigskip

We can also prove a weighted version of Poincar$\acute{\rm e}$'s lemmas. Note that the assumption \eqref{def3.1a4} is used in this lemma. Since the function $u$ is not involved in this lemma at all, all the implicit constants in this lemma are independent of $B$.

\lem{\label{lp2} For all $c\ll 1$ and all smooth functions $\phi=\phi(t,x)$ with $\phi(t)\in C^\infty_c(\R^3)$ for each $t\geq T_\eps$,  we have
\fm{
&\int_{\R^3}q_r^2\lra{q}^{-2}|\phi|^2w\cdot 1_{r/t\geq 1-c/2}\ dx\lesssim_{c} \int_{\R^3}(|\partial\phi|^2+r^{-2}|\phi|^21_{r/t\in[1-c,1-c/2]})w\ dx}
whenever $t\geq T_\eps$ with $\eps\ll1$.}\rm

\bigskip

Our proof relies on the asymptotics of the approximate optical function $q$.

\lem{\label{lp2.l1} Whenever  $t\geq T_\eps$, we have \eq{\label{lp2.l1c}&\partial_r(q_rw)\\
&\leq C\delta \lra{q}^{-1-\gamma_+}w|\mu|^21_{t/2<r<2t,\ q\geq 0}-\frac{\gamma_2}{16}c_0\eps\ln t\cdot \lra{q}^{-1-\gamma_2}w|\mu|^2-\frac{\gamma_1}{8}(1-q)^{-1}1_{q<0}w|\mu|^2.}
The constant $C$ in the first term is independent of the choice of $\delta\in(0,\delta_0)$.}
\begin{proof}
We first suppose that $r\leq t/2$, so $q=r-t$. In this case, $q_r=1$, so 
\fm{\partial_r(q_rw)&=(\frac{m'(r-t)}{m(r-t)}+c_0\eps\ln (t)\sigma'(r-t))w\\
&=(-\gamma_1(1-r+t)^{-1}1_{r<t}-c_0\eps\ln t\cdot \gamma_2(1+|r-t|)^{-1-\gamma_2})w.}
Since $|\mu|^2=4$ in this case, we obtain \eqref{lp2.l1c}.

Next suppose that $r>t/2$. The idea here is to express $\partial q$ and $\partial^2q$ in terms of $\mu$ and $\nu$. The terms involving $\nu$ are small, and we use \eqref{def3.1a4} to control the remaining terms. Since $q_r=(\nu-\mu)/2$, we have
\fm{\partial_r(q_rw)&=q_{rr}w+q_rw_r=\frac{\nu_q-\mu_q}{2}\cdot q_r w+(\frac{m'(q)}{m(q)}+c_0\eps\ln (t)\sigma'(q))q_r^2w\\
&=(\frac{1}{4}\mu\mu_q-\frac{1}{4}\mu\nu_q-\frac{1}{4}\mu_q\nu+\frac{1}{4}\nu\nu_q)\cdot w\\&\quad+(-\gamma_1(1-q)^{-1}1_{q<0}-c_0\eps\ln t\cdot \gamma_2(1+|q|)^{-1-\gamma_2})\cdot \frac{1}{4}(\nu-\mu)^2w\\
&=\frac{1}{4}\mu^2[\frac{\mu_q}{\mu}-\gamma_1(1-q)^{-1}1_{q<0}-c_0\eps\ln t\cdot \gamma_2(1+|q|)^{-1-\gamma_2}]w\\&\quad+O((|\nu\mu_q|+|\nu_q\mu|+|\nu\nu_q|+(\lra{q}^{-1}+c_0\eps\ln t\cdot \lra{q}^{-1-\gamma_2})|\nu|(|\mu|+|\nu|))w).}
Recall from \eqref{def3.1a4} that $|\mu_q|\lesssim \lra{q}^{-1-\gamma_{\sgn(q)}}|s\mu|$. Thus, as long as $c_0\gg1$, we have
\fm{\frac{\mu_q}{\mu}-c_0\eps\ln t\cdot \gamma_2(1+|q|)^{-1-\gamma_2}&\leq  C\lra{q}^{-1-\gamma_{\sgn(q)}}|\eps\ln t-\delta|-c_0\eps\ln t\cdot \gamma_2(1+|q|)^{-1-\gamma_2}\\
&\leq C\delta\lra{q}^{-1-\gamma_{\sgn(q)}}+(C-c_0\gamma_2)\eps\ln t\cdot (1+|q|)^{-1-\gamma_2}\\
&\leq C\delta\lra{q}^{-1-\gamma_{\sgn(q)}}-\frac{1}{2}c_0\eps\ln t\cdot \gamma_2(1+|q|)^{-1-\gamma_2}.}
We emphasize that the constant $C$ here is independent of the choice of $\delta$, because it comes from \eqref{def3.1a4}.
The second estimate holds because  $|\eps\ln t-\delta|\leq \delta+\eps\ln t$ and $\gamma_2<1/2<\min\{\gamma_+,\gamma_-\}$. The last estimate holds because  we can choose $c_0\gg1$ so that  $c_0 \gamma_2\geq 2C$.

To move forward, we apply Lemma \ref{lem4.2} and Lemma \ref{lem4.4} to control the terms involving $\nu$ and $\nu_q$. Note that these two lemmas are applicable whenever $t\geq T_\eps$ and $r\geq t/2$. We will control these terms in two cases: $r\leq 3t$, and $r> 2t$. If $r\leq 3t$, we have 
\fm{|\mu|\lesssim t^{C\eps},\ |\mu_q|\lesssim t^{C\eps}\lra{q}^{-1-\gamma_{\sgn(q)}},\ |\nu|\lesssim  \eps t^{-1+C\eps},\
|\nu_q|\lesssim \eps t^{-1+C\eps}\lra{q}^{-\gamma_{\sgn(q)}}+t^{-\gamma_-+C\eps}+\eps t^{-2+C\eps}.}
Since $\lra{q}\lesssim t^{1+C\eps}$, it follows that
\fm{|\nu\mu_q|+|\nu_q\mu|+|\nu\nu_q|&\lesssim
\eps t^{-1+C\eps}\lra{q}^{-\gamma_{\sgn(q)}}+t^{-\gamma_-+C\eps}+\eps t^{-2+C\eps}\\
&\lesssim \eps t^{\gamma_2-1+C\eps}(\lra{q}^{-\gamma_2-\gamma_{\sgn(q)}}+ \eps^{-1}t^{2-\gamma_-}\lra{q}^{-1-\gamma_2}+\lra{q}^{-1-\gamma_2})\\
&\lesssim  (\eps t^{\gamma_2-1+C\eps}+t^{\gamma_2+1-\gamma_-+C\eps})\lra{q}^{-1-\gamma_2}|\mu|^2;\\
(\lra{q}^{-1}+c_0\eps\ln t\cdot \lra{q}^{-1-\gamma_2})|\nu|(|\mu|+|\nu|)&\lesssim  \eps t^{-1+C\eps}\lra{q}^{-1}+c_0\eps^2\ln(t) t^{-1+C\eps}\lra{q}^{-1-\gamma_2}\\
&\lesssim (\eps t^{\gamma_2-1+C\eps}+c_0\eps^2\ln(t) t^{-1+C\eps})\lra{q}^{-1-\gamma_2}|\mu|^2.}
We have
\fm{&\eps t^{\gamma_2-1+C\eps}+t^{\gamma_2+1-\gamma_-+C\eps}+c_0\eps^2\ln(t) t^{-1+C\eps}\\&\lesssim (\frac{ t^{\gamma_2-1+C\eps}}{c_0\ln t}+\frac{t^{\gamma_2+1-\gamma_-+C\eps}}{c_0\eps\ln t}+\eps t^{-1+C\eps})\cdot c_0\eps \ln t\\&\lesssim ( t^{\gamma_2-1+C\eps}+t^{\gamma_2+2-\gamma_-+C\eps}+\eps )\cdot c_0\eps \ln t\lesssim \eps^{p}\cdot c_0\eps \ln t.}Here $p=p_{\gamma_2,\gamma_-}>0$ is a fixed constant depending only on $\gamma_2,\gamma_-$.
In the second last estimate, we notice that   $ \ln t\geq \ln \eps^{-1}\geq 1$. To obtain the last  estimate, we recall from \eqref{sec5wei4} that and $\gamma_2-1,\gamma_2-\gamma_-+2<0$, so the power of $t$ is less than a fixed negative constant as long as  $\eps\ll1$.   In summary, if $t/2\leq r\leq 3t$, we have
\fm{&\partial_r(q_rw)\\
&\leq C\delta\lra{q}^{-1-\gamma_{\sgn(q)}}w|\mu|^2+ (-\frac{\gamma_2}{8}+C\eps^{p})c_0\eps\ln t\cdot \lra{q}^{-1-\gamma_2}w|\mu|^2-\frac{\gamma_1}{4}(1-q)^{-1}1_{q<0}w|\mu|^2\\
&\leq C\delta\lra{q}^{-1-\gamma_{+}}w|\mu|^21_{q\geq 0}+ (-\frac{\gamma_2}{8}+C\eps^{p})c_0\eps\ln t\cdot \lra{q}^{-1-\gamma_2}w|\mu|^2-\frac{\gamma_1}{8}(1-q)^{-1}1_{q<0}w|\mu|^2.}
Recall that the constant $C$ in the first term is independent of the choice of $\delta$, so we choose  $\delta\ll_{\gamma_1} 1$ such that $C\delta\leq \gamma_1/8$. We also choose $\eps\ll1$ so that $C\eps^{p}<\gamma_2/16$.

Now suppose $r>2t$. In this case, we have  $0<q\sim r$ by Remark \ref{rmk4.1.1}. We have
\fm{&\mu=-2+O(r^{-\gamma_++C\eps}),\quad |\mu_q|\lesssim  r^{-1-\gamma_++C\eps} ,\\&|\nu|\lesssim  r^{-\gamma+C\eps}+\eps r^{-1+C\eps},\quad |\nu_q|\lesssim  r^{-\gamma_-+C\eps}+r^{-1-\gamma_++C\eps}+\eps r^{-2+C\eps}.}
Recall that $\gamma=\min\{\gamma_+,\gamma_-\}$.
It follows that
\fm{|\nu\mu_q|+|\nu_q\mu|+|\nu\nu_q|&\lesssim   r^{-\gamma_-+C\eps}+r^{-1-\gamma_++C\eps}+\eps r^{-2+C\eps};\\
(\lra{q}^{-1}+c_0\eps\ln t\cdot \lra{q}^{-1-\gamma_2})|\nu|(|\mu|+|\nu|)&\lesssim (r^{-1}+c_0 t^{C\eps}\cdot r^{-1-\gamma_2})\cdot(r^{-\gamma+C\eps}+\eps r^{-1+C\eps})\\&\lesssim r^{-1-\gamma_-+C\eps}+r^{-1-\gamma_++C\eps}+\eps r^{-2+C\eps}.}
By \eqref{sec5wei4}, we have
\fm{r^{-\gamma_-+C\eps}+r^{-1-\gamma_++C\eps}+\eps r^{-2+C\eps}&\lesssim \frac{ r^{1+\gamma_2-\gamma_-+C\eps}+r^{\gamma_2-\gamma_++C\eps}+\eps r^{\gamma_2-1+C\eps}}{c_0\eps\ln t}\cdot c_0\eps\ln t\cdot \lra{q}^{-1-\gamma_2}|\mu|^2\\
&\lesssim \frac{ r^{2+\gamma_2-\gamma_-+C\eps}+r^{1+\gamma_2-\gamma_++C\eps}+ r^{\gamma_2-1+C\eps}}{c_0\ln \eps^{-1}}\cdot c_0\eps\ln t\cdot \lra{q}^{-1-\gamma_2}|\mu|^2\\
&\lesssim (c_0\ln \eps^{-1})^{-1}\cdot c_0\eps\ln t\cdot \lra{q}^{-1-\gamma_2}|\mu|^2.}
To obtain the last estimate, we notice that $2+\gamma_2-\gamma_-,1+\gamma_2-\gamma_+,\gamma_2-1<0$.
We also notice that \fm{\delta\lra{q}^{-1-\gamma_+}\lesssim \delta r^{-\gamma_++\gamma_2}\lra{q}^{-1-\gamma_2}\leq \frac{\delta r^{1-\gamma_++\gamma_2}}{c_0 \ln t}\cdot c_0\eps\ln t \lra{q}^{-1-\gamma_2}\lesssim (c_0\ln\eps^{-1})^{-1}\cdot c_0\eps\ln t \lra{q}^{-1-\gamma_2}} whenever $r>2t$, $t\geq T_\eps$ and $\eps\ll 1$. Again, we use $1-\gamma_++\gamma_2<0$ here.
In summary, 
if $r>2t$, we have
\fm{\partial_r(q_rw)&\leq C\delta\lra{q}^{-1-\gamma_+}w|\mu|^2+ (-\frac{\gamma_2}{8}+C(c_0\ln\eps^{-1})^{-1})c_0\eps\ln t\cdot \lra{q}^{-1-\gamma_2}w|\mu|^2-\frac{\gamma_1}{4}(1-q)^{-1}1_{q<0}w|\mu|^2\\
&\leq(-\frac{\gamma_2}{8}+C(c_0\ln\eps^{-1})^{-1})c_0\eps\ln t\cdot \lra{q}^{-1-\gamma_2}w|\mu|^2-\frac{\gamma_1}{4}(1-q)^{-1}1_{q<0}w|\mu|^2.}
Finally we choose $\eps\ll1$ so that $C(c_0\ln\eps^{-1})^{-1}<\gamma_2/16$. 
\end{proof}\rm

We now return to the proof of Lemma \ref{lp2}.
\begin{proof}[Proof of Lemma \ref{lp2}] Recall that $\lra{q}=1+|q|$, so $\frac{d}{dq}(\lra{q})=\sgn(q)$. By Lemma \ref{lp2.l1}, we have
\fm{&\partial_r(r^2\phi^2q_r w\lra{q}^{-1})\\&=2r\phi^2q_rw\lra{q}^{-1}+2r^2\phi_r\phi q_rw\lra{q}^{-1}+r^2\phi^2\partial_r(q_rw)\lra{q}^{-1}-\sgn(q)r^2\phi^2q_r^2w\lra{q}^{-2}\\
&\leq 2r^2\phi(r^{-1}\phi+\phi_r)q_rw\lra{q}^{-1}-\sgn(q)r^2\phi^2q_r^2w\lra{q}^{-2}\\&\quad+[C\delta \lra{q}^{-\gamma_+}1_{t/2<r<2t,\ q\geq 0}-\frac{\gamma_2}{16}c_0\eps\ln t\cdot \lra{q}^{-\gamma_2}-\frac{\gamma_1}{8}1_{q<0}]w|\mu|^2r^2\phi^2\lra{q}^{-2}.}
Recall that $q_r=(\nu-\mu)/2=|\mu|/2+O(\eps(r+t)^{-1+C\eps}|\mu|)>0$, so \fm{|\mu|&\geq \frac{q_r}{1/2+C\eps (r+t)^{-1+C\eps}}\geq (2-C\eps (r+t)^{-1+C\eps})q_r,\\
|\mu|&\leq \frac{q_r}{1/2-C\eps (r+t)^{-1+C\eps}}\geq (2+C\eps (r+t)^{-1+C\eps})q_r.} Then,  we have
\fm{&\partial_r(r^2\phi^2q_r w\lra{q}^{-1})\\&\leq 2r^2\phi(r^{-1}\phi+\phi_r)q_rw\lra{q}^{-1} -\sgn(q)r^2\phi^2q_r^2w\lra{q}^{-2}\\
&\quad +(4+C\eps (r+t)^{-1+C\eps})C\delta\lra{q}^{-\gamma_+}1_{t/2<r<2t,\ q\geq 0}w|\mu|^2r^{2}\phi^2\lra{q}^{-2}\\&\quad-(4-C\eps (r+t)^{-1+C\eps})[\frac{\gamma_2c_0\eps\ln t}{16} \lra{q}^{-\gamma_2}+\frac{\gamma_1}{8}1_{q<0}]wq_r^2r^2\phi^2\lra{q}^{-2}\\
&\leq 2r^2\phi(r^{-1}\phi+\phi_r)q_rw\lra{q}^{-1}+r^2\phi^2q_r^2w\lra{q}^{-2}\\&\quad\cdot\kh{ -\sgn(q)+4C\delta\lra{q}^{-\gamma_+}1_{t/2<r<2t,\ q\geq 0}-\frac{\gamma_2c_0\eps\ln t}{4} \lra{q}^{-\gamma_2}-\frac{\gamma_1}{2}1_{q<0}+C\eps (r+t)^{-1+C\eps}}.}
If $q\geq 0$, then the coefficient of $wq_r^2r^2\phi^2\lra{q}^{-2}$ on the right hand side is not larger than
\fm{-1+4C\delta+C\eps(r+t)^{-1+C\eps}\leq -\frac{1}{2}.} 
To get the last estimate, we choose $\delta\ll1$  so that $4C\delta<1/8$.
If $q<0$, then the coefficient of $wq_r^2r^2\phi^2\lra{q}^{-2}$ on the right hand side is not larger than 
\fm{1-\gamma_1/2+C\eps \leq -\frac{\gamma_1-2}{4}<0.}
In the second last estimate, we use $\gamma_1>2$ and $\eps\ll1$. And since $\gamma_1<4$, we conclude that
\fm{&\partial_r(r^2\phi^2q_r w\lra{q}^{-1})\leq 2r^2\phi(r^{-1}\phi+\phi_r)q_rw\lra{q}^{-1}-\frac{\gamma_1-2}{4}r^2\phi^2q_r^2w\lra{q}^{-2}.}
Integrating this estimate with respect to $r\in[0,\infty)$ and $\omega\in\mathbb{S}^2$, we obtain
\fm{\int_{\R^3} \frac{\gamma_1-2}{4}wq_r^2\phi^2\lra{q}^{-2}\ dx&\leq \int_{\R^3}2\phi(r^{-1}\phi+\phi_r)q_rw\lra{q}^{-1}\ dx\\&\leq 2\kh{\int_{\R^3}\phi^2q_r^2w\lra{q}^{-2}\ dx}^{1/2}\kh{\int_{\R^3}(r^{-1}\phi+\phi_r)^2w\ dx}^{1/2}}
and thus
\eq{\label{lp2ff1}\int_{\R^3} wq_r^2\phi^2\lra{q}^{-2}\ dx\leq \frac{64}{(\gamma_1-2)^2}\int_{\R^3}(r^{-1}\phi+\phi_r)^2w\ dx.}

As in the proof of Lemma \ref{lp1}, we first consider the special case where $\phi\equiv 0$ unless $r>(1-c)t$ for some fixed small constant $0<c<1/2$.  We have
\fm{\partial_r(r\phi^2w)&=r^2(r^{-1}\phi+\phi_r)^2w-r^2\phi_r^2w+r\phi^2(\frac{m'(q)}{m(q)}+c_0\eps\ln t\cdot \sigma'(q))q_rw\\
&=r^2(r^{-1}\phi+\phi_r)^2w-r^2\phi_r^2w+r\phi^2(-\gamma_11_{q<0}-c_0\gamma_2\eps\ln t\cdot (1+|q|)^{-\gamma_2})\lra{q}^{-1}q_rw\\&=r^2(r^{-1}\phi+\phi_r)^2w-r^2\phi_r^2w+(rq_r)^{-1}(-\gamma_1\lra{q}1_{q<0}-c_0\gamma_2\eps\ln t\cdot \lra{q}^{1-\gamma_2})r^2\phi^2\lra{q}^{-2}q_r^2w.}
By \eqref{lem4.1c3}, we have $\lra{q}\lesssim \lra{r-t}(r+t)^{C\eps}\lesssim (r+t)^{1+C\eps}$. By Remark \ref{rmk4.2.1}, we have $q_r\geq C^{-1}t^{-C\eps}$. Since $0<\gamma_2<1$ and $r>(1-c)t$, we have $r+t<(1+\frac{1}{1-c})r<3r$ and
\fm{|\frac{c_0\gamma_2\eps\ln t\cdot\lra{q}^{1-\gamma_2}}{rq_r}|\lesssim c_0\gamma_2 t^{C\eps}(r+t)^{1-\gamma_2+C\eps} r^{-1}t^{C\eps}\lesssim c_0\gamma_2 r^{-\gamma_2+C\eps}.}
Since $t\geq T_\eps$, $r>(1-c)t$ and $\eps\ll1$, we have
\fm{\partial_r(r\phi^2w)&\geq r^2(r^{-1}\phi+\phi_r)^2w-r^2\phi_r^2w+[-(rq_r)^{-1} \gamma_1\lra{q}1_{q<0}-C\eps^{\gamma_2/2}]\cdot r^2\phi^2\lra{q}^{-2}q_r^2w.}
We now claim that for all $c\ll1$,  we have \fm{(rq_r)^{-1} \gamma_1\lra{q}1_{q<0}<\frac{(\gamma_1-2)^2}{128},\qquad\text{whenever } r\geq (1-c)t.} We  prove this claim by contradiction. If the claim is false, then there exists a point $(t,x)$ where $r\geq (1-c)t$, $q<0$ and \fm{\lra{q}\geq \frac{(\gamma_1-2)^2}{128\gamma_1}rq_r\geq C^{-1}rt^{-C\eps}\geq C^{-1}t^{1-C\eps}.}
The constant $C$ here could depend on $\gamma_1$. It then follows from Lemma \ref{lem4.2} and \eqref{def3.1a6} that
\fm{|q_r-1|=|\frac{\nu-\mu-2}{2}|\lesssim \eps(r+t)^{-1+C\eps}+t^{C\eps}\lra{q}^{-\gamma_-}\lesssim \eps+ t^{-1+C\eps},}
so we have $q_r\geq 1/2$ for $\eps\ll1$. Because of this lower bound for $q_r$, the lower bound for $\lra{q}$ above can be improved: we have $\lra{q}\geq C^{-1}r$. Since $q<0$ and $\gamma_->2$, by \eqref{lem4.1c1} we have
\fm{|q-r+t|\lesssim (r+t)^{C\eps}\lra{q}^{1-\gamma_-}\lesssim (r+\frac{r}{1-c})^{C\eps} r^{1-\gamma_-}\lesssim  r^{1-\gamma_-+C\eps}. }
As $r\geq t/2\geq T_\eps/2$, we have \fm{|r-t|\geq |q|-|q-r+t|\geq C^{-1}r-C r^{1-\gamma_-+C\eps}-1\geq (2C)^{-1}r.}
On the other hand, we have $r>(1-c)t$, so either $0\leq r-t\leq q+|q-r+t|<0+Cr^{1-\gamma_-+C\eps}$ or $0\leq t-r\leq cr/(1-c)\leq 2cr$. That is, we have
\fm{|r-t|\leq Cr^{1-\gamma_-+C\eps}+2cr=(2c+Cr^{-\gamma_-/2})r. }
By choosing $c\ll1$ and $\eps\ll1$ (note that $c$ is chosen before $\eps$ is chosen), we have $(2c+Cr^{-\gamma_-/2})r<(2C)^{-1}r$ which is a contradiction. It follows from the claim that
\fm{\int_{\R^3}(r^{-1}\phi+\phi_r)^2w\ dx\leq \int_{\R^3}\phi_r^2w+(\frac{(\gamma_1-2)^2}{128}+C\eps^{\gamma_2/2})\phi^2\lra{q}^{-2}q_r^2w\ dx.}
Combining this estimate with \eqref{lp2ff1}, we conclude that
\fm{\int_{\R^3}\phi^2\lra{q}^{-2}q_r^2w\ dx\leq \int_{\R^3}\frac{64}{(\gamma_1-2)^2}\phi_r^2w+(\frac{1}{2}+C\eps^{\gamma_2/2})\phi^2\lra{q}^{-2}q_r^2w\ dx.}
By choosing $\eps\ll1$ so that $1/2+C\eps^{\gamma_2/2}<3/4$, we conclude that \fm{
&\int_{\R^3}q_r^2\lra{q}^{-2}|\phi|^2w\ dx\lesssim_\eta \int_{\R^3}|\partial\phi|^2w\ dx}
in this special case.

In general, we choose the same cutoff function $\chi$ as in the proof of Lemma \ref{lp1}. Then,
\fm{&\int_{\R^3}q_r^2\lra{q}^{-2}|\chi(r/t)\phi|^2w\ dx\\
&\lesssim_{c} \int_{\R^3}|\partial(\chi(r/t)\phi)|^2w\ dx\lesssim \int_{\R^3}(|\partial\phi|^2+r^{-2}|\chi'(r/t)\phi|^2)w\ dx\\
&\lesssim \int_{\R^3}(|\partial\phi|^2+r^{-2}|\phi|^21_{r/t\in[1-c,1-c/2]})w\ dx.}
Again, we finish the proof by noticing that $\chi(r/t)\geq 1_{r/t\geq 1-c/2}$.
\end{proof}
\rmk{\label{rmk5.6.1}\rm  Here we have actually proved Lemma \ref{lp2} for infinitely many pairs of $(\gamma_1,\gamma_2)$ (with $\eps,c\ll_{\gamma_1,\gamma_2}1$ and $c_0\gg_{\gamma_1,\gamma_2}1$), as long as they satisfy \eqref{sec5wei4}. Thus, if $(\gamma_1,\gamma_2)$ satisfies \eqref{sec5wei4}, then Lemma \ref{lp2} with $(\gamma_1,\gamma_2)$ replaced by $(\frac{\gamma_1}{2}+1,\gamma_2)$ also holds. Note that we also change the weight $w$ accordingly. Since $2<\frac{\gamma_1}{2}+1<\gamma_1$, we have
\eq{\label{rmk5.6.1c}&\int_{\R^3}q_r^2(\lra{q}^{-2}1_{q\geq 0}+\lra{q}^{-\frac{\gamma_1+2}{2}}1_{q<0})\phi^2 1_{r/t\geq 1-c/2} w_{\gamma_1,\gamma_2} \ dx\\
&= \int_{\R^3}\kh{\lra{q}^{-2}   1_{q\geq 0}+\lra{q}^{-2+\frac{\gamma_1+2}{2}}  1_{q<0}} q_r^2\phi^2\exp(c_0\eps\ln t\cdot \sigma(q))1_{r/t\geq 1-c/2}\ dx\\
&= \int_{\R^3}\lra{q}^{-2}  q_r^2\phi^2w_{\frac{\gamma_1+2}{2},\gamma_2}1_{r/t\geq 1-c/2}\ dx\\
&\lesssim_c\int_{\R^3}(|\partial\phi|^2+r^{-2}|\phi|^21_{r/t\in[1-c,1-c/2]})w_{\frac{\gamma_1+2}{2},\gamma_2}\ dx\\
&\lesssim \int_{\R^3}|\partial\phi|^2w_{\gamma_1,\gamma_2}+r^{-\frac{\gamma_1+2}{2}}|\phi|^21_{r/t\in[1-c,1-c/2]}w_{\gamma_1,\gamma_2}\ dx.}
In the last estimate, we use that $w_{\frac{\gamma_1+2}{2},\gamma_2}=\lra{q}^{-\frac{\gamma_1-2}{2}}w_{\gamma_1,\gamma_2}$ for $q<0$ and  that $(-q)\sim r$ whenever $r/t\in[1-c,1-c/2]$.
}
\rm

\bigskip

We can now set the value of $c$ in Section \ref{sec4.3}. We choose $c\in(0,1/4)$ so that Lemma \ref{lp1} (with $\eta=\frac{\gamma_1}{2}$) and Lemma \ref{lp2} (with $\gamma_1$ replaced by $\frac{\gamma_1+2}{2}$) both hold with $c$ there replaced by $2c$.  Once we fix the value $c$, it can be treated as a universal constant.

We end this section with the following key lemma.  Note that assumption \eqref{def3.1a52} is used in this lemma. Again, the constant $B$ is not involved in this lemma, as we do not use the function $u$ at all.

\lem{\label{lp3}Suppose $\phi$ is smooth and $\phi(t)\in C_c^\infty(\R^3)$ for each $t\geq T_\eps$. Let $F=(F_{(J)}^{(K)}):=g_J^{\alpha\beta}\partial_\alpha\partial_\beta u^{(K)}_{app}$ where $u_{app}$ is defined in \eqref{uappdef}. Then for $t\geq T_\eps$ with $\eps\ll1$, we have \fm{\norm{\phi F}_{L^2(w)}\lesssim \eps t^{-1}\norm{\partial\phi}_{L^2(w)}+\eps t^{-1} \norm{r^{-\frac{\gamma_1+2}{4}}\phi 1_{r/t\in[1-2c,1-c]}}_{L^2(w)}.}}
\begin{proof}
Fix $1\leq J,K\leq M$ and write \fm{F_{(J)}^{(K)}=\eps (4r)^{-1}\psi(r/t)G_{2,J}(\omega)\mu\partial_q(\mu U_q^{(K)})+R_{(J)}^{(K)}.}
Here $\psi$ comes from the definition \eqref{uappdef} of $u_{app}$ in Section \ref{sec4.3}. By \eqref{lem4.11c2} in Lemma \ref{lem4.11}, we have 
\fm{|R_{(J)}^{(K)}|\lesssim \eps t^{-2+C\eps}\lra{q}^{-\gamma_{\sgn(q)}}+\eps t^{-3+C\eps}.}
Since $\psi(r/t)=0$ unless $r/t\in[1-c,1+c]$,  we have $R_{(J)}^{(K)}=0$ unless $r/t\in[1-c,1+c]$. Moreover, by Remark \ref{rmk4.1.1}, if we set $w_0(\rho):=1_{\rho\geq 0}+\lra{\rho}^{\gamma_1}1_{\rho<0}$, then whenever $r>t/2$ we have \fm{ w_0(r-t)\lesssim t^{C\eps}w_0(q)\lesssim t^{C\eps}w;\qquad
w\lesssim t^{C\eps}w_0(q)\lesssim t^{C\eps}w_0(r-t).}
Thus, by Lemma \ref{lp1} and Lemma \ref{lp2}, we have
\fm{\norm{\phi R^{(J)}_{(K)}}_{L^2(w)}&\lesssim \eps t^{-2+C\eps}\norm{\lra{q}^{-1}q_r\phi 1_{r/t\geq 1-c}}_{L^2(w)}+\eps t^{-3+C\eps}\norm{\phi 1_{r/t\geq 1-c}}_{L^2(w)}\\
&\lesssim \eps t^{-2+C\eps}(\norm{\partial\phi}_{L^2(w)}+\norm{r^{-1}\phi1_{r/t\in[1-2c,1-c]}}_{L^2(w)})\\
&\quad+\eps t^{-2+C\eps}\norm{w_0(r-t)^{1/2}\lra{r-t}^{-1}\phi}_{L^2(\R^3)}\\
&\lesssim \eps t^{-2+C\eps}(\norm{\partial\phi}_{L^2(w)}+\norm{r^{-1}\phi1_{r/t\in[1-2c,1-c]}}_{L^2(w)})\\
&\quad+\eps t^{-2+C\eps}(\norm{w_0(r-t)^{1/2}\partial\phi}_{L^2(\R^3)}+\norm{w_0(r-t)^{1/2}r^{-1}\phi1_{r/t\in[1-2c,1-c]}}_{L^2(\R^3)})\\
&\lesssim \eps t^{-2+C\eps}(\norm{\partial\phi}_{L^2(w)}+\norm{r^{-1}\phi1_{r/t\in[1-2c,1-c]}}_{L^2(w)})\\
&\lesssim \eps t^{-1}(\norm{\partial\phi}_{L^2(w)}+\norm{r^{-2+C\eps}\phi1_{r/t\in[1-2c,1-c]}}_{L^2(w)}).
}
Since $\gamma_1<4$, we have $(\gamma_1+2)/4<3/2$. By choosing $\eps\ll1$, we have $r^{-2+C\eps}\lesssim r^{-\frac{\gamma_1+2}{4}}$ whenever $r/t\in[1-2c,1-c]$ and $t\geq T_\eps$.

Recall from \eqref{def3.1a52} that $|G_{2,J}\partial_q(\mu U_q)|\lesssim \lra{q}^{-\gamma_{\sgn(q)}}$. Also recall that $|\mu|\sim q_r>0$. By \eqref{rmk5.6.1c} in Remark \ref{rmk5.6.1}, we obtain
\fm{\norm{\eps (4r)^{-1}\psi(r/t)G_{2,J}(\omega)\mu\partial_q(\mu U_q^{(K)})\phi}_{L^2(w)}&\lesssim \eps t^{-1}\norm{q_r\lra{q}^{-\gamma_{\sgn(q)}}\phi 1_{r/t\geq 1-c}}_{L^2(w)}\\
&\lesssim \eps t^{-1} \norm{r^{-\frac{\gamma_1+2}{4}}\phi 1_{r/t\in[1-2c,1-c]}}_{L^2(w)}+\eps t^{-1}\norm{\partial\phi}_{L^2(w)}.}
To obtain the last estimate, we notice that $\gamma_+>1$ and $\gamma_->(\gamma_1+2)/2$.
\end{proof}
\rm
\section{Continuity argument}\label{sec6}

\subsection{Setup}\label{sec6.1}
Fix $\chi(s)\in C_c^\infty(\R)$ such that $\chi\in[0,1]$ for all $s$, $\chi\equiv 1$ for $|s|\leq 1$ and $\chi\equiv 0$ for $|s|\geq 2$. Also fix a large time $T>T_\eps=1/\eps$, where $\eps\ll1$ are chosen so that all the results hold in the previous two sections. Consider the following equations for $v=(v^{(I)})_{I=1}^M$:
\eq{\label{eqn}
&g^{\alpha\beta}(u_{app}+v,\partial (u_{app}+v))\partial_\alpha\partial_\beta v^{(I)}\\&=-[ g^{\alpha\beta}(u_{app}+v,\partial (u_{app}+v))-g^{\alpha\beta}(u_{app},\partial u_{app})]\partial_\alpha\partial_\beta u_{app}^{(I)}\\
&\quad+[f^{(I)}(u_{app}+v,\partial (u_{app}+v))-f^{(I)}(u_{app},\partial u_{app})]\\
&\quad-\chi(t/T)[g^{\alpha\beta}(u_{app},\partial (u_{app}))\partial_\alpha\partial_\beta u_{app}^{(I)}-f^{(I)}(u_{app},\partial u_{app})]}
for all $t>T_\eps$, along with  data $v\equiv 0$ for $t\geq 2T$. We will add a superscript $T$ (i.e.\ $v^T$) if we hope to emphasize the dependence of $v$ on $T$.

We have the following results.
\begin{enumerate}[(a)]
\item By the local existence theory of quasilinear wave equations, we can find a local smooth solution to \eqref{eqn} near $t=2T$.
\item The solution on $[T_1,\infty)$ can be extended to $[T_1-\epsilon,\infty)$ for some small $\epsilon>0$ if  \[\norm{\partial^kv}_{L^\infty([T_1,\infty)\times\R^3)}<\infty,\qquad\text{ for all }k\leq 4.\]
\item The solution to \eqref{eqn} has a finite speed of propagation: $v^T(t,x)=0$ if $r+t>(4+2c)T$ or if $r-t<-2cT$, so $Z^I(t/T)=O(1)$ whenever $T/2\leq t\leq 2T$ and $v(t,x)\neq 0$.
\item If the solution exists for $t\leq T$,  we have  \eqref{qwe} holds for $t\leq T$ where  $u=u_{app}+v$.
\end{enumerate}

The proofs of these statements are standard. We refer to \cite{MR2455195} for the proofs of (a) and (b). For (c), we notice  that $Z^Ix_\alpha$ is a linear combination of $1$, $t$ and $x_i$ with constant coefficients, so $|Z^I(t/T)|\lesssim (1+r+t)/T\lesssim 1$ whenever $t\in[T,2T]$ and $v(t,x)\neq 0$ for each $I$. In this section, our goal is to prove the following proposition.

\prop{\label{prop6} Fix an integer $N\geq 6$,  a $(\gamma_+,\gamma_-)$-admissible global solution $(\mu,U)$ to the geometric reduced system \eqref{asy}, and a sufficiently small constant $\delta\in(0,\delta_0)$ where $\delta_0$ comes from Definition \ref{def3.1}. Choose $0<c\ll1$ as in Section \ref{sec5} and let $u_{app}$ be the corresponding approximate solution defined in Section \ref{sec4}. Then there exists  $\eps_{N}>0$,  depending on $(\mu,U),N,\delta,\gamma_\pm$, such that for any $0<\eps<\eps_{N}$, we can find a $C^{N+1}$ function\footnote{Since $u_{app},g^{**},f^{(*)}\in C^\infty$, one can show that both $u$ and $v$ are $C^\infty$. Here we apply \cite[Theorem I.4.3]{MR2455195}. However, the estimates \eqref{prop6c1} and \eqref{prop6c2} do not necessarily hold for $|I|>N$.} $u=(u^{(I)}(t,x))$ defined for all $t\geq 0$, such that the following statements are true:
\begin{enumerate}[\rm (a)]
\item If we set $v=u-u_{app}$, then $v$ is a solution to \eqref{eqn} for all $t\geq T_\eps$.
\item For $t\leq T$, $u$ is a solution to \eqref{qwe}.
\item For all $|I|\leq N$, we have
\begin{equation}\label{prop6c1}\norm{w_0(|x|-t)^{1/2}\partial Z^Iv(t)}_{L^2(\R^3)}\lesssim_I\eps t^{-1/2+C_{I}\eps},\qquad \forall t\geq T_\eps;\end{equation}
\begin{equation}\label{prop6c2}\norm{w_0(|x|-t)^{1/2}\partial Z^Iu(t)}_{L^2(\R^3)}\lesssim_I\eps , \qquad \forall 0\leq t\leq T_\eps.\end{equation}
\end{enumerate}
Here $w_0(\rho):=1_{\rho\geq 0}+\lra{\rho}^{\gamma_1}1_{\rho<0}$ for some $2<\gamma_1<\min\{2(\gamma_--1),4\}$.

Note that $\eps_N$ and all the constants in \eqref{prop6c1} and \eqref{prop6c2} are independent of the time $T$.
}
\rmk{\rm It should be pointed out that the $N$ in this proposition is different from the $N$ in the main theorem. 
}

\rm

\bigskip

We  use a continuity argument to prove this proposition. From now on we assume $\eps\ll1$, which means $\eps$ is arbitrary  in $(0,\eps_N)$ for some fixed small constant $\eps_N$ depending on $N$. First we  prove that a solution $u$ exists for  $t\geq T_\eps$. We start with a solution $v$ to \eqref{eqn} for $t\geq T_1$ such that for all $t\geq T_1> T_\eps$ and $k+i\leq N$,  \eq{\label{ca1}
E_{k,i}(t)
&:=\sum_{l\leq k,|I|\leq i}\left[E_u(\partial^lZ^Iv)(t)+\int_t^{\infty}E_q(\partial^lZ^Iv)(\tau)\ d\tau\right]\leq B_{k,i}
\eps^2 t^{-1+C_{k,i}\eps}}
and
\begin{equation}\label{ca2}
\sum_{|I|\leq 2}|Z^Iu|\leq B_0\eps t^{-1+2C_{0,N}\eps/2},\qquad |u|\leq B_1\eps t^{-1}.
\end{equation} 
Here $u=v+u_{app}$, $E_u(\cdot)$ is defined by \eqref{sec5ene}, and $E_q(\cdot)$ is defined by \eqref{sec5eneq}.  We remark that $C_{k,i},B_{k,i}$ depend on $k,i$ but not on $N$. 

Our goal is to prove that \eqref{ca1} and \eqref{ca2} hold with $B_{k,i},B_0,B_1$ replaced by  smaller constants $B_{k,i}^\prime,B_0^\prime,B_1^\prime$, and with $C_{k,i}$ unchanged, assuming that $\eps\ll 1$. To achieve this goal, we first induct on $i$, and then we induct on $k$ for each fixed $i$. For each $(k,i)$ with $k+i\leq N$, we want to prove the following inequality
\begin{equation}\label{sec6fff}\begin{aligned}
&\sum_{|I|\leq i}\norm{ g^{\alpha\beta}(u,\partial u)\partial_\alpha\partial_\beta \partial^kZ^Iv}_{L^2(w)}\\&\leq C\eps t^{-1}E_{k,i}(t)^{1/2}+C\eps^{1/2}(\ln\eps^{-1})^{-1/2} t^{-1+C\eps}\sum_{l\leq k\atop |I|\leq i}E_q(\partial^lZ^Iv)(t)^{1/2}\\&\quad+C\eps t^{-1+C\eps}(E_{k-1,i}(t)^{1/2}+E_{k+1,i-1}(t)^{1/2})+C\eps t^{-3/2+C\eps}.
\end{aligned}\end{equation}
Here we set $E_{-1,\cdot}=E_{\cdot,-1}=0$, and $C$ are constants independent of all the constants in \eqref{ca1} and \eqref{ca2}. This inequality will be proved in Section \ref{sec6.3}. We then combine \eqref{sec6fff} with the energy estimates \eqref{energymain} to derive an inequality on $E_{k,i}(t)$. 

After closing the continuity argument above, we obtain a solution $v$ to \eqref{eqn} for all $t\geq T_\eps$. As remarked above, we also obtain a solution $u=v+u_{app}$ to the original equations \eqref{qwe} for all $T_\eps\leq t\leq T$. In the second part of the proof, we shall prove that this $u$ can be extended to a solution to \eqref{qwe} for all $0\leq t\leq T$. The proof here is still based on a continuity argument, and we will go into its details in Section \ref{sec6.5}.

In the following computations,  let $C$ denote a universal constant  or a constant from the previous estimates for $q$ and $u_{app}$ (e.g.\ from Proposition \ref{mainprop4}). We also write  $A\lesssim B$ if the implicit constant factor is such a constant $C$. Here $C$ is allowed to depend on $(k,i)$ or $N$, but we will never write it as $C_{k,i}$ or $C_N$ because the value of $C$ is known before Section \ref{sec6}. We will choose the constants in the following order:
\eq{\label{sec6ord}C&\to C_{0,0},B_{0,0}\to C_{1,0},B_{1,0}\to\dots\to C_{N,0},B_{N,0}\\&\to C_{0,1},B_{0,1}\to\dots\to C_{N-1,1},B_{N-1,1}\\&\to C_{0,2},B_{0,2}\to\dots\to C_{N-2,2},B_{N-2,2}\\&\dots\\&\to C_{0,N},B_{0,N}\\&\to B_0,B_1\to C_N\to\eps_N.}
We emphasize that if a constant $K$ appears before a constant $K'$, then $K$ cannot depend on~$K'$. We also remark that $C_N$ can be chosen to be a polynomial of  all the constants appearing before it, and that the constant $B$ in Section \ref{sec5} can be replaced by this $C_N$. For simplicity, we allow  $C_N$ to vary from line to line, so it makes sense to write $C_N^2$, $C_N+C$, $CC_N$ as $C_N$.

In addition, since    $\eps \ll 1$ is chosen at the end, we can control terms like $C_N\eps$ and $C_N T_\eps^{-\kappa+C_N\eps}$ for $\kappa>0$ for any $(k,i)$ by a universal constant, e.g.\ $1$.

To end the setup, we derive a differential equation for $\partial^kZ^Iv$ from \eqref{eqn}. If we commute \eqref{eqn} with $\partial^kZ^I$, we have
\begin{equation}\label{eqn2}
\begin{aligned}
&g^{\alpha\beta}(u_{app}+v,\partial (u_{app}+v))\partial_\alpha\partial_\beta \partial^kZ^Iv^{(J)}\\
&=\partial^k[\Box,Z^I]v^{(J)}\\&\quad+[\left(g^{\alpha\beta}(u_{app}+v,\partial (u_{app}+v))-m^{\alpha\beta}\right),\partial^kZ^I]\partial_\alpha\partial_\beta v^{(J)}\\&\quad+\left(g^{\alpha\beta}(u_{app}+v,\partial (u_{app}+v))-m^{\alpha\beta}\right)\partial^k[\partial_\alpha\partial_\beta,Z^I]v^{(J)}\\
&\quad -\partial^kZ^I\left([ g^{\alpha\beta}(u_{app}+v,\partial (u_{app}+v))-g^{\alpha\beta}(u_{app},\partial u_{app})]\partial_\alpha\partial_\beta u_{app}^{(J)}\right)\\&\quad+\partial^kZ^I\left([f^{(J)}(u_{app}+v,\partial (u_{app}+v))-f^{(J)}(u_{app},\partial u_{app})]\right)\\
&\quad-\partial^kZ^I\left(\chi(t/T)[g^{\alpha\beta}(u_{app},\partial u_{app})\partial_\alpha\partial_\beta u_{app}^{(J)}-f^{(J)}(u_{app},\partial u_{app})]\right)\\
&=:\sum_{j=1}^6R_j^{(J),k,I}.
\end{aligned}
\end{equation} with $\partial^kZ^Iv\equiv 0$ for $t\geq 2T$.

\subsection{The pointwise bound \eqref{ca2}}\label{sec6.2}
In the next few subsections, we always assume $t\geq T_1> T_\eps$ and $\eps\ll1$. Define a function $w_0\in C^0(\R)$ by $w_0(\rho):=1_{\rho\geq 0}+\lra{\rho}^{\gamma_1}1_{\rho<0}$. By \eqref{rmk4.1.1c2} in Remark \ref{rmk4.1.1} and since $q=r-t$ whenever $|x|<t/2$, we have \fm{w_0(r-t)\lesssim t^{C\eps}w_0(q)\lesssim t^{C\eps}w,\qquad
w\lesssim t^{C\eps}w_0(q)\lesssim t^{C\eps}w_0(r-t),\qquad\forall t\geq \exp(\delta/\eps).} Here the implicit constants depend only on $c_0$ in \eqref{sec5wei}. By \eqref{ca2} and \eqref{sec5ene} we have
\eq{\label{equivnorm}
&\norm{\partial\phi}_{L^2(w_0)}\lesssim t^{C\eps}\norm{\partial\phi}_{L^2(w)}\sim t^{C\eps}E_u(\phi)^{1/2};\\
&\norm{\partial\phi}_{L^2(w)}\sim E_u(\phi)^{1/2}\lesssim t^{C\eps}\norm{\partial\phi}_{L^2(w_0)}.} 
Here all the norms are taken for fixed $t$, and $\norm{f}_{L^2(w_0)}^2:=\int |f|^2w_0(r-t)\ dx$. 
We can choose  $\eps\ll1$ so that all the implicit constants in this inequality are universal and do not depend on $N$. Then,
\fm{\norm{w_0(|\cdot|-t)^{1/2}\partial Z^Iv(t)}_{L^2(\R^3)}^2&\lesssim t^{C\eps}E_u(Z^Iv)(t)\lesssim B_{0,i} \eps^2t^{-1+(C_{0,i}+C)\eps}}
for all $|I|=i\leq N$.
By the weighted Klainerman-Sobolev inequality  \eqref{corklc}, we have 
\eq{\label{ptb1}|\partial Z^Iv(t)|&\lesssim B_{0,i+2}^{1/2}\eps t^{-1/2+(C_{0,i+2}+C)\eps/2}\lra{r+t}^{-1}\lra{t-r}^{-1/2}w_0(r-t)^{-1/2}\\
&\lesssim B_{0,i+2}^{1/2}\eps t^{-1/2+(C_{0,i+2}+C)\eps/2}\lra{r+t}^{-1}(\lra{t-r}^{-1/2}1_{r-t\geq 0}+\lra{t-r}^{-(1+\gamma_1)/2}1_{r-t<0})\\
&\lesssim B_{0,i+2}^{1/2}\eps t^{-3/2+(C_{0,i+2}+C)\eps/2}}
whenever  $|I|=i\leq N-2$ and $t\in[T_1,2T]$.  We also claim that  
\eq{\label{ptb2} &|Z^Iv(t,x)|\\
&\lesssim B_{0,i+2}^{1/2}\eps t^{-1/2+(C_{0,i+2}+C)\eps/2} (t^{-1}\lra{r-t}^{-(\gamma_1-1)/2}1_{r< t}+t^{-1}\lra{r-t}^{1/2}1_{t\leq r\leq 2t}+r^{-1/2}1_{r>2t})\\
&\lesssim B_{0,i+2}^{1/2}\eps t^{-1+(C_{0,i+2}+C)\eps/2}}
whenever $|I|=i\leq N-2$ and $t\in[T_1,2T]$. To see this, we first recall from Section \ref{sec6.1} that  $v\equiv 0$ whenever $t\geq 2T$, or whenever $r+t\geq (4+2c)T$ and $T_1\leq t\leq 2T$. Then \eqref{ptb2} follows from the following lemma.
\lem{\label{lemptb2} Fix $\gamma_1>2$. Suppose that $\phi=\phi(t,x)$ is a $C^1$ function  in $[T_1,T_2]\times\R^3$ with $0\leq T_1\leq T_2<\infty$. Suppose that $\phi(T_2,x)=0$  for all $x\in\R^3$, that $\lim_{|x|\to\infty}\phi(t,x)=0$ for each $t\in[T_1,T_2]$,  and that \eq{\label{lemptb2a}|\partial\phi(t,x)|\leq A_0(t)\lra{r+t}^{-1}\lra{r-t}^{-1/2}(1_{r\geq t}+\lra{r-t}^{-\gamma_1/2}1_{r<t}),\qquad\forall(t,x)\in[T_1,T_2]\times\R^3.}
Here $A_0(t)\geq 0$ is a nonincreasing function of $t$. Then, we have
\eq{\label{lemptb2c}|\phi(t,x)|\leq CA_0(t)(\lra{t}^{-1}\lra{r-t}^{-(\gamma_1-1)/2}1_{r< t}+\lra{t}^{-1}\lra{r-t}^{1/2}1_{t\leq r\leq  2t}+\lra{r}^{-1/2}1_{r>2t})}
everywhere in $[T_1,T_2]\times\R^3$. Here $C>1$ is a  constant  independent of $T_1,T_2$ and $A_0$.
}
\begin{proof}
We first prove this lemma for $t\in[T_1,T_2]$ and $r=|x|<t$. By the fundamental theorem of calculus, we have
\fm{|\phi(t,x)|&\leq |\phi(T_2,x)|+\int_t^{T_2}|\partial_\tau\phi(\tau,x)|\ d\tau\leq 0+\int_t^{T_2}A_0(\tau) \lra{\tau+r}^{-1}\lra{r-\tau}^{-(1+\gamma_1)/2}\ d\tau\\&\leq A_0(t)\int_t^\infty \lra{\tau+r}^{-1}\lra{r-\tau}^{-(1+\gamma_1)/2}\ d\tau.}
Note that $\tau\geq t> r$ in the integrals. Besides, we have
\fm{&\int_t^\infty \lra{\tau+r}^{-1}\lra{r-\tau}^{-(1+\gamma_1)/2}\ d\tau\\
&\leq \int_t^{2(t+r)} \lra{t+r}^{-1}\lra{r-\tau}^{-(1+\gamma_1)/2}\ d\tau+\int_{2(t+r)}^\infty \lra{\tau}^{-1}\lra{\tau/2}^{-(1+\gamma_1)/2}\ d\tau\\
&\leq \lra{t+r}^{-1}\int_{-\infty}^{r-t} \lra{\rho}^{-(1+\gamma_1)/2}\ d\rho+\int_{2(t+r)}^\infty \lra{\tau/2}^{-(3+\gamma_1)/2}\ d\tau\\&\leq C \lra{t}^{-1}\lra{r-t}^{(1-\gamma_1)/2}+C\lra{r+t}^{-(1+\gamma_1)/2}\leq C \lra{t}^{-1}\lra{r-t}^{(1-\gamma_1)/2}. }
This gives us \eqref{lemptb2c}.
Here the constant $C>1$ is determined only by $\gamma_1>2$.

Next we consider the case when $t\in[T_1,T_2]$ and $r/t\in[1,2]$. It is known from the proof above that $|\phi(t,t\omega)|\leq CA_0(t)\lra{t}^{-1}$.  Here $\omega=x/r\in\mathbb{S}^2$. It thus follows that
\fm{|\phi(t,x)|&\leq |\phi(t,t\omega)|+\int_t^r|\partial_r\phi(t,\rho\omega)|\ d\rho\\
&\leq CA_0(t)\lra{t}^{-1}+\int_t^rA_0(t)\lra{\rho+t}^{-1}\lra{\rho-t}^{-1/2}\ d\rho\\&\leq CA_0(t)\lra{t}^{-1}+A_0(t)\lra{t}^{-1}\int_0^{r-t}\lra{\rho'}^{-1/2}\ d\rho'\leq CA_0(t)\lra{t}^{-1}\lra{r-t}^{1/2}.}
Again, the constant $C>1$ is determined only by $\gamma_1>2$.

Finally, we assume $r>2t$. Since $\lim_{|x|\to\infty}\phi(t,x)=0$, we have
\fm{|\phi(t,x)|&\leq \int_r^\infty |\partial_r\phi(t,\rho\omega)|\ d\rho\leq \int_r^\infty A_0(t)\lra{\rho+t}^{-1}\lra{\rho-t}^{-1/2}\ d\rho\leq \int_r^\infty A_0(t)\lra{\rho/2}^{-3/2}\ d\rho\\
&\leq CA_0(t)\lra{r}^{-1/2}.}
Here we use $\rho-t\geq \rho/2$ whenever $\rho\geq r>2t$. Again, the constant $C>1$ is determined only by $\gamma_1>2$.
\end{proof}\rm

A direct corollary of \eqref{ptb2} is that
\eq{\label{ptb2cor}\norm{v(t)\cdot 1_{r/t\in[1-2c,1-c/2]}}_{L^2(w)}&\lesssim t^{C\eps}\norm{\lra{r-t}^{\gamma_1/2}v(t)\cdot 1_{r/t\in[1-2c,1-c/2]}}_{L^2(\R^3)}\\
&\lesssim t^{\gamma_1/2+C\eps}\cdot t^{3/2}\cdot C_N\eps t^{-3/2-(\gamma_1-1)/2+C_N\eps}\lesssim C_N\eps t^{1/2+C_N\eps}.}
Here we notice that $0<-q\sim r$ and $|r-t|\sim r$ whenever $r/t\in[1-2c,1-c/2]$. For completeness, we present the following estimate which is closely related  to \eqref{ptb2cor}. Since $|Z\phi|\lesssim \lra{r+t}|\partial\phi|$, whenever $0<|I|\leq N+1$ we have
\eq{\label{ptb2cor2} \norm{Z^Iv(t)\cdot 1_{r/t\in[1-2c,1-c/2]}}_{L^2(w)}&\lesssim t\sum_{|J|=|I|-1}\norm{\partial Z^Jv(t)\cdot 1_{r/t\in[1-2c,1-c/2]}}_{L^2(w)}\lesssim tE_{0,|I|-1}(t)^{1/2}.}

Since $u=v+u_{app}$, by Proposition \ref{mainprop4}, we conclude that
\fm{\sum_{|I|\leq 2}|Z^Iu(t,x)|&\leq C\eps t^{-1+C\eps}+CB_{0,4}^{1/2}\eps t^{-1+(C_{0,4}+C)\eps/2},\\
|\partial u(t,x)|&\leq C\eps t^{-1}+CB_{0,3}^{1/2}\eps t^{-3/2+(C_{0,3}+C)\eps/2}.}
This allows us to replace $B_0,B_1$ with $B_0/2,B_1/2$ in \eqref{ca2} as long as we choose  $B_0,B_1\gg1$ and $\eps\ll 1$. Again we emphasize that $B_0$ and $B_1$ are chosen after all the $C_{k,i}$ and $B_{k,i}$ are chosen.

We  have three remarks before we end this subsection. First, if $N\geq 6$, \eqref{ptb1} and \eqref{ptb2} allow us to extend the solution to \eqref{eqn} below $t=T_1$, by the local existence theory of quasilinear wave equations. Moreover, these two pointwise bounds, together with Proposition \ref{mainprop4}, allow us to use Lemma \ref{l24} freely to $g^{\alpha\beta}(\cdot,\partial(\cdot))$ and $f^{(I)}(\cdot,\partial(\cdot))$. Here we need to replace the $u$ and $v$ in Lemma \ref{l24} by $(u_{app},\partial  u_{app} )$ and $(v,\partial v)$, respectively.  Finally, the assumptions \eqref{sec5a1} and \eqref{sec5a2} are satisfied because of \eqref{ca2}, \eqref{ptb1}, and \eqref{ptb2}, so  we can apply the energy estimates in Section \ref{sec5} freely for $t\geq T_1$.

\subsection{The $L^2$ estimates for $R_j^{(J),k,I}$}\label{sec6.3} We set
\eq{\mcl{R}_{j}^{k,i}:=\sum_{|I|=i}\sum_{J=1}^N\norm{R_j^{(J),k,I}}_{L^2(w)}.}
Here the $R^{(J),k,I}_j$'s come from \eqref{eqn2}. Our goal in this subsection is to  derive several estimates for $R_j^{(J),k,I}$ when we assume that \eqref{ca1} and \eqref{ca2} (and therefore \eqref{ptb1} and \eqref{ptb2}) hold. We seek to prove the following proposition.

\prop{\label{propsec6.3} For each $k,i\geq0$ with $k+i\leq N$, we have
\eq{\label{propsec6.3c}\sum_{j=1}^6\mcl{R}_{j}^{k,i}&\leq C\eps t^{-1}E_{k,i}(t)^{1/2}+C\eps^{1/2}(\ln\eps^{-1})^{-1/2} t^{-1+C\eps}\sum_{l\leq k\atop |I|\leq i}E_q(\partial^lZ^Iv)(t)^{1/2}\\&\quad+C\eps t^{-1+C\eps}(E_{k-1,i}(t)^{1/2}+E_{k+1,i-1}(t)^{1/2})+C\eps t^{-3/2+C\eps}.}In other words, the inequality \eqref{sec6fff} holds. Note that all the  constants here are chosen before any of $C_{k,i},B_{k,i},B_0,B_1$ is chosen.}\rm\bigskip

We start with the easiest estimate.
\lem{ For each $k,i\geq 0$ with $i+k\leq N$, we have $\mcl{R}_6^{k,i}\lesssim \eps t^{-3/2+C\eps}$. The constants here are independent of the choice of $C_{k,i},\ B_{k,i}$, etc. }
\begin{proof}
First we note the following fact.  If $P(t,x)=a_0+c^\alpha x_\alpha$ where the $a_0$ and $c^\alpha$ are constants, then for all multiindices $J$ with $|J|>0$, $Z^JP$ is  of the form $Z^JP=c_{J}^{\alpha} x_\alpha$ where the  $c_J^\alpha$ are constants. As a result, we have
\fm{Z^{J}(t/T)=O((1+t+r)/T),\qquad \partial Z^{J}(t/T)=O(1/T),\qquad \partial^2 Z^{J}(t/T)=0.}
In summary, whenever $1\lesssim r\lesssim T$ and $1\lesssim t\sim T$, we have
\eq{\label{sec6.3l6f1}|\partial^lZ^J (t/T)|\lesssim t^{-l}1_{l\leq 1}.  }

Now fix $k\geq 0$ and a multtindenx $I$. By the chain rule, we can express $\partial^kZ^I\chi(t/T)$ as a linear combination of terms of the form
\fm{\chi^{(m)}(t/T)\cdot \prod_{s=1}^m \partial^{k_s}Z^{I_s}(t/T),\qquad \text{whenever }m\geq0,\ \sum k_*=k,\ \sum |I_s|=|I|,\ k_*+|I_*|>0.}
By \eqref{sec6.3l6f1} we may also require that $k_s\leq 1$, so $m\geq k$. It follows that 
\fm{|\partial^kZ^I\chi(t/T)|\lesssim\sum_{m=k}^{k+|I|}\prod_{k_1+\dots+k_m=k}t^{-k_s}\lesssim t^{-k}.}
It then follows by Leibniz's rule, Proposition \ref{mainprop4} and \eqref{equivnorm} that
\fm{\mcl{R}_{6}^{k,i}&\lesssim\sum_{k_1+k_2=k\atop |I_1|+|I_2|=i}\norm{\partial^{k_1}Z^{I_1}\chi(t/T)\cdot \partial^{k_2}Z^{I_2}(\eps S^{-3,1-\gamma_-}_-+\eps S^{-3,0}_+)1_{|r-t|\leq c t}}_{L^2(w)}\\
&\lesssim\sum_{k_1+k_2=k}\norm{t^{-k_1}\cdot \lra{r-t}^{-k_2}\eps t^{-3+C\eps}(1_{q\geq 0}+\lra{r-t}^{1-\gamma_-}1_{q<0})w_0(r-t)^{1/2}1_{|r-t|\leq ct}}_{L^2(\R^3)}\cdot t^{C\eps}\\
&\lesssim\norm{ \lra{r-t}^{-k}\eps t^{-3+C\eps}(1_{q\geq 0}+\lra{r-t}^{1-\gamma_-+\gamma_1/2}1_{q<0})1_{|r-t|\leq ct}}_{L^2(\R^3)}\cdot t^{C\eps}\\
&\lesssim \eps t^{-3+C\eps}\cdot t^{3/2}\lesssim \eps t^{-3/2+C\eps}.}
Note that to obtain the second last estimate, we use $1-\gamma_-+\gamma_1/2<0$.
\end{proof}\rm

\bigskip

Next we estimate $\mcl{R}_4^{k,i}$ and $\mcl{R}_5^{k,i}$. Here we need to apply Lemma \ref{l24}.
\lem{\label{sec6.3l4} For each $k,i$ with $i+k\leq N$, we have\eq{\label{sec6.3l4cc}\mcl{R}^{k,i}_4&\lesssim \eps t^{-1}E_{k,i}^{1/2}+ \eps^{1/2}(\ln\eps^{-1})^{-1/2} t^{-1+C\eps}\sum_{|I|\leq i}E_q(\partial^kZ^Iv)(t)^{1/2}\\&\quad+\eps t^{-1+C\eps}(E_{k,i-1}^{1/2}+E_{k-1,i}^{1/2})+\eps t^{-3/2}.}}
\begin{proof}By applying  Leibniz's rule, we can expand $R^{(*),k,I}_4$ as a linear combination of terms of the form
\eq{\label{sec6.3l4f1} \partial^{k_1}Z^{I_1}[ g^{\alpha\beta}(u_{app}+v,\partial (u_{app}+v))-g^{\alpha\beta}(u_{app},\partial u_{app})]\cdot \partial^{k_2} Z^{I_2}\partial_\alpha\partial_\beta u_{app}} for $k_1+k_2=k,\ |I_1|+|I_2|=|I|=i$. By Lemma \ref{l24}, we can control (with a constant factor $C$ independent of $B_{*,*}$ and $C_{*,*}$) \eqref{sec6.3l4f1}  by the sum of
\eq{\label{sec6.3l4f20}\abs{(g^{\alpha\beta}_J \partial^{k}Z^{I}v^{(J)} +g^{\alpha\beta\lambda}_J \partial^{k}Z^{I}\partial_\lambda v^{(J)})\cdot \partial_\alpha\partial_\beta u_{app}},}
\eq{\label{sec6.3l4f2}&\abs{g^{\alpha\beta}_J \partial^{k_1}Z^{I_1}v^{(J)} +g^{\alpha\beta\lambda}_J \partial^{k_1}Z^{I_1}\partial_\lambda v^{(J)}}\cdot|\partial^{k_2}Z^{I_2}\partial^2u_{app}|,\\&\qquad\qquad\qquad k_1+k_2=k,\ |I_1|+|I_2|=i,\ k_1+|I_1|<k+i,}
\eq{\label{sec6.3l4f3}&\sum_{k'+k''\leq k_1\atop |I'|+|I''|\leq |I_1|}(|\partial^{k'}Z^{I'}(v,\partial v)|(|\partial^{k''}Z^{I''}(u_{app},\partial u_{app})|+|\partial^{k''}Z^{I''}(v,\partial v)|))\cdot|\partial^{k_2}Z^{I_2}\partial^2u_{app}|,\\&\qquad\qquad\qquad\qquad k_1+k_2=k,\ |I_1|+|I_2|=i.}

We start with \eqref{sec6.3l4f20}. Since $v(t)
\in C_c^\infty$ for each $t$, by  Lemma  \ref{lp3} we have 
\fm{&\norm{g^{\alpha\beta}_J \partial^{k}Z^{I}v^{(J)} \cdot \partial_\alpha\partial_\beta u_{app}}_{L^2(w)}\\&\lesssim \eps t^{-1}\norm{\partial^{k+1}Z^Iv(t)}_{L^2(w)}+\eps t^{-1}\norm{r^{-\frac{\gamma_1+2}{4}}\partial^kZ^Iv(t)1_{r/t\in[1-2c,1-c]}}_{L^2(w)}\\
&\lesssim \eps t^{-1}\norm{\partial^{k+1}Z^Iv(t)}_{L^2(w)}+\eps t^{-1}\norm{\partial^{k}Z^Iv(t)}_{L^2(w)}\cdot 1_{k>0}\\&\quad+\eps t^{-1}\sum_{|J|=|I|-1}\norm{r^{-\frac{\gamma_1+2}{4}} Z^Jv(t)1_{r/t\in[1-2c,1-c]}}_{L^2(w)}\cdot 1_{k=0}\\
&\lesssim \eps t^{-1}E_{k,i}^{1/2}+\eps t^{-\frac{\gamma_1+2}{4}}E_{0,i-1}^{1/2}+\eps t^{-\frac{\gamma_1+4}{4}+C_N\eps}\lesssim \eps t^{-1}E_{k,i}^{1/2}+\eps t^{-3/2}.
}
To obtain the second last estimate, we apply \eqref{ptb2cor} and \eqref{ptb2cor2}. To obtain the last one, we recall that $\gamma_1>2$. Moreover, since $[\partial,Z]=C\partial$, by \eqref{lem4.10c0}  we have
\fm{&|g^{\alpha\beta\lambda}_J \partial^{k}Z^{I}\partial_\lambda v^{(J)}\cdot \partial_\alpha\partial_\beta u_{app}|\\
&\lesssim |g^{\alpha\beta\lambda}_J \partial_\lambda\partial^{k}Z^{I} v^{(J)}\cdot \partial_\alpha\partial_\beta u_{app}|+\sum_{|J|<|I|}|\partial^{k+1}Z^Jv^{(J)}|\cdot | \partial^2 u_{app}|\\
&\lesssim \eps t^{-1}|\partial^{k+1}Z^{I} v^{(J)}|+\eps t^{-1+C\eps}|\wt{T}\partial^{k}Z^{I} v^{(J)}|\lra{q}^{-1-\gamma_{\sgn(q)}}+\sum_{|J|<|I|}|\partial^{k+1}Z^Jv^{(J)}|\cdot | \partial^2 u_{app}|.}
Here recall that $\wt{T}_\alpha:=q_t\partial_\alpha -q_\alpha\partial_t$. By \eqref{sec5eneq}, we have \fm{\norm{\lra{q}^{-1-\gamma_{\sgn(q)}}\wt{T}\partial^kZ^Iv(t)}_{L^2(w)}^2&\lesssim  t^{C\eps}\int_{\R^3}|q_t|^{-1}\lra{q}^{-1-\gamma_2}|\wt{T}\partial^kZ^Iv(t)|^2\ dx\\&\lesssim t^{C\eps}(\eps\ln\eps^{-1})^{-1}E_q(\partial^kZ^Iv)(t).} Also note  that the $L^2(w)$ norm of the sum involving $|\partial^2u_{app}|$ is controlled by 
\fm{& \eps t^{-1+C\eps}\sum_{|J|<i}\norm{\lra{r-t}^{-2}(1_{q\geq 0}+\lra{r-t}^{1-\gamma_-}1_{q<0})\partial^{k+1}Z^Jv(t)}_{L^2(w)}\\
&\lesssim \eps t^{-1+C\eps}\sum_{|J|<i}\norm{\lra{r-t}^{-2}(1_{q\geq 0}+\lra{r-t}^{1-\gamma_-+\gamma_1/2}1_{q<0})\partial^{k+1}Z^Jv(t)}_{L^2(\R^3)}\\
&\lesssim \eps t^{-1+C\eps}E_{k,i-1}(t)^{1/2}.}
Here we use $1-\gamma_-+\gamma_1/2<0$. In summary, the $L^2(w)$ norm of \eqref{sec6.3l4f20} can be controlled by the right side of \eqref{sec6.3l4cc}.

Next we estimate \eqref{sec6.3l4f2}. By Proposition \ref{mainprop4} and Lemma \ref{c1l2.1}, we have
\fm{|\partial^{k_2}Z^{I_2}\partial^2u_{app}|\lesssim \eps t^{-1+C\eps}\lra{r-t}^{-2-k_2}(1_{q\geq 0}+\lra{r-t}^{1-\gamma_-}1_{q<0})1_{r/t\in[1-c,1+c]}.}
By \eqref{rmk4.1.1c2}  and \eqref{lem4.1c3}, we can replace $1_{q\geq 0}$ and $1_{q<0}$ with $1_{r\geq t}$ and $1_{r<t}$ respectively. Because of \eqref{equivnorm}, we have
\fm{&\norm{\abs{g^{\alpha\beta}_J \partial^{k_1}Z^{I_1}v^{(J)} +g^{\alpha\beta\lambda}_J \partial^{k_1}Z^{I_1}\partial_\lambda v^{(J)}}\cdot|\partial^{k_2}Z^{I_2}\partial^2u_{app}|}_{L^2(w)}\\
&\lesssim\norm{\partial^{k_1}Z^{I_1}v \cdot\partial^{k_2}Z^{I_2}\partial^2u_{app}}_{L^2(w)}+\sum_{|J|\leq |I_1|}\norm{\partial^{k_1+1}Z^{J} v \cdot \partial^{k_2}Z^{I_2}\partial^2u_{app}}_{L^2(w)}\\
&\lesssim \eps t^{-1+C\eps}\norm{\partial^{k_1}Z^{I_1}v \cdot\lra{r-t}^{-2}(1_{r\geq t}+\lra{r-t}^{1-\gamma_-+\gamma_1/2}1_{r<t})1_{r/t\geq 1-c}}_{L^2(\R^3)}+\eps t^{-1+C\eps}E_{k_1,|I_1|}^{1/2}\\
&\lesssim\eps t^{-1+C\eps}\norm{\partial^{k_1}Z^{I_1}v \cdot\lra{r-t}^{-2}1_{r/t\geq 1-c}}_{L^2(\R^3)}+\eps t^{-1+C\eps}E_{k_1,|I_1|}^{1/2} \\
&\lesssim \eps t^{-1+C\eps}E_{k,i-1}^{1/2}+\eps t^{-1+C\eps}E_{k-1,i}^{1/2}+\eps t^{-1+C\eps}\norm{(1_{r\geq t}+\lra{r-t}^{\frac{\gamma_1+2}{4}}1_{r<t})\partial^{k_1+1}Z^{I_1}v }_{L^2(\R^3)}\\&\quad+\eps t^{-1+C\eps}\norm{\lra{r-t}^{\frac{\gamma_1+2}{4}}r^{-1}\partial^{k_1}Z^{I_1}v1_{r/t\in[1-2c,1-c]}}_{L^2(\R^3)}\\
&\lesssim \eps t^{-1+C\eps}E_{k,i-1}^{1/2}+\eps t^{-1+C\eps}E_{k-1,i}^{1/2}+\eps t^{-3/2}.}
In the first estimate, we apply the triangle inequality and use $[Z,\partial]=C\partial$ to commute $Z^{I_1}$ and $\partial_{\lambda}$.
In the third one, we use $1-\gamma_-+\gamma_1/2<0$. In the fourth one, we apply Lemma \ref{lp1} with $\eta=\frac{\gamma_1}{2}$ and $\phi=\partial^{k_1}Z^{I_1}v$. The last estimate follows from  \eqref{equivnorm}, \eqref{ptb2cor} and \eqref{ptb2cor2}. In summary, the $L^2(w)$ norm of \eqref{sec6.3l4f2} can be controlled by the right side of \eqref{sec6.3l4cc}.

Finally we estimate \eqref{sec6.3l4f3}. Note that $Z^{J}(u_{app},\partial u_{app})=O(\eps t^{-1+C\eps}1_{r/t\in[1-c,1+c]})$ for each $J$, that $Z^J(v,\partial v)=O(C_N\eps t^{-1+C_N\eps})$ whenever $|J|\leq N-2$ by \eqref{ptb2}, and that $\min\{k'+1+|I''|,k''+|I''|+1\}\leq (k+i+2)/2\leq (N+2)/2\leq N-2$. So \eqref{sec6.3l4f3} is controlled by
\fm{
&C_N\eps t^{-1+C_N\eps}\sum_{k_1+k_2\leq k\atop |I_1|+|I_2|\leq i}|\partial^{k_1}Z^{I_1}v| |\partial^{k_2}Z^{I_2}\partial^2 u_{app}|\\
&\lesssim C_N\eps t^{-2+C_N\eps}\sum_{k_1\leq k\atop |I_1|\leq i}|\partial^{k_1}Z^{I_1}v| \lra{r-t}^{-2}(1_{q\geq 0}+\lra{r-t}^{1-\gamma_-}1_{q<0})\cdot 1_{r/t\geq 1-c}.}
By Lemma \ref{lp1} with $\eta=\gamma_1/2$, \eqref{equivnorm} and \eqref{rmk4.1.1c2}, we have
\fm{&C_N\eps t^{-2+C_N\eps}\sum_{k_1\leq k\atop |I_1|\leq i}\norm{|\partial^{k_1}Z^{I_1}v| \lra{r-t}^{-2}(1_{q\geq 0}+\lra{r-t}^{1-\gamma_-}1_{q<0})\cdot 1_{r/t\geq 1-c}}_{L^2(w)}\\
&\lesssim C_N\eps t^{-2+C_N\eps}\sum_{k_1\leq k\atop |I_1|\leq i}\norm{|\partial^{k_1}Z^{I_1}v| \lra{r-t}^{-2}(1_{r\geq t}+\lra{r-t}^{1-\gamma_-+\gamma_1/2}1_{r<t})\cdot 1_{r/t\geq 1-c}}_{L^2(\R^3)}\\
&\lesssim C_N\eps t^{-2+C_N\eps}\\&\cdot\sum_{k_1\leq k\atop |I_1|\leq i}(\norm{(1_{r\geq t}+\lra{r-t}^{\frac{\gamma_1+2}{4}}1_{r<t})\partial^{k_1+1}Z^{I_1}v}_{L^2(\R^3)}+\norm{r^{\frac{\gamma_1-2}{4}}\partial^{k_1}Z^{I_1}v\cdot 1_{r/t\in[1-2c,1-c]} }_{L^2(\R^3)})\\
&\lesssim C_N\eps t^{-2+C_N\eps}E_{k,i}^{1/2}+C_N\eps t^{-2+C_N\eps}\cdot C_N\eps t^{-\gamma_1/4+C_N\eps}.}
By choosing $\eps\ll1$,  we conclude that the $L^2(w)$ norm of \eqref{sec6.3l4f3} is controlled by the right side of \eqref{sec6.3l4cc}.
\end{proof}\rm

\lem{\label{sec6.3l5} For each $k,i$ with $k+i\leq N$, we have\eq{\label{sec6.3l5c}\mcl{R}_5^{k,i}\lesssim \eps t^{-1}E_{k,i}^{1/2}+\eps t^{-1+C\eps}(E^{1/2}_{k,i-1}+E^{1/2}_{k-1,i})+\eps t^{-3/2}.}}
\begin{proof}
Recall the Taylor expansion \eqref{sec3tay} for $f$. Let us first consider the quadratic terms. These terms are of the form
\fm{&\partial^kZ^I(f_{JK}^{*,\alpha\beta}(\partial_\alpha u_{app}^{(J)}+\partial_\alpha v^{(J)})(\partial_\beta u_{app}^{(K)}+\partial_\beta v^{(K)})-f_{JK}^{*,\alpha\beta} \partial_\alpha u_{app}^{(J)}  \partial_\beta u_{app}^{(K)})\\
&=f_{JK}^{*,\alpha\beta}\partial^kZ^I( \partial_\alpha u_{app}^{(J)} \partial_\beta v^{(K)}+\partial_\alpha v^{(J)} \partial_\beta u_{app}^{(K)}+ \partial_\alpha v^{(J)}  \partial_\beta v^{(K)}),}
so they are controlled pointwisely by a linear combination of terms of the form
\fm{&|\partial^{k_1}Z^{I_1}\partial u_{app}\cdot\partial^{k_2}Z^{I_2}\partial v|+|\partial^{k_1}Z^{I_1}\partial v\cdot\partial^{k_2}Z^{I_2}\partial v|,}where $k_1+k_2=k,\ |I_1|+|I_2|=i$. Here we have $\min\{k_1+|I_1|,k_2+|I_2|\}\leq (k+|I|)/2\leq N/2$. By applying \eqref{ptb1} to $\partial^{k_*}Z^{I_*}\partial v$ with smaller $k_*+|I_*|$, we notice that the $L^2(w)$ norm of these terms are controlled by
\fm{&\eps t^{-1}\norm{\partial^kZ^I\partial v}_{L^2(w)}+\eps t^{-1+C\eps}\sum_{k_2\leq k,\ |I_2|\leq i\atop k_2+|I_2|<k+|I|}\norm{\partial^{k_2}Z^{I_2}\partial v}_{L^2(w)}+C_N\eps t^{-3/2+C_N\eps}\sum_{k_2\leq k\atop |I_2|\leq i}\norm{\partial^{k_2}Z^{I_2}\partial v}_{L^2(w)}\\
&\lesssim \eps t^{-1}\norm{\partial^{k+1}Z^Iv}_{L^2(w)}+\eps t^{-1+C\eps}\sum_{k_2\leq k,\ |I_2|\leq i\atop k_2+|I_2|<k+|I|}\norm{ \partial^{k_2+1}Z^{I_2} v}_{L^2(w)}+C_N\eps t^{-3/2+C_N\eps}\sum_{k_2\leq k\atop |I_2|\leq i}\norm{\partial^{k_2+1}Z^{I_2} v}_{L^2(w)}\\
&\lesssim \eps t^{-1}E^{1/2}_{k,i}+\eps t^{-1+C\eps}(E^{1/2}_{k,i-1}+E^{1/2}_{k-1,i}).}
To get the first estimate, we note that $\partial u_{app}=O(\eps t^{-1})$. In the last estimate, we use \eqref{equivnorm} and $C_N\eps t^{-3/2+C_N\eps}\leq C_N\eps t^{-5/4}\leq \eps t^{-1}$ for $\eps \ll1$ and $t\geq T_\eps$.

Next we consider the cubic terms. These terms are linear combinations (with constant coefficients) of 
\eq{\label{sec6.3l5f1}\partial^{k}Z^{I}(\partial^{p_1}u_{app}\partial^{p_2}u_{app}\partial^{p_3}v+\partial^{p_1}u_{app}\partial^{p_2}v\partial^{p_3}v+\partial^{p_1}v\partial^{p_2}v\partial^{p_3}v),\qquad p_*\leq1,\ \sum p_*>0.}
We note that the constraint $\sum p_*>0$ results from $f(u,0)=O(|u|^4)$. By  Leibniz's rule, we expand \eqref{sec6.3l5f1} as a linear combination of  terms of the forms
\eq{\label{sec6.3l5f2}\partial^{k_1}Z^{I_1}\partial^{p_1} u_{app} \partial^{k_2}Z^{I_2}\partial^{p_2} u_{app} \partial^{k_3}Z^{I_3}\partial^{p_3}v;&\\\partial^{k_1}Z^{I_1}\partial^{p_1} u_{app}\partial^{k_2}Z^{I_2}\partial^{p_2}v\partial^{k_3}Z^{I_3}\partial^{p_3}v,&\qquad k_2+|I_2|+p_2\leq k_3+|I_3|+p_3;\\\partial^{k_1}Z^{I_1}\partial^{p_1}v\partial^{k_2}Z^{I_2}\partial^{p_2}v\partial^{k_3}Z^{I_3}\partial^{p_3}v,&\qquad k_1+|I_1|+p_1\leq k_2+|I_2|+p_2\leq k_3+|I_3|+p_3.}
Here $\sum k_*=k,\ \sum |I_*|=i,\ p_*\leq1,\ \sum p_*>0$.  Moreover, by Proposition \ref{mainprop4}, \eqref{ptb1}, and \eqref{ptb2}, we have \fm{\sum_{|J|\leq N-2}(|Z^Ju_{app}|+|Z^Jv|)&\lesssim C_N\eps t^{-1+C_N\eps}(1_{r\geq t}+\lra{r-t}^{-\gamma_1/2}1_{r<t});\\
\sum_{|J|\leq N-2}(|Z^J\partial u_{app}|+|Z^J\partial v|)&\lesssim C_N\eps t^{-1+C_N\eps}\lra{r-t}^{-1}(1_{r\geq t}+\lra{r-t}^{-\gamma_1/2}1_{r<t}).}
For terms in \eqref{sec6.3l5f2} with $k_3+p_3>0$, their $L^2(w)$ norms are controlled by \fm{&C_N\eps^2 t^{-2+C_N\eps}\norm{(1_{r\geq t}+\lra{r-t}^{-\gamma_1}1_{r<t})|\partial^{k_3}Z^{I_3}\partial^{p_3}v(t)|}_{L^2(w)}\\
&\lesssim C_N\eps^2 t^{-2+C_N\eps}\norm{\partial^{k_3}Z^{I_3}\partial^{p_3}v(t)}_{L^2(\R^3)}t^{C\eps}\lesssim\eps t^{-1}E_{k,i}^{1/2}.}
For terms with $k_3+p_3=0$, we have $p_1+p_2>0$, so the $L^2(w)$ norms of those terms are controlled by
\fm{&C_N\eps^2 t^{-2+C_N\eps}\norm{\lra{r-t}^{-1}(1_{r\geq t}+\lra{r-t}^{-\gamma_1}1_{r<t})|Z^{I_3}v(t)|}_{L^2(w)}\\
&\lesssim C_N\eps^2 t^{-2+C_N\eps}\norm{\lra{r-t}^{-1}Z^{I_3}v(t)}_{L^2(\R^3)}\cdot t^{C\eps}
\\
&\lesssim C_N\eps^2 t^{-2+C_N\eps}\norm{\lra{r-t}^{-1}Z^{I_3}v(t)1_{r/t\geq 1-c}}_{L^2(\R^3)}+C_N\eps^2 t^{-2+C_N\eps}\norm{r^{-1}Z^{I_3}v(t)1_{r/t< 1-c}}_{L^2(\R^3)}\\
&\lesssim C_N\eps^2 t^{-2+C_N\eps}\norm{(1_{r\geq t}+\lra{r-t}^{\frac{\gamma_1+2}{4}}1_{r<t})\partial Z^{I_3}v(t)}_{L^2(\R^3)}\\&\quad+C_N\eps^2 t^{-2+C_N\eps}\norm{\lra{r-t}^{\frac{\gamma_1+2}{4}}r^{-1}Z^{I_3}v(t)1_{r/t\in[1-2c,1-c]}}_{L^2(\R^3)}+ C_N\eps^2 t^{-2+C_N\eps}\norm{\partial Z^{I_3}v(t)}_{L^2(\R^3)}\\
&\lesssim \eps  t^{-1}E_{k,i}^{1/2}+\eps t^{-3/2}.}
Here we apply  Hardy's inequality and Lemma \ref{lp1} with $\eta=\gamma_1/2$. The last step has been proved in the proof of Lemma \ref{sec6.3l4}.

Finally, let us consider the terms of order larger than $3$. By Lemma \ref{l24}, those terms are controlled by a linear combination of terms of the form
\eq{\label{sec6.3l5f3}|\partial^{k_1}Z^{I_1}v|\prod_{j=2}^4 |\partial^{k_j}Z^{I_j}(v,\partial v,u_{app},\partial u_{app})|.}
Here $k_*\leq k$, $|I_*|\leq i$ and $k_j+|I_j|>(k+i)/2$ for at most one index $1\leq j\leq 4$. Using the pointwise bounds for $v$ and $u_{app}$ and the support of $u_{app}$, we control \eqref{sec6.3l5f3} by 
\fm{&C_N\eps^3 t^{-3+C_N\eps}(1_{r\geq t}+\lra{r-t}^{-3\gamma_1/2}1_{r<t})\sum_{l\leq k\atop |J|\leq i}(|\partial^lZ^Jv|+|\partial^lZ^J\partial v|)1_{r/t\in[1-c,1+c]}\\&+C_N\eps^3 t^{-3+C_N\eps}(1_{r\geq t}+\lra{r-t}^{-3\gamma_1/2}1_{r<t})\sum_{1\leq l\leq k+1\atop |J|\leq i}|\partial^lZ^Jv|\\&+C_N\eps^3 t^{-3+C_N\eps}\lra{r-t}^{-1}(1_{r\geq t}+\lra{r-t}^{-3\gamma_1/2}1_{r<t})\sum_{ |J|\leq i}|Z^Jv|\\&+C_N\eps t^{-1+C_N\eps}(1_{r\geq t}+\lra{r-t}^{-\gamma_1/2}1_{r<t})\sum_{ |J|\leq i}|Z^Jv|^3.}
Since $w^{1/2}\lesssim t^{C\eps}(1_{r\geq t}+\lra{r-t}^{\gamma_1/2}1_{r<t})$ and since $\lra{r-t}\lesssim t$ whenever $r\sim t$, the $L^2(w)$ norm of the formula above is controlled by the $L^2(\R^3)$ norm of
\fm{&C_N\eps^3 t^{-3+C_N\eps}(1_{r\geq t}+\lra{r-t}^{-\gamma_1}1_{r<t})\sum_{1\leq l\leq k+1\atop |J|\leq i}|\partial^lZ^Jv|\\&+C_N\eps^3 t^{-3+C_N\eps}\lra{r-t}^{-1}(1_{r\geq t}+\lra{r-t}^{-\gamma_1}1_{r<t})\sum_{ |J|\leq i}|Z^Jv|\\&+C_N\eps t^{-1+C_N\eps}\sum_{ |J|\leq i}|Z^Jv|^3.}
The $L^2(\R^3)$ norm of the first row is controlled by 
$C_N\eps^3t^{-3+C_N\eps}E_{k,i}^{1/2}$ because of \eqref{equivnorm}. To control the second row, we apply Lemma \ref{lp1} with $\eta=\gamma_1/2$ to obtain an upper bound
\fm{&C_N\eps^3t^{-3+C_N\eps}\sum_{|J|\leq i}\norm{(1_{r\geq t}+\lra{r-t}^{\frac{\gamma_1+2}{4}}1_{r<t})\partial Z^Jv}_{L^2(\R^3)}\\
&+C_N\eps^3t^{-3+C_N\eps}\sum_{|J|\leq i}\norm{\lra{r-t}^{\frac{\gamma_1+2}{4}}r^{-1} Z^Jv1_{r/t\in[1-2c,1-c]}}_{L^2(\R^3)}\\
&\lesssim C_N\eps^3t^{-3+C_N\eps} E_{k,i}^{1/2}+\eps t^{-3/2}.}
By the Sobolev embedding $L^6(\R^3)\hookrightarrow \dot{H}^1(\R^3)$, we control the $L^2(\R^3)$ norm of the last row by
\fm{C_N\eps t^{-1+C_N\eps}\sum_{|J|\leq i}\norm{\partial Z^Jv}_{L^2(\R^3)}^3\lesssim C_N\eps^3 t^{-2+C_N\eps}E_{k,i}^{1/2}.}
We end the proof by choosing $\eps\ll1$.
\end{proof}

\rmk{\rm If we include terms of the form $u\cdot u\cdot u$ in the Taylor expansions of $f^{*}(u,\partial u)$, then using the proof here we need to replace $\eps t^{-1}E_{k,i}^{1/2}$ with  $C_N\eps t^{-1+C_N\eps}E_{k,i}^{1/2}$ in \eqref{sec6.3l5c}. Such a weaker bound is not sufficient for our proof. It is unclear to the author how to improve these cubic terms, so we need to assume  \eqref{sec3tay}.}\rm

\bigskip

We continue with $\mcl{R}_2^{k,i}$ and $\mcl{R}_3^{k,i}$.
\lem{\label{sec6.3l2}For each $k,i$ with $k+i\leq N$, we have
\eq{\mcl{R}_{2}^{k,i}&\lesssim \eps t^{-1}E_{k,i}^{1/2}+\eps t^{-1+C \eps}(E_{k+1,i-1}^{1/2}+E_{k-1,i}^{1/2})\\&\quad+\eps^{1/2}(\ln\eps^{-1})^{-1/2} t^{-1+C\eps}\sum_{|I'|\leq i}E_q(\partial^{k}Z^{I'}v)(t)^{1/2}+\eps t^{-3/2}.}
\begin{proof}
We write $R_{2}^{(*),k,I}$ as a linear combination of terms of the form
\fm{\partial^{k_1}Z^{I_1}g^{\alpha\beta}(u_{app}+v,\partial(u_{app}+v))\cdot \partial^{k_2}Z^{I_2}\partial_\alpha\partial_\beta v^{(*)},\qquad k_1+k_2\leq k,\ |I_1|+|I_2|\leq i,\ k_1+|I_1|>0.}
We apply Lemma \ref{l24}  to obtain an upper bound
\eq{\label{sec6.3l2f1}&\sum_{k_1+k_2\leq k,\ |I_1|+|I_2|\leq i\atop k_1+|I_1|>0}\kh{|\partial^{k_1}Z^{I_1}(u_{app}+v)\partial^{k_2}Z^{I_2}\partial^2v|+|g^{\alpha\beta\lambda}_J\partial^{k_1}Z^{I_1}\partial_\lambda(u^{(J)}_{app}+v^{(J)})\partial^{k_2}Z^{I_2}\partial_\alpha\partial_\beta v|}\\&\quad+\sum_{k'+k''+k_2\leq k,\ |I'|+|I''|+|I_2|\leq i\atop k_2+|I_2|<k+i }\sum_{l',l''\leq 1}|\partial^{k'}Z^{I'}\partial^{l'}(u_{app}+v)||\partial^{k''}Z^{I''}\partial^{l''}(u_{app}+v)|\cdot|\partial^{k_2}Z^{I_2}\partial^2v|.}

By setting $(k_2,|I_2|)=(k-1,i)$ in the first sum,  we have $(k_1,|I_1|)=(1,0)$ and 
\fm{&\sum_{|I'|=i}\kh{|\partial (u_{app}+v)\partial^{k-1}Z^{I'}\partial^2v|+|g^{\alpha\beta\lambda}_J\partial \partial_\lambda(u^{(J)}_{app}+v^{(J)})\partial^{k-1}Z^{I'}\partial_\alpha\partial_\beta v|}\\
&\lesssim (|\partial u_{app}|+|\partial v|)\sum_{|I'|\leq i} |\partial^{k+1}Z^{I'}v|+\sum_{|I'|<i}|\partial^2(u^{(J)}_{app}+v^{(J)})\partial^{k+1} Z^{I'} v|\\
&\quad+\sum_{|I'|=i} (|g^{\alpha\beta\lambda}_J\partial_\lambda\partial  u^{(J)}_{app} \partial_\alpha\partial_\beta\partial^{k-1}Z^{I'} v|+|\partial^2 v \partial^{k+1}Z^{I'} v|) \\
&\lesssim (\eps t^{-1}+C_N\eps t^{-3/2+C_N\eps})\sum_{|I'|\leq i} |\partial^{k+1}Z^{I'}v|+(\eps t^{-1+C\eps}+C_N\eps t^{-3/2+C_N\eps})\sum_{|I'|< i} |\partial^{k+1}Z^{I'}v|\\
&\quad+\sum_{|I'|=i} |g^{\alpha\beta\lambda}_J\partial_\lambda\partial  u^{(J)}_{app}\partial_\alpha\partial_\beta\partial^{k-1}Z^{I'} v| .}
We then apply \eqref{lem4.10c1} in Lemma \ref{lem4.10} to obtain 
\fm{&\sum_{|I'|=i} |g^{\alpha\beta\lambda}_J\partial_\lambda\partial  u^{(J)}_{app}\partial_\alpha\partial_\beta\partial^{k-1}Z^{I'} v|&\lesssim \sum_{|I'|=i} (\eps t^{-1}|\partial^{k+1}Z^{I'} v|+\eps t^{-1+C\eps}\lra{q}^{-1-\gamma_{\sgn(q)}}|\wt{T}\partial^{k}Z^{I'}v|).}
Following the proof in Lemma \ref{sec6.3l4}, we have the $L^2(w)$ norm of this part is controlled by 
\fm{\eps t^{-1}E_{k,i}^{1/2}+\eps^{1/2}(\ln \eps^{-1})^{-1/2} t^{-1+C\eps}\sum_{|I'|=i}E_q(\partial^kZ^{I'}v)^{1/2}.}

Next, we estimate
\eq{\label{sec6.3l2f2}&\sum_{k_1+k_2\leq k,\ |I_1|+|I_2|\leq i\atop k_1+|I_1|>0,\ (k_2,|I_2|)\neq (k-1,i)}\kh{|\partial^{k_1}Z^{I_1}(u_{app}+v)\partial^{k_2}Z^{I_2}\partial^2v|+|g^{\alpha\beta\lambda}_J\partial^{k_1}Z^{I_1}\partial_\lambda(u^{(J)}_{app}+v^{(J)})\partial^{k_2}Z^{I_2}\partial_\alpha\partial_\beta v|}\\
&\lesssim \sum_{k_2\leq k,\ |I_2|\leq i\atop k_2+|I_2|<k+i,\ (k_2,|I_2|)\neq (k-1,i)}\eps t^{-1+C\eps}|\partial^{k_2+2}Z^{I_2}v|+\sum_{k_1+k_2\leq k,\ |I_1|+|I_2|\leq i\atop k_1+|I_1|>0,\ (k_2,|I_2|)\neq (k-1,i)}|\partial^{k_1+1}Z^{I_1} v\partial^{k_2+2}Z^{I_2}v|\\&\quad+\sum_{k_1\geq 1,\ k_1+k_2\leq k,\ |I_1|+|I_2|\leq i\atop (k_2,|I_2|)\neq (k-1,i)}|\partial^{k_1}Z^{I_1}v\partial^{k_2+2}Z^{I_2}v|+\sum_{k_2\leq k,\ |I_1|+|I_2|\leq i\atop |I_1|>0}|Z^{I_1}v\partial^{k_2+2}Z^{I_2}v|.}
Since $\min\{k_1+|I_1|,k_2+|I_2|+1\}\leq N/2+1\leq N-2$ whenever $N\geq 6$, by \eqref{ptb1} and \eqref{ptb2} we have
\fm{|\partial^{k_1}Z^{I_1}v|\lesssim C_N\eps t^{-1-\frac{1}{2}1_{k_1>0}+C_N\eps}\qquad\text{or}\qquad |\partial^{k_2+1} ZZ^{I_2}v|\lesssim C_N\eps t^{-3/2+C_N\eps};}
\fm{|\partial^{k_1+1}Z^{I_1}v|\lesssim C_N\eps t^{-3/2+C_N\eps}\qquad\text{or}\qquad |\partial^{k_2+1}Z Z^{I_2}v|\lesssim C_N\eps t^{-3/2+C_N\eps}.}
Also note that if $|I_1|\leq N-2$, by \eqref{ptb2} we have 
\fm{\lra{r-t}^{-1}|Z^{I_1}v|&\lesssim C_N\eps t^{-3/2+C_N\eps}(\lra{r-t}^{-(\gamma_1+1)/2}1_{r<t}+\lra{r-t}^{-1/2}1_{t\leq r\leq 2t}+t^{1/2}\lra{r-t}^{-1}1_{r>2t})\\
&\lesssim C_N\eps t^{-3/2+C_N\eps};}
if $|I_1|\geq N-1$, then we must have $k_2+|I_2|\leq 1$ in the last sum in \eqref{sec6.3l2f2}, so  there we have
\fm{w^{1/2}|\partial^{k_2+2}Z^Iv|\lesssim C_N\eps t^{-3/2+C_N\eps}\lra{r-t}^{-3/2}.}
It follows that \eqref{sec6.3l2f2} can be controlled by
\fm{&\sum_{k_2\leq k,\ |I_2|\leq i\atop k_2+|I_2|<k+i,\ (k_2,|I_2|)\neq (k-1,i)}\eps t^{-1+C\eps}|\partial^{k_2+2}Z^{I_2}v|\\&+\sum_{k_1\leq k,\ |I_1|\leq i\atop k_1+|I_1|>0}C_N\eps t^{-3/2+C_N\eps}|\partial^{k_1+1}Z^{I_1} v|+\sum_{k_2\leq k,\ |I_2|\leq i\atop k_2+|I_2|<k+i,\ (k_2,|I_2|)\neq (k-1,i)}C_N\eps t^{-3/2+C_N\eps}|\partial^{k_2+2}Z^{I_2}v|\\&+\sum_{1\leq k_1\leq k,\ |I_1|\leq i}C_N\eps t^{-3/2+C_N\eps}|\partial^{k_1}Z^{I_1}v|+\sum_{k_2\leq k-1,\ |I_2|\leq i\atop (k_2,|I_2|)\neq (k-1,i)}C_N\eps t^{-3/2+C_N\eps}|\partial^{k_2+2}Z^{I_2}v|\\
&+\sum_{k_2\leq k,\ |I_2|<i}C_N\eps t^{-3/2+C_N\eps}|\partial^{k_2+1}ZZ^{I_2}v|+\sum_{0<|I_1|\leq i}C_N\eps t^{-3/2+C_N\eps}|w^{-1/2}\lra{r-t}^{-3/2}Z^{I_1}v|.}
Each row corresponds to one term in \eqref{sec6.3l2f2}. Now we take the $L^2(w)$ norm and let $\eps\ll1$. We conclude that the $L^2(w)$ norm of \eqref{sec6.3l2f2} is controlled by 
\fm{\eps t^{-1+C\eps}(E_{k+1,i-1}^{1/2}+E_{k-1,i}^{1/2})+\eps t^{-1}E_{k,i}^{1/2}+\eps t^{-3/2}.}
Here we only explain how to estimate the $L^2(w)$  norm of $|w^{-1/2}\lra{r-t}^{-3/2}Z^{I_1}v|$. In fact, by Lemma \ref{lp1} with $\eta=\gamma_1/2$, we have
\fm{&\norm{w^{-1/2}\lra{r-t}^{-3/2}Z^{I_1}v}_{L^2(w)}\lesssim \norm{\lra{r-t}^{-1}Z^{I_1}v}_{L^2(\R^3)}\\&\lesssim \norm{(1_{r\geq t}+\lra{r-t}^{\frac{\gamma_1+2}{4}}1_{r<t})\partial Z^{I_1}v}_{L^2(\R^3)}+\norm{r^{-1}\lra{r-t}^{\frac{\gamma_1+2}{4}} Z^{I_1}v1_{r/t\in[1-2c,1-c]}}_{L^2(\R^3)}\\
&\quad+\norm{r^{-1}Z^{I_1}v1_{r/t\leq 1-c}}_{L^2(\R^3)}\\
&\lesssim t^{C\eps}E_{k,i}^{1/2}+C_N\eps t^{-\gamma_1/4+C_N\eps}.
}
In the last estimate apply  Hardy's inequality, \eqref{equivnorm},  \eqref{ptb2cor} and \eqref{ptb2cor2}.

Finally, we estimate the second sum in \eqref{sec6.3l2f1}. Without loss of generality, we assume $k'+l'+|I'|\leq k''+l''+|I''|$, so $k'+l'+|I'|\leq N/2+1\leq N-2$. We can thus replace $|\partial^{k'}Z^{I'}\partial^{l'}(u_{app}+v)|$ with $C_N\eps t^{-3/2+C_N\eps}$. If $l''=0$, then the product is controlled by the first sum in \eqref{sec6.3l2f1}. If $l''=1$, then the product is controlled by
\fm{C_N\eps t^{-3/2+C_N\eps}\sum_{k''+k_2\leq k,\ |I''|+|I_2|\leq i\atop k_2+|I_2|<k+i}|\partial^{k''}Z^{I''}\partial (u_{app}+v)||\partial^{k_2}Z^{I_2}\partial^2v |.}
Since $\min\{k''+|I''|,k_2+|I_2|+1\}\leq N/2+1$, we can control the sum above by
\fm{C_N\eps^2 t^{-5/2+C_N\eps}\sum_{k_2\leq k,\ |I_2|\leq i\atop k_2+|I_2|<k+i}|\partial^{k_2}Z^{I_2}\partial^2v |+C_N\eps^2 t^{-5/2+C_N\eps}\sum_{k''\leq k,\ |I''|\leq i}|\partial^{k''}Z^{I''}\partial v|.}
The $L^2(w)$ norm of this formula is controlled by $\eps t^{-1}(E_{k,i}^{1/2}+E_{k+1,i-1}^{1/2})$.
\end{proof}

\lem{\label{sec6.3l3}For each $k,i$ with $k+i\leq N$, we have
\eq{\mcl{R}_{3}^{k,i}\lesssim \eps t^{-1}E_{k,i}^{1/2}+\eps t^{-1+C\eps}E_{k+1,i-1}^{1/2}.}}
\begin{proof}
We can write $R_{3}^{(*),k,I}$ as a linear combination of terms of the form
\fm{&(g^{\alpha\beta}(u_{app}+v,\partial(u_{app}+v))-m^{\alpha\beta})\partial^{k+2}Z^{I'}v\\
&=[g_J^{\alpha\beta}(u^{(J)}_{app}+v^{(J)})+g_{J}^{\alpha\beta\lambda}\partial_\lambda (u^{(J)}_{app}+v^{(J)})+O((|u_{app}|+|v|+|\partial u_{app}|+|\partial v|)^2)]\partial^{k+2}Z^{I'}v}for $|I'|<i$.
This is because $[\partial,Z]=C\cdot \partial$ for some constant $C$. Since $\partial (u_{app},v)=O(\eps t^{-1})$, $u_{app}=O(\eps t^{-1+C\eps})$ and $v=O(C_N\eps t^{-1+C_N\eps})$, we have
\fm{&|(g^{\alpha\beta}(u_{app}+v,\partial(u_{app}+v))-m^{\alpha\beta})\partial^{k+2}Z^{I'}v|\\
&\lesssim |v||\partial^{k+2}Z^{I'}v|+(\eps t^{-1+C\eps}+(C_N\eps t^{-1+C_N\eps})^2)|\partial^{k+2}Z^{I'}v|\\
&\lesssim \lra{r-t}^{-1}|v||Z\partial^{k+1}Z^{I'}v|+\eps t^{-1+C\eps}|\partial^{k+2}Z^{I'}v|.}
In the previous lemma, we have shown that $\lra{r-t}^{-1}|v|\lesssim C_N\eps t^{-3/2+C_N\eps}$. Thus, by taking the $L^2(w)$ norm, we conclude that
\fm{\mcl{R}_3^{k,i}\lesssim \eps t^{-1}E_{k,i}^{1/2}+\eps t^{-1+C\eps}E_{k+1,i-1}^{1/2}.}
\end{proof}\rm

Finally, let us estimate $\mcl{R}_{1}^{k,i}$.
\lem{For each $k,i$ with $k+i\leq N$, we have
\eq{\mcl{R}_{1}^{k,i}&\lesssim \eps t^{-1}E_{k,i}^{1/2}+\eps t^{-1+C \eps}(E_{k+1,i-1}^{1/2}+E_{k-1,i}^{1/2})\\&\quad+\eps^{1/2}(\ln\eps^{-1})^{-1/2} t^{-1+C\eps}\sum_{l\leq k\atop|I'|\leq i}E_q(\partial^lZ^{I'}v)^{1/2}+\eps t^{-3/2+C\eps}.}}
\begin{proof}Note that $R_{1}^{(*),k,0}=0$ and 
\fm{R_{1}^{(*),k,I}&=\sum_{(I_1,I_2,I_3)=I,\ |I_2|=1}\partial^kZ^{I_1}[\Box,Z^{I_2}]Z^{I_3}v^{(*)}=\sum_{|I_1|+|I_3|=i-1}C\cdot\partial^kZ^{I_1}\Box Z^{I_3}v^{(*)}\\
&=\sum_{|I_1|+|I_3|=i-1}C\cdot\Box \partial^kZ^{I_1} Z^{I_3}v^{(*)}+\sum_{|I_1|+|I_3|=i-1}C\cdot\partial^k [Z^{I_1},\Box] Z^{I_3}v^{(*)}.}
Here we recall that $[\Box,Z]=C\cdot \Box$. Similarly, we can keep expanding $[Z^{I_1},\Box]$ and then move $\Box$ to the front by adding an extra term with commutators. Note that the order of derivatives in the term with commutator is strictly decreasing in each step, so after finitely many steps, we get no extra term with commutators. As a result, we can write $R_{1}^{(*),k,I}$ as a linear combination of terms of the form
\fm{\Box\partial^kZ^{I'}v^{(*)},\qquad |I'|<i.}

We now induct on $|I|$. When $|I|=0$, there is nothing to prove. If $|I|=i>0$, by making use of \eqref{eqn2}, we have
\fm{R_{1}^{(*),k,I}=\sum_{|I'|<i}C\cdot\Box\partial^kZ^{I'}v^{(*)}&=\sum_{|I'|<i}(C\cdot (g^{\alpha\beta}-m^{\alpha\beta})\partial_\alpha\partial_\beta \partial^kZ^{I'}v^{(*)}+\sum_{j=1}^6C\cdot R_{j}^{(*),k,I'}).}
The $L^2(w)$ norm of the first sum is controlled by \fm{&\sum_{|I'|<i}\norm{(|u_{app}|+|v|+|\partial u_{app}|+|\partial v|)|\partial^{k+2}Z^{I'}v|}_{L^2(w)}\\
&\lesssim \eps t^{-1+C\eps}E_{k+1,i-1}^{1/2}+\sum_{|I'|<i}\norm{\lra{r-t}^{-1}|v||Z\partial^{k+1}Z^{I'}v|}_{L^2(w)}\\
&\lesssim \eps t^{-1+C\eps}E_{k+1,i-1}^{1/2}+C_N\eps t^{-3/2+C_N\eps}E_{k,i}^{1/2}.}
In the second sum, if $j\neq 1$, we can estimate $\mcl{R}_{j}^{(*),k,i'}$ with $i'<i$ by the previous five lemmas. As for $j=1$, we use the induction hypotheses. This finishes the proof.
\end{proof}\rm

\subsection{Improving of the energy estimate \eqref{ca1}}\label{sec6.4}

Fix $T_1\leq t\leq 2T$ and a pair $(k,i)$ with $k,i\geq 0$ and $k+i\leq N$. For each fixed multiindex $I$ with $|I|=i$, by \eqref{energymain} in Proposition \ref{prop5.1} we have
\fm{&E_u(\partial^kZ^Iv)(t)+\int_t^{2T} E_q(\partial^kZ^Iv)(\tau)\ d\tau\\
&\lesssim \int_t^\infty\norm{g^{\alpha\beta}(u,\partial u)\partial_\alpha\partial_\beta\partial^kZ^I(\tau)}_{L^2(w)}\norm{\partial\partial^kZ^Iv(\tau)}_{L^2(w)}+(C_N\eps \tau^{-1}+\tau^{-17/16})\norm{\partial\partial^kZ^Iv}^2_{L^2(w)}\ d\tau.}
This inequality holds for all $T_1\leq t\leq 2T$.  Here we recall that $C_N$ in this inequality depends on the constant $B$ in \eqref{sec5a1} and \eqref{sec5a2}. We also recall that $v\equiv 0$ whenever $t\geq 2T$. 

Now by applying \eqref{equivnorm} and \eqref{propsec6.3c} in Proposition \ref{propsec6.3}, we have
\fm{&\sum_{|I|=i}E_u(\partial^kZ^Iv)(t)+\int_t^{2T} \sum_{|I|=i}E_q(\partial^kZ^Iv)(\tau)\ d\tau\\
&\lesssim \sum_{|I|=i}\kh{\int_t^{2T}\sum_{j=1}^6\mcl{R}_{j}^{k,i}E_u(\partial^kZ^Iv)(\tau)^{1/2}+(C_N\eps \tau^{-1}+\tau^{-17/16})E_u(\partial^kZ^Iv)(\tau)\ d\tau}\\
&\lesssim   \int_t^{2T}(C_N\eps\tau^{-1}+\tau^{-17/16})E_{k,i}(\tau)+\eps^{1/2}(\ln\eps^{-1})^{-1/2} \tau^{-1+C\eps}\sum_{l\leq k\atop |I|\leq i}E_q(\partial^lZ^Iv)(\tau)^{1/2}E_{k,i}(\tau)^{1/2}\\&\quad+\eps \tau^{-1+C\eps}(E_{k-1,i}(\tau)^{1/2}+E_{k+1,i-1}(\tau)^{1/2})E_{k,i}(\tau)^{1/2}+\eps t^{-3/2+C\eps}E_{k,i}(\tau)^{1/2}\ d\tau.}
Also note that 
\fm{&\int_t^{2T} \tau^{-1+C\eps}\sum_{l\leq k\atop|I|\leq i}E_q(\partial^lZ^Iv)(\tau)^{1/2}E_{k,i}(\tau)^{1/2}\ d\tau\\
&\lesssim  \kh{\int_t^{2T}  \sum_{l\leq k\atop |I|\leq i}E_q(\partial^lZ^Iv)(\tau)\ d\tau}^{1/2}\kh{\int_t^{2T}  \tau^{-2+C\eps}E_{k,i}(\tau)\ d\tau}^{1/2}\\
&\lesssim   E_{k,i}(t)^{1/2}\kh{\int_t^{2T} \tau^{-2+C\eps}E_{k,i}(\tau)\ d\tau}^{1/2}.}
In the last step we use   the definition of the energy \eqref{ca1}.

Now we apply \eqref{ca1}. Our goal is to prove \eqref{ca1} with  $B_{k,i}$ replaced by $B_{k,i}'<B_{k,i}$. Note that
\eq{\label{sec6.4me}&\sum_{|I|=i}E_u(\partial^kZ^Iv)(t)+\int_t^{2T} \sum_{|I|=i}E_q(\partial^kZ^Iv)(\tau)\ d\tau\\
&\leq  \int_t^{2T}C_N\eps^3 B_{k,i}\tau^{-2+C_{k,i}\eps}+CB_{k,i}\eps^2\tau^{-33/16+C_{k,i}\eps}+CB_{k,i}^{1/2}\eps^2 t^{-2+(C_{k,i}/2+C)\eps}\\&\qquad\qquad+C_N\eps^3\tau^{-1+C\eps}(B_{k,i}^{1/2}B^{1/2}_{k-1,i}\tau^{-1+(C_{k,i}+C_{k-1,i})\eps/2}+B_{k,i}^{1/2}B^{1/2}_{k+1,i-1}\tau^{-1+(C_{k,i}+C_{k-1,i+1})\eps/2})\ d\tau\\
&\quad+CB_{k,i}^{1/2}\eps^{3/2}(\ln\eps^{-1})^{-1/2} t^{-1/2+C_{k,i}\eps/2}\kh{\int_t^{2T} B_{k,i}\eps^2\tau^{-3+(C+C_{k,i})\eps}\ d\tau}^{1/2}\\
&\leq C_N\eps^3t^{-1+C_{k,i}\eps}+C_N\eps^{2+1/16}t^{-1+C_{k,i}\eps}+CB^{1/2}_{k,i}\eps^2 t^{-1+(C_{k,i}/2+C)\eps}\\&\quad+C_N\eps^3 t^{-1+(C_{k,i}+C_{k-1,i})\eps/2+C\eps}+C_N\eps^3 t^{-1+(C_{k,i}+C_{k+1,i-1})\eps/2+C\eps}+C_N\eps^{5/2}(\ln\eps^{-1})^{-1/2} t^{-3/2+C_N\eps}.}
For simplicity, we set $B_{\cdot,-1}=B_{-1,\cdot}=C_{\cdot,-1}=C_{-1,\cdot}=0$.  In the last step here, we use the estimate $t^{-1/16}\leq \eps^{1/16}$ whenever $t\geq T_\eps$.
We then choose the constants in the order given in \eqref{sec6ord}. For each $k,i$ with $k+i\leq N$, suppose we have chosen $C_{k',i'}$ and $B_{k',i'}$ for all $i'<i$ and $k'\leq N-i'$, or for all $i'=i$ and $k'<k$. We now choose a  sufficiently  large constant $B_{k,i}$ such that $CB_{k,i}^{1/2}\leq B_{k,i}/10$ and $B_{k,i-1}+B_{k-1,i}\leq B_{k,i}/10$. We also choose a sufficiently large constant $C_{k,i}$ such that $C_{k,i}/2+C\leq C_{k,i}$, $(C_{k,i}+C_{k-1,i})/2+C\leq C_{k,i}$ and $(C_{k,i}+C_{k-1,i+1})/2+C\leq C_{k,i}$; in other words, the exponents of all terms on the right side of \eqref{sec6.4me} except the last term are no smaller than $-1+C_{k,i}\eps$.
It follows from \eqref{sec6.4me} that
\fm{E_{k,i}(t)&\leq E_{k,i-1}(t)+E_{k-1,i}(t)+\sum_{|I|=i}E_u(\partial^kZ^Iv)(t)+\int_t^{2T} \sum_{|I|=i}E_q(\partial^kZ^Iv)(\tau)\ d\tau\\
&\leq C_N\eps^{2+1/16}t^{-1+C_{k,i}\eps}+\frac{1}{5}B_{k,i}\eps^2 t^{-1+ C_{k,i}\eps}+C_N\eps^{5/2}(\ln\eps^{-1})^{-1/2} t^{-3/2+C_N\eps}.}
After we have chosen all the $B_{*,*}$ and $C_{*,*}$, we can also choose the value of $C_N$. We can take $C_N$ to be a polynomial of all the constants appearing before it in \eqref{sec6ord}.
Finally we choose  $\eps\ll_N1$ such that 
\fm{C_N\eps^{1/16}\leq 1/10,\qquad C_N\eps^{1/2}(\ln\eps^{-1})^{-1/2}\leq 1,\qquad C_N\eps\leq  1/4.}
We  conclude that $E_{k,i}(t)\leq \frac{1}{2}B_{k,i}\eps^2t^{-1+C_{k,i}\eps}$. The continuity argument thus ends.

\label{sec6.5}\subsection{Existence of the solutions for $0\leq t\leq T_\eps$}
So far we have proved that there exists a $C^{N+1}$ function $u=(u^{(I)})$ defined for $t\geq T_\eps$, such that $u-u_{app}$ is a solution to \eqref{eqn} for $t\geq T_\eps$. We now extend $u$ to the whole upper half space $[0,\infty)\times\R^3$ by solving  a backward Cauchy problem \eqref{qwe} with initial data at  $T_\eps$ (or at any time in $[T_\eps,T]$).

It seems that the existence of such a  $u$ for $t\geq 0$ follows from the almost global existence result for \eqref{qwe}. However, we remark that our $u$ depends on the time $T$ chosen at the beginning of the current section, so to obtain the almost global existence result, we also need to let $\eps\ll1$ depend on $T$. That is unsatisfactory. Instead we hope to prove the existence result for all $\eps\ll1$, regardless of what the value of $T$ is.

Our main result in this subsection is the following proposition.
\prop{\label{prop6.9} Fix an integer $N\geq 6$.  Let $u$ be a $C^{N+1}$ solution to \eqref{qwe} for $T_\eps\leq t\leq T_\eps+1$. Suppose that $u(t)$ has compact support\footnote{Later we will make use of Poincar$\acute{\rm e}$'s lemmas in the proof. In Section \ref{sec5.3}, we always assume that $\phi(t)$ has compact support for each $t$, so we also assume our $u(t)$ has compact support in this proposition.} for each $T_\eps\leq t\leq T_\eps+1$, that \eq{\label{prop6.9a0}\sum_{|I|\leq N}|Z^Iu(T_\eps,0)|\leq B\eps T_\eps^{-(\gamma_1+1)/2},} and that \eq{\label{prop6.9a}\sum_{|I|\leq N}\norm{\partial Z^Iu(t)}_{L^2(w_0)}&\leq B\eps,\qquad T_\eps\leq t\leq T_\eps+1.}
Then, for $\eps\ll_B1$, $u$ can be extended to a $C^{N+1}$ solution to \eqref{qwe} for all $0\leq t\leq T_\eps+1$, such that  \eqref{prop6.9a} (with $B$ replaced by a possibly larger constant $B'$) holds for all $0\leq t\leq T_\eps+1$. We emphasize that the choice of $\eps\ll_B1$ does not depend on  $\supp u$; that is, we can choose $\eps_0>0$ depending on $B$ but not on $\supp u$, so that all the conclusions hold for each $\eps\in(0,\eps_0)$.}\rm

\rmk{\rm Let $u$ be a function constructed in the previous subsections. Let us show that the assumptions of Proposition \ref{prop6.9} hold for this $u$.

We start with the support condition on $u$. This is obvious because $\supp u_{app}\subset\{t\geq T_\eps,\ r/t\in[1-c,1+c]\}$ and the result (c) at the beginning of Section \ref{sec6.1}.

Next we show \eqref{prop6.9a0}. This is a direct corollary of \eqref{ptb2}. Since we are using a continuity argument in Section \ref{sec6.2}, we obtain \eqref{ptb2} for each $t\geq T_\eps$ once  this continuity argument is closed. We also notice that $T_\eps^{C\eps}=\eps^{-C\eps}\leq 2$ for $\eps\ll1$.

Finally we show \eqref{prop6.9a}. Let $w_0=w_0(\rho)$ and $\norm{f}_{L^2(w_0)}$ be defined in Proposition  \ref{prop6}. Note that by \eqref{ca1} and Proposition \ref{mainprop4}, for $T_\eps\leq t\leq T_\eps+1$ we have
\fm{ &\sum_{|I|\leq N}\norm{\partial Z^Iu(t)}_{L^2(w_0)}\\&\leq \sum_{|I|\leq N}\norm{\partial Z^Iv(t)}_{L^2(w_0)}+\sum_{|I|\leq N}\norm{\partial Z^Iu_{app}(t)}_{L^2(w_0)}\\
&\leq CE_{0,N}(t)^{1/2}+C\norm{\eps\lra{r-t}^{-1}w_0(r-t)^{1/2}t^{-1+C\eps}(1_{r\geq t}+\lra{r-t}^{1-\gamma_-}1_{r<t})1_{|r-t|\leq ct}}_{L^2(\R^3)}\\
&\leq CE_{0,N}(t)^{1/2}+C\norm{\eps\lra{r-t}^{-1}t^{-1+C\eps}1_{|r-t|\leq ct}}_{L^2(\R^3)}\\
&\leq CC_N\eps t^{-1/2+C_N\eps}+C\eps t^{-1+C\eps}\kh{\int_{(1-c)t}^{(1+c)t}\rho^2\lra{\rho-t}^{-2}\ d\rho}^{1/2}\\
&\leq CC_N\eps t^{-1/2+C_N\eps}+C\eps t^{C\eps}\kh{\int_{0}^{ct}\lra{\rho}^{-2}\ d\rho}^{1/2}\\
&\leq CC_N\eps t^{-1/2+C_N\eps}+C\eps t^{C\eps}\leq CC_N\eps T_\eps^{-1/4}+C\eps T_\eps^{C\eps}\leq C\eps.}
To obtain the last two estimates, we choose  $\eps \ll1$.  So we obtain \eqref{prop6.9a}. 

In conclusion, we are allowed to apply Proposition \ref{prop6.9}. Note that both $\eps\ll 1$ and the constant $C$ in the last estimate are independent of $T$, which is what we want. We thus end the proof of Proposition \ref{prop6}.}

\rm

\bigskip

We use a continuity argument to show the existence of  $u$ in Proposition \ref{prop6.9}. Suppose that $u$ exists for all $t\geq T_1$ for some $0< T_1\leq T_\eps$ such that 
\eq{\label{ca3} \sum_{|I|\leq N}\norm{\partial Z^Iu(t)}_{L^2(w_0)}\leq B_N\eps,\qquad T_1\leq t\leq T_\eps.}
Here $B_N$ is a constant depending on $N$. In particular, we assume that $B_N\geq 2B$ where $B$ is the constant in \eqref{prop6.9a}. By the weighted Klainerman-Sobolev inequality \eqref{corklc}, we have
\eq{\label{ca4}\sum_{|I|\leq N-2}|\partial Z^Iu(t,x)|\leq CB_N\eps\lra{r+t}^{-1}\lra{r-t}^{-1/2}w_0(r-t)^{-1/2},\qquad T_1\leq t\leq T_\eps.}
It then follows from \eqref{prop6.9a0} that
\eq{\label{ca5}&\sum_{|I|\leq N-2}|Z^Iu(t,x)|\\&\leq\sum_{|I|\leq N-2}|Z^Iu(T_\eps,0)|+\int_t^{T_\eps} \sum_{|I|\leq N-2}|\partial_tZ^Iu(\tau,0)|\ d\tau+\int_0^r\sum_{|I|\leq N-2}|\partial_rZ^Iu(t,\rho\omega)|\ d\rho\\
&\leq CB\eps T_\eps^{-(\gamma_1+1)/2}+\int_t^{T_\eps} CB_N\eps\lra{\tau}^{-(3+\gamma_1)/2}\ d\tau \\&\quad+\int_0^r CB_N\eps \lra{\rho+t}^{-1}(\lra{\rho-t}^{-1/2}1_{\rho\geq t}+\lra{\rho-t}^{-(\gamma_1+1)/2}1_{\rho<t})\ d\rho\\
&\leq  CB_N\eps\lra{t}^{-1}\lra{r-t}^{(1-\gamma_1)/2}1_{r<t}+CB_N\eps\lra{t}^{-1}\lra{r-t}^{1/2}1_{ t\leq r\leq 2t}+CB_N\eps\lra{t}^{-1/2}1_{r>2t}}
whenever $T_1\leq t\leq T_\eps$. Thus, we can choose $\eps\ll 1$ so that the right hand sides of \eqref{ca4} and \eqref{ca5} are no more than $1/2$. This allows us to apply Lemma \ref{l24} freely. Moreover, since the assumption \eqref{sec5a3} is satisfied,  we can freely use the energy estimate \eqref{energymain2} in Section \ref{sec5.1}. It also follows from \eqref{ca5} that
\eq{\label{ca6}&\lra{r-t}^{-1}\sum_{|I|\leq N-2}|Z^Iu(t,x)|\\&\leq CB_N\eps\lra{t}^{-1}\lra{r-t}^{-(1+\gamma_1)/2}1_{r<t}+CB_N\eps\lra{t}^{-1}\lra{r-t}^{-1/2}1_{ r\geq t},\qquad T_1\leq t\leq T_\eps.}

By applying $Z^I$ with $|I|\leq N$ to \eqref{qwe}, we obtain
\eq{\label{eqn3}g^{\alpha\beta}(u,\partial u)\partial_\alpha\partial_\beta Z^Iu&=[\Box,Z^I]u+[g^{\alpha\beta}(u,\partial u)-m^{\alpha\beta},Z^I]\partial_\alpha\partial_\beta u\\&\quad+(g^{\alpha\beta}(u,\partial u)-m^{\alpha\beta})[\partial_\alpha\partial_\beta, Z^I]u+Z^If(u,\partial u).}
Since $f$ contains no terms of the form $u\cdot u\cdot u$, it follows from Lemma \ref{l24}  that
\fm{&|g^{\alpha\beta}(u,\partial u)\partial_\alpha\partial_\beta Z^Iu|\\
&\lesssim \sum_{|I_1|+|I_2|\leq |I|\atop |I_2|<|I|}|Z^{I_1}(g^{\alpha\beta}(u,\partial u)-m^{\alpha\beta})\partial^2Z^{I_2}u|+\sum_{|I'|\leq |I|}|Z^{I'}f(u,\partial u)|\\
&\lesssim \sum_{|I_1|+|I_2|\leq |I|\atop |I_2|<|I|}(|Z^{I_1}u\partial^2Z^{I_2}u|+|\partial Z^{I_1}u\partial^2Z^{I_2}u|)+\sum_{|I_1|+|I_2|\leq |I|}|\partial Z^{I_1}u\partial Z^{I_2}u|\\&\quad+\sum_{|I_1|+|I_2|+|I_3|\leq |I|\atop k_1,k_2,k_3\leq 1,\ k_1+k_2+k_3>0}\prod_{l=1}^3|\partial^{k_l}Z^{I_l}u|+\sum_{|I_1|+|I_2|+|I_3|+|I_4|\leq |I|\atop k_1,k_2,k_3,k_4\leq 1}\prod_{l=1}^4|\partial^{k_l}Z^{I_l}u|.}
Also notice that $|Z^{I_1}u\partial^2Z^{I_2}u|\lesssim |\lra{r-t}^{-1}Z^{I_1}uZ\partial Z^{I_2}u|$. In the first two sums, since $\min\{|I_1|,|I_2|+1\}\leq (N+1)/2\leq N-2$, we can apply \eqref{ca4} and \eqref{ca6} to the lower order terms. This gives us an upper bound
\fm{&\sum_{|I'|\leq |I|}CB_N\eps\lra{t}^{-1}\lra{r-t}^{-1/2}w_0(r-t)^{-1/2}(|\lra{r-t}^{-1}Z^{I'}u|+|\partial Z^{I'}u|)\\
&\leq \sum_{|I'|\leq |I|}CB_N\eps\lra{t}^{-1}w_0(r-t)^{-1/2}(|\lra{r-t}^{-1}Z^{I'}u|+|\partial Z^{I'}u|).}
In the third sum, we can assume without loss of generality that $|I_3|$ is the largest among $|I_*|$. Thus,
$|I_1|+1,|I_2|+1\leq N/2+1\leq N-2$. If $k_3=1$, then we can apply \eqref{ca4} and \eqref{ca5} to control $|\partial^{k_1}Z^{I_1}u\partial^{k_2}Z^{I_2}u|$. If $k_3=0$, then we must have $k_1=1$ or $k_2=1$, so we obtain an additional factor $\lra{r-t}^{-1}$ by applying  Lemma \ref{c1l2.1}. This gives us an upper bound 
\fm{&\sum_{|I_3|\leq |I|}CB_N^2\eps^2\lra{t}^{-1}(\lra{t}^{-1}\lra{r-t}^{1-\gamma_1}1_{r<t}+1_{r\geq t})(|\partial Z^{I_3}u|+|\lra{r-t}^{-1}Z^{I_3}u|)\\
&\leq \sum_{|I_3|\leq |I|}CB_N^2\eps^2\lra{t}^{-1}w_0(r-t)^{-1}(|\partial Z^{I_3}u|+|\lra{r-t}^{-1}Z^{I_3}u|)\\
&\leq \sum_{|I'|\leq |I|}\eps\lra{t}^{-1}w_0(r-t)^{-1/2}(|\lra{r-t}^{-1}Z^{I'}u|+|\partial Z^{I'}u|),\qquad \text{whenever }\eps\ll1.}
In the fourth sum, we apply \eqref{ca5}, so we obtain an upper bound
\fm{CB_N\eps(\lra{t}^{-1}\lra{r-t}^{(1-\gamma_1)/2}1_{r<t}+\lra{t}^{-1/2}1_{r\geq t})\sum_{|I'|\leq |I|}|Z^{I'}u|^3.}

Our goal is to estimate the $L^2(w_0)$ norm of \eqref{eqn3}. This is achieved  in the following two lemmas.

\lem{For each $T_1\leq t\leq T_\eps$, we have
\fm{\sum_{|I'|\leq N}\norm{w_0(r-t)^{-1/2}\lra{r-t}^{-1}|Z^{I'}u|}_{L^2(w_0)}\lesssim \sum_{|I'|\leq N}\lra{t}^{\frac{-\gamma_1+6}{4}}\norm{\partial Z^{I'}u}_{L^2(w_0)}+B_N\eps\lra{t}^{\frac{-\gamma_1+6}{4}}.}}
\begin{proof}
We have
\fm{\sum_{|I'|\leq N}\norm{w_0(r-t)^{-1/2}\lra{r-t}^{-1}|Z^{I'}u|}_{L^2(w_0)}=\sum_{|I'|\leq N}\norm{\lra{r-t}^{-1}|Z^{I'}u|}_{L^2(\R^3)}.}

By Lemma \ref{lp1} with $\eta=\gamma_1/2$, we have
\fm{&\norm{\lra{r-t}^{-1}|Z^{I'}u|}_{L^2(\R^3)}\\
&\lesssim \norm{(1_{r\geq t}+\lra{r-t}^{\frac{\gamma_1+2}{4}}1_{r<t})|\partial Z^{I'}u|}_{L^2(\R^3)}+\norm{\lra{r-t}^{\frac{\gamma_1+2}{4}}| Z^{I'}u|1_{r/t\in[1-c,1-c/2]}}_{L^2(\R^3)}\\
&\qquad+\norm{\lra{r-t}^{-1}|Z^{I'}u|1_{r/t<1-c/2}}_{L^2(\R^3)}\\
&\lesssim \norm{(1_{r\geq t}+\lra{r-t}^{\frac{-\gamma_1+2}{4}}1_{r<t})|\partial Z^{I'}u|}_{L^2(w_0)}+\lra{t}^{\frac{\gamma_1+2}{4}}\norm{| Z^{I'}u|1_{r/t\in[1-c,1-c/2]}}_{L^2(\R^3)}\\
&\qquad+\norm{r^{-1}|Z^{I'}u|1_{r/t<1-c/2}}_{L^2(\R^3)}\\
&\lesssim \norm{\partial Z^{I'}u}_{L^2(w_0)}+\lra{t}^{\frac{\gamma_1+2}{4}}\norm{| Z^{I'}u|1_{r/t\in[1-c,1-c/2]}}_{L^2(\R^3)}} for each $|I'|\leq N$
whenever $t\geq \max\{T_1,100(\eta+1)^2/\eta^2\}$. Here $0<c\ll1$ is a sufficiently small constant chosen in Lemma \ref{lp1}.  In the last step, we use Hardy's inequality and $w_0\geq 1$. If $|I'|>0$, we have \fm{\lra{t}^{\frac{\gamma_1+2}{4}}|Z^{I'}u|1_{r/t\in[1-c,1-c/2]}&\lesssim \lra{t}^{\frac{\gamma_1+6}{4}}\sum_{|J|=|I'|-1}|\partial Z^Ju|1_{r/t\in[1-c,1-c/2]}\\&
\lesssim w_0(r-t)^{1/2}\lra{t}^{\frac{-\gamma_1+6}{4}}\sum_{|J|=|I'|-1}|\partial Z^Ju|1_{r/t\in[1-c,1-c/2]}.}
The $L^2(\R^3)$ norm of this formula is controlled by $\lra{t}^{\frac{-\gamma_1+6}{4}}\sum_{|J|=|I'|-1}\norm{\partial Z^{J}u}_{L^2(w_0)}$. Recall that $\gamma_1<4$ by \eqref{sec5wei4}, so $\lra{t}^{\frac{-\gamma_1+6}{4}}\geq 1$. If $|I'|=0$, then we apply \eqref{ca5} to obtain
\fm{&\lra{t}^{\frac{\gamma_1+2}{4}}\norm{| u|1_{r/t\in[1-c,1-c/2]}}_{L^2(\R^3)}\\
&\lesssim B_N\eps\lra{t}^{\frac{\gamma_1-2}{4}}\norm{\lra{r-t}^{(1-\gamma_1)/2}1_{r/t\in[1-c,1-c/2]}}_{L^2(\R^3)}\\
&\lesssim B_N\eps\lra{t}^{\frac{\gamma_1-2}{4}}\cdot \lra{t}^{(1-\gamma_1)/2}\cdot \lra{t}^{3/2}\lesssim B_N\eps \lra{t}^{3/2-\gamma_1/4}.}
This finishes the proof when  $t\geq \max\{T_1,100(\eta+1)^2/\eta^2\}$.

If $T_1\leq t<100(\eta+1)^2/\eta^2$ (such $t$ may not exist), we apply Lemma \ref{lp1st} to get 
\fm{\norm{\lra{r-t}^{-1}|Z^{I'}u|}_{L^2(\R^3)}\lesssim\norm{\partial Z^{I'}u}_{L^2(\R^3)}\lesssim\norm{\partial Z^{I'}u}_{L^2(w_0)}.}
This finishes our proof immediately.
\end{proof}

\lem{For $T_1\leq t\leq T_\eps$, we have\fm{&\norm{CB_N\eps(\lra{t}^{-1}\lra{r-t}^{(1-\gamma_1)/2}1_{r<t}+\lra{t}^{-1/2}1_{r\geq t})\sum_{|I'|\leq |N}|Z^{I'}u|^3}_{L^2(w_0)}\\
&\lesssim B_N^3\eps^3\lra{t}^{-1/2}\sum_{|I'|\leq N}\norm{\partial Z^{I'}u}_{L^2(w_0)}. }}
\begin{proof}
We have
\fm{&\norm{CB_N\eps(\lra{t}^{-1}\lra{r-t}^{(1-\gamma_1)/2}1_{r<t}+\lra{t}^{-1/2}1_{r\geq t})\sum_{|I'|\leq |N}|Z^{I'}u|^3}_{L^2(w_0)}\\
&=CB_N\eps\norm{(\lra{t}^{-1}\lra{r-t}^{1/2}1_{r<t}+\lra{t}^{-1/2}1_{r\geq t})\sum_{|I'|\leq N}|Z^{I'}u|^3}_{L^2(\R^3)}\\
&\lesssim B_N\eps\lra{t}^{-1/2}\sum_{|I'|\leq N}\norm{Z^{I'}u}_{L^6(\R^3)}^{3}\lesssim B_N\eps\lra{t}^{-1/2}\sum_{|I'|\leq N}\norm{\partial Z^{I'}u}_{L^2(\R^3)}^{3} \\
&\lesssim B_N^3\eps^3\lra{t}^{-1/2}\sum_{|I'|\leq N}\norm{\partial Z^{I'}u}_{L^2(w_0)}. }
In this proof, we use the homogeneous Sobolev embedding $L^6(\R^3)\hookrightarrow \dot{H}^1(\R^3)$ and the assumption \eqref{ca3}.
\end{proof}

\rm

\bigskip

We conclude that for each $T_1\leq t\leq T_\eps$,
\eq{&
\sum_{|I|\leq N}\norm{g^{\alpha\beta}(u,\partial u)\partial_\alpha\partial_\beta Z^Iu(t)}_{L^2(w_0)}\\
&\lesssim B_N\eps\lra{t}^{\frac{-\gamma_1+2}{4}}\sum_{|I|\leq N}\norm{\partial Z^{I}u}_{L^2(w_0)}+B_N^2\eps^2\lra{t}^{\frac{-\gamma_1+2}{4}}+B_N^3\eps^3\lra{t}^{-1/2}\sum_{|I'|\leq N}\norm{\partial Z^{I'}u}_{L^2(w_0)}\\
&\lesssim B_N\eps\lra{t}^{\frac{-\gamma_1+2}{4}}\sum_{|I|\leq N}\norm{\partial Z^{I}u}_{L^2(w_0)}+B_N^2\eps^2\lra{t}^{\frac{-\gamma_1+2}{4}}.}
In the last estimate, we use $\gamma_1<4$ and $\eps\ll1$.

Now we can apply the energy estimate \eqref{energymain2} in Section \ref{sec5.1}. Note that the $m=m_{\gamma_1}(q)$ in that section is the same as $w_0$ in this section. We thus have
\eq{&\sum_{|I|\leq N}\norm{\partial Z^Iu(t)}_{L^2(w_0)}^2\lesssim\sum_{|I|\leq N}\norm{\partial Z^Iu(T_\eps)}_{L^2(w_0)}^2\\&\qquad+\sum_{|I|\leq N}\int_t^{T_\eps}(C_{B_N}\eps\lra{\tau}^{-1/2}+B_N\eps\lra{\tau}^{\frac{-\gamma_1+2}{4}})\norm{\partial Z^{I}u(\tau)}^2_{L^2(w_0)}+B_N^2\eps^2\lra{\tau}^{\frac{-\gamma_1+2}{4}}\norm{\partial Z^{I}u(\tau)}_{L^2(w_0)}\ d\tau\\
&\lesssim B^2\eps^2+\int_t^{T_\eps}B^2\eps^2(C_{B_N}\eps\lra{\tau}^{-1/2}+B_N\eps\lra{\tau}^{\frac{-\gamma_1+2}{4}})+BB_N^2\eps^3\lra{\tau}^{\frac{-\gamma_1+2}{4}}\ d\tau\\
&\lesssim B^2\eps^2+B^2\eps^2(C_{B_N}\eps\lra{T_\eps}^{1/2}+B_N\eps\lra{T_\eps}^{\frac{-\gamma_1+6}{4}})+BB_N^2\eps^3\lra{T_\eps}^{\frac{-\gamma_1+6}{4}}\\
&\lesssim B^2\eps^2+B^2(C_{B_N}\eps^{5/2}+B_N\eps^{\frac{\gamma_1+6}{4}})+BB_N^2\eps^{\frac{\gamma_1+6}{4}}\lesssim B^2\eps^2.}
for each $T_1\leq t\leq T_\eps$. In the second estimate, we use \eqref{ca3} and \eqref{prop6.9a}. In the last two estimates, we notice that $T_\eps=1/\eps$ and $2<\gamma_1<4$.  All the implicit constants are all known before we choose the value of $B_N$, so by choosing $B_N\gg_B1$ we conclude that
\fm{\sup_{t\in[T_1,T_\eps]}\sum_{|I|\leq N}\norm{\partial Z^Iu(t)}_{L^2(w_0)}\leq CB\eps\leq \frac{1}{2}B_N\eps.}
This ends the continuity argument so we obtain Proposition \ref{prop6.9}.

\section{Limit as $T\to\infty$}\label{sec7}

Our goal for this section is to prove the following proposition. 

\prop{\label{prop7} Fix $N\geq 6$. Then for the same $\eps_{N}$ in Proposition \ref{prop6} and for  $0<\eps<\eps_{N}$, there is a solution $u$ to \eqref{qwe}  in $C^{N-4}$ for all $t\geq 0$, such that for all  $|I|\leq N-5$
\eq{\label{p7f1}
\norm{w_0(|x|-t)^{1/2}\partial Z^I(u-u_{app})(t)}_{L^2(\R^3)}&\lesssim_I\eps t^{-1/2+C_I\eps},\qquad t\geq T_\eps;\\
\norm{w_0(|x|-t)^{1/2}\partial Z^Iu(t)}_{L^2(\R^3)}&\lesssim_I\eps ,\qquad 0\leq t\leq T_\eps.
}
Here $w_0$ is defined in Proposition \ref{prop6}.

Besides, for all $|I|\leq N-5$ and $t\geq T_\eps$,
\begin{equation}\label{p7f2}
|\partial Z^I(u-u_{app})(t,x)|\lesssim_I \eps t^{-1/2+C_I\eps}\lra{r+t}^{-1}\lra{t-r}^{-1/2}w_0(r-t)^{-1/2},
\end{equation}
\begin{equation}\label{p7f3}
|Z^I(u-u_{app})(t,x)|\lesssim_I\eps t^{-1/2+C_I\eps/2} (t^{-1}\lra{r-t}^{-(\gamma_1-1)/2}1_{r< t}+t^{-1}\lra{r-t}^{1/2}1_{t\leq r\leq 2t}+r^{-1/2}1_{r>2t}).
\end{equation}
}
\rmk{\rm It should be pointed out that the value of ``$N$'' in the main theorem is equal to $N-4$ for the $N$ in this proposition.}

\rmk{\label{rmk712}\rm The solution $u$ is unique in the following sense. Fix $N\geq 11$. Suppose that $u_1$ and $u_2$ are two $C^{N-4}$ solutions to \eqref{qwe} for $t\geq 0$ corresponding to the same $u_{app}$. Suppose that \eqref{p7f1}, \eqref{p7f2}, and \eqref{p7f3} hold for both $u=u_1$ and $u=u_2$. Then, if $\eps\ll_N1$ is sufficiently small, we have $u_1\equiv u_2$ everywhere. We will prove this uniqueness result at the end of this section.

This uniqueness statement above does imply the one in Remark \ref{rmkthmunique}. This is because the lemma in the next remark (Remark \ref{rmk7.1.3Uuapp}), we show that the estimates in the main theorem imply \eqref{p7f1} -- \eqref{p7f3}.}

\rmk{\rm\label{rmk7.1.3Uuapp} We make the following important claims about $u_{app}$ and $\eps r^{-1}U$. In \eqref{p7f1} and \eqref{p7f2}, we can replace $\partial Z^I(u-u_{app})$ on the left sides with \fm{\partial Z^I(u-\eps r^{-1}U)\cdot 1_{r/t\in[1/2,3/2]}+\partial Z^Iu\cdot 1_{r/t\notin[1/2,3/2]}.} Similarly, in \eqref{p7f3} we can replace $Z^I(u-u_{app})$ with\fm{Z^I(u-\eps r^{-1}U)\cdot 1_{r/t\in[1/2,3/2]}+Z^Iu\cdot 1_{r/t\notin[1/2,3/2]}.}
Note that $u_{app}=\eps r^{-1}U$ when $|r/t-1|< c/2$ and that $u_{app}=0$ when $|r/t-1|>1/2$. Thus, to prove the claims, we need to estimate $\partial Z^I(u_{app}-\eps r^{-1}U)\cdot 1_{|r/t-1|\in[c/2,1/2]}$ and $Z^I(u_{app}-\eps r^{-1}U)\cdot 1_{|r/t-1|\in[c/2,1/2]}$.

\lem{For all $(t,x)$ with $t\geq T_\eps$ and $||x|/t-1|\in[c/2,1/2]$, for each $I$ we have
\fm{|Z^I(\eps r^{-1}U)|+|Z^Iu_{app}|+t(|\partial Z^I(\eps r^{-1}U)|+|\partial Z^Iu_{app}|)&\lesssim_{c,I} \eps t^{-\gamma_-+C\eps}1_{r\leq t}+\eps t^{-1+C_I\eps}1_{r>t},}
\fm{\norm{w_0(|x|-t)^{1/2}|\partial Z^I(u_{app}-\eps r^{-1}U)|1_{|r/t-1|\in[c/2,1/2]}}_{L^2(\R^3)}\lesssim \eps t^{-1/2+C\eps}.}}
\begin{proof}
By part (b) of Lemma \ref{lem4.5}, we have $\eps r^{-1}U\in \eps S^{-1,1-\gamma_-}_-\cap \eps S^{-1,0}_+$. In other words, we have 
\fm{|Z^I(\eps r^{-1}U)|\lesssim_I \eps t^{-1+C_I\eps}(\lra{r-t}^{1-\gamma_-}1_{r\leq t}+1_{r>t}),\qquad \forall t\geq 1/\eps, ||x|-t|\leq t/2,\ \forall I.}
Moreover, by Proposition \ref{mainprop4}, we obtain the same estimate with $\eps r^{-1}U$ replaced by $u_{app}$. By Lemma \ref{c1l2.1}, for all $t\geq 1/\eps$ and $||x|/t-1|\in[c/2,1/2]$, we have
\fm{|\partial Z^Iu_{app}|+|\partial Z^I(\eps r^{-1}U)| \lesssim_{c,I} \eps t^{-1-\gamma_-+C\eps}1_{r\leq t}+\eps t^{-2+C_I\eps}1_{r>t},\qquad \forall I.}
It follows that 
\fm{&\norm{w_0(|x|-t)^{1/2}|\partial Z^I(u_{app}-\eps r^{-1}U)|1_{|r/t-1|\in[c/2,1/2]}}_{L^2(\R^3)}\\&\lesssim t^{\gamma_1/2}\cdot \eps t^{-1-\gamma_-+C
\eps}\cdot t^{3/2}+\eps t^{-2+C
\eps}\cdot t^{3/2}\lesssim\eps t^{-1/2+C\eps}.}
Here we use $\gamma_1/2<\gamma_--1$.
\end{proof}\rm

Whenever $|r/t-1|\in[c/2,1/2]$, we have $r\sim t\sim\lra{r-t}$ and thus
\fm{\eps t^{-1/2+C\eps}\lra{r+t}^{-1}\lra{r-t}^{-1/2}w_0(r-t)^{-1/2}&\sim\eps t^{-2+C\eps}(t^{-\gamma_1/2}1_{r\leq t}+1_{r> t}),\\ \eps t^{-1/2+C\eps/2} (t^{-1}\lra{r-t}^{-(\gamma_1-1)/2}1_{r< t}+t^{-1}\lra{r-t}^{1/2}1_{t\leq r\leq 2t})&\sim\eps t^{-1+C\eps/2} (t^{-\gamma_1/2}1_{r< t}+1_{t\leq r}).}
Since $\gamma_1/2<\gamma_--1$, we have $\eps t^{-1-\gamma_1/2+C\eps}\geq \eps t^{-\gamma_-}$. In other words, the right sides of the estimates in this lemma can be absorbed by the right sides of \eqref{p7f1} -- \eqref{p7f3}, so the claims are proved.

In summary, we obtain the first two estimates in part (ii) of Theorem \ref{mthm}. It also explains why the $u_{app}$ depends on a parameter $c$, while we do not include such a parameter in Theorem~\ref{mthm}. To obtain the last estimate, we notice that whenever $t\geq 1/\eps$ and $r/t\in[1/2,3/2]$,  \fm{(\partial_t-\partial_r)(\eps r^{-1}U)&=-\eps r^{-2}U+\eps r^{-1}(\eps t^{-1}U_s+\mu U_q)\\
&=\eps r^{-1}\mu U_q+O(\eps t^{-2+C\eps}\lra{\max\{0,-q\}}^{1-\gamma_-})\\
&=\eps r^{-1}\mu U_q+O(\eps t^{-2+C\eps}(\lra{r-t}^{1-\gamma_-}1_{r<t}+1_{r\geq t})).}
}

\rm\bigskip

This section is structured as follows. In Section \ref{sec7.1}, we prove Proposition \ref{prop7} assuming that a key lemma, Lemma \ref{lem7.2}, holds. Suppose that we start with two different times $T_1$ and $T_2$ in Section \ref{sec6}, obtain two corresponding solutions $u_1$ and $u_2$ by applying Proposition \ref{prop6}. The key lemma above gives us an estimate for the energy of $u_1-u_2$. Then, in Section \ref{sec7.2}, we prove  Lemma \ref{lem7.2}. The proof there is similar to that of Proposition \ref{prop6}, and our main tool is the energy estimates proved in Section \ref{sec5}. Finally, in Section \ref{s7uni} and Section \ref{sec7.4pfprop1.1}, we prove the results in Remark \ref{rmkthmunique} and Remark \ref{mthmrmk3}, respectively.

From now on, we will use the constant $C$ to denote all the constants in the previous sections including $C_{k,i},B_{k,i},C_N$, etc. As usual, the constant $C$ cannot depend on $\eps,T_1,T_2,T$, etc.

\subsection{A key lemma}\label{sec7.1}

Fix $N\geq 6$ and $T_2>T_1\gg 1$.  By Proposition \ref{prop6}, for each $0<\eps<\eps_{N}$ (note that $\eps_N$ is independent of $T_1$ and $T_2$), we  get two corresponding solutions $u_1=u^{T_1}$ and $u_2=u^{T_2}$, where $u^{T}$ denotes the solution constructed in Proposition \ref{prop6} corresponding to $T$. Both of them are defined  for all $t\geq 0$. Our goal now is to prove that $u_1-u_2$ tends to $0$ in some Banach space as $T_2>T_1\to\infty$. 

The key lemma in this subsection is the following.

\lem{\label{lem7.2}Fix a constant $\lambda_0\in(0,1/2)$. Then, as long as $\eps\ll_{\lambda_0}1$, for each $t\geq 0$, we have
\eq{\label{lem7.2c1}\sum_{|I|\leq N-3}\norm{w_0(|x|-t)^{1/2}\partial Z^I(u_1-u_2)(t)}_{L^2(\R^3)}\leq C_{\lambda_0}\eps T_1^{-1/2+\lambda_0}(t^{-\lambda_0+C\eps}1_{t\geq 1/\eps}+1_{t<1/\eps}).}}
\rm

\bigskip

We  postpone the proof of Lemma \ref{lem7.2} to the next subsection. Here we explain how to derive Proposition \ref{prop7} from this key lemma. By applying the weighted Klainerman-Sobolev inequality \eqref{corklc}, we have\eq{\label{lem7.2c2}&\sum_{|I|\leq N-5}|\partial Z^I(u_1-u_2)(t,x)|\\
&\lesssim_{\lambda_0}\eps T_1^{-1/2+\lambda_0}(t^{-\lambda_0+C\eps}1_{t\geq 1/\eps}+1_{t<1/\eps})\lra{r+t}^{-1}\lra{r-t}^{-1/2}w_0(r-t)^{-1/2}}
for all $t\geq 0$ and $x\in\R^3$. Recall from Section \ref{sec6} that $u_1=u_2=u_{app}$ whenever $t\geq T_2$, so by Lemma \ref{lemptb2} we have
\eq{\label{lem7.2c3}&\sum_{|I|\leq N-5}|Z^I(u_1-u_2)(t,x)|\\
&\lesssim_{\lambda_0} \eps T_1^{-1/2+\lambda_0}(t^{-\lambda_0+C\eps}1_{t\geq 1/\eps}+1_{t<1/\eps})\\&\qquad\qquad\cdot(\lra{t}^{-1}\lra{r-t}^{-(\gamma_1-1)/2}1_{r<t}+\lra{t}^{-1}\lra{r-t}^{1/2}1_{t\leq r\leq 2t}+\lra{r}^{-1/2}1_{r>2t})}
for all $t\geq 0$ and $x\in\R^3$. That is,
\fm{\sup_{(t,x)\in[0,\infty)\times\R^3} \sum_{|I|\leq N-5}(|Z^I(u^{T_1}-u^{T_2})(t,x)|+|\partial Z^I(u^{T_1}-u^{T_2})(t,x)|)\leq C_{\lambda_0}\eps T_1^{-1/2+\lambda_0}\to 0}
as $T_2>T_1\to\infty$. Then, there is $u^\infty\in C^{N-4}(\{(t,x):\ t\geq 0\})$, such that $\partial Z^Iu^T\to \partial Z^Iu^\infty$ and $Z^Iu^T \to Z^Iu^\infty$ uniformly for all $(t,x)\in[0,\infty)\times\R^3$ as $T\to\infty$, for each $|I|\leq N-5$. Now \eqref{p7f2} and \eqref{p7f3} follows from taking limit  in \eqref{ptb1} and \eqref{ptb2}. By Fatou's lemma and Proposition \ref{prop6},  we have
\eq{\label{vinfty}&\sum_{|I|\leq N-5}\norm{w_0(|x|-t)^{1/2}\partial Z^I(u^\infty-u_{app})(t)}_{L^2(\R^3)}\\&\leq\sum_{|I|\leq N-5} \liminf_{T\to\infty}\norm{w_0(|x|-t)^{1/2}\partial Z^I(u^T-u_{app})(t)}_{L^2(\R^3)}\lesssim\eps t^{-1/2+C\eps},\qquad t\geq T_\eps;}
\eq{\label{vinfty2}&\sum_{|I|\leq N-5}\norm{w_0(|x|-t)^{1/2}\partial Z^Iu^\infty(t)}_{L^2(\R^3)}\\&\leq\sum_{|I|\leq N-5} \liminf_{T\to\infty}\norm{w_0(|x|-t)^{1/2}\partial Z^Iu^T(t)}_{L^2(\R^3)}\lesssim\eps ,\qquad 0\leq t\leq T_\eps.}
Meanwhile, if $N\geq 6$, then by taking $T\to\infty$ in \fm{&g^{\alpha\beta}(u^T,\partial u^T)\partial_\alpha\partial_\beta (u^T-u_{app})^{(I)}\\&=-[ g^{\alpha\beta}(u^T,\partial u^T)-g^{\alpha\beta}(u_{app},\partial u_{app})]\partial_\alpha\partial_\beta u_{app}^{(I)}+[f^{(I)}(u^T,\partial u^T)-f^{(I)}(u_{app},\partial u_{app})]\\
&\quad-\chi(t/T)[g^{\alpha\beta}(u_{app},\partial (u_{app}))\partial_\alpha\partial_\beta u_{app}^{(I)}-f^{(I)}(u_{app},\partial u_{app})],} we conclude that $u^\infty$ is a solution to \eqref{qwe} for $t\geq T_\eps$. Here we note that $\lim_{T\to\infty}\chi(t/T)=1$ for all $t\geq T_\eps$. Also notice that $u^T$ solves \eqref{qwe} whenever $t\in[0,T_\eps]$ for all $T$, so by taking $T\to\infty$, we conclude that $u^\infty$ solves \eqref{qwe} whenever $t\in[0,T_\eps]$. We thus conclude the proof of Proposition \ref{prop7}.

\subsection{Proof of Lemma \ref{lem7.2}} \label{sec7.2}

In this subsection, we seek to prove Lemma \ref{lem7.2}. The proof is in fact very similar to that of Proposition \ref{prop6}. 

\subsubsection{The case $t\geq T_1$}
Set $v_j:=u^{T_j}-u_{app}$ for $j=1,2$. Since $u_1-u_2=v_1-v_2$, by \eqref{prop6c1} we have
\fm{\sum_{|I|\leq N-5}\norm{\partial Z^I(u_1-u_2)(t)}_{L^2(w_0)}\lesssim \sum_{|I|\leq N-5}[\norm{\partial Z^Iv_1(t)}_{L^2(w_0)}+\norm{\partial Z^Iv_2(t)}_{L^2(w_0)}]\lesssim \eps t^{-1/2+C\eps}.}
For all $t\geq T_1$, we have $\eps t^{-1+C\eps}\leq \eps T_1^{-1/2+\lambda_0} t^{-\lambda_0+C\eps}$, so we obtain \eqref{lem7.2c1} for all $t\geq T_1$.

\subsubsection{The case $T_\eps\leq t\leq T_1$}
Suppose that $t\in[T_\eps,T_1]$  with $T_\eps=1/\eps$ and $\eps\ll1$. Recall that $\eps$ is independent of the choice of $T$. Now we have  $\chi(t/T_1)=\chi(t/T_2)=1$, so  for each $k+|I|\leq N$, 
\begin{equation}\label{teqn}\begin{aligned}&g^{\alpha\beta}(u_1,\partial u_1)\partial_\alpha\partial_\beta \partial^kZ^I\wt{v}\\
&=\partial^k[\Box,Z^I]\wt{v}-\partial^kZ^I\left([ g^{\alpha\beta}(u_2,\partial u_2)-g^{\alpha\beta}(u_1,\partial u_1)]\partial_\alpha\partial_\beta u_{app}\right)\\&\quad+\partial^kZ^I\left(f (u_2,\partial u_2)-f (u_1,\partial u_1)\right)-(g^{\alpha\beta}(u_2,\partial u_2)-g^{\alpha\beta}(u_1,\partial u_1))\partial_\alpha\partial_\beta\partial^k Z^Iv_2\\&\quad+[\left(g^{\alpha\beta}(u_2,\partial  u_2)-m^{\alpha\beta}\right),\partial^kZ^I]\partial_\alpha\partial_\beta \wt{v}+[\left(g^{\alpha\beta}(u_2,\partial u_2)-g^{\alpha\beta}(u_1,\partial u_1)\right),\partial^kZ^I]\partial_\alpha\partial_\beta v_1\\&\quad+\left(g^{\alpha\beta}(u_2,\partial u_2)-m^{\alpha\beta}\right)\partial^k[\partial_\alpha\partial_\beta,Z^I]\wt{v}+\left(g^{\alpha\beta}(u_2,\partial u_2)-g^{\alpha\beta}(u_1,\partial u_1)\right)\partial^k[\partial_\alpha\partial_\beta,Z^I]v_1.\end{aligned}\end{equation}
Here we set $\wt{v}=v_1-v_2$ and make use of \eqref{eqn2}.

Define a new energy \fm{\wt{E}_{k,i}(t):=\sum_{l\leq k,|I|\leq i}\kh{ E_{u_1}(\partial^lZ^I\wt{v})(t)+\int_t^\infty E_q(\partial^lZ^I\wt{v})(\tau)\ d\tau}.}Here $E_{u_1}$ is defined in \eqref{sec5ene} with $u$ replaced by $u_1$ and $E_q$ is defined by \eqref{sec5eneq}. Using the pointwise bounds in Section \ref{sec6}, we can show that \eqref{equivnorm} still holds, so by the weighted Klainerman-Sobolev inequality, we have
\eq{\label{s7ptb1} |\partial Z^I\wt{v}(t,x)|&\lesssim \eps t^{C\eps}\wt{E}_{0,|I|+2}^{1/2}(t)\lra{r+t}^{-1}\lra{r-t}^{-1/2}w_0(r-t)^{-1/2}}
for all $t\in[T_\eps,T_2]$, $x\in\R^3$ and $|I|\leq N-5$. Since $u_1=u_2=u_{app}$ for $t\geq T_2$, we can apply Lemma \ref{lemptb2} and obtain
\eq{\label{s7ptb2}&|Z^I\wt{v}(t,x)|\\
&\lesssim \eps\sup_{\tau\in[t,T_2]} (\tau^{C\eps}\wt{E}_{0,|I|+2}^{1/2}(\tau))\cdot (\lra{t}^{-1}\lra{r-t}^{-(\gamma_1-1)/2}1_{r< t}+\lra{t}^{-1}\lra{r-t}^{1/2}1_{t\leq r\leq  2t}+\lra{r}^{-1/2}1_{r>2t})}for all $t\in[T_\eps,T_2]$, $x\in\R^3$ and $|I|\leq N-5$. Similar to \eqref{ptb2cor} and \eqref{ptb2cor2}, here we have
\eq{\label{s7ptb2cor}\norm{\wt{v}(t)\cdot 1_{r/t\in[1-2c,1-c]}}_{L^2(w)}\lesssim \lra{t}^{1+C\eps}\cdot \sup_{\tau\in[t,T_2]} (\tau^{C\eps}\wt{E}_{0,2}(\tau)^{1/2}),}
and if $|I|>0$,
\eq{\label{s7ptb2cor2}\norm{Z^I\wt{v}(t)\cdot 1_{r/t\in[1-2c,1-c/2]}}_{L^2(w)}&\lesssim t\sum_{|J|=|I|-1}\norm{\partial Z^I\wt{v}(t)\cdot 1_{r/t\in[1-2c,1-c/2]}}_{L^2(w)}\lesssim tE_{0,|I|-1}(t)^{1/2}}
whenever $t\in[T_\eps,T_2]$.

For $k+i\leq N-3$ with $|I|=i$, and for $t\in[T_\eps,T_1]$, we have 
\eq{\label{f71}&\norm{g^{\alpha\beta}(u_1,\partial u_1)\partial_\alpha\partial_\beta \partial^kZ^I\wt{v}}_{L^2(w)}\\&\leq C\eps t^{-1}\wt{E}_{k,i}(t)^{1/2}+C\eps^{1/2}(\ln\eps^{-1})^{-1/2} t^{-1+C\eps}\sum_{|I|\leq i}E_q(\partial^lZ^I\wt{v})(t)^{1/2}\\&\quad+C\eps t^{-1+C\eps}(\wt{E}_{k-1,i}(t)^{1/2}+\wt{E}_{k+1,i-1}(t)^{1/2})+C\lra{t}^{-\frac{\gamma_1+2}{4}+C\eps}\cdot \sup_{\tau\in[t,T_2]} (\tau^{C\eps}\wt{E}_{0,2}(\tau)^{1/2})}with $\wt{E}_{-1,\cdot}=\wt{E}_{\cdot,-1}=0$. This is a simple application of Lemma \ref{c1l2.1}, Lemma \ref{l24} and the estimates for $u_1,v_1,u_2,v_2$.  We skip the detail of the proof here, since it is very similar to that of \eqref{sec6fff}. We only make three remarks here.
\begin{enumerate}
\item  We should always put $L^2(w)$ norm on the terms involving $\wt{v}$ and put $L^\infty$ norm on the terms involving $u_1,u_2,v_1,v_2$. The pointwise bounds only holds for $|I|\leq N-2$, as seen in \eqref{ptb1} and \eqref{ptb2}, so we need to assume $k+i\leq N-3$ instead of $k+i\leq N$ above.
\item There is no term like $R_6^{(*),k,I}$ in the previous section, so we expect $\wt{E}_{k,i}$ to have a better decay than $E_{k,i}$.
\item The last term in \eqref{f71} comes from \eqref{s7ptb2cor}.
\end{enumerate}

Since \eqref{ptb1} and \eqref{ptb2} hold for $v_1$, we can apply energy estimate \eqref{energymain}  for  $E_{u_1}$.  Thus, for all $T_\eps\leq t\leq T_1$ and for $k+i\leq N-3$, 
\fm{\wt{E}_{k,i}(t)&\leq \wt{E}_{k,i-1}(t)+\wt{E}_{k-1,i}(t)+ C\eps^2 T_1^{-1+C\eps}+C\int_t^{T_1}(\eps\tau^{-1}+\tau^{-17/16})\wt{E}_{k,i}(\tau)\ d\tau\\&\quad+C\eps\int_t^{T_1}\tau^{-1+C\eps}(\wt{E}_{k-1,i}(\tau)^{1/2}+\wt{E}_{k+1,i-1}(\tau)^{1/2})\wt{E}_{k,i}(\tau)^{1/2}\ d\tau\\
&\quad+C\eps^{1/2}(\ln\eps^{-1})^{-1/2}\int_t^{T_1} \tau^{-1+C\eps}\sum_{l\leq k,\ |I|\leq i}E_q(\partial^lZ^I\wt{v})(\tau)^{1/2}\wt{E}_{k,i}(\tau)^{1/2}\ d\tau\\
&\quad+C\int_t^{T_1}\tau^{-(\gamma_1+2)/4+C\eps}\sup_{\tau'\in[\tau,T_2]}((\tau')^{C\eps}\wt{E}_{0,2}(\tau')^{1/2})\wt{E}_{k,i}(\tau)^{1/2}\ d\tau.}
Note that
\fm{&C\eps^{1/2}(\ln\eps^{-1})^{-1/2}\int_t^{T_1} \tau^{-1+C\eps}\sum_{l\leq k,\ |I|\leq i}E_q(\partial^lZ^I\wt{v})(\tau)^{1/2}\wt{E}_{k,i}(\tau)^{1/2}\ d\tau\\
&\leq C\eps^{1/2}(\ln\eps^{-1})^{-1/2}\int_t^{T_1} \kh{\sum_{l\leq k,\ |I|\leq i}E_q(\partial^lZ^I\wt{v})(\tau)+\tau^{-2+2C\eps}\wt{E}_{k,i}(\tau)}\ d\tau\\
&\leq C\eps^{1/2}(\ln\eps^{-1})^{-1/2}(\wt{E}_{k,i}(t)+\int_t^{T_1} \tau^{-2+2C\eps}\wt{E}_{k,i}(\tau)\ d\tau)}
By choosing $\eps\ll 1$, we have
\eq{\label{sec7ff6}\wt{E}_{k,i}(t)&\leq \wt{E}_{k,i-1}(t)+\wt{E}_{k-1,i}(t)+ C\eps^2 T_1^{-1+C\eps}+C\int_t^{T_1}(\eps\tau^{-1}+\tau^{-17/16})\wt{E}_{k,i}(\tau)\ d\tau\\&\quad+C\eps\int_t^{T_1}\tau^{-1+C\eps}(\wt{E}_{k-1,i}(\tau)^{1/2}+\wt{E}_{k+1,i-1}(\tau)^{1/2})\wt{E}_{k,i}(\tau)^{1/2}\ d\tau\\
&\quad+C\int_t^{T_1}\tau^{-(\gamma_1+2)/4+C\eps}\sup_{\tau'\in[\tau,T_2]}((\tau')^{C\eps}\wt{E}_{0,2}(\tau')^{1/2})\wt{E}_{k,i}(\tau)^{1/2}\ d\tau.}
We emphasize that to derive \eqref{sec7ff6}, we do not need any assumption on the decay rate of $\wt{E}_{k,i}(t)$.

We now use a continuity argument to prove \eqref{lem7.2c1}. Suppose that for some $T_c\in(T_\eps,T_1]$, we have proved
\eq{\label{s7.2.1ca} \wt{E}_{k,i}(t)\leq \wt{B}_{k,i}\eps^2T_1^{-1+2\lambda_0}t^{-2\lambda_0+\wt{C}_{k,i}\eps},\qquad \forall t\in[T_c,T_1].}
Here $\wt{B}_{k,i}$ and $\wt{C}_{k,i}$ are large constants to be chosen later. We remark that they play the same role as in the $B_{k,i}$ and $C_{k,i}$ in Section \ref{sec6}. Note that we have
\eq{\label{s7.2.1ca2}\wt{E}_{k,i}(t)\leq C\eps^2 t^{-1+C\eps}\leq C\eps^2 T_1^{-1+2\lambda_0}t^{-2\lambda_0+C\eps},\qquad\forall t\in[T_1,T_2].}  This is because we have  proved not only \eqref{prop6c1} but also $E_{k,i}(t)\lesssim \eps^2 t^{-1+C\eps}$ in Section \ref{sec6}.

Given \eqref{s7.2.1ca} and \eqref{s7.2.1ca2}, we have
\fm{\sup_{\tau\in[t,T_2]}(\tau^{C\eps}\wt{E}_{0,2}(\tau)^{1/2})&\leq \sup_{\tau\in[t,T_1]}\wt{B}_{0,2}^{1/2}\eps T_1^{-1/2+\lambda_0}t^{-\lambda_0+\wt{C}_{0,2}\eps/2+C\eps}+\sup_{\tau\in[T_1,T_2]}C\eps \tau^{-1/2+C\eps}\\
&\leq C \wt{B}_{0,2}^{1/2}\eps T_1^{-1/2+\lambda_0}t^{-\lambda_0+\wt{C}_{0,2}\eps/2+C\eps}.}
whenever $t\in[T_\eps,T_1]$. So the last integral in \eqref{sec7ff6} is controlled by 
\fm{&C\int_t^{T_1}\tau^{-(\gamma_1+2)/4+C\eps}\cdot C \wt{B}_{0,2}^{1/2}\eps T_1^{-1/2+\lambda_0}\tau^{-\lambda_0+\wt{C}_{0,2}\eps/2+C\eps}\cdot \wt{B}_{k,i}^{1/2}\eps T_1^{-1/2+\lambda_0}\tau^{-\lambda_0+\wt{C}_{k,i}\eps/2}\ d\tau\\
&\leq \wt{C}_N\eps^2T_1^{-1+2\lambda_0}\int_t^{T_1}\tau^{-(\gamma_1+2)/4-2\lambda_0+\wt{C}_N\eps}\ d\tau\leq \wt{C}_N\eps^2T_1^{-1+2\lambda_0} t^{(2-\gamma_1)/4-2\lambda_0+\wt{C}_N\eps}\leq \eps^2T_1^{-1+2\lambda_0} t^{-2\lambda_0}.}
As in Section \ref{sec6}, here $\wt{C}_N$ denotes a sufficiently large constant depending on all  the $\wt{B}_{*,*}$ and $\wt{C}_{*,*}$. Since $\gamma_1>2$, by choosing $\eps\ll1$, we have $(2-\gamma_1)/4+\wt{C}_N\eps\leq (2-\gamma_1)/8$ and thus $\wt{C}_Nt^{(2-\gamma_1)/4+\wt{C}_N}\leq \wt{C}_NT_\eps^{(2-\gamma_1)/8}\leq 1$. This is why we have the last estimate above.

We now seek to prove \eqref{s7.2.1ca} with $\wt{B}_{k,i}$ replaced by $\wt{B}_{k,i}/2$. To achieve this goal, we first induct on $i=0,1,\dots,N$ and then on $k=0,\dots,N-3-i$ for each fixed $i$. For convenience we set $\wt{C}_{-1,\cdot}=\wt{C}_{\cdot,-1}=\wt{B}_{-1,\cdot}=\wt{B}_{\cdot,-1}=0$. Fix $(k,i)$ and suppose  we have chosen the values of  $\wt{B}_{k',i'},\wt{C}_{k',i'}$  for all $i'<i$, or for all $i'=i$ and $k'<k$. Then, it follows from \eqref{sec7ff6} that
\fm{\wt{E}_{k,i}(t)&\leq \wt{B}_{k,i-1}\eps^2T_1^{-1+2\lambda_0} t^{-2\lambda_0+\wt{C}_{k,i-1}\eps}+\wt{B}_{k-1,i}\eps^2T_1^{-1+2\lambda_0} t^{-2\lambda_0+\wt{C}_{k-1,i}\eps}+ C\eps^2 T_1^{-1+C\eps}\\&\quad+C\int_t^{T_1}(\eps\tau^{-1}+\tau^{-17/16})\wt{B}_{k,i}\eps^2T_1^{-1+2\lambda_0}\tau^{-2\lambda_0+\wt{C}_{k,i}\eps}\ d\tau\\&\quad+C\eps\int_t^{T_1}\tau^{-1+C\eps}\wt{B}_{k-1,i}^{1/2}\wt{B}_{k+1,i-1}^{1/2}\eps^2T_1^{-1+2\lambda_0}\cdot \tau^{-2\lambda_0+(\wt{C}_{k,i}+\wt{C}_{k-1,i}+\wt{C}_{k+1,i-1})\eps/2} \ d\tau\\
&\quad+\eps^2T_1^{-1+2\lambda_0} t^{-2\lambda_0}\\
&\leq \wt{B}_{k,i-1}\eps^2T_1^{-1+2\lambda_0} t^{-2\lambda_0+\wt{C}_{k,i-1}\eps}+\wt{B}_{k-1,i}\eps^2T_1^{-1+2\lambda_0} t^{-2\lambda_0+\wt{C}_{k-1,i}\eps}\\&\quad+ C\eps^2 T_1^{-1+C\eps}+\wt{C}_N(\eps +T_\eps^{-1/16})\eps^2T_1^{-1+2\lambda_0}t^{-2\lambda_0+\wt{C}_{k,i}\eps} \\&\quad+\wt{C}_N\eps\tau^{-1+C\eps}\eps^2T_1^{-1+2\lambda_0}\cdot t^{-2\lambda_0+(\wt{C}_{k,i}+\wt{C}_{k-1,i}+\wt{C}_{k+1,i-1}+C)\eps/2}+\eps^2T_1^{-1+2\lambda_0} t^{-2\lambda_0}.}
If we choose $\wt{C}_{k,i}\gg1$ so that $\wt{C}_{k-1,i}+\wt{C}_{k+1,i-1}+C\leq \wt{C}_{k,i}$, then
$\wt{E}_{k,i}$ has an upper bound
\fm{\eps^2T_1^{-1+2\lambda_0}t^{-2\lambda_0+\wt{C}_{k,i}\eps}\cdot\kh{\wt{B}_{k,i-1}+\wt{B}_{k-1,i}+C+\wt{C}_N\eps^{1/16}}.}
If we choose $\wt{B}_{k,i}\gg1$ so that $\wt{B}_{k,i-1}+\wt{B}_{k-1,i}+C\leq \wt{B}_{k,i}/4$, then we have \fm{\wt{E}_{k,i}(t)\leq \frac{1}{2}\wt{B}_{k,i}\eps^2T_1^{-1+2\lambda_0}t^{-2\lambda_0+\wt{C}_{k,i}\eps},\qquad t\in[T_c,T_1]} as long as $\eps\ll1$. This finishes the proof of \eqref{lem7.2c1} for $t\in[T_\eps,T_1]$.

\subsubsection{The case $0\leq t\leq T_{\eps}$}
For $0\leq t\leq T_{\eps}$, we use the equation
\eq{\label{diffeqn}&g^{\alpha\beta}(u_1,\partial u_1)\partial_\alpha\partial_\beta \partial^kZ^I\wt{v}\\&=\partial^k[\Box,Z^I]\wt{v}+[g^{\alpha\beta}(u_1,\partial u_1)-m^{\alpha\beta},\partial^kZ^I]\partial_\alpha\partial_\beta \wt{v}+(g^{\alpha\beta}(u_1,\partial u_1)-m^{\alpha\beta})\partial^k[\partial_\alpha\partial_\beta, Z^I]\wt{v}\\&\quad-\partial^kZ^I(g^{\alpha\beta}(u_2,\partial u_2)-g^{\alpha\beta}(u_1,\partial u_1))\partial_\alpha\partial_\beta u_2+\partial^kZ^I(f(u_2,\partial u_2)-f(u_1,\partial u_1))} for each $k\geq 0$ and $I$ with $k+|I|\leq N-3$. Here we still set $\wt{v}=u_1-u_2$.  The proof of \eqref{lem7.2c1} is now similar to that of Proposition \ref{prop6.9}. We thus skip the details here. Note that to estimate the $L^2(w_0)$ norm of \eqref{diffeqn}, we always put $L^2(w_0)$ norm on the terms involving $\wt{v}$ and put $L^\infty$ norm on the terms involving $u_1$ and $u_2$. This again explains why we only have $|I|\leq N-3$ instead of $|I|\leq N$ in \eqref{lem7.2c1}.

\subsection{Uniqueness}\label{s7uni}
We now briefly explain how we show the uniqueness results in Remark \ref{rmkthmunique}. Note that it is sufficient to prove the uniqueness result in Remark \ref{rmk712}. In fact, as suggested by the computations in Remark \ref{rmk7.1.3Uuapp}, the estimates \eqref{p7f1}, \eqref{p7f2}, and \eqref{p7f3} are essentially equivalent to those estimates in Theorem \ref{mthm}, even if we replace $(c,c')$ in the main theorem by other $(d,d')$ with $d,d'\in(0,1/4)$. Let $\wt{v}$ be the difference of two solutions, and we notice that $\wt{v}$ satisfies the equations \eqref{diffeqn} such that for all $T\geq T_\eps= 1/\eps$,
\fm{\sum_{|I|\leq N-3}\norm{\partial Z^I \wt{v}(T)}_{L^2(w_0)}\lesssim\eps T^{-1/2+C\eps}.}
This estimate follows from \eqref{p7f1}. Now, we apply the energy estimates \eqref{energymain} on the time interval $[t,T]$ with $t\geq T_\eps$. Following the proof in Section \ref{sec7.1}, we can show that 
\fm{\sum_{|I|\leq N-3}\norm{\partial Z^I \wt{v}(t)}_{L^2(w_0)}\lesssim\eps T^{-1/2+\lambda_0}t^{-\lambda_0+C\eps},\qquad t\geq T_\eps.}
Here the implicit constants are all independent of $T$, so by taking $\eps\ll 1$ and $T\to\infty$, we conclude that $\wt{v}\equiv 0$ for $t\geq T_\eps$. It is then straightforward to show that $\wt{v}\equiv 0$ for all $t\geq 0$. We skip the detail in this proof.

\subsection{Proof of the results in Remark  \ref{mthmrmk3}}
\label{sec7.4pfprop1.1}

We now return to Remark  \ref{mthmrmk3}. We will prove the following proposition.

\prop{\label{mthmrmk3prop}Let $(\mu,U)$ be an admissible global solution to the geometric reduced system. For  $\eps\ll1$, let $u$ be the global solution constructed in Theorem \ref{mthm}. 

\begin{enumerate}[\rm (a)]
\item Suppose that there exist $(s^0,q^0,\omega^0)\in(1,\infty)\times\R\times\mathbb{S}^2$ such that \fm{|U(s,q^0,\omega^0)|+|\mu U_q(s,q^0,\omega^0)|\gtrsim \exp(-Cs),\qquad\forall s>s^0.} Then, for $\eps\ll1$, we have $u\not\equiv 0$.
\item Suppose that there exist $\kappa\geq 0$, $M>0$, and  $(s^0,q^0,\omega^0)\in(1,\infty)\times\R\times\mathbb{S}^2$ such that \fm{|U(s,q,\omega^0)|\geq M \lra{q}^{-\kappa}\exp(-Cs),\qquad \forall s>s^0,\ q>q^0.}Then, for each $\lambda\in(0,\frac{1}{2\kappa+1})$ and $\eps\ll_{\kappa,\lambda,q^0,M}1$, we have
\fm{|u(t,r\omega^0)|\gtrsim_{\kappa}M\eps t^{-1-C_\kappa \eps}\lra{r-t}^{-\kappa}}
for all $t\geq e^{(s^0+\delta)/\eps}$ and $r>0$ such that $|r-t|<t^\lambda$ and $q(t,r\omega^0)>q^0$. Thus, $\supp(u,u_t)|_{t=0}$ is not compact.
\item Suppose that there exist $\kappa\geq 1$, $M>0$, and  $(s^0,q^0,\omega^0)\in(1,\infty)\times\R\times\mathbb{S}^2$ such that \fm{|\mu U_q(s,q,\omega^0)|\geq M \lra{q}^{-\kappa}\exp(-Cs),\qquad \forall s>s^0,\ q>q^0.}Then, for each $\lambda\in(0,\frac{1}{2\kappa-1})$ and $\eps\ll_{\kappa,\lambda,q^0,M}1$, we have
\fm{|(u_t-u_r)(t,r\omega^0)|\gtrsim_\kappa M\eps t^{-1-C_\kappa \eps}\lra{r-t}^{-\kappa}}
for all $t\geq e^{(s^0+\delta)/\eps}$ and $r>0$ such that $|r-t|<t^\lambda$ and $q(t,r\omega^0)>q^0$. Thus, $\supp(u,u_t)|_{t=0}$ is not compact.
\end{enumerate}
Here the choice of $\eps$ depends on the admissible solution $(\mu,U)$ and $M$. The constants in these estimates, however, do not depend on $M$. 
}
\rm
\bigskip

Now we prove Proposition \ref{mthmrmk3prop}. By part (ii) of Theorem \ref{mthm}, we have
\eq{\label{keyestsec7.4}(u_t-u_r)(t,x)&=\eps r^{-1}(\mu U_q)(t,x)+O(\eps t^{-3/2+C\eps}\lra{r-t}^{-1/2}),\\
u(t,x)&=\eps r^{-1}U(t,x)+O(\eps t^{-3/2+C\eps}\lra{r-t}^{1/2})} for all 
$t\geq 1/\eps$ and $|x|/t\in[1/2,3/2]$.
Recall that $U$ and $\mu U_q$ can be viewed as a function of both  $(t,x)$ and $(s,q,\omega)$ at the same time.

Fix $q^0$ and $\omega^0$. Fix an arbitrary time $t\geq e^{(s^0+\delta)/\eps}$. We first claim that there exists a unique $r(t)>0$ such that $r(t)/t\in[1/2,3/2]$ and that $q(t,r(t)\omega^0)=q^0$, as long as $\eps\ll_{q^0}1$. Here recall that $q$ was defined in Section \ref{sec4.1} as the solution to the transport equation \eqref{qeqn}. Since $s^0>1$, we have $t>e^{1/\eps}$. Then, we have $t>\max\{4|q^0|,1/\eps\}$ for all $\eps\ll_{q^0}1$. The existence of such an $r(t)$ was proved in Remark \ref{rmk4.2.1}. By \eqref{lem4.1c1} in Lemma \ref{lem4.1}, we have
\fm{|q(t,\frac{3}{2}t\omega^0)-t/2|+|q(t,\frac{1}{2}t\omega^0)+t/2|\lesssim t^{C\eps}.}
Since $Ct^{C\eps}<t/4$ and $|q^0|<t/4$ for $t\geq e^{(s^0+\delta)/\eps}$ and $\eps\ll1$, we have 
\fm{q(t,\frac{1}{2}t\omega^0)<-t/4<q^0<t/4<q(t,\frac{3}{2}t\omega^0).}
By Remark \ref{rmk4.2.1}, we have $q_r>0$ everywhere, so the map $\rho\mapsto q(t,\rho\omega^0)$ is strictly increasing. This gives the uniqueness of $r(t)$ and forces $r(t)/t\in[1/2,3/2]$.

Let us first prove part (a). By \eqref{keyestsec7.4}, for all $t\geq e^{(s^0+\delta)/\eps}$ and $\eps\ll_{q^0}1$, we have
\fm{&|u(t,r(t)\omega^0)|+|(u_t-u_r)(t,r(t)\omega^0)|\\&\geq |\eps r(t)^{-1}U(\eps\ln t-\delta,q^0,\omega^0)|+|\eps r(t)^{-1}(\mu U_q)(\eps\ln t-\delta,q^0,\omega^0)|-C\eps t^{-3/2+C\eps}\lra{r(t)-t}^{1/2}\\
&\geq C^{-1}\eps t^{-1-C\eps}-C_{q^0}\eps t^{-3/2+C\eps}.}
In the last inequality, we use $t^{-C\eps}\lesssim \lra{q(t,x)}/\lra{r-t}\lesssim t^{C\eps}$ whenever $t\geq 1/\eps$ and $|r-t|\leq t/2$; see \eqref{lem4.1c3} in Lemma \ref{lem4.1}. By choosing $\eps\ll1$, we can make $-1-C\eps>-3/2+C\eps$. Thus, for $t\geq e^{(s^0+\delta)/\eps}$ and $\eps\ll_{q^0}1$, we have $|u(t,r(t)\omega^0)|+|(u_t-u_r)(t,r(t)\omega^0)|>0$. That is, $u\not\equiv 0$.

Next, we prove part (b). Fix $\kappa\geq 0$ and $\lambda\in(0,\frac{1}{2\kappa+1})$. By the assumptions on $|U|$, \eqref{keyestsec7.4}, and \eqref{lem4.1c3}, we have
\fm{|u(t,r\omega^0)|&\geq C^{-1}M\eps t^{-1-C\eps}\lra{q}^{-\kappa}-C\eps t^{-3/2+C\eps}\lra{r-t}^{1/2}\\
&\geq C_\kappa^{-1}M\eps t^{-1-C_\kappa\eps}\lra{r-t}^{-\kappa}-C\eps t^{-3/2+C\eps}\lra{r-t}^{1/2}}
whenever $t\geq e^{(s^0+\delta)/\eps}$, $r/t\in[1/2,3/2]$, and $q(t,r\omega^0)>q^0$. If we further assume that $|r-t|<t^\lambda$, then we have
\fm{\frac{C_\kappa^{-1}M\eps t^{-1-C_\kappa\eps}\lra{r-t}^{-\kappa}}{C\eps t^{-3/2+C\eps}\lra{r-t}^{1/2}}=C_\kappa^{-1}Mt^{1/2-C_\kappa\eps}\lra{r-t}^{-1/2-\kappa}\geq C_\kappa^{-1}Mt^{1/2-C_\kappa\eps-\lambda(\kappa+1/2)}.}
Since $0<\lambda (\kappa+1/2)<1/2$, the power of $t$ on the right side is positive as long as $\eps\ll_{\lambda,\kappa} 1$. Since $t\geq e^{1/\eps}$, by taking $\eps\ll_{\lambda,\kappa,q^0,M}1$, we can make the quotient above larger than $2$. We thus conclude that \fm{|u(t,r\omega^0)|\gtrsim_\kappa M \eps t^{-1-C_\kappa\eps}\lra{r-t}^{-\kappa},\qquad \forall t\geq e^{(s^0+\delta)/\eps}, |r-t|<t^\lambda,q(t,r\omega^0)>q^0.}
Note that if $r=t+t^\lambda/2$, we have
\fm{q(t,r\omega^0)\geq (r-t)-|q(t,r\omega^0)-(r-t)|\geq t^\lambda/2-Ct^{C\eps}>t^\lambda/4,\quad \forall t\geq e^{(s^0+\delta)/\eps},\eps\ll_{\kappa,\lambda,q^0}1.}
In other words, we have $q(t,(t+t^\lambda/2)\omega^0)>q^0$ for $t\geq e^{(s^0+\delta)/\eps}$, so $u(t,(t+t^\lambda/2)\omega^0)\neq 0$ for all sufficiently large $t$. On the other hand, if the data of $u$ have compact support, then  $u(t,r\omega^0)=0$ whenever $t\geq 0$ and $r-t\gtrsim 1$. A contradiction. Thus, the data of $u$ cannot have compact support.

Finally, the proof of part (c) is essentially the same as that of part (b).
\section{Examples}\label{sec8}

In this section, we apply our main theorem to several examples.  

\subsection{A scalar  quasilinear wave equation satisfying the weak null condition}\label{sec8.1}
Consider the scalar quasilinear wave equation 
\eq{\label{exmeqn1} g^{\alpha\beta}(u,\partial u)\partial_\alpha\partial_\beta u=f(u,\partial u).}
Here the unknown $u$ is $\R$-valued, and the coefficients $g^{\alpha\beta}$ and $f$ satisfy the same assumptions as in \eqref{qwe}. We also assume that they satisfy a ``partial'' null condition.  That is, if 
\eq{g^{\alpha\beta}(u,\partial u)&=m^{\alpha\beta}+g^{\alpha\beta}_0 u+g^{\alpha\beta\lambda}_0\partial_\lambda u+O(|u|^2+|\partial u|^2),}
\eq{f(u,\partial u)&=f_0^{\alpha\beta}\partial_\alpha u\partial_\beta u+O(|u|^4+|\partial u|(|u|^2+|\partial u|^2)),}
then we assume that for all $\omega\in\mathbb{S}^2$,
\eq{g^{\alpha\beta\lambda}_0\wh{\omega}_\alpha\wh{\omega}_\beta\wh{\omega}_\gamma\equiv 0,}
\eq{f^{\alpha\beta\lambda}_0\wh{\omega}_\alpha\wh{\omega}_\beta\equiv 0.}
One example of \eqref{exmeqn1} is $\Box u=u\Delta u$.

By setting $G(\omega):=g^{\alpha\beta}_0\wh{\omega}_\alpha\wh{\omega}_\beta$, we obtain the geometric reduced system for \eqref{exmeqn1}:
\eq{\label{exmeqn1rs}\left\{\begin{array}{l}
\displaystyle \partial_s(\mu U_q)=0,\\
\displaystyle \partial_s \mu=\frac{1}{4}G(\omega)\mu^2U_q.
\end{array}\right.}
This geometric reduced system is exactly the same as that of the quasilinear wave equation 
\fm{g^{\alpha\beta}(u)\partial_\alpha\partial_\beta u=0} which was studied in  \cite{MR4232783,MR4315017}. By following the proof there and applying part (a) of Remark \ref{mthmrmk3}, we have the following result.

\cor{Fix $\gamma_+>1$ and $\gamma_->2$. For each function $A\in C^\infty(\R\times\mathbb{S}^2)$ such that   $\partial_q^a\partial_\omega^cA=O_{a,c}(\lra{q}^{-a-\gamma_{\sgn(q)}})$ for each $a,c\geq 0$, we set
\fm{\mu(s,q,\omega):=-2\exp(-\frac{1}{2}G(\omega)A(q,\omega)s)}
and
\fm{U(s,q,\omega)=\int_{-\infty}^qA(p,\omega)\exp(\frac{1}{2}G(\omega)A(p,\omega)s)\ dp.} Then, $(\mu,U)$ is a $(\gamma_+,\gamma_-)$-admissible global solution to the geometric reduced system. Moreover, for each $\eps\ll1$, there exists a matching global solution $u$ to the system \eqref{exmeqn1} of size $\eps$. If furthermore we assume $A\not\equiv 0$, the global solution $u$ is nontrivial.}\rm
\bigskip

We remark that this corollary  extends Theorem 1.1 in \cite{MR4232783} or Theorem 3.1 in \cite{MR4315017}. Here we relax the assumptions on the scattering data; see Remark \ref{rmkthmext}. Moreover, we do not need to assume that  the coefficients $g$ are independent of $\partial u$ and that the semilinear term $f$ is vanishing.

\subsection{A scalar  semilinear wave equation violating the weak null condition}\label{sec8.2}

Consider the scalar semilinear wave equation 
\eq{\label{exmeqn2} \Box u=f(u,\partial u).}
Here the unknown $u$ is $\R$-valued, and  $f$ is a smooth function such $f$ satisfies \eqref{qwefi}.  One example of \eqref{exmeqn2} is John's counterexample $\Box u=u_t^2$.

We now set $F(\omega):=f_0^{\alpha\beta}\wh{\omega}_\alpha\wh{\omega}_\beta$ for $\omega\in\mathbb{S}^2$. Here $f_0^{\alpha\beta}$ comes from the Taylor expansion
\eq{f(u,\partial u)=f_0^{\alpha\beta}\partial_\alpha u\partial_\beta u+O(|u|^4+|\partial u|(|u|^2+|\partial u|^2)).} Then, the geometric reduced system for \eqref{exmeqn2} is
\eq{\label{exmeqn2rs}\left\{\begin{array}{l}
\displaystyle \partial_s(\mu U_q)=-\frac{1}{4}F(\omega)(\mu U_q)^2,\\
\displaystyle \partial_s \mu=0.
\end{array}\right.}
For simplicity we set $\mu|_{s=0}\equiv -2$, so we have $\mu\equiv -2$ for all $s$. The first equation then reduces to
\fm{\partial_s  U_q=\frac{1}{2}F(\omega)  U_q ^2.}
Note that this equation is the same as H\"{o}rmander's asymptotic equation.

We now obtain the following lemma.
\lem{\label{lemexm2.1}Fix $A\in C^\infty(\R\times\mathbb{S}^2)$, and fix two constants $\gamma_+>1$, $\gamma_->2$. Suppose that $F(\omega)A(q,\omega)\leq 0$ for all $(q,\omega)$ and that $\partial_q^a\partial_\omega^cA=O_{a,c}(\lra{q}^{-a-\gamma_{\sgn(q)}})$ for each $a,c\geq 0$. If we set
\fm{U(s,q,\omega)&:=\int_{-\infty}^q \frac{2A(p,\omega)}{-F(\omega)A(p,\omega)s+2}\ dp,}
then $(-2,U)$ is a  $(\gamma_+,\gamma_-)$-admissible global solution to the geometric reduced system \eqref{exmeqn2rs}.

Moreover, if $F(\omega^0)A(q^0,\omega^0)<0$ for some $(q^0,\omega^0)\in\R\times\mathbb{S}^2$, then there exists a large constant $\wt{C}>1$ such that 
\eq{\label{lemexm2.1:c}|U(s,q,\omega)|\gtrsim  |A(q^0,\omega^0)|^2\lra{s}^{-1},\quad \forall s>0,\ q>q^0,\ |\omega-\omega^0|\leq \wt{C}^{-1}\min\{|F(\omega^0)|,|A(q^0,\omega^0)|\}.}
Both $\wt{C}$ and the implicit constant are independent of $(q^0,\omega^0)$.}
\begin{proof}It is easy to check that $(-2,U)$ is a solution to \eqref{exmeqn2rs}. Moreover, by choosing sufficiently small $\delta_0>0$, we have $0\leq -F(\omega)A(q,\omega)\delta_0\leq 1$ for all $(q,\omega)$. This implies that $F(\omega)A(q,\omega)s+2\geq 1$ whenever $s\geq -\delta_0$. As a result, $(-2,U)$ is indeed a global solution for all $s\geq -\delta_0$.

It remains to check  \eqref{def1.1a3}--\eqref{def1.1a8} (or  \eqref{def3.1a3}--\eqref{def3.1a8} in Definition \ref{def3.1}). Most of these estimates hold trivially, and we only need to check \eqref{def1.1a7} and \eqref{def1.1a8}. Note that
$\partial_s^a\partial_q^b\partial_\omega^cU_q$ is equal to a linear combination of terms of the form
\fm{(-FAs+2)^{-1-m}\partial_{q}^{b_0}\partial_\omega^{c_0}A\cdot\prod_{j=1}^m\partial_s^{a_j}\partial_q^{b_j}\partial_\omega^{c_j}(FAs).}
Here $m\geq 0$, $\sum a_*=a$, $\sum b_*=b$, $\sum c_*=c$ and $a_j+b_j+c_j>0$ for each $j$. Since $-FAs+2\geq 1$ for each $s\geq -\delta_0$, it is clear that
\fm{|(-FAs+2)^{-1-m}\partial_{q}^{b_0}\partial_\omega^{c_0}A\cdot\prod_{j=1}^m\partial_s^{a_j}\partial_q^{b_j}\partial_\omega^{c_j}(FAs)|\lesssim \lra{q}^{-b-(m+1)\gamma_{\sgn(q)}}\lra{s}^m\lesssim \lra{q}^{-b-\gamma_{\sgn(q)}}\exp(C s ).}
This gives us \eqref{def1.1a8}. By integrating \eqref{def1.1a8}, we obtain \eqref{def1.1a7}.

Next, suppose that $F(\omega^0)A(q^0,\omega^0)<0$ for some $(q^0,\omega^0)$. Fix a large constant $\wt{C}>1$ independent of the choice of $(q^0,\omega^0)$. We will choose its value later in the proof. Since $|\partial_\omega F(\omega)|\lesssim 1$, we have
$|F(\omega)-F(\omega^0)|\lesssim |\omega-\omega^0|$. If $\wt{C}/2$ is larger than the implicit constant here (depending only on $F$ but not on $\omega^0$), then whenever $|\omega-\omega^0|\leq \wt{C}^{-1}|F(\omega^0)|$, we have
\fm{|F(\omega)-F(\omega^0)|\leq \frac{1}{2}|F(\omega^0)|\Longrightarrow \sgn(F(\omega))=\sgn(F(\omega^0))\neq 0.}
Since $F(\omega)A(q,\omega)\leq 0$ everywhere, we conclude that whenever $q\in\R$ and $|\omega-\omega^0|\leq \wt{C}^{-1}|F(\omega^0)|$, we have
\fm{\sgn(F(\omega^0))\cdot A(q,\omega)\leq 0}
and thus
\fm{|U(s,q,\omega)|&=\int_{-\infty}^{q}\frac{2|A(p,\omega)|}{-F(\omega)A(p,\omega)s+2}\ dp\geq \int_{-\infty}^{q^0}\frac{2|A(p,\omega)|}{-F(\omega^0)A(p,\omega)s+2}\ dp\\
&\geq \frac{2\int_{-\infty}^{q^0}|A(p,\omega)|\ dp}{\sup_{(q,\omega)\in\R\times\mathbb{S}^2}|F(\omega)\cdot A(q,\omega)|s+2}\gtrsim \lra{s}^{-1}\int_{-\infty}^{q^0}|A(p,\omega)|\ dp.}
The implicit constant in the last estimate is independent of $(q^0,\omega^0)$. To continue, we notice that $|\partial_q A|+|\partial_\omega A|\lesssim 1$, so $|A(q,\omega)-A(q^0,\omega^0)|\lesssim |q-q^0|+|\omega-\omega^0|$. If $\wt{C}/2$ is larger than the implicit constant here (depending only on $A$ but not on $(q^0,\omega^0)$), then whenever $|\omega-\omega^0|\leq \wt{C}^{-1}|F(\omega^0)|$ and $|q-q^0|+|\omega-\omega^0|\leq \wt{C}^{-1}|A(q^0,\omega^0)|$, we have
\fm{|A(q,\omega)-A(q^0,\omega^0)|\leq \frac{1}{2}|A(q^0,\omega^0)|\Longrightarrow |A(q,\omega)|\geq \frac{1}{2}|A(q^0,\omega^0)|}
and therefore
\fm{\int_{-\infty}^{q^0}|A(p,\omega)|\ dp\geq \int_{q^0-\wt{C}^{-1}|A(q^0,\omega^0)|}^{q^0} \frac{1}{2}|A(q^0,\omega^0)|\ dp=\wt{C}^{-1}|A(q^0,\omega^0)|\cdot\frac{1}{2}|A(q^0,\omega^0)|. }
This finishes the proof.
\end{proof}\rm

\bigskip

We can now apply Theorem \ref{mthm} and Remark \ref{mthmrmk3} (or Proposition \ref{mthmrmk3prop}). This gives us the following corollary.

\cor{\label{corexm2}For each function $A=A(q,\omega)$ as in Lemma \ref{lemexm2.1}, we let $(-2,U)$ be the corresponding admissible global solution to the geometric reduced system. Then, for each $\eps\ll1$, there exists a global solution $u$  for $t\geq 0$ to the system \eqref{exmeqn2} matching $(-2,U)$ at infinite time of size $\eps$. 

Moreover, the solution $u$  is nontrivial as long as $A\not\equiv 0$. If furthermore $F\cdot A\not\equiv 0$, then the data of $u$ are not compactly supported. 

In summary, there exists a family of nontrivial global solutions to \eqref{exmeqn2} for $t\geq 0$. Moreover, in the case when the null condition is not satisfied (i.e.\ $F\not\equiv 0$), then there exists a nonempty subset of the family of global solutions above whose initial data are not compactly supported. In the case when $F\neq 0$ everywhere, then the initial data of any of the solutions constructed above are not compactly supported. }

\rmk{\rm For $\Box u=u_t^2$, we mentioned in Remarks \ref{rmkfutureglobal} and \ref{rmkfutureglobal:john} that all nontrivial global solutions constructed here are future global but not past global.  }

\rmk{\label{corexm2:rmk}\rm In the case when $F(\omega^0)A(q^0,\omega^0)<0$ for some $(q^0,\omega^0)\in\R\times\mathbb{S}^2$, we can prove a lower bound for $\norm{u(t)}_{L^2(\R^3)}$:
\eq{\norm{u(t)}_{L^2(\R^3)}\gtrsim_\lambda |A(q^0,\omega^0)|^2 \eps t^{\lambda/2-C\eps},\qquad \forall t\geq e^{(1+\delta)/\eps},\ \eps\ll_{\lambda}1,\ \lambda\in(0,1).} 
Let us prove this estimate. By the second half of Lemma \ref{lemexm2.1}, we obtain a large constant $\wt{C}>0$. By Lemma \ref{lem4.1}, we have $|q-(r-t)|\lesssim t^{C\eps}$. By enlarging $\wt{C}$ if necessary, whenever $t\geq e^{(1+\delta)/\eps}$ and $q^0+\wt{C}t^{\wt{C}\eps}<r-t<t/2$, we have
\fm{q(t,x)\geq r-t-|q-(r-t)|>q^0.}
Then, by part (b) of Proposition \ref{mthmrmk3prop} (with $\kappa=0$ and $M=|A(q^0,\omega^0)|^2>0$), for all $0<\lambda<1/2$ and $\eps\ll_{\lambda,q^0,\omega^0}1$, we have
\fm{|u(t,r\omega)|\gtrsim |A(q^0,\omega^0)|^2 \cdot \eps t^{-1-C\eps}}
whenever $t\geq e^{(1+\delta)/\eps}$, $|\omega-\omega^0|\leq \wt{C}^{-1}\min\{|F(\omega^0)|,|A(q^0,\omega^0)|\}$, and $q^0+\wt{C}t^{\wt{C}\eps}< r-t< t^\lambda$. From its proof, we notice that the implicit constant is independent of $(t,r,\omega)$. As a result, we have
\fm{\norm{u(t)}_{L^2(\R^3)}&\gtrsim|A(q^0,\omega^0)|^2 \cdot\eps t^{-1-C\eps}\cdot |\{x\in\R^3:\ q^0+\wt{C}t^{\wt{C}\eps}< |x|-t< t^\lambda\}|^{1/2}\\
&\gtrsim|A(q^0,\omega^0)|^2 \cdot \eps t^{-1-C\eps}\cdot (t^{\lambda+2})^{1/2}\gtrsim|A(q^0,\omega^0)|^2 \cdot\eps t^{\lambda/2-C\eps},\qquad \forall t\geq e^{(1+\delta)/\eps}.}
In the second estimate we use $q^0+\wt{C}t^{\wt{C}\eps}<t^\lambda/2$ as long as $t\geq e^{(1+\delta)/\eps}$ and $\eps\ll_{q^0,\omega^0,\wt{C},\lambda}1$. 
}
\rm

\subsection{A scalar  quasilinear wave equation violating the weak null condition}\label{sec8.3}
Consider the scalar quasilinear wave equation
\eq{\label{exmeqn3} g^{\alpha\beta}(\partial u)\partial_\alpha\partial_\beta u=0.}
Here the unknown $u$ is $\R$-valued, and the coefficients $g^{\alpha\beta}$ satisfy the same assumptions as in \eqref{qwe}. One example of \eqref{exmeqn3} is John's counterexample $\Box u=u_t u_{tt}$. 

For convenience, we set  $G(\omega):=g^{\alpha\beta\lambda}_0\wh{\omega}_\alpha\wh{\omega}_\beta\wh{\omega}_\lambda$ for all $\omega\in\mathbb{S}^2$. The $g^{\alpha\beta\lambda}_0$ come from the Taylor expansions \eq{\label{exmeqn3tay}g^{\alpha\beta}(\partial u)=m^{\alpha\beta}+g^{\alpha\beta\lambda}_0\partial_\lambda u+O(|\partial u|^2).}

Unlike the previous two examples, here we do not apply Theorem \ref{mthm} directly to \eqref{exmeqn3}. This is because we cannot obtain an interesting global solution if we do so.  To see this, we first notice that the geometric reduced system for \eqref{exmeqn3} is
\eq{\label{exmeqn3original}\left\{\begin{array}{l}
\displaystyle \partial_s(\mu U_q)=0,\\
\displaystyle \partial_s \mu=-\frac{1}{8}G(\omega)\mu^2\partial_q(\mu U_q).
\end{array}\right.}
Set $\mu|_{s=0}\equiv- 2$. Then,  this system has a global solution for $s\geq 0$ if and only if $G(\omega)  U_{qq}|_{s=0}\geq 0$. However, by \eqref{def3.1a8} in Definition \ref{def3.1}, we have $\lim_{q\to\infty}U_q(s,q,\omega)=0$, so we need \fm{ G(\omega)U(0,q,\omega)\equiv 0,\qquad\forall (q,\omega)\in\R\times\mathbb{S}^2.}
In the most interesting case, we may have $G(\omega)\neq 0$ everywhere on $\mathbb{S}^2$, so the reduced system above only has a trivial admissible solution. As a result, applying Theorem \ref{mthm} directly to \eqref{exmeqn3} is not satisfactory.

In this subsection, we need to use an alternative way to construct nontrivial global solutions. Our construction is based on the simple observation that  $u$ solves \eqref{exmeqn3} if and only if $(v_{(\alpha)})_{\alpha=0,1,2,3}=(\partial_\alpha u)_{\alpha=0,1,2,3}$ solves the following system of quasilinear wave equations
\eq{\label{exmeqn3var}g^{\alpha\beta}(v)\partial_\alpha\partial_\beta v_{(\sigma)}&=-\frac{1}{4}(\partial_{v^\lambda}g^{\alpha\beta})(v)\cdot (\partial_\sigma v_{(\lambda)}+\partial_\lambda v_{(\sigma)})(\partial_\alpha v_{(\beta)}+\partial_\beta v_{(\alpha)}),\qquad \sigma=0,1,2,3}
coupled to a system of constraint equations
\eq{\label{exmeqn3ct} \partial_\alpha v_{(\beta)}=\partial_\beta v_{(\alpha)},\qquad \alpha,\beta=0,1,2,3.}
Note that $g^{\alpha\beta}$ is a given function of $v=(v^\alpha)\in\R^4$, so $\partial_{v^\lambda}g^{\alpha\beta}$ is also a given function defined in $\R^4$. Moreover, because of \eqref{exmeqn3tay}, we have $\partial_{v^\lambda}g^{\alpha\beta}(0)=g_0^{\alpha\beta\lambda}$. Thus, the geometric reduced system for \eqref{exmeqn3var} is
\eq{\label{exmeqn3rs}
\left\{
\begin{array}{l}
\displaystyle \partial_s(\mu \partial_qU_{(\sigma)})=\frac{1}{16}g^{\alpha\beta\lambda}_0 (\wh{\omega}_\sigma \partial_q U_{(\lambda)}+\wh{\omega}_\lambda \partial_q U_{(\sigma)})(\wh{\omega}_\alpha \partial_q U_{(\beta)}+\wh{\omega}_\beta \partial_q U_{(\alpha)})\mu^2,\quad \sigma=0,1,2,3;\\[1em]
\displaystyle \partial_s\mu=\frac{1}{4}g^{\alpha\beta\sigma}_0\wh{\omega}_\alpha\wh{\omega}_\beta\mu^2 \partial_q U_{(\sigma)}. 
\end{array}
\right.}
Since we expect $\partial_\alpha v_{(\beta)}\approx -\frac{\eps}{2r}\wh{\omega}_\alpha\mu \partial_qU_{(\beta)}$ for $r\approx t\gg1$, from \eqref{exmeqn3ct} we obtain an additional equation
\eq{\label{exmeqn3ctvar} \wh{\omega}_\alpha\partial_qU_{(\beta)}=\wh{\omega}_\beta\partial_qU_{(\alpha)},\qquad \alpha,\beta=0,1,2,3.}
Since $\wh{\omega}_0\equiv-1$, we have $\partial_qU_{(\alpha)}=-\wh{\omega}_\alpha\partial_qU_{(0)}$. Then, the reduced system \eqref{exmeqn3rs} can be reduced to
\eq{\left\{
\begin{array}{l}
\displaystyle \partial_s(\mu \partial_qU_{(0)})=-\frac{1}{4}G(\omega) (\mu\partial_q U_{(0)}  )^2;\\[1em]
\displaystyle \partial_s\mu=-\frac{1}{4}G(\omega)\mu^2 \partial_q U_{(0)}. 
\end{array}
\right.}
By solving this system explicitly,  we obtain the following lemma.

\lem{\label{lemexm3.1}Fix $A\in C^\infty(\R\times\mathbb{S}^2)$, and fix two constants $\gamma_+>1,\gamma_->2$. Suppose that $G(\omega)A(q,\omega)\leq 0$ for all $(q,\omega)$ and that $\partial_q^a\partial_\omega^cA=O_{a,c}(\lra{q}^{-a-\gamma_{
\sgn(q)}})$ for each $a,c\geq 0$. If we set
\fm{U_{(\alpha)}(s,q,\omega)&:=\int_{-\infty}^q -\wh{\omega}_{\alpha}A(p,\omega)\ dp,\qquad \alpha=0,1,2,3}and
\fm{\mu(s,q,\omega)&:=\frac{4}{G(\omega)A(q,\omega)s-2},}
then $(\mu,(U_{(\alpha)}))$ is a $(\gamma_+,\gamma_-)$-admissible global solution to the geometric reduced system \eqref{exmeqn3rs}.

Moreover, if $G(\omega^0)A(q^0,\omega^0)<0$ for some $(q^0,\omega^0)\in\R\times\mathbb{S}^2$, then there exists a large constant $\wt{C}>1$ such that
\eq{\label{lemexm3.1:c}|U(s,q,\omega)|\gtrsim  |A(q^0,\omega^0)|^2,\quad \forall s>0,\ q>q^0,\ |\omega-\omega^0|\leq \wt{C}^{-1}\min\{|G(\omega^0)|,|A(q^0,\omega^0)|\}.}
Both $\wt{C}$ and the implicit constant are independent of $(q^0,\omega^0)$. }
\begin{proof}
It is easy to check that $(\mu,(U_{(\alpha)}))$ is a solution to \eqref{exmeqn3rs}. Moreover, by choosing sufficiently small $\delta_0>0$, we have $0\leq -G(\omega)A(q,\omega)\delta_0\leq 1$ for all $(q,\omega)$. This implies that $G(\omega)A(q,\omega)s-2\leq -1$ whenever $s\geq -\delta_0$. As a result, $(\mu,(U_{(\alpha)}))$ is indeed a global solution for all $s\geq -\delta_0$.

It remains to check  \eqref{def1.1a3}--\eqref{def1.1a8} (or  \eqref{def3.1a3}--\eqref{def3.1a8} in Definition \ref{def3.1}). The  proofs rely heavily on the estimate $G(\omega)A(q,\omega)s-2\leq -1$ whenever $s\geq -\delta_0$. Here we shall only check \eqref{def1.1a52} and leave the proofs of other assumptions as exercise. Note that 
\fm{\partial_q(\mu\partial_qU_{(0)})&=\partial_q(\frac{4A}{GAs-2})=\frac{4A_q}{GAs-2}-\frac{4GAA_qs}{(GAs-2)^2}.}
Since $A_q=O(\lra{q}^{-1-\gamma_{\sgn(q)}})$ and $GAs-2\leq -1$, we have 
\fm{|\partial_q(\mu\partial_qU_{(0)})|&\lesssim |A_q|+|A_q|\cdot\frac{|GAs|}{|GAs-2|^2}\lesssim \lra{q}^{-1-\gamma_{\sgn(q)}}(1+\frac{|GAs|}{2-GAs}).}
If $|GAs|\leq 2$, then we have $\frac{2}{2-GAs}\leq 2$. If $|GAs|>2$, then we have $s>0$ and thus $2-GAs\geq -GAs\geq 2$ and thus $\frac{|GAs|}{2-GAs}\leq \frac{|GAs|}{-GAs}=1$. Thus \eqref{def1.1a52} follows.

The proof of \eqref{lemexm3.1:c} is essentially the same as that of \eqref{lemexm2.1:c}. Here we do not have a factor $\lra{s}^{-1}$ on the right side because of the different definitions of $U$.
\end{proof}
\rmk{\rm The global solution $(\mu,(U_{(\alpha)}))$ constructed in this lemma induces a global solution to the geometric reduced system for \eqref{exmeqn3}. In fact, if we set $W(s,q,\omega)$ by  \fm{W(s,q,\omega)=\int_{-\infty}^q 2(\mu^{-1}U_{(0)})(s,p,\omega)\ dp} then this $(\mu,U)$ is a global solution to \eqref{exmeqn3original}. However,  this  $(\mu,W)$ is not an admissible solution. It can be easily checked that $\mu W_q=2U_{(0)}$ does not decay as $q\to\infty$ unless $A\equiv 0$. This again explains why we need to study the differentiated equations \eqref{exmeqn3var}. See Remark~\ref{rmk1.1.1}.
}
\rm

\bigskip 
By applying Theorem \ref{mthm} and Remark \ref{mthmrmk3}, we obtain the following corollary.

\cor{\label{cor8.5intermediate}For each function $A=A(q,\omega)$ as in Lemma \ref{lemexm3.1}, we let $(\mu,U)$ be the corresponding admissible global solution to the geometric reduced system. Then, for each $\eps\ll1$, there exists a global solution $v=(v_{(\alpha)})_{\alpha=0,1,2,3}$ for $t\geq 0$ to the system \eqref{exmeqn3var} matching $(\mu,U)$ at infinite time of size $\eps$.

Moreover, the solution $v$  is nontrivial as long as $A\not\equiv 0$. If we have $G\cdot A\not\equiv 0$, then the data of $v$ are not compactly supported.}
\rmk{\label{corexm3:rmk}\rm By following the proof in Remark \ref{corexm2:rmk}, in the case when $GA\not\equiv 0$, we can also prove \eq{\norm{v(t)}_{L^2(\R^3)}\gtrsim_\lambda |A(q^0,\omega^0)|^2 \cdot\eps t^{\lambda/2-C\eps},\qquad \forall t\geq e^{(1+\delta)/\eps},\ \eps\ll_{\lambda}1,\ \lambda\in(0,1).}
Here we choose $(q^0,\omega^0)$ so that $G(\omega^0)A(q^0,\omega^0)< 0$.
Later we will check that $\partial u= v$, so this estimate gives a lower bound of $\norm{\partial u(t)}_{L^2(\R^3)}$.}\rm

\bigskip

Our last step is to show that there does exist some solution $u$ to \eqref{exmeqn3} such that $v_{(\alpha)}=\partial_\alpha u$. It suffices to check the constraint equations \eqref{exmeqn3ct}. For each $\alpha,\beta=0,1,2,3$, we set
\eq{w_{\alpha\beta}:=\partial_\alpha v_{(\beta)}-\partial_\beta v_{(\alpha)}.}
In order to show $w\equiv 0$, we make use of the following two lemmas. Lemma \ref{lemexm3.2} states  that $w$ is a solution to some general linear wave equations, and Lemma \ref{lemexm3.3} states that the energy of $w$ tends to $0$ as $t\to\infty$. 

\lem{\label{lemexm3.2}We have 
\fm{g^{\alpha\beta}(v)\partial_\alpha\partial_\beta w_{\sigma\gamma}&=h(v,\partial v)\cdot \partial w\cdot \partial v+h(v,\partial v)\cdot w\cdot \partial^2v+h(v,\partial v)\cdot w\cdot \partial v\cdot \partial v.}Here $h(v,\partial v)\cdot \partial w\cdot \partial v$ means a linear combination of terms of the form
\fm{h^{\alpha'\beta'\sigma'\gamma'\lambda'}(v,\partial v)\partial_{\lambda'}w_{\alpha'\beta'}\partial_{\sigma'} v_{(\gamma')}}
where $h(v,\partial v)$ is a smooth function of $(v,\partial v)$. The expressions $h(v,\partial v)\cdot w\cdot \partial^2v$ and $h(v,\partial v)\cdot w\cdot \partial v\cdot \partial v$ have  similar meanings.}
\begin{proof}
For $\sigma,\gamma=0,1,2,3$,
\fm{&g^{\alpha\beta}(v)\partial_\alpha\partial_\beta\partial_\sigma v_{(\gamma)} =\partial_\sigma(g^{\alpha\beta}(v)\partial_\alpha\partial_\beta v_{(\gamma)})-\partial_\sigma(g^{\alpha\beta}(v))\cdot \partial_\alpha\partial_\beta v_{(\gamma)}\\
&=\partial_\sigma(-\frac{1}{4}(\partial_{v^\lambda}g^{\alpha\beta})(v)\cdot (\partial_\gamma v_{(\lambda)}+\partial_\lambda v_{(\gamma)})(\partial_\alpha v_{(\beta)}+\partial_\beta v_{(\alpha)}))-(\partial_{v^\lambda}g^{\alpha\beta})(v)\cdot\partial_\sigma v_{(\lambda)}\partial_\alpha\partial_\beta v_{(\gamma)}\\
&=-\frac{1}{4}(\partial_{v^{\lambda'}}\partial_{v^\lambda}g^{\alpha\beta})(v)\cdot \partial_\sigma v_{(\lambda')}\cdot (\partial_\gamma v_{(\lambda)}+\partial_\lambda v_{(\gamma)})(\partial_\alpha v_{(\beta)}+\partial_\beta v_{(\alpha)})\\
&\quad -\frac{1}{4}(\partial_{v^\lambda}g^{\alpha\beta})(v)\cdot (\partial_\sigma\partial_\gamma v_{(\lambda)}+\partial_\sigma\partial_\lambda v_{(\gamma)})(\partial_\alpha v_{(\beta)}+\partial_\beta v_{(\alpha)})\\
&\quad -\frac{1}{4}(\partial_{v^\lambda}g^{\alpha\beta})(v)\cdot (\partial_\gamma v_{(\lambda)}+\partial_\lambda v_{(\gamma)})(\partial_\sigma\partial_\alpha v_{(\beta)}+\partial_\sigma\partial_\beta v_{(\alpha)})\\
&\quad-(\partial_{v^\lambda}g^{\alpha\beta})(v)\cdot\partial_\sigma v_{(\lambda)}\partial_\alpha\partial_\beta v_{(\gamma)}\\
&=:R^{1}_{\sigma\gamma}+R^{2}_{\sigma\gamma}+R^{3}_{\sigma\gamma}+R^{4}_{\sigma\gamma}.}
In the second row, we use the fact that $v$ solves \eqref{exmeqn3var}. In the third identity, we apply the product rule and the chain rule. Similarly, we can show that $g^{\alpha\beta}(v)\partial_\alpha\partial_\beta \partial_\gamma v_{(\sigma)}=\sum_{j=1}^4 R_{\gamma\sigma}^j$. Here, to obtain $R^{j}_{\gamma\sigma}$, we simply interchange the roles of $\sigma$ and $\gamma$ in the definition of $R^j_{\sigma\gamma}$. To finish the proof, we now need to compute $\sum_{j=1}^4 R_{\sigma\gamma}^j-\sum_{j=1}^4 R_{\gamma\sigma}^j$.

Let us start with $R^{1}_{\sigma\gamma}-R^1_{\gamma\sigma}$. By the definition of $w$, we have $\partial_\lambda v_{(\gamma)}=w_{\lambda\gamma}+\partial_\gamma v_{(\lambda)}$. As a result, we have 
\fm{R^{1}_{\sigma\gamma}&=-\frac{1}{4}(\partial_{v^{\lambda'}}\partial_{v^\lambda}g^{\alpha\beta})(v)\cdot (\partial_\alpha v_{(\beta)}+\partial_\beta v_{(\alpha)})\cdot \partial_\sigma v_{(\lambda')}\cdot (2\partial_\gamma v_{(\lambda)}+w_{\lambda\gamma})\\
&=-\frac{1}{2}(\partial_{v^{\lambda'}}\partial_{v^\lambda}g^{\alpha\beta})(v)\cdot (\partial_\alpha v_{(\beta)}+\partial_\beta v_{(\alpha)})\cdot \partial_\sigma v_{(\lambda')}\cdot \partial_\gamma v_{(\lambda)}+h(v)\cdot \partial v\cdot\partial v\cdot w.}
Here $h(v)$ denotes a smooth function of $v$.  Since $\partial_{v^{\lambda'}}\partial_{v^\lambda}g^{\alpha\beta}=\partial_{v^\lambda}\partial_{v^{\lambda'}}g^{\alpha\beta}$, the term $(\partial_{v^{\lambda'}}\partial_{v^\lambda}g^{\alpha\beta})(v)\cdot \partial_\sigma v_{(\lambda')}\cdot \partial_\gamma v_{(\lambda)}$ is symmetric with respect to the indices $\sigma$ and $\gamma$. Thus,
\fm{R^{1}_{\sigma\gamma}-R^1_{\gamma\sigma}
&=h(v)\cdot w\cdot \partial v\cdot\partial v.}

We move on to $R^2_{\sigma\gamma}-R^2_{\gamma\sigma}$. From the definition of $w$, we have $\partial_\sigma\partial_\lambda v_{(\gamma)}=\partial_\sigma w_{\lambda\gamma}+\partial_\sigma\partial_\gamma v_{(\lambda)}$ and thus
\fm{R^2_{\sigma\gamma}&=-\frac{1}{4}(\partial_{v^\lambda}g^{\alpha\beta})(v)\cdot (2\partial_\sigma\partial_\gamma v_{(\lambda)}+\partial_\sigma w_{\lambda\gamma})(\partial_\alpha v_{(\beta)}+\partial_\beta v_{(\alpha)})\\
&=-\frac{1}{2}(\partial_{v^\lambda}g^{\alpha\beta})(v)\cdot  \partial_\sigma\partial_\gamma v_{(\lambda)}(\partial_\alpha v_{(\beta)}+\partial_\beta v_{(\alpha)})+h(v)\cdot \partial v\cdot \partial w.}
Since the first term on the right hand side is symmetric with respect to $\sigma$ and $\gamma$, we have
\fm{R^2_{\sigma\gamma}-R^2_{\gamma\sigma}=h(v)\cdot \partial v\cdot \partial w.}

Next, we compute $R^3_{\sigma\gamma}-R^3_{\gamma\sigma}+R^4_{\sigma\gamma}-R^4_{\gamma\sigma}$. By the definition of $w$, we have
\fm{R^3_{\sigma\gamma}&=-\frac{1}{4}(\partial_{v^\lambda}g^{\alpha\beta})(v)\cdot (\partial_\gamma v_{(\lambda)}+\partial_\lambda v_{(\gamma)})(\partial_\alpha \partial_\beta v_{(\sigma)}+\partial_\alpha w_{\sigma\beta}+\partial_\beta\partial_\alpha v_{(\sigma)}+\partial_\beta w_{\sigma\alpha})\\
&=-\frac{1}{2}(\partial_{v^\lambda}g^{\alpha\beta})(v)\cdot (2\partial_\gamma v_{(\lambda)}+w_{\lambda\gamma})\partial_\alpha \partial_\beta v_{(\sigma)}+h(v)\cdot \partial v\cdot\partial w\\
&=-(\partial_{v^\lambda}g^{\alpha\beta})(v)\cdot \partial_\gamma v_{(\lambda)}\cdot \partial_\alpha \partial_\beta v_{(\sigma)}+h(v)\cdot \partial v\cdot\partial w+h(v)\cdot \partial^2v\cdot w.}
It turns out that the first term on the right hand side is $R^4_{\gamma\sigma}$. In other words, we have
\fm{(R^3_{\sigma\gamma}-R^4_{\gamma\sigma})-(R^3_{\gamma\sigma}-R^4_{\sigma\gamma})&=h(v)\cdot \partial v\cdot\partial w+h(v)\cdot \partial^2v\cdot w.}
This finishes our proof.
\end{proof}

\lem{\label{lemexm3.3} For $t\geq 1/\eps$, we have 
\fm{\sum_{\alpha,\beta=0,1,2,3}\norm{w_0^{1/2}\partial w_{\alpha\beta}(t)}_{L^2(\R^3)}\lesssim \eps t^{-1/2+C\eps}.}
Here $w_0$ is defined in Theorem \ref{mthm}.}
\begin{proof}
By Theorem \ref{mthm}, we have
\eq{\label{lemexm3.3f1}\norm{w_0^{1/2}\partial \partial_\alpha(v_{(\beta)}-\eps r^{-1}U_{(\beta)})(t)}_{L^2(\{|x|/t\in[1-c',1+c'])}+\norm{w_0^{1/2} \partial\partial_\alpha v_{(\beta)}(t)}_{L^2(\{|x|/t\notin[1-c',1+c']\})}&\lesssim \eps t^{-1/2+C\eps}.}
Moreover, recall that $\partial_qU_{(\alpha)}=-\wh{\omega}_\alpha\partial_qU_{(0)}$ and that $U_{(*)}= 0$ as  $q\to-\infty$. Thus, we have $U_{(\alpha)}=-\wh{\omega}_\alpha U_{(0)}$. This implies that whenever $ |x|/t\in[1-c',1+c']$, 
\fm{\partial_\alpha (\eps r^{-1}U_{(\beta)})&=\partial_\alpha (-\eps r^{-1}\wh{\omega}_\beta U_{(0)})\\
&=\eps r^{-1}\wh{\omega}_\beta q_\alpha A  +\partial (\eps r^{-1}\wh{\omega}_\beta) \cdot  U_{(0)} +\eps r^{-1}\sum_{a+c=1} \partial_s^a\partial_\omega^cU_{(0)} \cdot(\eps t^{-1},\partial \omega)\\
&=\eps r^{-1}\wh{\omega}_\beta q_\alpha A  +\eps S^{-2,1-\gamma_-}_-\cap \eps S^{-2,0}_+.}
In the last estimate, we use the assumptions in Definition \ref{def3.1} and the notation $S_\pm^{s,p}$ introduced in Definition \ref{c1defn1.5}. By Lemma \ref{c1lem6} and Lemma \ref{lem4.5},  we have $q\in S^{0,1}$ and thus $ \wh{\omega}_\alpha q_\beta-\wh{\omega}_\beta q_\alpha\in S^{-1,1}$. Since $\partial_q^a\partial_\omega^cA=O(\lra{q}^{-a-\gamma_{\sgn(q)}})=O(\lra{r-t}^{-a-\gamma_{\sgn(q)}}t^{C\eps})$ by Lemma \ref{lem4.1c3},  we have $A\in S^{0,-\gamma_-}_-\cap S^{0,-\gamma_+}_+$ by following the proof of Lemma \ref{lem4.5}. As a result, we have
\fm{&\partial_\alpha (\eps r^{-1}U_{(\beta)})-\partial_\beta (\eps r^{-1}U_{(\alpha)})=\eps r^{-1}(\wh{\omega}_\beta q_\alpha-\wh{\omega}_\alpha q_\beta) A+\eps S^{-2,1-\gamma_-}_-\cap \eps S^{-2,0}_+\\&\in \eps S^{-2,1}\cdot (S^{0,-\gamma_-}_-\cap S^{0,-\gamma_+}_+)+\eps S^{-2,1-\gamma_-}_-\cap \eps S^{-2,0}_+\subset \eps S^{-2,1-\gamma_-}_-\cap \eps S^{-2,0}_+.}
The last inclusion follows from $\gamma_+>1$. Since 
\fm{ \norm{w_0^{1/2}\cdot \eps t^{-2+C\eps}(1_{r\geq t}+\lra{r-t}^{1-\gamma_-}1_{r<t})}_{L^2(\{|x|/t\leq [1-c',1+c']\})}\lesssim \eps t^{-1/2+C\eps},}
and since $|\partial w|\lesssim|\partial^2v|$,
by \eqref{lemexm3.3f1} we conclude that
\fm{\sum_{\alpha,\beta=0,1,2,3}\norm{w_0^{1/2}\partial w_{\alpha\beta}(t)}_{L^2(\R^3)}\lesssim \eps t^{-1/2+C\eps}.}
\end{proof}\rm

We can now prove $w\equiv 0$ by applying the energy estimates.
\lem{\label{lemexm3.4}We have $w\equiv 0$ for all $t\geq 0$.}
\begin{proof}
The proof here is essentially the same as that in Section \ref{sec7}. We shall combine the equations for $w_{**}$ with the energy estimates.  Details are skipped in this proof. 
\end{proof}\rm

\bigskip

We can now state our final result for this example. Note that the second half of this result follows from the second half of Corollary \ref{cor8.5intermediate}.
\cor{\label{corexm3} There exists a family of nontrivial global solutions to \eqref{exmeqn3} for $t\geq 0$.

Moreover, in the case when the null condition is not satisfied (i.e.\ $G\not\equiv 0$), there exists a nonempty subset of the family of global solutions above whose initial data are not compactly supported. In the case when $G\neq 0$ everywhere, then the initial data of any of the solutions constructed above are not compactly supported. }\rm

\bigskip

We remind our readers that if $u$ is a solution to \eqref{exmeqn3}, then it is $\partial u$ instead of $u$ that satisfies the estimates in Theorem \ref{mthm}.

We end this subsection with a result on the asymptotic behaviors of $u$ itself. The only way to estimate $u$ is to integrate the pointwise bounds in the main theorem.

\prop{\label{propasymujohn}For each function $A=A(q,\omega)$ as in Lemma \ref{lemexm3.1}, we let $(\mu,U)$ be the corresponding $(\gamma_+,\gamma_-)$-admissible solution for some $\gamma_+>1$ and $\gamma_->2$. Also fix $\gamma_1\in(2,\min\{4,2(\gamma_--1)\})$. For $\eps\ll_{A,\gamma_\pm,\gamma_1}1$, we let $v=(v_{(\alpha)})_{\alpha=0,1,2,3}$ be the global solution  to \eqref{exmeqn3var} for $t\geq 0$ of size $\eps$ constructed in Corollary \ref{cor8.5intermediate}, and let $u$ be a global solution  to \eqref{exmeqn3} for  $t\geq 0$ such that $v=\partial u$. Define
\fm{W(s,q,\omega)&=\int_{-\infty}^q 2(\mu^{-1}U_{(0)})(s,p,\omega)\ dp;\\ W(t,x)&=W(\eps\ln t-\delta,q(t,x),\omega),\quad t\geq 1/\eps,\ ||x|-t|\leq t/2.}
Then, the limit $u_{\infty,0}:=\lim_{t\to\infty}u(t,0)$ exists, and for all $(t,x)$ with $t\geq 1/\eps$, we have
\fm{&|u(t,x)-u_{\infty,0}-\eps|x|^{-1}W(t,x) 1_{|x|/t\in[1/2,3/2]}|\lesssim \eps t^{-1/2+C\eps}K(\gamma_1,t,r)+\eps t^{-\gamma_1/2+C\eps}.}
Here we define \fm{K(\gamma_1,t,r):= \left\{
\begin{array}{ll}
   t^{(1-\gamma_1)/2}+t^{-1}\lra{r-t}^{(3-\gamma_1)/2},  & r\leq t,\ \gamma_1\neq 3;  \\
   t^{-1}\ln(\frac{1+t}{1+t-r}),  & r\leq t,\ \gamma_1=3;\\
   t^{(1-\gamma_1)/2}+t^{-1}\lra{r-t}^{3/2}, &t< r\leq 2t,\ \gamma_1\neq 3;\\
   t^{-1}\ln t+t^{-1}\lra{r-t}^{3/2}, &t< r\leq 2t,\ \gamma_1= 3;\\
   r^{1/2},&r>2t.
\end{array}
\right.}
}
\rmk{\rm Following the same way, one can derive the asymptotics of $Z^Iu$. The formulas will be much more complicated, and we omit them for simplicity.}
\rmk{\rm Recall that any polynomial of $(t,x)$ of degree at most $1$ is a global solution to $\Box u=u_tu_{tt}$. A corollary of this proposition is that the solution $u$ here cannot be  a nonzero polynomial of $(t,x)$ of degree at most $1$. Note that $|u(t,x)-u_{\infty,0}|\lesssim\eps t^{-1-}$ for $|x|<t/2$. This clearly cannot hold if $u$ is a nonzero linear polynomial.}
\rmk{\rm Note that \fm{\eps t^{-1/2+C\eps}K(\kappa_1,t,r)+\eps t^{-\gamma_1/2+C\eps}\lesssim\eps t^{-1-}+\eps t^{-3/2+C\eps}\lra{r-t}^{3/2}1_{r/t\in(1,2]}+\eps t^{-1/2+C\eps}r^{1/2}1_{r>2t}.} Thus, when $|x|-t<t^{1/3-}$, the right side of the estimate in the proposition is $O(\eps t^{-1-})$, which indicates that we obtain a relatively good approximation of $u$ in that region. On the other hand, if $|x|/t>1+c$ for some fixed constant $c>0$, then the remainder is $O(\eps t^{-1/2+C\eps} r^{1/2})$. This remainder is too large and tends to infinity when $r\to\infty$, so we do not have a good approximation of $u$ there. }
\begin{proof}[Proof of Proposition \ref{propasymujohn}]
We first show that $u_{\infty,0}=\lim_{t\to\infty}u(t,0)$ exists. By part (ii) of Theorem \ref{mthm}, for $t\geq 1/\eps$ and $\alpha=0,1,2,3$, we have
\fm{&|(u_\alpha-\eps r^{-1}U_{(\alpha)})(t,x)|\cdot1_{|x|/t\in[1/2,3/2]}+|u_\alpha(t,x)|\cdot1_{|x|/t\notin[1/2,3/2]}\\
&\lesssim \eps t^{-1/2+C\eps}(t^{-1}\lra{r-t}^{-(\gamma_1-1)/2}1_{r<t}+t^{-1}\lra{r-t}^{1/2}1_{r/t\in[1,2]}+r^{-1/2}1_{r>2t}).}
By setting $x=0$, we have
$|u_t(t,0)|\lesssim \eps t^{-1-\gamma_1/2+C\eps}$, so $u_t(\cdot,0)\in L^1(1/\eps,\infty)$.
And since $\partial u$ is at least continuous in $[0,\infty)\times\R^3$, we have $u_t(\cdot,0)\in L^1(0,\infty)$. Thus, the limit
\fm{u_{\infty,0}=\lim_{t\to\infty}u(t,0)=u(0,0)+\int_0^\infty u_\tau(\tau,0)\ d\tau}
exists and is a finite real number.

Next, we fix $(t,x)$ with $t\geq 1/\eps$. By the fundamental theorem of calculus, we have
\eq{\label{asymujohnest1}&|u(t,x)-u_{\infty,0}+\int_{0}^{|x|}\eps\rho^{-1} U_{(0)}(t,\rho\omega) 1_{\rho/t\in[1/2,3/2]}\ d\rho|\\
&\leq\int_{t}^\infty |u_t(t,0)|\ dt+\int_{0}^{|x|}|\omega\cdot \nabla_xu(t,\rho\omega)-\eps \rho^{-1} \sum_{j=1}^3\omega_j\cdot U_{(j)}(t,\rho\omega) 1_{\rho/t\in[1/2,3/2]}|\ d\rho\\
&\lesssim \eps t^{-\gamma_1/2+C\eps}+\eps t^{-1/2+C\eps}\int_{0}^{|x|}(t^{-1}\lra{\rho-t}^{-(\gamma_1-1)/2}1_{\rho<t}+t^{-1}\lra{\rho-t}^{1/2}1_{\rho/t\in[1,2]}+\rho^{-1/2}1_{\rho>2t})\ d\rho.}
Here we recall from Lemma \ref{lemexm3.1} that $U_{(0)}=-\sum_{j=1}^3\omega_jU_{(j)}$. To continue, we compute the integral on the right side of \eqref{asymujohnest1}.

\lem{For $t\geq 1/\eps$ and $r\geq 0$, we have
\fm{&\int_{0}^{r}(t^{-1}\lra{\rho-t}^{-(\gamma_1-1)/2}1_{\rho<t}+t^{-1}\lra{\rho-t}^{1/2}1_{\rho/t\in[1,2]}+\rho^{-1/2}1_{\rho>2t})\ d\rho\lesssim K(\gamma_1,t,r).}
Recall that $K(\gamma_1,t,r)$ was defined in the statement of Proposition \ref{propasymujohn}.}
\begin{proof}
We first assume $r\leq t$. In this case, the integral is
\fm{\int_{0}^{r}t^{-1}\lra{\rho-t}^{-(\gamma_1-1)/2}\ d\rho\lesssim t^{-1}\int_{0}^{r}(1+t-\rho)^{-(\gamma_1-1)/2}\ d\rho=t^{-1}\int_{t-r}^{t}(1+\rho)^{-(\gamma_1-1)/2}\ d\rho.}
Recall that $\gamma_1\in(2,\min\{4,2(\gamma_--1)\})$. If $\gamma_1\neq 3$, we obtain an upper bound $t^{-1} (t^{(3-\gamma_1)/2}+\lra{r-t}^{(3-\gamma_1)/2})$ for the integral. If $\gamma_1=3$, we obtain an upper bound $t^{-1}\ln(\frac{1+t}{1+t-r})$.

Next, suppose $t< r\leq 2t$. In the previous paragraph, we have computed $\int_{0}^{t}t^{-1}\lra{\rho-t}^{-(\gamma_1-1)/2}\ d\rho$. Note that
\fm{\int_{t}^r t^{-1}\lra{\rho-t}^{1/2}\ d\rho&=\int_{0}^{r-t} t^{-1}(\rho+1)^{1/2}\ d\rho\lesssim t^{-1}\lra{r-t}^{3/2}.} Thus, for the integral in the statement of the lemma, we have an upper bound $t^{(1-\gamma_1)/2}+t^{-1}+t^{-1}\lra{r-t}^{3/2}\lesssim t^{(1-\gamma_1)/2}+t^{-1}\lra{r-t}^{3/2}$ if $\gamma_1\neq 3$, and an upper bound $t^{-1}\ln t+t^{-1}\lra{r-t}^{3/2}$ if $\gamma_1= 3$.

Finally, suppose that $r>2t$. In the previous  two paragraphs, we have computed $\int_{0}^{t}t^{-1}\lra{\rho-t}^{-(\gamma_1-1)/2}\ d\rho+\int_t^{2t} t^{-1}\lra{\rho-t}^{1/2} \ d \rho$. Since $\int_{2t}^r \rho^{-1/2}\ d\rho\lesssim r^{1/2}$, the integral in the statement of the lemma has an upper bound
$t^{(1-\gamma_1)/2}+t^{-1}\cdot t^{3/2}+r^{1/2}\lesssim r^{1/2}$ if $\gamma_1\neq 3$, and an upper bound $t^{-1}\ln t+t^{-1} \cdot t^{3/2}+r^{1/2}\lesssim r^{1/2}$ if $\gamma_1=3$.
\end{proof}\rm

By this lemma, we can replace the right side of \eqref{asymujohnest1} by $\eps t^{-\gamma_1/2+C\eps}+\eps t^{-1/2+C\eps}K(\gamma_1,t,|x|)$. We also notice that this upper bound coincides with the right side of the estimate in Proposition \ref{propasymujohn}.

To finish the proof of the proposition, we seek to replace the integral involving $U_{(0)}$ on the left side of \eqref{asymujohnest1} by $\eps |x|^{-1}W\cdot 1_{|x|/t\in[1/2,3/2]}$. We now claim that
\eq{\label{asymujohnest2}&|\int_{0}^{|x|}\eps\rho^{-1} U_{(0)}(t,\rho\omega) 1_{\rho/t\in[1/2,3/2]}\ d\rho+\eps |x|^{-1}W(t,x)1_{|x|/t\in[1/2,3/2]}|\\
&\lesssim \eps t^{-\gamma_1/2+C\eps}+\eps t^{-1+C\eps}K(\gamma_1,t,|x|).}
There is nothing to prove when $|x|<t/2$. When $|x|>3t/2$, we notice that \fm{|U_{(0)}(s,q,\omega)|\lesssim \int_{-\infty}^q |A(p,\omega)|\ dp\lesssim \int_{-\infty}^q \lra{p}^{-\gamma_{\sgn(p)}}\ dp\lesssim \lra{\max\{0,-q\}}^{1-\gamma_-}.}
It follows that
\fm{|\int_{0}^{|x|} \eps \rho^{-1}U_{(0)}(t,\rho\omega)\cdot 1_{\rho/t\in[1/2,3/2]}\ d\rho|\lesssim \eps t^{-1+C\eps}\int_{t/2}^{3t/2}(\lra{\rho-t}^{1-\gamma_-}1_{\rho\leq t}+1_{\rho>t})\ d\rho\lesssim\eps t^{C\eps}.}
In the first estimate, we use  \eqref{rmk4.1.1c2} in Remark \ref{rmk4.1.1} to convert powers of $\lra{q}$ to powers of $\lra{r-t}$ (up to a factor $t^{C\eps}$). Meanwhile, we notice that $\eps t^{-1/2+C\eps}K(\gamma_1,t,|x|)\gtrsim \eps t^{C\eps}$ whenever $|x|>3t/2$, so we obtain \eqref{asymujohnest2} when $|x|/t>3/2$. It remains to check \eqref{asymujohnest2} when $|x|/t\in[1/2,3/2]$.

We first notice the following estimate connecting $W$ and $U_{(0)}$.
\lem{For all $t\geq 1/\eps$ and $|x|/t\in[1/2,3/2]$, we have 
\fm{\abs{\eps |x|^{-1}\int_{t/2}^{r}U_{(0)}(t,\rho\omega)\ d\rho+\eps |x|^{-1}W(t,x)}&\lesssim \eps t^{-\gamma_1/2+C\eps}+\eps t^{-1/2+C\eps}K(\gamma_1,t,|x|).}}
\begin{proof}
Fix $(t,x)$ with $t\geq 1/\eps$ and $|x|/t\in[1/2,3/2]$. First, since $q_r=(\nu-\mu)/2$, we have
\fm{&\int_{t/2}^{r}U_{(0)}(t,\rho\omega)\ d\rho=\int_{t/2}^{r}-2(\mu^{-1} U_{(0)})(t,\rho\omega)\ q_r(t,\rho\omega)d\rho+\int_{t/2}^{r}(\nu\mu^{-1} U_{(0)} )(t,\rho\omega)\ d\rho.}
We first control the  integral involving $\nu$. Since $\gamma_->2$, we have
\fm{|\mu^{-1}U_{(0)}(s,q,\omega)|&\leq|(GAs-2)\int_{-\infty}^q A(p,\omega)\ dp |\lesssim \lra{s}\int_{-\infty}^q|A(p,\omega)|\ dp\lesssim e^{Cs}\lra{\max\{0,-q\}}^{1-\gamma_-}.}
By \eqref{rmk4.1.1c2}, we conclude that $|\mu^{-1}U_{(0)}(t,x)|\lesssim t^{C\eps}(\lra{r-t}^{1-\gamma_-}1_{r<t}+1_{r\geq t})$ for all $t\geq 1/\eps$ and $|x|/t\in[1/2,3/2]$.
By Lemma \ref{lem4.2} and \eqref{rmk4.1.1c2}, we have $\nu=O(t^{-1+C\eps}(\lra{r-t}^{1-\gamma_-}1_{r\leq t}+1_{r> t}))$  for $(t,x)$ in the same range. As a result,
\fm{&|\int_{t/2}^{r}(\nu\mu^{-1} U_{(0)} )(t,\rho\omega)\ d\rho|\lesssim \int_{t/2}^r t^{-1+C\eps}(\lra{\rho-t}^{2-2\gamma_-}1_{\rho\leq t}+1_{\rho> t})\ d\rho\\
& \lesssim t^{-1+C\eps}(\lra{r-t}^{3-2\gamma_-}1_{r\leq t}+\lra{r-t}1_{r> t}).}
In the last step, we use $3-2\gamma_-<-1$.

Next, we compute the first integral. By the change of variables $q=q(t,\rho\omega)$, we have 
\fm{&\int_{t/2}^{r}-2(\mu^{-1}U_{(0)})(\eps \ln t-\delta,q(t,\rho\omega),\omega)q_r(t,\rho\omega)\ d\rho\\
&=\int_{q(t,t\omega/2)}^{q(t,x)}-2(\mu^{-1}U_{(0)})(\eps \ln t-\delta,q,\omega)\ dq=-W(t,x)+W(t,t\omega/2).}
Since $|\mu^{-1}U_{(0)}(s,q,\omega)|\lesssim e^{Cs}\lra{\max\{0,-q\}}^{1-\gamma_-}$, we have $|W(s,q,\omega)|\lesssim e^{Cs}(\lra{q}^{2-\gamma_-}1_{q<0}+\lra{q}1_{q\geq 0})$ by integration. As a result, $W(t,t\omega/2)=W(\eps\ln t-\delta,q(t,t\omega/2),\omega)=O(t^{2-\gamma_-+C\eps})$. Here we recall the boundary condition $q|_{|x|=t/2}=-t/2$ in Section \ref{sec4}. In summary, we obtain 
\fm{\eps r^{-1}|\int_{t/2}^{r}U_{(0)}(t,\rho\omega)\ d\rho+W(t,x)|\lesssim \eps t^{1-\gamma_-+C\eps}+\eps t^{-2+C\eps}(\lra{r-t}^{3-2\gamma_-}1_{r\leq t}+\lra{r-t}1_{r> t}).}
We now need to show that whenever $t\geq 1/\eps$ and $r/t\in[1/2,3/2]$, 
\fm{&\eps t^{1-\gamma_-+C\eps}+\eps t^{-2+C\eps}(\lra{r-t}^{3-2\gamma_-}1_{r\leq t}+\lra{r-t}1_{r> t})\lesssim \eps t^{-\gamma_1/2+C\eps}+\eps t^{-1/2+C\eps}K(\gamma_1,t,r).}
Since $3-2\gamma_-<-1$ and since $\gamma_1/2<\min\{\gamma_--1,2\}$, for all $(t,r)$ in the range above  we have
\fm{\eps t^{1-\gamma_-+C\eps}+\eps t^{-2+C\eps} \lra{r-t}^{3-2\gamma_-}\lesssim \eps t^{1-\gamma_-+C\eps}+\eps t^{-2+C\eps}\lesssim \eps t^{-\gamma_1/2+C\eps}.} Moreover, we have $K(\gamma_1,t,r)\geq t^{-1}\lra{r-t}^{3/2}$ when $t< r\leq 3t/2$. Since $t\geq 1/\eps$ and $\lra{r-t}\geq 1$, we have \fm{\eps t^{-2+C\eps}\lra{r-t}1_{r>t}&\leq \eps 
 t^{-3/2+C\eps}\lra{r-t}^{3/2}\leq \eps t^{-1/2+C\eps} K(\gamma_1,t,r).} 
\end{proof}
\rm

We compare this lemma with \eqref{asymujohnest2}. It now remains to prove the next lemma.

\lem{For all $t\geq 1/\eps$ and $r\in[t/2,3t/2]$, we have\fm{&|\int_{t/2}^r \eps \rho^{-1}U_{(0)}(t,\rho\omega)\ d\rho-\int_{t/2}^r \eps r^{-1}U_{(0)}(t,\rho\omega)\ d\rho|\lesssim \eps t^{-\gamma_1/2+C\eps}+\eps t^{-1/2+C\eps}K(\gamma_1,t,r).}}
\begin{proof}
Fix $t\geq 1/\eps$ and . Since $\rho^{-1}-r^{-1}=(\rho r)^{-1}(r-\rho)$, we have
\eq{\label{lem813est1}&|\int_{t/2}^r \eps \rho^{-1}U_{(0)}(t,\rho\omega)\ d\rho-\int_{t/2}^r \eps r^{-1}U_{(0)}(t,\rho\omega)\ d\rho|\leq \int_{t/2}^r \eps (\rho r)^{-1}|r-\rho|\cdot |U_{(0)}(t,\rho\omega)|\ d\rho\\
&\lesssim \eps t^{-2+C\eps}\int_{t/2}^r  (r-\rho)\cdot (\lra{\rho-t}^{1-\gamma_-}1_{\rho\leq t}+1_{\rho>t})\ d\rho.}
Here we use the bounds for $U_{(0)}$ and \eqref{rmk4.1.1c2} to replace $q$ with $r-t$.

We first consider the case when $r\leq t$. In this case, the right side of \eqref{lem813est1} is bounded by 
\fm{\int_{t/2}^{r}  (r-\rho)\cdot (1+t-\rho)^{1-\gamma_-}\ d\rho=\int_{0}^{r-t/2}  \xi\cdot (\xi+t-r+1)^{1-\gamma_-}\ d\xi.}
Since $1-\gamma_-<-1$, we  integrate by parts to get
\fm{&\left.\frac{\xi(\xi+t-r+1)^{2-\gamma_-}}{2-\gamma_-}\right]_{\xi=0}^{\xi=r-t/2}-\int_{0}^{r-t/2}  \frac{ (\xi+t-r+1)^{2-\gamma_-}}{2-\gamma_-}\ d\xi.}
If $\gamma_-\neq 3$, the sum is 
\fm{0<-\frac{(r-t/2)(t/2+1)^{2-\gamma_-}}{\gamma_--2}+\frac{(t/2+1)^{3-\gamma_-}-(t-r+1)^{3-\gamma_-}}{(\gamma_--2)(3-\gamma_-)}\lesssim t^{3-\gamma_-}+\lra{r-t}^{3-\gamma_-}.}
To relate the right side with $K(\gamma_1,t,r)$, we notice that 
\fm{&\eps t^{-2+C\eps}(t^{3-\gamma_-}+\lra{r-t}^{3-\gamma_-})\lesssim\eps t^{1-\gamma_-+C\eps}+\eps t^{-3/2+C\eps} \lra{r-t}^{5/2-\gamma_-}\\&\lesssim \eps t^{-\gamma_1/2+C\eps}+\eps t^{-3/2+C\eps}\lra{r-t}^{3/2-\gamma_1/2}\lesssim\eps t^{-\gamma_1/2+C\eps}+ \eps t^{-1/2+C\eps} K(\gamma_1,t,r).}
In the second estimate, we use $\gamma_1/2<\gamma_--1$. In the last one, we use the definition of $K(\gamma_1,t,r)$ for $r\leq t$. If $\gamma_-=3$, the sum is
\fm{0<-(r-t/2)(t/2+1)^{-1}+\ln\frac{t/2+1}{t-r+1}\leq \ln \frac{1+t}{1+t-r}=tK(\gamma_1,t,r).}
Multiply both sides by $\eps t^{-2+C\eps}$ and we obtain an upper bound $\eps t^{-1+C\eps}K(\gamma_1,t,r)$. In summary, we finish the proof in the case when $r\leq t$.

Next, suppose that $r>t$. Since $\int_{t}^r(r-\rho)\ d\rho=(r-t)^2/2$, the integral on the right side of \eqref{lem813est1} is bounded by \fm{&\eps t^{-\gamma_1/2+C\eps}+\eps t^{-1/2+C\eps}K(\gamma_1,t,t)+\eps t^{-2+C\eps}\lra{r-t}^2.}
If $\gamma_1\neq 3$, then \fm{&\eps t^{-\gamma_1/2+C\eps}+\eps t^{-1/2+C\eps}K(\gamma_1,t,t)+\eps t^{-2+C\eps}\lra{r-t}^2\\&\lesssim \eps t^{-\gamma_1/2+C\eps}+\eps t^{-3/2+C\eps}+\eps t^{-2+C\eps}\lra{r-t}^2\lesssim \eps t^{-\gamma_1/2+C\eps}+\eps t^{-3/2+C\eps}\lra{r-t}^{3/2}\\
&\lesssim  \eps t^{-\gamma_1/2+C\eps}+\eps t^{-1/2+C\eps}K(\gamma_1,t,r).} Here we use $\lra{r-t}\lesssim t$ whenever $r\sim t$. If $\gamma_1= 3$, then \fm{&\eps t^{-\gamma_1/2+C\eps}+\eps t^{-1/2+C\eps}K(\gamma_1,t,t)+\eps t^{-2+C\eps}\lra{r-t}^2\\
&\lesssim \eps t^{-\gamma_1/2+C\eps}+\eps t^{-3/2+C\eps}\ln t+\eps t^{-3/2+C\eps}\lra{r-t}^{3/2}\lesssim\eps t^{-\gamma_1/2+C\eps}+\eps t^{-1/2+C\eps}K(\gamma_1,t,r).}
\end{proof}

\rm

We have finished the proof of Proposition \ref{propasymujohn}.
\end{proof}
\rm

\subsection{The 3D compressible Euler equations with no vorticity}\label{sec8.4}
Consider the  compressible Euler equations in three space dimensions without vorticity
\eq{\label{exmeqn4o}\left\{
\begin{array}{l}
\rho_t+ \nabla\cdot(\rho v)=0,\\[1em]
\rho(v_t+v\cdot\nabla v)+\nabla p=0,\\[1em]
\partial_av^{b}=\partial_bv^a,\qquad a,b\in\{1,2,3\}.
\end{array}
\right.}
Here the unknowns are the velocity $v=(v^a)_{a=1,2,3}:\R^{1+3}\to\R^{3}$ and the density $\rho:\R^{1+3}\to (0,\infty)$. The pressure $p$ in the equations is a function of the density $p=p(\rho)$ where the map $\rho\mapsto p(\rho)$ is given.  Here we  assume that $dp/d\rho>0$ everywhere. 

We now seek to formulate the 3D compressible Euler equations with no vorticity as a system of covariant quasilinear wave equations.  Since $\rho>0$ everywhere, we  can set $\varrho:=\ln(\rho/\rho_0)$ where $\rho_0>0$ is a constant background density. Set $c_{\rm s}:=\sqrt{dp/d\rho}$ and view it as a function of $\varrho$, so $c_{\rm s}=c_{\rm s}(\varrho)$ and $c_{\rm s}'(\varrho)=\frac{d}{d\varrho}c_{\rm s}(\varrho)$. By setting $B=\partial_t+v^a\partial_a$,  we reduce \eqref{exmeqn4o} to 
\eq{\label{exmeqn4o2}\left\{
\begin{array}{l}
B\varrho+ \partial_c v^c=0,\\[1em]
Bv+c_{\rm s}^2(\varrho)\nabla\varrho=0,\\[1em]
\partial_av^{b}=\partial_bv^a,\qquad a,b\in\{1,2,3\}.
\end{array}
\right.}
Now, the acoustical metric $g=(g_{\alpha\beta})$ is defined by 
\eq{g:=-dt\otimes dt+c_{\rm s}^{-2}\sum_{a=1}^3 (dx^a-v^adt)\otimes (dx^a-v^adt).}
The inverse acoustical metric $g^{-1}=(g^{\alpha\beta})$ is defined by
\eq{g^{-1}:=-B\otimes B+c_{\rm s}^2\sum_{a=1}^3 \partial_a\otimes \partial_a.}
It is obvious that both $g$ and $g^{-1}$ are functions of $(v,\varrho)$. Without loss of generality, we also assume that $c_{\rm s}(0)=1$, so $g$ and $g^{-1}$ are both the Minkowski metrics when $(v,\varrho)=0$.

With all these notations, it has been proved in \cite{MR4011696,MR4109292} that if $(v,\rho)$ solves \eqref{exmeqn4o}, then $(v,\varrho)$ solves the following system:
\eq{\label{exmeqn4}\left\{
\begin{array}{l}
\displaystyle \Box_g v^a=-(1+c_{\rm s}^{-1}c_{\rm s}')g^{\alpha\beta}\partial_\alpha\varrho\partial_\beta v^a,\qquad a=1,2,3;\\[1em]
\displaystyle \Box_g \varrho=-3c_{\rm s}^{-1}c_{\rm s}'g^{\alpha\beta}\partial_\alpha\varrho\partial_\beta \varrho+2\sum_{1\leq a<b\leq 3}(\partial_av^a\partial_bv^b-\partial_av^b\partial_bv^a).
\end{array}
\right.}
Here $\Box_g$ is a covariant wave operator defined by 
\eq{\Box_g\phi:=\frac{1}{\sqrt{|\det g|}}\partial_\alpha(\sqrt{|\det g|}g^{\alpha\beta}\partial_\beta \phi).}

In order to apply Theorem \ref{mthm} to \eqref{exmeqn4}, we first give the corresponding geometric reduced system. For simplicity, we write $u=(u^{ \alpha })_{\alpha=0,1,2,3}:=(\varrho,v)$. In other words, we have $u^0=\varrho$ and $u^a=v^a$, $a=1,2,3$. Besides, we have 
\fm{\Box_g\phi&=g^{\alpha\beta}\partial_\alpha\partial_\beta \phi+(\partial_\alpha g^{\alpha\beta}-\frac{1}{2}g_{\lambda\mu}g^{\alpha\beta} \partial_\alpha g^{\lambda\mu})\partial_\beta \phi\\
&=\Box\phi+g^{\alpha\beta}_\sigma u^\sigma\partial_\alpha\partial_\beta \phi+( g^{\alpha\beta}_\sigma \partial_\alpha u^\sigma-\frac{1}{2}m_{\lambda\mu}m^{\alpha\beta}  g^{\lambda\mu}_\sigma \partial_\alpha u^\sigma)\partial_\beta \phi+\text{cubic terms}.}
Here the $g^{\alpha\beta}_\sigma$ are linear coefficients in the Taylor expansion of $g^{\alpha\beta}$: $g^{\alpha\beta}=m^{\alpha\beta}+g^{\alpha\beta}_\sigma u^\sigma+O(|u|^2)$. In fact, since
\fm{g^{-1}&=-\partial_t\otimes\partial_t-\sum_{a=1}^3v^a\partial_t\otimes\partial_a-\sum_{a=1}^3v^a\partial_a\otimes\partial_t-\sum_{a,b=1}^3v^av^b\partial_a\otimes\partial_b+c_{\rm s}^2\sum_{a=1}^3\partial_a\otimes\partial_a\\
&=-\partial_t\otimes\partial_t-\sum_{a=1}^3v^a\partial_t\otimes\partial_a-\sum_{a=1}^3v^a\partial_a\otimes\partial_t+(1+2c_{\rm s}'(0)\varrho)\sum_{a=1}^3\partial_a\otimes\partial_a+\text{quadratic terms},}
we have
\fm{g^{0i}_i=g^{i0}_i=-1,\quad  g^{ii}_0=2c_{\rm s}'(0),\qquad i=1,2,3;\qquad g^{**}_*=0\quad \text{in all other cases}.}
It follows that 
\fm{\Box_g\phi&=\Box\phi-2u^a\partial_t\partial_a\phi+2c_{\rm s}'(0) u^0\Delta \phi-\partial_tu^a\partial_a\phi-\partial_au^a\partial_t\phi+2c_{\rm s}'(0)\partial_au^0\partial_a\phi-3c_{\rm s}'(0)m^{\alpha\beta}   \partial_\alpha u^0\partial_\beta \phi\\&\quad+\text{cubic terms}.}
As a result, we can rewrite \eqref{exmeqn4} as
\fm{ &\Box u^a-2u^b\partial_t\partial_b u^a+2c_{\rm s}'(0) u^0\Delta u^a\\
&=\partial_tu^b\partial_b u^a+\partial_bu^b\partial_tu^a-(2c_{\rm s}'(0)-1)\partial_t u^0\partial_t u^a-\partial_b u^0\partial_b u^a+\text{cubic terms},\qquad a=1,2,3;\\
& \Box u^0-2u^a\partial_t\partial_a u^0+2c_{\rm s}'(0) u^0\Delta  u^0\\
&=\partial_tu^a\partial_a u^0+\partial_au^a\partial_t u^0-2c_{\rm s}'(0)\partial_au^0\partial_a u^0+2\sum_{1\leq a<b\leq 3}(\partial_au^a\partial_bu^b-\partial_au^b\partial_bu^a)+\text{cubic terms}.
}
Also note that the cubic terms here do not contain $u^*\cdot u^*\cdot u^*$. Then, the corresponding geometric reduced system is
\eq{\label{exmeqn4rs}\left\{
\begin{array}{l}
\displaystyle \partial_s(\mu U_q^{a})=\frac{1}{2}(\omega_b U_q^b+c_{\rm s}'(0) U_q^0)\mu^2U_q^a,\qquad a=1,2,3;\\[1em]
\displaystyle \partial_s(\mu U_q^{0})=\frac{1}{2}(\omega_aU_q^a+c_{\rm s}'(0) U_q^0 )\mu^2U_q^0;\\[1em]
\displaystyle \partial_s \mu=\frac{1}{2}(\omega_a U^a_q+c_{\rm s}'(0) U_q^0)\mu^2.
\end{array}
\right.}
Since $\partial_\alpha u^\sigma\approx -\frac{\eps}{2r}\wh{\omega}_\alpha\mu \partial_qU^{\sigma}$ and since $c_{\rm s}^2(0)=1$, from \eqref{exmeqn4o2} we obtain  additional constraint equations for $U$:
\eq{\label{exmeqn4rs2}\left\{
\begin{array}{l}
\displaystyle U_q^0=\omega_cU_q^c;\\[1em]
\displaystyle U_q^a=\omega_aU_q^0,\qquad a=1,2,3;\\[1em]
\omega_aU_q^b=\omega_bU_q^a,\qquad a,b=1,2,3.
\end{array}
\right.}
Then, the reduced system \eqref{exmeqn4rs} can be reduced to 
\eq{\left\{
\begin{array}{l}
\displaystyle \partial_s(\mu U_q^{0})=\frac{1}{2}(1+c_{\rm s}'(0)  )(\mu U_q^0)^2;\\[1em]
\displaystyle \partial_s \mu=\frac{1}{2}(1+c_{\rm s}'(0))\mu^2 U_q^0.
\end{array}
\right.}
By solving this system explicitly, we obtain the following lemma.

\lem{\label{lemexm4.1} Fix $A\in C^\infty(\R\times\mathbb{S}^2)$, and fix two constants $\gamma_+>1,\gamma_->1$. Suppose that $(1+c_{\rm s}'(0))A(q,\omega)\geq 0$  for all $(q,\omega)$  and that $\partial_q^a\partial^b_\omega A=O_{a,c}(\lra{q}^{-a-\gamma_{\sgn(q)}})$ for each $a,c\geq 0$. If we set
\fm{U^{a}(s,q,\omega)&:=\int_{-\infty}^q \omega_aA(p,\omega)\ dp,\qquad
U^{0}(s,q,\omega):=\int_{-\infty}^q A(p,\omega)\ dp,}
\fm{\mu(s,q,\omega)&:=\frac{-2}{(1+c_{\rm s}'(0))A(q,\omega)s+1},}
then $(\mu,(U^\alpha))$ is a $(\gamma_+,\gamma_-)$-admissible global solution to the geometric reduced system \eqref{exmeqn4rs}.

Moreover, if $(1+c'_{\rm s}(0))A(q^0,\omega^0)>0$ for some $(q^0,\omega^0)\in\R\times\mathbb{S}^2$, then there exists a large constant $\wt{C}>1$ such that
\eq{\label{lemexm4.1:c}|U(s,q,\omega)|\gtrsim  |A(q^0,\omega^0)|^2,\quad \forall s>0,\ q>q^0,\ |\omega-\omega^0|\leq \wt{C}^{-1} |A(q^0,\omega^0)|.}
Both $\wt{C}$ and the implicit constant are independent of $(q^0,\omega^0)$. }
\begin{proof}The solution in this lemma has a similar form to that in Lemma \ref{lemexm3.1}. As a result, the proof is essentially the same as that of Lemma \ref{lemexm3.1}. The details are thus omitted. 
\end{proof}
\rm

By applying Theorem \ref{mthm} and Remark \ref{mthmrmk3}, we obtain the following corollary. 

\cor{For each function $A=A(q,\omega)$ as in Lemma \ref{lemexm4.1}, we let $(\mu,U)$ be the corresponding admissible global solution to the geometric reduced system. Then, for each $\eps\ll1$, there exists a global solution $u=(u^{\alpha})_{\alpha=0,1,2,3}$ for $t\geq 0$ to the system \eqref{exmeqn4} matching $(\mu,U)$ at infinite time of size $\eps$.

Moreover, the solution $u$  is nontrivial as long as $A\not\equiv  0$. If furthermore we have $c'_{\rm s}(0)\neq -1$, then the initial data of $u$ are not compactly supported.}
\rmk{\label{corexm4:rmk}\rm By following the proof in Remark \ref{corexm2:rmk}, in the case when $(1+c'_{\rm s}(0))A\not\equiv 0$, we can also prove \eq{\norm{u(t)}_{L^2(\R^3)}\gtrsim_\lambda|A(q^0,\omega^0)|^2 \cdot \eps t^{\lambda/2-C\eps},\qquad \forall t\geq e^{(1+\delta)/\eps},\ \eps\ll_{\lambda}1,\ \lambda\in(0,1).}
Here we choose $(q^0,\omega^0)$ so that $(1+c'_{\rm s}(0))A(q^0,\omega^0)> 0$.
}\rm

\bigskip
It remains to show that the original equations \eqref{exmeqn4o} and \eqref{exmeqn4o2}  hold. The idea of the proof is similar to that in Section \ref{sec8.3}. Recall that we write $u=(\varrho,v)$. For each $a,b=1,2,3$, we set
\eq{w_{ab}:=\partial_au^b-\partial_bu^a.}
For each $a=1,2,3$, we  set
\eq{w_{0a}=-w_{a0}:=c_{\rm s}^2\partial_au^0+Bu^a.}
We also set 
\eq{w_{00}:=Bu^0+\partial_cu^c.}
Recall that $B:=\partial_t+v^a\partial_a$. Note that these definitions are motivated by the original equations \eqref{exmeqn4o2}.

In order to prove $w\equiv 0$, we prove two lemmas. The first one  plays a similar role as Lemma \ref{lemexm3.2} in the previous subsection.

\lem{\label{lemexm4.2} We have
\fm{g^{\alpha\beta}\partial_\alpha\partial_\beta w_{\sigma\gamma}=h(u,\partial u)\cdot \partial w\cdot \partial u+h(u,\partial u)\cdot w\cdot \partial^2u+h(u,\partial u)\cdot w\cdot \partial u\cdot \partial u.}
See Lemma \ref{lemexm3.2} for the meaning of the right hand side.}
\rm

\bigskip

The proof of Lemma \ref{lemexm4.2} is much more complicated than that of Lemma \ref{lemexm3.2}. We will return to its proof at the end of this subsection.

The next lemma plays a similar role as  Lemma \ref{lemexm3.3} in the previous subsection.

\lem{\label{lemexm4.3} For $t\geq 1/\eps$, we have 
\fm{\sum_{\alpha,\beta=0,1,2,3}\norm{w_0^{1/2}\partial w_{\alpha\beta}(t)}_{L^2(\R^3)}\lesssim \eps t^{-1/2+C\eps}.}}
\begin{proof}
By Theorem \ref{mthm}, we have
\eq{\label{lemexm4.3f1}\norm{w_0^{1/2}\partial^l\partial_\alpha (u^\sigma-\eps r^{-1}U^\sigma)(t)}_{L^2(\{x/t\in[1-c',1+c'] \})}+\norm{w_0^{1/2}\partial^l\partial_\alpha u^\sigma(t)}_{L^2(\{|x|/t\notin[1-c',1+c']\})}\lesssim \eps t^{-1/2+C\eps}}
for each $\alpha,\sigma=0,1,2,3$ and $l\leq 1$. See Lemma \ref{lemexm4.1} for the definitions of $U^*$ and $\mu$. Thus, for $ |x|/t\in[1-c',1+c']$, we have
\fm{\partial_\alpha(\eps r^{-1}U^0,\eps r^{-1}U^a)&=(\eps r^{-1}Aq_\alpha,\eps \omega_ar^{-1}Aq_\alpha)+\eps S^{-2,1-\gamma_-}_-\cap\eps S^{-2,0}_+.
}
Here we use the assumptions in Definition \ref{def3.1} and the notation $S^{s,p}_\pm$ introduced in Definition \ref{c1defn1.5}.  As discussed in the proof of Lemma \ref{lemexm3.3},  we have $q\in S^{0,1}$, $\wh{\omega}_\alpha q_\beta-\wh{\omega}_\beta q_\alpha\in S^{-1,1}$, and $A\in S^{0,-\gamma_-}_-\cap S^{0,-\gamma_+}_+$.  It follows that
\fm{\partial_a(\eps r^{-1}U^b)-\partial_b(\eps r^{-1}U^a)&=\eps r^{-1}A(\omega_a q_b-\omega_bq_a)+\eps S^{-2,1-\gamma_-}_-\cap\eps S^{-2,0}_+=\eps S^{-2,1-\gamma_-}_-\cap\eps S^{-2,0}_+,\\
\partial_a(\eps r^{-1}U^0)+\partial_t(\eps r^{-1}U^a)&=\eps r^{-1}A(q_a+\omega_a q_t)+\eps S^{-2,1-\gamma_-}_-\cap\eps S^{-2,0}_+=\eps S^{-2,1-\gamma_-}_-\cap\eps S^{-2,0}_+,\\
\partial_t(\eps r^{-1}U^0)+\partial_c(\eps r^{-1}U^c)&=\eps r^{-1}A (q_t+q_r)+\eps S^{-2,1-\gamma_-}_-\cap\eps S^{-2,0}_+=\eps S^{-2,1-\gamma_-}_-\cap\eps S^{-2,0}_+.}
Whenever $ |x|/t\in[1-c',1+c']$, we have $ u^\sigma=\eps r^{-1}U^\sigma$, so
\fm{\partial w_{ab}&=\partial_a(\eps r^{-1}U^b)-\partial_b(\eps r^{-1}U^a),\\
\partial w_{0a}&=\partial[c_{\rm s}^2(u^0)\partial_a(\eps r^{-1}U^0)+\partial_t(\eps r^{-1}U^a)+\eps r^{-1}U^c\partial_c (\eps r^{-1}U^a)]\\
&=\partial[\partial_a(\eps r^{-1}U^0)+\partial_t(\eps r^{-1}U^a)]+O(\eps t^{-2+C\eps}(1_{r\geq t}+\lra{r-t}^{1-\gamma_-}1_{r<t})),\\
\partial w_{00}&=\partial [\partial_t(\eps r^{-1}U^0)+\eps r^{-1}U^c\partial_c (\eps r^{-1}U^0)+\partial_c(\eps r^{-1}U^c)]\\
&=\partial[\partial_t(\eps r^{-1}U^0)+\partial_c(\eps r^{-1}U^c)]+O(\eps t^{-2+C\eps}(1_{r\geq t}+\lra{r-t}^{1-\gamma_-}1_{r<t})).}
Thus, $\partial w=O(\eps t^{-2+C\eps}(1_{r\geq t}+\lra{r-t}^{1-\gamma_-}1_{r<t}))$ whenever $ |x|/t\in[1-c',1+c']$. Since $|\partial w|\lesssim |\partial^2 u|+|\partial u|^2$, we  conclude from \eqref{lemexm4.3f1} that $\norm{w_0^{1/2}\partial w(t)}_{L^2(\R^3)}\lesssim \eps t^{-1/2+C\eps}$ for $t\geq1/\eps$.
\end{proof}\rm

By combining Lemma \ref{lemexm4.2} and \eqref{lemexm4.3}, we can use an energy estimate to prove that $w\equiv 0$. In other words, we obtain a global solution to the 3D compressible Euler with no vorticity. The proof here is essentially the same as that of Lemma \ref{lemexm3.4}, so we skip the proof here. Our final conclusion is as follows.

\cor{\label{corexm4} There exists a family of nontrivial global solutions to the 3D compressible Euler equations with no vorticity for $t\geq 0$. If moreover we have $c_{\rm s}'(0)\neq -1$, then the initial data of all the solutions constructed here are not compactly supported.}\rm

\subsubsection{Proof of Lemma \ref{lemexm4.2}}It remains to prove Lemma \ref{lemexm4.2}. In the rest of this subsection,  we use $\mathcal{R}$ to denote an arbitrary sum of terms of the form $h(u,\partial u)\cdot\partial u\cdot  \partial w$, $h(u,\partial u)\cdot \partial^2 u\cdot  w$ or $h(u,\partial u)\cdot \partial  u\cdot\partial  u\cdot  w$. Our goal is to show $g^{\alpha\beta}\partial_\alpha\partial_\beta w_{\sigma\gamma}=\mcl{R}$.

We start with some useful identities related to the acoustical metric. Recall that $B=\partial_t+u^a\partial_a$, 
\fm{g_{00}=-1+c_{\rm s}^{-2}\sum_{a=1}^3(u^a)^2,\qquad g_{0a}=g_{a0}=-c_{\rm s}^{-2}u^a,\qquad g_{ab}=c_{\rm s}^{-2}\delta_{ab};}
\fm{g^{00}=-1,\qquad g^{0a}=g^{a0}=-u^a,\qquad g^{ab}=-u^au^b+c_{\rm s}^2\delta^{ab}.} Also recall that  
\fm{\Box_g\phi&=|\det g|^{-1/2}\partial_\alpha (g^{\alpha\beta}|\det g|^{1/2}\partial_\beta \phi)\\&=g^{\alpha\beta}\partial_\alpha\partial_\beta \phi+(\partial_\alpha g^{\alpha\beta}-\frac{1}{2}g_{\lambda\mu}g^{\alpha\beta} \partial_\alpha g^{\lambda\mu})\partial_\beta \phi.} It follows from the chain rule  that $\partial g^{**}=h(u)\partial u$, so $\Box_gw_{**}=g^{\alpha\beta}\partial_\alpha\partial_\beta w_{**}+\mathcal{R}$. It suffices to prove $\Box_gw_{**}=\mathcal{R}$.

We start our proof with the following lemma.\lem{\label{l4am0} We have\eq{\label{l4am0c}\Box_g\phi&=-BB\phi+c_{\rm s}^2\Delta \phi+2c_{\rm s}^{-1}c_{\rm s}'Bu^0\cdot B\phi-\partial_cu^c B\phi-c_{\rm s}^{-1}c_{\rm s}'g^{\alpha\beta}\partial_\alpha u^0\partial_\beta\phi.}
Here the right hand side is expressed in Cartesian coordinates.}
\begin{proof}
The identity \eqref{l4am0c} comes from   Luk--Speck \cite[Lemma 5.2]{MR4109292}. We  skip its proof, but we notice that we only need to use the definitions of the acoustic metric $g$, its inverse $g^{-1}$, and the vector field $B$. We do not need to assume that $u=(u^\alpha)$ satisfies the equations \eqref{exmeqn4o} or \eqref{exmeqn4}.
\end{proof}\rm

A direct corollary is that \fm{\Box_gw_{**}=-BBw_{**}+c_{\rm s}^2\Delta w_{**}+\mathcal{R}.}
Thus it suffices to compute $-BBw_{**}+c_{\rm s}^2\Delta w_{**}$. Moreover, by \eqref{l4am0c} and the equations \eqref{exmeqn4} for $u=(\varrho,v)$, we have
\fm{-BB u^a+c_{\rm s}^2\Delta u^a+2c_{\rm s}^{-1}c_{\rm s}'Bu^0\cdot Bu^a-\partial_cu^c Bu^a-c_{\rm s}^{-1}c_{\rm s}'g^{\alpha\beta}\partial_\alpha u^0\partial_\beta u^a&=-(1+c_{\rm s}^{-1}c_{\rm s}')g^{\alpha\beta}\partial_\alpha u^0\partial_\beta u^a}
for each $a=1,2,3$, and
\fm{&-BB u^0+c_{\rm s}^2\Delta u^0+2c_{\rm s}^{-1}c_{\rm s}'Bu^0\cdot Bu^0-\partial_cu^c Bu^0-c_{\rm s}^{-1}c_{\rm s}'g^{\alpha\beta}\partial_\alpha u^0\partial_\beta u^0\\
&=-3c_{\rm s}^{-1}c_{\rm s}'g^{\alpha\beta}\partial_\alpha u^0\partial_\beta u^0+2\sum_{1\leq a<b\leq 3}(\partial_au^a\partial_bu^b-\partial_au^b\partial_bu^a).}
Since $g^{-1}=-B\otimes B+c_{\rm s}^2\sum_{a=1}^3\partial_a\otimes\partial_a$, we have $g^{\alpha\beta}\partial_\alpha \phi\partial_\beta \psi=-B\phi\cdot B\psi+c_{\rm s}^2\sum_{a=1}^3\phi_a\psi_a$ for any pair of functions $(\phi,\psi)$. By the definition of $w$ and by easy computations, we obtain
\eq{\label{sec8eqnuaBB} -BB u^a+c_{\rm s}^2\Delta u^a=-2c_{\rm s}^{-1}c_{\rm s}'Bu^0\cdot Bu^a+w_{00} Bu^a-c_{\rm s}^2\delta^{bc}\partial_bu^0\partial_cu^a,\qquad a=1,2,3;}
\eq{\label{sec8eqnu0BB}-BB u^0+c_{\rm s}^2\Delta u^0=w_{00}\partial_cu^c -2c_{\rm s}c_{\rm s}'\delta^{ab} \partial_au^0\partial_bu^0 -\partial_au^b\partial_bu^a.}
Later, we will mainly use \eqref{sec8eqnuaBB} and \eqref{sec8eqnu0BB} as well as the definitions of $w$ to compute $-BBw+c_{\rm s}^2\Delta w$. For example, since $w_{00}=Bu^0+\partial_cu^c$, we  rewrite
\eq{\label{for8.36}&-BBw_{00}+c_{\rm s}^2\Delta w_{00}=-BBBu^0-BB\partial_cu^c+c_{\rm s}^2\Delta Bu^0+c_{\rm s}^2\Delta \partial_cu^c\\
&=B(-BBu^0+c_{\rm s}^2\Delta u^0)+\partial_c(-BBu^c+c_{\rm s}^2\Delta u^c)+[-BB+c_{\rm s}^2\Delta,\partial_c]u^c+[-BB+c_{\rm s}^2\Delta,B]u^0.}We now apply \eqref{sec8eqnuaBB} and \eqref{sec8eqnu0BB} to handle the first two terms on the right side. For the remaining terms, we compute the commutators directly. Following the same idea, since $w_{0a}=Bu^a+c_{\rm s}^2\partial_au^0$ and $w_{ab}=\partial_au^b-\partial_bu^a$, for $a,b=1,2,3$ we have
\eq{\label{for8.37}&-BBw_{0a}+c_{\rm s}^2\Delta w_{0a}=-BBBu^a-BB(c_{\rm s}^2\partial_au^0)+c_{\rm s}^2\Delta (Bu^a)+c_{\rm s}^2\Delta(c_{\rm s}^2\partial_au^0)\\
&=B(-BBu^a+c_{\rm s}^2\Delta u^a)+c_{\rm s}^2\partial_a(-BBu^0+c_{\rm s}^2\Delta u^0)+[-BB+c_{\rm s}^2\Delta,B]u^a+[-BB+c_{\rm s}^2\Delta,c_{\rm s}^2\partial_a]u^0,}
\eq{\label{for8.38}&-BBw_{ab}+c_{\rm s}^2\Delta w_{ab}=-BB\partial_au^b+BB\partial_bu^a+c_{\rm s}^2\Delta \partial_au^b-c_{\rm s}^2\Delta\partial_bu^a\\
&=\partial_a(-BBu^b+c_{\rm s}^2\Delta u^b)-\partial_b(-BBu^a+c_{\rm s}^2\Delta u^a)+[-BB+c_{\rm s}^2\Delta,\partial_a]u^b-[-BB+c_{\rm s}^2\Delta,\partial_b]u^a.}

Motivated by \eqref{for8.36}--\eqref{for8.38}, we now compute each of the following terms:
\fm{B(-BBu^*+c_{\rm s}^2\Delta u^*),\ \partial_a(-BBu^*+c_{\rm s}^2\Delta u^*);\\ [-BB+c_{\rm s}^2\Delta,B],\ [-BB+c_{\rm s}^2\Delta,\partial_a],\ [-BB+c_{\rm s}^2\Delta,c_{\rm s}^2\partial_a].}
The computations from now on will be complicated and tedious, but the main idea is clear. Whenever we see a second derivative of $u^*$, we try to convert it to $BBu^*$ or $\Delta u^*$ with the help of $w$. For example, from the definitions of $w$, we have identities
\fm{B\partial_au^0&=B(c_{\rm s}^{-2}w_{0a}-c_{\rm s}^{-2}Bu^a)=c_{\rm s}^{-2}BBu^a+c_{\rm s}^{-2}Bw_{0a}+\text{lower order terms}.}
If it is unclear how to change a term (such as $B\partial_au^b$ or $\partial_a\partial_bu^c$), we can temporarily keep them unchanged. 
Then we use \eqref{sec8eqnuaBB} and \eqref{sec8eqnu0BB} to handle terms involving $BBu^*$ or $\Delta u^*$. As for the terms not containing $\partial^2u$, we combine them together and try to express the sum in terms of $w$.

We first compute $B(-BBu^*+c_{\rm s}^2\Delta u^*)$ and $\partial_a(-BBu^*+c_{\rm s}^2\Delta u^*)$. The result is summarized in the next lemma.
\lem{\label{techlemm8.16}We have
\fm{B(-BBu^a+c_{\rm s}^2\Delta u^a)&=-2(c_{\rm s}^{-1}c_{\rm s}''-c_{\rm s}^{-2}(c_{\rm s}')^2)(Bu^0)^2 Bu^a-2c_{\rm s}^{-1}c_{\rm s}'(BBu^0  Bu^a+Bu^0BBu^a)\\
&\quad+ BBu^c\partial_cu^a+Bu^cB\partial_cu^a +\mcl{R},}
\fm{B(-BBu^0+c_{\rm s}^2\Delta u^0)&=-2( c_{\rm s} c_{\rm s}''-3(c_{\rm s}')^2)Bu^0\cdot \delta^{ab}\partial_au^0\partial_bu^0+4 c_{\rm s}^{-1}c_{\rm s}'\partial_au^0\cdot B Bu^a\\
&\quad-2\partial_au^bB\partial_bu^a+\mcl{R},}
\fm{\partial_b(-BBu^a+c_{\rm s}^2\Delta u^a)&=-2(c_{\rm s}^{-1}c_{\rm s}''-3c_{\rm s}^{-2}(c_{\rm s}')^2)\partial_bu^0 Bu^0Bu^a-2c_{\rm s}^{-1}c_{\rm s}' Bu^0\partial_bBu^a+\partial_b Bu^c \partial_cu^a\\
&\quad+ Bu^c \partial_b \partial_cu^a+2c_{\rm s}^{-3}c_{\rm s}'B Bu^b Bu^a-2c_{\rm s}^{-1}c_{\rm s}'\partial_bu^c\partial_cu^0Bu^a+\mcl{R},}
\fm{\partial_a(-BBu^0+c_{\rm s}^2\Delta u^0)&=-2\delta^{bc}(c_{\rm s}c_{\rm s}''-3(c_{\rm s}')^2) \partial_au^0\partial_bu^0\partial_cu^0+4c_{\rm s}^{-1}c_{\rm s}'\partial_bu^0 \partial_a Bu^b -2\partial_bu^c\partial_a\partial_cu^b+\mcl{R}.}
Here $\mcl{R}$ denotes an arbitrary function which can be written as $h(u)\cdot\partial w\cdot \partial u+h(u)\cdot w\cdot \partial u\cdot \partial u+h(u)\cdot w\cdot \partial^2u$.}
\begin{proof}
By \eqref{sec8eqnuaBB}, we have
\fm{&B(-BBu^a+c_{\rm s}^2\Delta u^a)=B(-2c_{\rm s}^{-1}c_{\rm s}'Bu^0\cdot Bu^a+w_{00} Bu^a-c_{\rm s}^2\delta^{cd}
\partial_cu^0\partial_du^a)\\
&=-2(c_{\rm s}^{-1}c_{\rm s}''-c_{\rm s}^{-2}(c_{\rm s}')^2)Bu^0 Bu^0Bu^a-2c_{\rm s}^{-1}c_{\rm s}'(BBu^0  Bu^a+Bu^0BBu^a)\\
&\quad +B(w_{00}Bu^a)-
B((\delta^{cd}w_{0c}-Bu^d)\partial_du^a)\\
&=-2(c_{\rm s}^{-1}c_{\rm s}''-c_{\rm s}^{-2}(c_{\rm s}')^2)Bu^0 Bu^0Bu^a-2c_{\rm s}^{-1}c_{\rm s}'(BBu^0  Bu^a+Bu^0BBu^a)\\
&\quad +BBu^c\partial_cu^a+Bu^cB\partial_cu^a+\mcl{R}.}
In the second identity, we apply the product rule, the chain rule, and the identity $c_{\rm s}^2\partial_au^0=w_{0a}-Bu^
a$. To get the last identity, we note that $B(w\cdot Bu)=Bw\cdot Bu+w B(B^\alpha)\partial_\alpha u+wB^{\alpha\beta}\partial_\alpha\partial_\beta u$. Since $B^0=1$ and $B^a=u^a$, we can write $B(w\cdot Bu)=\mcl{R}$.  Similarly, we have $B(w\cdot \partial u)=\mcl{R}$.

Next, by \eqref{sec8eqnu0BB}, we have
\fm{&B(-BBu^0+c_{\rm s}^2\Delta u^0)=B(w_{00}\partial_cu^c -2\delta^{ab}c_{\rm s}c_{\rm s}' \partial_au^0\partial_bu^0-\partial_au^b\partial_bu^a)\\
&=B( -2\delta^{ab}c_{\rm s}^{-3} c_{\rm s}'\cdot c_{\rm s}^2\partial_au^0\cdot c_{\rm s}^2\partial_bu^0)-2\partial_au^bB\partial_bu^a+\mcl{R}\\
&=-2(c_{\rm s}^{-3} c_{\rm s}')' Bu^0\cdot  c_{\rm s}^4\delta^{ab}\partial_au^0\partial_bu^0-4c_{\rm s}^{-1}c_{\rm s}'\partial_au^0\cdot B(\delta^{ab}w_{0b}-Bu^a)-2\partial_au^bB\partial_bu^a+\mcl{R}\\
&=-2( c_{\rm s} c_{\rm s}''-3(c_{\rm s}')^2)Bu^0\cdot \delta^{ab}\partial_au^0\partial_bu^0+4 c_{\rm s}^{-1}c_{\rm s}'\partial_au^0\cdot B Bu^a-2\partial_au^bB\partial_bu^a+\mcl{R}.}
Here we notice that $B(\partial_au^b\partial_bu^a)=2\partial_au^bB\partial_bu^a$ because the roles of $a,b$ in the sum can be interchanged here. Similarly, we have $B( -2\delta^{ab}c_{\rm s}^{-3} c_{\rm s}'\cdot c_{\rm s}^2\partial_au^0\cdot c_{\rm s}^2\partial_bu^0)=B( -2\delta^{ab}c_{\rm s}^{-3} c_{\rm s}')\cdot c_{\rm s}^2\partial_au^0\cdot c_{\rm s}^2\partial_bu^0-4 \delta^{ab}c_{\rm s}^{-3} c_{\rm s}'\cdot c_{\rm s}^2\partial_au^0\cdot B(c_{\rm s}^2\partial_bu^0)$.

Next, by \eqref{sec8eqnuaBB}, we have
\fm{&\partial_b(-BBu^a+c_{\rm s}^2\Delta u^a)=\partial_b(-2c_{\rm s}^{-1}c_{\rm s}'Bu^0\cdot Bu^a+w_{00} Bu^a-c_{\rm s}^2\delta^{cd}
\partial_cu^0\partial_du^a)\\
&=-2(c_{\rm s}^{-1}c_{\rm s}''-c_{\rm s}^{-2}(c_{\rm s}')^2)\partial_bu^0\cdot Bu^0Bu^a-2c_{\rm s}^{-1}c_{\rm s}'(\partial_bBu^0Bu^a+Bu^0\partial_bBu^a)\\
&\quad-\partial_b(
(\delta^{cd}w_{0c}-Bu^d)\partial_du^a)+\mcl{R}\\
&=-2(c_{\rm s}^{-1}c_{\rm s}''-c_{\rm s}^{-2}(c_{\rm s}')^2)\partial_bu^0\cdot Bu^0Bu^a-2c_{\rm s}^{-1}c_{\rm s}'(\partial_bBu^0Bu^a+Bu^0\partial_bBu^a)\\
&\quad+\partial_b Bu^c \partial_cu^a+ Bu^c \partial_b \partial_cu^a+\mcl{R}.}
At this moment it is unclear how to handle $\partial_b Bu^c,\partial_bBu^a,\partial_b \partial_cu^a$, so we keep them unchanged. However, we notice that 
\fm{&-2c_{\rm s}^{-1}c_{\rm s}'\partial_bBu^0Bu^a=-2c_{\rm s}^{-3}c_{\rm s}'B(c_{\rm s}^2\partial_bu^0)Bu^a+2c_{\rm s}^{-3}c_{\rm s}'[B,c_{\rm s}^2\partial_b]u^0Bu^a\\
&=-2c_{\rm s}^{-3}c_{\rm s}'B(w_{0b}-Bu^b)Bu^a+2c_{\rm s}^{-3}c_{\rm s}'( 2c_{\rm s}c_{\rm s}'Bu^0\partial_bu^0-c_{\rm s}^2\partial_bu^c\partial_cu^0)Bu^a\\
&=2c_{\rm s}^{-3}c_{\rm s}'B Bu^b Bu^a+4c_{\rm s}^{-2}(c_{\rm s}')^2Bu^0\partial_bu^0Bu^a-2c_{\rm s}^{-1}c_{\rm s}'\partial_bu^c\partial_cu^0Bu^a+\mcl{R}.}
In summary,
\fm{\partial_b(-BBu^a+c_{\rm s}^2\Delta u^a)&=-2(c_{\rm s}^{-1}c_{\rm s}''-3c_{\rm s}^{-2}(c_{\rm s}')^2)\partial_bu^0\cdot Bu^0Bu^a-2c_{\rm s}^{-1}c_{\rm s}' Bu^0\partial_bBu^a+\partial_b Bu^c \partial_cu^a\\
&\quad+ Bu^c \partial_b \partial_cu^a+2c_{\rm s}^{-3}c_{\rm s}'B Bu^b Bu^a-2c_{\rm s}^{-1}c_{\rm s}'\partial_bu^c\partial_cu^0Bu^a+\mcl{R}.}

Finally, by \eqref{sec8eqnu0BB} we have 
\fm{&\partial_a(-BBu^0+c_{\rm s}^2\Delta u^0)=\partial_a(w_{00}\partial_cu^c -2\delta^{bc}c_{\rm s}c_{\rm s}' \partial_bu^0\partial_cu^0-\partial_bu^c\partial_cu^b)\\
&=\partial_a(-2\delta^{bc}c_{\rm s}^{-3}c_{\rm s}' \cdot c_{\rm s}^2\partial_bu^0\cdot c_{\rm s}^2\partial_cu^0)-2\partial_bu^c\partial_a\partial_cu^b+\mcl{R}
\\
&=-2\delta^{bc}(c_{\rm s}c_{\rm s}''-3(c_{\rm s}')^2) \partial_au^0\partial_bu^0\partial_cu^0-4c_{\rm s}^{-1}c_{\rm s}'\partial_bu^0\cdot \partial_a( \delta^{bc}w_{0c}-Bu^b)-2\partial_bu^c\partial_a\partial_cu^b+\mcl{R}\\
&=-2\delta^{bc}(c_{\rm s}c_{\rm s}''-3(c_{\rm s}')^2) \partial_au^0\partial_bu^0\partial_cu^0+4c_{\rm s}^{-1}c_{\rm s}'\partial_bu^0  \partial_a Bu^b -2\partial_bu^c\partial_a\partial_cu^b+\mcl{R}.}
\end{proof}\rm

Next, we consider the commutator terms.
\lem{\label{techlemm8.17}Given a function $\phi$, we have \fm{\ [-BB+c_{\rm s}^2\Delta,B]\phi&=c_{\rm s}^2\Delta u^c\partial_c\phi+\delta^{ab}c_{\rm s}^2\partial_bu^c\partial_a\partial_c\phi+\delta^{ab}c_{\rm s}^2 \partial_au^c\partial_c\partial_b\phi-2c_{\rm s}c'_{\rm s} Bu^0\Delta \phi,\\ [-BB+c_{\rm s}^2\Delta,\partial_a]\phi&=B(\partial_au^b\partial_b\phi)+\partial_au^b\partial_bB\phi-2c_{\rm s}c_{\rm s}'\partial_au^0\Delta\phi,\\
[-BB+c_{\rm s}^2\Delta,c_{\rm s}^2\partial_a]\phi&=c_{\rm s}^2B(\partial_au^b\partial_b\phi)+c_{\rm s}^2\partial_au^b\partial_bB\phi-2c_{\rm s}^3c_{\rm s}'\partial_au^0\Delta\phi -B(2c_{\rm s}c_{\rm s}' Bu^0\partial_a\phi)\\
&\quad-2c_{\rm s}c_{\rm s}' Bu^0B\partial_a\phi+\delta^{cd}c_{\rm s}^2\partial_c(2c_{\rm s}c_{\rm s}' \partial_du^0\partial_a\phi)+2\delta^{cd}c_{\rm s}^3c_{\rm s}' \partial_cu^0\partial_d\partial_a\phi.}}
\begin{proof}
We start with $[-BB+c_{\rm s}^2\Delta,B]=[c_{\rm s}^2\Delta,B]$. We have
\fm{\ [c_{\rm s}^2\Delta,B]\phi&=\delta^{ab}[c_{\rm s}^2\partial_a\partial_b,B]\phi=\delta^{ab}c_{\rm s}^2\partial_a[\partial_b,B]\phi+\delta^{ab}c_{\rm s}^2[\partial_a,B]\partial_b\phi+[c_{\rm s}^2,B]\Delta\phi\\
&=\delta^{ab}c_{\rm s}^2\partial_a(\partial_bu^c\partial_c\phi)+\delta^{ab}c_{\rm s}^2 \partial_au^c\partial_c\partial_b\phi-2c_{\rm s}c'_{\rm s} Bu^0\Delta \phi\\
&=c_{\rm s}^2\Delta u^c\partial_c\phi+\delta^{ab}c_{\rm s}^2\partial_bu^c\partial_a\partial_c\phi+\delta^{ab}c_{\rm s}^2 \partial_au^c\partial_c\partial_b\phi-2c_{\rm s}c'_{\rm s} Bu^0\Delta \phi.}

Next, we compute $[-BB+c_{\rm s}^2\Delta,\partial_a]$. We have
\fm{\ [-BB+c_{\rm s}^2\Delta,\partial_a]\phi&=-B[B,\partial_a]\phi-[B,\partial_a]B\phi+[c_{\rm s}^2,\partial_a]\Delta\phi\\
&=B(\partial_au^b\partial_b\phi)+\partial_au^b\partial_bB\phi-2c_{\rm s}c_{\rm s}'\partial_au^0\Delta\phi.}

Finally, we compute $[-BB+c_{\rm s}^2\Delta,c_{\rm s}^2\partial_a]=[-BB+c_{\rm s}^2\Delta,c_{\rm s}^2]\partial_a+c_{\rm s}^2[-BB+c_{\rm s}^2\Delta,\partial_a]$. It remains to compute $[-BB+c_{\rm s}^2\Delta,c_{\rm s}^2]\partial_a$. Note that
\fm{&[-BB+c_{\rm s}^2\Delta,c_{\rm s}^2]\partial_a\phi=-B[B,c_{\rm s}^2]\partial_a\phi-[B,c_{\rm s}^2]B\partial_a\phi+\delta^{cd}c_{\rm s}^2\partial_c[\partial_d,c_{\rm s}^2]\partial_a\phi+\delta^{cd}c_{\rm s}^2[\partial_c,c_{\rm s}^2]\partial_d\partial_a\phi\\
&=-B(2c_{\rm s}c_{\rm s}' Bu^0\partial_a\phi)-2c_{\rm s}c_{\rm s}' Bu^0B\partial_a\phi+\delta^{cd}c_{\rm s}^2\partial_c(2c_{\rm s}c_{\rm s}' \partial_du^0\partial_a\phi)+2\delta^{cd}c_{\rm s}^3c_{\rm s}' \partial_cu^0\partial_d\partial_a\phi.}
\end{proof}\rm 

We now show $-BBw+c_{\rm s}^2\Delta w=\mcl{R}$ by combining both the previous lemmas and \eqref{for8.36}--\eqref{for8.38}. Let us start with $w_{00}$. 

\lem{We have $-BBw_{00}+c_{\rm s}^2\Delta w_{00}=\mcl{R}$.}
\begin{proof}
By \eqref{for8.36}, we write $-BBw_{00}+c_{\rm s}^2\Delta w_{00}$ as  
\fm{B(-BBu^0+c_{\rm s}^2\Delta u^0)+\partial_c(-BBu^c+c_{\rm s}^2\Delta u^c)+[-BB+c_{\rm s}^2\Delta,\partial_c]u^c+[-BB+c_{\rm s}^2\Delta,B]u^0.}
These four terms equal, respectively,
\fm{\underline{-2( c_{\rm s} c_{\rm s}''-3(c_{\rm s}')^2)Bu^0\cdot \delta^{ab}\partial_au^0\partial_bu^0}\dashuline{+4 c_{\rm s}^{-1}c_{\rm s}'\partial_au^0\cdot B Bu^a}\uuline{-2\partial_au^bB\partial_bu^a}+\mcl{R},}
\fm{&\underline{-2(c_{\rm s}^{-1}c_{\rm s}''-3c_{\rm s}^{-2}(c_{\rm s}')^2)\partial_cu^0 Bu^0Bu^c-2c_{\rm s}^{-1}c_{\rm s}' Bu^0\partial_cBu^c}+\partial_b Bu^c \partial_cu^b\\
&\quad\dashuline{+ Bu^c \partial_b \partial_cu^b+2c_{\rm s}^{-3}c_{\rm s}'\delta_{ab}B Bu^a Bu^b-2c_{\rm s}^{-1}c_{\rm s}'\partial_bu^c\partial_cu^0Bu^b}+\mcl{R},}
\fm{\uuline{B(\partial_cu^b\partial_b u^c)}+\partial_cu^b\partial_bBu^c\dashuline{-2c_{\rm s}c_{\rm s}'\partial_cu^0\Delta u^c},}
\fm{\dashuline{c_{\rm s}^2\Delta u^c\partial_c u^0}+\delta^{ab}c_{\rm s}^2\partial_bu^c\partial_a\partial_cu^0+\delta^{ab}c_{\rm s}^2 \partial_au^c\partial_c\partial_bu^0\underline{-2c_{\rm s}c'_{\rm s} Bu^0\Delta u^0}.}
Here we apply the previous two lemmas. It is easy to see that the terms with a double underline cancel with each other. Moreover, the sum of the terms with a factor $Bu^0$ (i.e. the underlined terms) is
\fm{&Bu^0\kh{-2( c_{\rm s} c_{\rm s}''-3(c_{\rm s}')^2)\delta^{ab}\partial_au^0\partial_bu^0-2(c_{\rm s}^{-1}c_{\rm s}''-3c_{\rm s}^{-2}(c_{\rm s}')^2)\partial_cu^0 Bu^c-2c_{\rm s}^{-1}c_{\rm s}' \partial_cBu^c-2c_{\rm s}c'_{\rm s} \Delta u^0}\\
&=Bu^0\kh{-2(c_{\rm s}^{-1}c_{\rm s}''-3c_{\rm s}^{-2}(c_{\rm s}')^2)\delta^{ab}w_{0b}\partial_au^0 -2c_{\rm s}^{-1}c_{\rm s}' B\partial_cu^c-2c_{\rm s}^{-1}c_{\rm s}' \partial_cu^b\partial_bu^c-2c_{\rm s}c'_{\rm s} \Delta u^0}\\
&=2c_{\rm s}^{-1}c_{\rm s}'Bu^0\kh{ BBu^0-c_{\rm s}^2 \Delta u^0- \partial_cu^b\partial_bu^c}+\mcl{R}\\
&=2c_{\rm s}^{-1}c_{\rm s}'Bu^0\kh{-w_{00}\partial_cu^c+2c_{\rm s}c_{\rm s}'\delta^{ab}\partial_au^0\partial_bu^0+\partial_au^b\partial_bu^a- \partial_cu^b\partial_bu^c}+\mcl{R}\\
&=4(c_{\rm s}')^2Bu^0 \delta^{ab}\partial_au^0\partial_bu^0+\mcl{R}.}
Here we keep using the definitions of $w$. For example,  we combine the first two terms on the first row as follows: $c_{\rm s}^2\delta^{ab}\partial_au^0\partial_bu^0+\partial_cu^0Bu^c=\partial_au^0(Bu^a+c_{\rm s}^2\delta^{ab}\partial_bu^0)=\partial_au^0 w_{0a}$. Moreover, when we compute $\partial_cBu^c$, we use $\partial_cBu^c=B\partial_cu^c+[\partial_c,B]u^c=B(w_{00}-Bu^0)+[\partial_c,B]u^c$. In the second last identity, we apply \eqref{sec8eqnu0BB}.

Next, the sum of the terms with a factor $\partial_a u^0$ or $Bu^a$ (i.e. the terms with a dashed line below) is 
\fm{&4 c_{\rm s}^{-1}c_{\rm s}'\partial_au^0 B Bu^a+ Bu^c \partial_b \partial_cu^b+2c_{\rm s}^{-3}c_{\rm s}'\delta_{ab}B Bu^a Bu^b\\
&-2c_{\rm s}^{-1}c_{\rm s}'\partial_bu^c\partial_cu^0Bu^b-2c_{\rm s}c_{\rm s}'\partial_cu^0\Delta u^c+c_{\rm s}^2\Delta u^c\partial_c u^0\\
&=4 c_{\rm s}^{-3}c_{\rm s}'(w_{0a}-\delta_{ab}Bu^b)  B Bu^a+ Bu^c \partial_c(w_{00}-Bu^0)+2c_{\rm s}^{-3}c_{\rm s}'\delta_{ab}B Bu^a Bu^b\\
&\quad-2c_{\rm s}^{-1}c_{\rm s}'\partial_bu^c\partial_cu^0Bu^b+(1-2c_{\rm s}^{-1}c_{\rm s}')\Delta u^c(w_{0c}-\delta_{cd}Bu^d)\\
&= -Bu^c \partial_cBu^0-2c_{\rm s}^{-3}c_{\rm s}'\delta_{ab}B Bu^a Bu^b -2c_{\rm s}^{-1}c_{\rm s}'\partial_bu^c\partial_cu^0Bu^b-(1-2c_{\rm s}^{-1}c_{\rm s}')\Delta u^c \delta_{cd}Bu^d +\mcl{R}.}
To continue, we notice that
\fm{&-Bu^c \partial_cBu^0=-Bu^c B\partial_cu^0-Bu^c\partial_cu^b\partial_bu^0 =-Bu^c B(c_{\rm s}^{-2}(w_{0c}-\delta_{bc}Bu^b))-Bu^c\partial_cu^b\partial_bu^0\\
&=-2 \delta_{bc} c_{\rm s}^{-3}c_{\rm s}' Bu^cBu^0Bu^b+ c_{\rm s}^{-2}\delta_{bc}BBu^bBu^c -Bu^c\partial_cu^b\partial_bu^0+\mcl{R}.}
As a result, the sum of the terms with a dashed line below is
\fm{&-2 \delta_{bc} c_{\rm s}^{-3}c_{\rm s}' Bu^cBu^0Bu^b+ (c_{\rm s}^{-2}-2c_{\rm s}^{-3}c_{\rm s}')\delta_{bc}(BBu^b-c_{\rm s}^2\Delta u^b)Bu^c \\&-(1+2c_{\rm s}^{-1} c_{\rm s}')Bu^b\partial_bu^c\partial_cu^0+\mcl{R}\\
&=-2 \delta_{bc} c_{\rm s}^{-3}c_{\rm s}' Bu^cBu^0Bu^b+ (c_{\rm s}^{-2}-2c_{\rm s}^{-3}c_{\rm s}')\delta_{bc}(2c_{\rm s}^{-1}c_{\rm s}'Bu^0Bu^b+c_{\rm s}^2\delta^{de}\partial_du^0\partial_eu^b)Bu^c \\&\quad-(1+2c_{\rm s}^{-1} c_{\rm s}')Bu^b\partial_bu^c\partial_cu^0+\mcl{R}\\
&=-4c_{\rm s}^{-4}(c_{\rm s}')^2 \delta_{bc}  Bu^cBu^0Bu^b+ (1-2c_{\rm s}^{-1}c_{\rm s}')\delta_{bc} \delta^{de}\partial_du^0\partial_eu^bBu^c -(1+2c_{\rm s}^{-1} c_{\rm s}')Bu^b\partial_bu^c\partial_cu^0+\mcl{R}\\
&=-4c_{\rm s}^{-4}(c_{\rm s}')^2 \delta_{bc}  Bu^cBu^0Bu^b-4c_{\rm s}^{-1}c_{\rm s}'\partial_cu^0\partial_bu^cBu^b +\mcl{R}.}
Here we use \eqref{sec8eqnuaBB}. In the last identity, we use $\partial_eu^b=\partial_bu^e+w_{eb}$.

Then, we consider the remaining terms. The sum is 
\fm{&2\partial_bBu^c\partial_cu^b+2\delta^{ab}c_{\rm s}^2\partial_bu^c\partial_a\partial_cu^0=2\partial_bBu^c\partial_cu^b+2c_{\rm s}^2\partial_bu^c\partial_c(c_{\rm s}^{-2}w_{0a}\delta^{ab}-c_{\rm s}^{-2}Bu^b)\\
&=2\partial_bBu^c\partial_cu^b+4c_{\rm s}^{-1}c_{\rm s}'\partial_bu^c\partial_cu^0Bu^b-2\partial_bu^c\partial_cBu^b+\mcl{R}=4c_{\rm s}^{-1}c_{\rm s}'\partial_bu^c\partial_cu^0Bu^b+\mcl{R}.}

We now add all the terms together. In conclusion, $-BBw_{00}+c_{\rm s}^2\Delta w_{00}$ equals
\fm{&4(c_{\rm s}')^2Bu^0\delta^{ab}\partial_au^0\partial_bu^0-4c_{\rm s}^{-4}(c_{\rm s}')^2 \delta_{bc}  Bu^cBu^0Bu^b\\
&+4c_{\rm s}^{-1}c_{\rm s}'\partial_bu^c\partial_cu^0Bu^b - 4c_{\rm s}^{-1} c_{\rm s}'Bu^b\partial_bu^c\partial_cu^0+\mcl{R}\\
&=4c_{\rm s}^{-2}(c_{\rm s}')^2Bu^0(\delta^{ab}w_{0a}-Bu^b)\partial_bu^0-4c_{\rm s}^{-4}(c_{\rm s}')^2 \delta_{ab}  Bu^aBu^0Bu^b+0+\mcl{R}\\
&=-4c_{\rm s}^{-4}(c_{\rm s}')^2Bu^0Bu^b (w_{0b}-\delta_{ab}Bu^a)-4c_{\rm s}^{-4}(c_{\rm s}')^2 \delta_{ab}  Bu^aBu^0Bu^b+0+\mcl{R}=\mcl{R}.}
This finishes the proof. Again, here we use $w_{0a}=Bu^a+c_{\rm s}^2\partial_au^0$.
\end{proof}\rm

Next, we consider $w_{0a}$.
\lem{We have $-BBw_{0a}+c_{\rm s}^2\Delta w_{0a}=\mcl{R}$ for $a=1,2,3$.}
\begin{proof}
By \eqref{for8.37}, we write $-BBw_{0a}+c_{\rm s}^2\Delta w_{0a}$ as
\fm{B(-BBu^a+c_{\rm s}^2\Delta u^a)+c_{\rm s}^2\partial_a(-BBu^0+c_{\rm s}^2\Delta u^0)+[-BB+c_{\rm s}^2\Delta,B]u^a+[-BB+c_{\rm s}^2\Delta,c_{\rm s}^2\partial_a]u^0.}
By Lemmas \ref{techlemm8.16} and \ref{techlemm8.17}, these four terms equal, respectively,
\fm{&-2(c_{\rm s}^{-1}c_{\rm s}''-c_{\rm s}^{-2}(c_{\rm s}')^2)(Bu^0)^2 Bu^a\underline{-2c_{\rm s}^{-1}c_{\rm s}' BBu^0  Bu^a} \uuline{-2c_{\rm s}^{-1}c_{\rm s}' Bu^0BBu^a}\\
&+ \uuline{BBu^c\partial_cu^a}+\uwave{Bu^cB\partial_cu^a} +\mcl{R},}
\fm{&-2\delta^{bc}(c_{\rm s}^3c_{\rm s}''-3c_{\rm s}^2(c_{\rm s}')^2) \partial_au^0\partial_bu^0\partial_cu^0+\uwave{4 c_{\rm s}c_{\rm s}'\partial_bu^0 \partial_a Bu^b -2c_{\rm s}^2\partial_bu^c\partial_a\partial_cu^b}+\mcl{R},}
\fm{\uuline{c_{\rm s}^2\Delta u^c\partial_c u^a}+\uwave{2\delta^{bd}c_{\rm s}^2\partial_bu^c\partial_d\partial_cu^a}\uuline{-2c_{\rm s}c'_{\rm s} Bu^0\Delta u^a},}
\fm{&\uwave{c_{\rm s}^2B\partial_au^b\partial_bu^0}+\uuline{c_{\rm s}^2\partial_au^bB\partial_bu^0+c_{\rm s}^2\partial_au^b\partial_bBu^0}\underline{-2c_{\rm s}^3c_{\rm s}'\partial_au^0\Delta u^0 -2c_{\rm s}c_{\rm s}' BBu^0\partial_au^0}\\
&-2Bu^0 (c_{\rm s}c_{\rm s}')'Bu^0 \partial_au^0\uuline{-4c_{\rm s}c_{\rm s}' Bu^0B\partial_au^0}+\underline{2c_{\rm s}^3c_{\rm s}' \Delta u^0\partial_au^0}\\&+2\delta^{cd}c_{\rm s}^2\partial_du^0(c_{\rm s}c_{\rm s}')'\partial_cu^0\partial_au^0+\uwave{4c_{\rm s}^3c_{\rm s}' \delta^{cd}\partial_du^0\partial_c\partial_au^0}.}

Let us start with  the underlined terms:
\fm{&-2c_{\rm s}^{-1}c_{\rm s}'BBu^0Bu^a-2c_{\rm s}^3c_{\rm s}'\partial_au^0\Delta u^0-2c_{\rm s}c_{\rm s}' BBu^0\partial_au^0+2c_{\rm s}^3c_{\rm s}' \Delta u^0\partial_au^0\\
&=-2c_{\rm s}^{-1}c_{\rm s}'BBu^0(Bu^a+c_{\rm s}^2\partial_au^0)+0=-2c_{\rm s}^{-1}c_{\rm s}'BBu^0 \cdot w_{0a}=\mcl{R}.}
In the last step, we use the product rule to get $BBu^0=h(u)\cdot \partial^2u+h(u)\cdot \partial u\cdot \partial u$.

We move on to the terms with a double underline:
\fm{&-2c_{\rm s}^{-1}c_{\rm s}' Bu^0BBu^a+BBu^c\partial_cu^a+c_{\rm s}^2 \Delta u^c\partial_cu^a-2c_{\rm s}c_{\rm s}'Bu^0\Delta u^a\\
&\quad+c_{\rm s}^2\partial_au^bB\partial_bu^0+c_{\rm s}^2\partial_au^b\partial_bBu^0-4c_{\rm s}c_{\rm s}' Bu^0B\partial_au^0.}
Note that $\partial_bBu^0=B\partial_bu^0+\partial_bu^c\partial_cu^0$ and
\fm{B\partial_bu^0&=B(c_{\rm s}^{-2}(w_{0b}-Bu^b))=-2c_{\rm s}^{-3}c_{\rm s}'Bu^0(w_{0b}-Bu^b)+c_{\rm s}^{-2}(Bw_{0b}-BBu^b).}
Thus, the sum above is equal to
\fm{&-2c_{\rm s}^{-1}c_{\rm s}' Bu^0BBu^a+BBu^c\partial_cu^a+c_{\rm s}^2 \Delta u^c\partial_cu^a-2c_{\rm s}c_{\rm s}'Bu^0\Delta u^a\\
&+2c_{\rm s}^2\partial_au^bB\partial_bu^0-4c_{\rm s}c_{\rm s}' Bu^0B\partial_au^0+c_{\rm s}^2\partial_au^b\partial_bu^c\partial_cu^0\\
&=-2c_{\rm s}^{-1}c_{\rm s}' Bu^0BBu^a+BBu^c\partial_cu^a+c_{\rm s}^2 \Delta u^c\partial_cu^a-2c_{\rm s}c_{\rm s}'Bu^0\Delta u^a\\
&\quad+4c_{\rm s}^{-1}c_{\rm s}'\delta_{bc}\partial_au^bBu^0 Bu^c -2\delta_{bc}\partial_au^bBBu^c-8c_{\rm s}^{-2}(c_{\rm s}')^2(Bu^0)^2 Bu^a \\
&\quad+4c_{\rm s}^{-1}c_{\rm s}' Bu^0 BBu^a +c_{\rm s}^2\partial_au^b\partial_bu^c\partial_cu^0+\mcl{R}\\
&=2c_{\rm s}^{-1}c_{\rm s}' Bu^0(BBu^a-c_{\rm s}^2\Delta u^a)+2BBu^cw_{ca}+(-BBu^c+c_{\rm s}^2 \Delta u^c)\partial_cu^a\\
&\quad+4c_{\rm s}^{-1}c_{\rm s}'\delta_{bc}\partial_au^bBu^0 Bu^c -8c_{\rm s}^{-2}(c_{\rm s}')^2(Bu^0)^2 Bu^a +c_{\rm s}^2\partial_au^b\partial_bu^c\partial_cu^0+\mcl{R}\\
&=2c_{\rm s}^{-1}c_{\rm s}' Bu^0(2c_{\rm s}^{-1}c_{\rm s}'Bu^0Bu^a+c_{\rm s}^2\delta^{bc}\partial_bu^0\partial_cu^a)+\partial_cu^a(-2c_{\rm s}^{-1}c_{\rm s}'Bu^0Bu^c-c_{\rm s}^2\delta^{bd}\partial_bu^0\partial_du^c)\\
&\quad+4c_{\rm s}^{-1}c_{\rm s}'\delta_{bc}\partial_au^bBu^0 Bu^c -8c_{\rm s}^{-2}(c_{\rm s}')^2(Bu^0)^2 Bu^a +c_{\rm s}^2\partial_au^b\partial_bu^c\partial_cu^0+\mcl{R}\\
&=4c_{\rm s}^{-1}c_{\rm s}'\delta_{bc}\partial_au^bBu^0 Bu^c +2c_{\rm s} c_{\rm s}' Bu^0\delta^{bc}\partial_bu^0\partial_cu^a-2c_{\rm s}^{-1}c_{\rm s}'\partial_cu^a Bu^0Bu^c\\
&\quad-c_{\rm s}^2\delta^{bd}\partial_bu^0\partial_du^c\partial_cu^a-4c_{\rm s}^{-2}(c_{\rm s}')^2(Bu^0)^2 Bu^a +c_{\rm s}^2\partial_au^b\partial_bu^c\partial_cu^0+\mcl{R}\\
&=-4c_{\rm s}^{-2}(c_{\rm s}')^2(Bu^0)^2 Bu^a +\mcl{R}.}
In the first identity, we use the formula for $B\partial_bu^0$ above and put all the terms involving $w$ in $\mcl{R}$. In the third one, we use \eqref{sec8eqnuaBB}. To obtain the last line, we notice that $2c_{\rm s} c_{\rm s}' Bu^0\delta^{bc}\partial_bu^0\partial_cu^a=-2c_{\rm s}^{-1} c_{\rm s}' Bu^0 Bu^b\partial_bu^a+\mcl{R}$, $-c_{\rm s}^2\delta^{bd}\partial_bu^0\partial_du^c\partial_cu^a=-c_{\rm s}^2\partial_bu^0\partial_cu^b\partial_au^c+\mcl{R}$. Here we use the definitions of $w_{0a}$ and $w_{ab}$.

Next, we consider the terms with a tilde below:
\fm{&(Bu^cB\partial_cu^a+c_{\rm s}^2B\partial_au^b\partial_bu^0) +(-2c_{\rm s}^2\partial_bu^c\partial_a\partial_cu^b+2\delta^{bd}c_{\rm s}^2\partial_bu^c\partial_d\partial_cu^a)\\&+(4 c_{\rm s}c_{\rm s}'\partial_bu^0 \partial_a Bu^b+4c_{\rm s}^3c_{\rm s}' \delta^{cd}\partial_du^0\partial_c\partial_au^0).}
In the first bracket, we have $c_{\rm s}^2B\partial_au^b\partial_bu^0=B\partial_au^b(w_{0b}-\delta_{bc}Bu^b)+\mcl{R}=-\delta_{bc}B\partial_au^bBu^c+\mcl{R}=-B\partial_cu^aBu^c+\mcl{R}$. Similarly, we have  $2\delta^{bd}c_{\rm s}^2\partial_bu^c\partial_d\partial_cu^a=2c_{\rm s}^2\partial_bu^c\partial_a\partial_cu^b+\mcl{R}$ in the second one. For the third one, we have
\fm{&4 c_{\rm s}c_{\rm s}'\partial_bu^0 \partial_a Bu^b+4c_{\rm s}^3c_{\rm s}' \delta^{cd}\partial_du^0\partial_c\partial_au^0=4 c_{\rm s}c_{\rm s}'\partial_bu^0 \partial_a Bu^b+4c_{\rm s}^3c_{\rm s}' \partial_du^0\partial_a(c_{\rm s}^{-2}(\delta^{cd}w_{0c}-Bu^d))\\
&=4 c_{\rm s}c_{\rm s}'\partial_bu^0 \partial_a Bu^b-4c_{\rm s}c_{\rm s}' \partial_du^0\partial_a Bu^d +8(c_{\rm s}')^2 \partial_du^0 \partial_au^0Bu^d+\mcl{R}=8(c_{\rm s}')^2 \partial_du^0 \partial_au^0Bu^d+\mcl{R}.}
In summary, the sum of the terms with a tilde below is $8(c_{\rm s}')^2 \partial_bu^0 \partial_au^0Bu^b+\mcl{R}$.

Finally, we combine the sum of the underlined terms, of the terms with a double underline, of the terms with a tilde, and of the remaining unhighlighted terms. We have $-BBw_{0a}+c_{\rm s}^2\Delta w_{0a}$ equals
\fm{&-4c_{\rm s}^{-2}(c_{\rm s}')^2(Bu^0)^2 Bu^a+8(c_{\rm s}')^2 \partial_bu^0 \partial_au^0Bu^b\\
&-2(c_{\rm s}^{-1}c_{\rm s}''-c_{\rm s}^{-2}(c_{\rm s}')^2)(Bu^0)^2 Bu^a-2\delta^{bc}(c_{\rm s}^3c_{\rm s}''-3c_{\rm s}^2(c_{\rm s}')^2) \partial_au^0\partial_bu^0\partial_cu^0\\
&-2((c_{\rm s}')^2+c_{\rm s}c_{\rm s}'')(Bu^0)^2\partial_au^0+2\delta^{cd}\partial_du^0(c_{\rm s}^3c_{\rm s}''+c_{\rm s}^2(c_{\rm s}')^2)\partial_cu^0\partial_au^0+\mcl{R}.}
We now use the following identities: $8(c_{\rm s}')^2 \partial_bu^0 \partial_au^0Bu^b=-8c_{\rm s}^2(c_{\rm s}')^2 \delta^{bc}\partial_bu^0 \partial_au^0\partial_cu^0+\mcl{R}$ and $-2((c_{\rm s}')^2+c_{\rm s}c_{\rm s}'')(Bu^0)^2\partial_au^0=2(c_{\rm s}^{-2}(c_{\rm s}')^2+c_{\rm s}^{-1}c_{\rm s}'')(Bu^0)^2 Bu^a+\mcl{R}$. Thus, the sum above equals
\fm{&-4c_{\rm s}^{-2}(c_{\rm s}')^2(Bu^0)^2 Bu^a-8c_{\rm s}^2(c_{\rm s}')^2 \delta^{bc}\partial_bu^0 \partial_au^0\partial_cu^0\\
&-2(c_{\rm s}^{-1}c_{\rm s}''-c_{\rm s}^{-2}(c_{\rm s}')^2)(Bu^0)^2 Bu^a-2\delta^{bc}(c_{\rm s}^3c_{\rm s}''-3c_{\rm s}^2(c_{\rm s}')^2) \partial_au^0\partial_bu^0\partial_cu^0\\
&+2(c_{\rm s}^{-2}(c_{\rm s}')^2+c_{\rm s}^{-1}c_{\rm s}'')(Bu^0)^2 Bu^a+2\delta^{cd}(c_{\rm s}^3c_{\rm s}''+c_{\rm s}^2(c_{\rm s}')^2)\partial_du^0\partial_cu^0\partial_au^0+\mcl{R}.}
By computing the coefficients of $(Bu^0)^2Bu^a$ and of $\delta^{bc}\partial_au^0\partial_bu^0\partial_cu^0$, we conclude that this sum is $\mcl{R}$.
\end{proof}\rm 

Finally, we consider $w_{ab}$.
\lem{We have $-BBw_{ab}+c_{\rm s}^2\Delta w_{ab}=\mcl{R}$ for $a,b=1,2,3$.}
\begin{proof}
By \eqref{for8.38}, we write $-BBw_{ab}+c_{\rm s}^2\Delta w_{ab}$ as 
\fm{\partial_a(-BBu^b+c_{\rm s}^2\Delta u^b)-\partial_b(-BBu^a+c_{\rm s}^2\Delta u^a)+[-BB+c_{\rm s}^2\Delta,\partial_a]u^b-[-BB+c_{\rm s}^2\Delta,\partial_b]u^a.}
By Lemmas \ref{techlemm8.16} and \ref{techlemm8.17}, we have
\fm{&\partial_a(-BBu^b+c_{\rm s}^2\Delta u^b)+[-BB+c_{\rm s}^2\Delta,\partial_a]u^b\\
&=-2(c_{\rm s}^{-1}c_{\rm s}''-3c_{\rm s}^{-2}(c_{\rm s}')^2)\partial_au^0 Bu^0Bu^b\\
&\quad+B(\partial_au^c\partial_cu^b)+\partial_au^c\partial_cBu^b+ Bu^c \partial_a \partial_cu^b+\partial_a Bu^c \partial_cu^b-2c_{\rm s}^{-1}c_{\rm s}' Bu^0\partial_aBu^b\\
&\quad +2c_{\rm s}^{-3}c_{\rm s}'B Bu^a Bu^b-2c_{\rm s}c_{\rm s}'\partial_au^0\Delta u^b-2c_{\rm s}^{-1}c_{\rm s}'\partial_au^c\partial_cu^0Bu^b+\mcl{R}.}
Later, we will interchange the roles of $a$ and $b$, and then take the difference.

For the first row on the right side, we note that $-2(c_{\rm s}^{-1}c_{\rm s}''-3c_{\rm s}^{-2}(c_{\rm s}')^2)\partial_au^0 Bu^0Bu^b=2(c_{\rm s}^{-1}c_{\rm s}''-3c_{\rm s}^{-2}(c_{\rm s}')^2)Bu^a Bu^0Bu^b+\mcl{R}$ is symmetric with respect to $a,b$ modulus $\mcl{R}$. If we interchange the roles of $a$ and $b$ and take the difference, we simply get $\mcl{R}$. 

For the second row, we have 
\fm{&B(\partial_au^c\partial_cu^b)+\partial_au^c\partial_cBu^b+ Bu^c \partial_a \partial_cu^b+\partial_a Bu^c \partial_cu^b-2c_{\rm s}^{-1}c_{\rm s}' Bu^0\partial_aBu^b\\
&=B(\partial_au^c\partial_cu^b)+\partial_au^cB\partial_cu^b+\partial_au^c\partial_cu^d\partial_du^b+ Bu^c \partial_c\partial_au^b\\
&\quad+B\partial_a u^c \partial_cu^b+\partial_a u^d\partial_du^c \partial_cu^b-2c_{\rm s}^{-1}c_{\rm s}' Bu^0B\partial_au^b-2c_{\rm s}^{-1}c_{\rm s}' Bu^0\partial_au^c\partial_cu^b\\
&=2B(\partial_au^c\partial_cu^b)+2\partial_au^c\partial_cu^d\partial_du^b+ Bu^c \partial_c \partial_au^b-2c_{\rm s}^{-1}c_{\rm s}' Bu^0B\partial_au^b-2c_{\rm s}^{-1}c_{\rm s}' Bu^0\partial_au^c\partial_cu^b.}
Recall that we need to interchange the roles of $a$ and $b$ and take the difference. Note that
\fm{\partial_au^c\partial_cu^b-\partial_bu^c\partial_cu^a&=w_{ac}\partial_cu^b+w_{db} \partial_cu^a,}
so we have $2B(\partial_au^c\partial_cu^b)-2B(\partial_bu^c\partial_cu^a), -2c_{\rm s}^{-1}c_{\rm s}' Bu^0(\partial_au^c\partial_cu^b-\partial_bu^c\partial_cu^a)=\mcl{R}$. Similarly, we can show that $2\partial_au^c\partial_cu^d\partial_du^b-2\partial_bu^c\partial_cu^d\partial_du^a=\mcl{R}$. Moreover, $Bu^c\partial_c\partial_au^b-Bu^c\partial_c\partial_bu^a=Bu^c\partial_cw_{ab}=\mcl{R}$ and $-2c_{\rm s}^{-1}c_{\rm s}' Bu^0B(\partial_au^b-\partial_bu^a)=-2c_{\rm s}^{-1}c_{\rm s}' Bu^0Bw_{ab}=\mcl{R}$. In summary, the second row also generates $\mcl{R}$ only.

For the third row, we have
\fm{&2c_{\rm s}^{-3}c_{\rm s}'B Bu^a Bu^b-2c_{\rm s}c_{\rm s}'\partial_au^0\Delta u^b-2c_{\rm s}^{-1}c_{\rm s}'\partial_au^c\partial_cu^0Bu^b\\
&=2c_{\rm s}^{-3}c_{\rm s}'B Bu^a Bu^b+2c_{\rm s}^{-1}c_{\rm s}'Bu^a\Delta u^b-2c_{\rm s}^{-1}c_{\rm s}'\partial_au^c\partial_cu^0Bu^b+\mcl{R}
\\
&=2c_{\rm s}^{-3}c_{\rm s}'B Bu^a Bu^b+2c_{\rm s}^{-3}c_{\rm s}'Bu^a(BBu^b-2c_{\rm s}^{-1}c_{\rm s}'Bu^0Bu^b-c_{\rm s}^2\delta^{cd}\partial_cu^0\partial_du^b)\\
&\quad-2c_{\rm s}^{-1}c_{\rm s}'\partial_au^c\partial_cu^0Bu^b+\mcl{R}\\
&=2c_{\rm s}^{-3}c_{\rm s}'(B Bu^a Bu^b+BBu^bBu^a)-4c_{\rm s}^{-4}(c_{\rm s}')^2Bu^0Bu^bBu^a\\
&\quad-2c_{\rm s}^{-1}c_{\rm s}'\partial_cu^0(\partial_bu^cBu^a+\partial_au^c Bu^b)+\mcl{R}.}
In the second identity, we use \eqref{sec8eqnuaBB}. It is clear that every term here is symmetric with respect to $a,b$ modulus $\mcl{R}$. In summary, we have $-BBw_{ab}+c_{\rm s}^2\Delta w_{ab}=\mcl{R}$.
\end{proof}\rm

\bibliography{paperjce}{}
\bibliographystyle{plain}

\end{document}